\documentclass[11pt,a4paper]{amsart}
\usepackage[latin9]{inputenc}
\usepackage{mathtools}
\usepackage{amscd}
\usepackage{amstext}
\usepackage{amsthm}
\usepackage{amssymb}
\usepackage[unicode=true,
 bookmarks=false,
 breaklinks=false,pdfborder={0 0 1},backref=false,colorlinks=true,linkcolor=blue,citecolor=blue]
 {hyperref}
\usepackage{cleveref}
\usepackage{soul}
\usepackage{xcolor}
\usepackage{comment}

\makeatletter

\newcommand*\LyXbar{\rule[0.585ex]{1.2em}{0.25pt}}
\pdfpageheight\paperheight
\pdfpagewidth\paperwidth

\numberwithin{equation}{section}
\numberwithin{figure}{section}
\theoremstyle{plain}
\newtheorem{thm}{\protect\theoremname}[section]
  \theoremstyle{definition}
  
 \theoremstyle{definition}
 \newtheorem{quesintro}{}
 
  \theoremstyle{definition}
  \newtheorem{defn}[thm]{\protect\definitionname}
 \theoremstyle{definition}
 \newtheorem{stateintro}{}
 
  \theoremstyle{definition}
  \newtheorem{example}[thm]{\protect\examplename}
  \theoremstyle{definition}
  \newtheorem{notation}[thm]{\protect{Notation}}
  \theoremstyle{plain}
  \newtheorem{prop}[thm]{\protect\propositionname}
  \theoremstyle{plain}
  \newtheorem{cor}[thm]{\protect\corollaryname}
  \theoremstyle{plain}
  \newtheorem{lem}[thm]{\protect\lemmaname}
 \theoremstyle{definition}
 \newtheorem{rem}[thm]{\protect\remarkname}

\@ifundefined{date}{}{\date{}}

\usepackage{tikz}
\usepackage{amsthm}
\usepackage[hang]{footmisc}
\usepackage[all]{xypic}
\usepackage{color}
\usepackage[inline]{enumitem}
\usepackage{mathrsfs}

\allowdisplaybreaks

\makeatother

  \providecommand{\corollaryname}{Corollary}
  \providecommand{\definitionname}{Definition}
  \providecommand{\examplename}{Example}
  \providecommand{\lemmaname}{Lemma}
  \providecommand{\problemname}{Problem}
  \providecommand{\propositionname}{Proposition}
 \providecommand{\remarkname}{Remark}
\providecommand{\theoremname}{Theorem}

\begin{document}

\global\long\def\floorstar#1{\lfloor#1\rfloor}
\global\long\def\ceilstar#1{\lceil#1\rceil}

\global\long\def\B{B}
 \global\long\def\A{A}
 \global\long\def\J{J}
 \global\long\def\K{\mathcal{K}}
 \global\long\def\D{D}
 \global\long\def\Ch{D}
 \global\long\def\Zh{\mathcal{Z}}
 \global\long\def\E{E}
 \global\long\def\Oh{\mathcal{O}}

\global\long\def\T{{\mathbb{T}}}
 \global\long\def\BR{{\mathbb{R}}}
 \global\long\def\N{{\mathbb{N}}}
 \global\long\def\Z{{\mathbb{Z}}}
 \global\long\def\C{{\mathbb{C}}}
 \global\long\def\Q{{\mathbb{Q}}}

\global\long\def\aut{\mathrm{Aut}}
 \global\long\def\supp{\mathrm{supp}}

\global\long\def\eps{\varepsilon}

\global\long\def\id{\mathrm{id}}

\global\long\def\halpha{\widehat{\alpha}}
 \global\long\def\calpha{\widehat{\alpha}}

\global\long\def\tih{\widetilde{h}}

\global\long\def\opFol{\operatorname{F{\o}l}}

\global\long\def\opRange{\operatorname{Range}}

\global\long\def\opIso{\operatorname{Iso}}

\global\long\def\dimnuc{\dim_{\operatorname{nuc}}}

\global\long\def\set#1{\left\{  #1\right\}  }


\global\long\def\mset#1{\left\{  \!\!\left\{  #1\right\}  \!\!\right\}  }



\global\long\def\IA{\mathbb{A}}
 \global\long\def\IB{\mathbb{B}}
 \global\long\def\IC{\mathbb{C}}
 \global\long\def\ID{\mathbb{D}}
 \global\long\def\IE{\mathbb{E}}
 \global\long\def\IF{\mathbb{F}}
 \global\long\def\IG{\mathbb{G}}
 \global\long\def\IH{\mathbb{H}}
 \global\long\def\II{\mathbb{I}}
 \global\long\def\IJ{\mathbb{J}}
 \global\long\def\IK{\mathbb{K}}
 \global\long\def\IL{\mathbb{L}}
 \global\long\def\IM{\mathbb{M}}
 \global\long\def\IN{\mathbb{N}}
 \global\long\def\IO{\mathbb{O}}
 \global\long\def\IP{\mathbb{P}}
 \global\long\def\IQ{\mathbb{Q}}
 \global\long\def\IR{\mathbb{R}}
 \global\long\def\IS{\mathbb{S}}
 \global\long\def\IT{\mathbb{T}}
 \global\long\def\IU{\mathbb{U}}
 \global\long\def\IV{\mathbb{V}}
 \global\long\def\IW{\mathbb{W}}
 \global\long\def\IX{\mathbb{X}}
 \global\long\def\IY{\mathbb{Y}}
 \global\long\def\IZ{\mathbb{Z}}


\global\long\def\CA{\mathcal{A}}
 \global\long\def\CB{\mathcal{B}}
 \global\long\def\CC{\mathcal{C}}
 \global\long\def\CalD{\mathcal{D}}
 \global\long\def\CD{\mathcal{D}}
 \global\long\def\CE{\mathcal{E}}
 \global\long\def\CF{\mathcal{F}}
 \global\long\def\CG{\mathcal{G}}
 \global\long\def\CH{\mathcal{H}}
 \global\long\def\CI{\mathcal{I}}
 \global\long\def\CJ{\mathcal{J}}
 \global\long\def\CK{\mathcal{K}}
 \global\long\def\CL{\mathcal{L}}
 \global\long\def\CM{\mathcal{M}}
 \global\long\def\CN{\mathcal{N}}
 \global\long\def\CO{\mathcal{O}}
 \global\long\def\CP{\mathcal{P}}
 \global\long\def\CQ{\mathcal{Q}}
 \global\long\def\CR{\mathcal{R}}
 \global\long\def\CS{\mathcal{S}}
 \global\long\def\CT{\mathcal{T}}
 \global\long\def\CU{\mathcal{U}}
 \global\long\def\CV{\mathcal{V}}
 \global\long\def\CW{\mathcal{W}}
 \global\long\def\CX{\mathcal{X}}
 \global\long\def\CY{\mathcal{Y}}
 \global\long\def\CZ{\mathcal{Z}}


\global\long\def\FA{\mathfrak{A}}
 \global\long\def\FB{\mathfrak{B}}
 \global\long\def\FC{\mathfrak{C}}
 \global\long\def\FD{\mathfrak{D}}
 \global\long\def\FE{\mathfrak{E}}
 \global\long\def\FF{\mathfrak{F}}
 \global\long\def\FG{\mathfrak{G}}
 \global\long\def\FH{\mathfrak{H}}
 \global\long\def\FI{\mathfrak{I}}
 \global\long\def\FJ{\mathfrak{J}}
 \global\long\def\FK{\mathfrak{K}}
 \global\long\def\FL{\mathfrak{L}}
 \global\long\def\FM{\mathfrak{M}}
 \global\long\def\FN{\mathfrak{N}}
 \global\long\def\FO{\mathfrak{O}}
 \global\long\def\FP{\mathfrak{P}}
 \global\long\def\FQ{\mathfrak{Q}}
 \global\long\def\FR{\mathfrak{R}}
 \global\long\def\FS{\mathfrak{S}}
 \global\long\def\FT{\mathfrak{T}}
 \global\long\def\FU{\mathfrak{U}}
 \global\long\def\FV{\mathfrak{V}}
 \global\long\def\FW{\mathfrak{W}}
 \global\long\def\FX{\mathfrak{X}}
 \global\long\def\FY{\mathfrak{Y}}
 \global\long\def\FZ{\mathfrak{Z}}


\global\long\def\Ra{\Rightarrow}
 \global\long\def\La{\Leftarrow}
 \global\long\def\LRa{\Leftrightarrow}

\global\long\def\quer{\overline{}}
 \global\long\def\eins{\mathbf{1}}
 \global\long\def\diag{\operatorname{diag}}
 \global\long\def\ad{\operatorname{Ad}}
 \global\long\def\ev{\operatorname{ev}}
 \global\long\def\fin{{\subset\!\!\!\subset}}
 \global\long\def\diam{\operatorname{diam}}
 \global\long\def\Hom{\operatorname{Hom}}
 \global\long\def\dst{{\displaystyle }}
 \global\long\def\spp{\operatorname{supp}}
 \global\long\def\spo{\operatorname{supp}_{o}}
 \global\long\def\del{\partial}
 \global\long\def\lsc{\operatorname{Lsc}}
 \global\long\def\GU{\CG^{(0)}}
 \global\long\def\HU{\CH^{(0)}}
 \global\long\def\AU{\CA^{(0)}}
 \global\long\def\BU{\CB^{(0)}}
 \global\long\def\CUU{\CC^{(0)}}
 \global\long\def\DU{\CD^{(0)}}
 \global\long\def\CUUU{\CC'{}^{(0)}}

\global\long\def\AUl{(\CA^{l})^{(0)}}
\global\long\def\BUl{(B^{l})^{(0)}}
\global\long\def\HUp{(\CH^{p})^{(0)}}

\global\long\def\properlength{proper}

\global\long\def\interior#1{#1^{\operatorname{o}}}


\title[Almost elementariness]{Almost elementary \'{e}tale groupoids
}

\author{Xin~Ma}

\address{X.~Ma: Institute for Advanced Study in Mathematics, Harbin Institute of Technology, Harbin, 150001, China}

\email{xma17@hit.edu.cn}

\author{Jianchao~Wu}

\address{J.~Wu: Shanghai Center for Mathematical Sciences, Fudan University,
	Shanghai 200438, China}

\email{jianchao\_wu@fudan.edu.cn}



\begin{abstract}
Motivated by Matui and Kerr's work on almost finiteness, we introduce a new finite approximation property for (possibly non-ample) \'{e}tale groupoids called {almost elementariness}, which unifies and generalizes both almost finiteness and pure infiniteness. 
This property serves as a dynamical analogue of regularity properties of $C^*$-algebras. In support of this view, we prove that minimal almost elementary groupoids yield tracially $\mathcal{Z}$-stable reduced groupoid $C^*$-algebras. 
Consequently, we obtain as a corollary that the reduced $C^*$-algebras of all minimal amenable second countable almost finite groupoids in Matui's sense are $\mathcal{Z}$-stable and thus classifiable by the Elliott invariants.

Two basic ingredients underlying the definition of almost elementariness are castles in groupoids and groupoid subequivalence, both of which are developed extensively in this work.  
Notably, building on our flexible framework of castles, we introduce the technique of nesting of castles, which play a key role in the proof of our main theorem. 

Besides our main theorem on tracial $\mathcal{Z}$-stability, we also discuss in depth the relations between almost elementariness and other properties for groupoids such as effectiveness, the groupoid small boundary property, groupoid strict comparison and Matui's and Kerr's notions of almost finiteness. 
\end{abstract}

\maketitle
\tableofcontents{}

\section{Introduction}

Groupoids are a generalization of groups where the multiplication
operation is allowed to be only partially defined. The study of topological
groupoids lies at the crossroads of group theory, dynamics, geometry,
topology, mathematical physics, and operator algebras, largely thanks
to their Swiss-knife-like ability to handle many mathematical objects,
such as groups, group actions, equivalence relations on topological
spaces, nonperiodic tilings, etc., through a unifying framework. 

A recurring theme in many of these mathematical topics is the various
ways in which infinite structures may be approximated by finite structures.
Analysis often enters the picture this way. An influential poster
child of this theme is the notion of amenability in group theory,
together with its many reincarnations in other fields, such as injectivity
in von Neumann algebra theory, nuclearity in $C^{*}$-algebraic theory,
topological amenability in topological dynamics, and metric amenability
in coarse geometry, etc. For topological groupoids \LyXbar \LyXbar{}
usually assumed to be locally compact, $\sigma$-compact and Hausdorff,
and sometimes also \emph{\'{e}tale}, which is a groupoid generalization
of having discretized time, as opposed to continuous time, in a classical
topological dynamical system \LyXbar \LyXbar{} the richness of their
structures is reflected in a number of different yet intricately related
approximation properties that realize this theme. 

Among the strongest approximation properties for topological groupoids
is the notion of \emph{AF groupoids} (\cite{Renault1980groupoid}).
The terminology AF was borrowed from operator algebras and was originally
short for \emph{approximately finite}. This property applies to \emph{ample}
(i.e., totally disconnected) \'{e}tale groupoids and demands that every
compact subset of a topological groupoid is contained in a subgroupoid
that is \emph{elementary }\LyXbar \LyXbar{} namely, it is isomorphic
to a principal (i.e., not containing nontrivial subgroups) finite
groupoid, typically ``fattened up'' with topological spaces. 
The way these elementary groupoids embed reminds one of Kakutani-Rokhlin
towers, a fundamental tool in measure-theoretic and topological dynamics.
Thanks to the transparent and rigid structure of elementary groupoids,
AF groupoids are well studied and classified (\cite{Krieger1979/80dimension,GiordanoPutnamSkau2004Affable,GiordanoMatuiPutnamSkau2008absorption}),
though the notion is relatively restrictive. 

Inspired by the F{\o}lner set approach to (group-theoretic) amenability,
Matui \cite{Matui2012Homology} introduced a more general notion called
\emph{almost finite groupoids}\footnote{Almost finiteness should not be abbreviated as AF or confused with
\emph{almost AF} \emph{groupoids} (\cite{Phillips2005Crossed})! }, which, like AF groupoids, also applies to ample \'{e}tale groupoids,
but it only demands that every compact subset of a topological groupoid
is \emph{almost }contained in an elementary subgroupoid in a F{\o}lner-like
sense. Almost finiteness strikes a remarkable balance between applicability
and utility: it was shown to be enjoyed by \emph{transformation groupoids
}arising from free actions on the Cantor set by $\Z^{n}$ (later generalized
in \cite{KerrSzabo2020Almost}; see below), as well as those arising
from aperiodic tilings (\cite{ItoWhittakerZacharias-tiling}); on
the other hand, this notion has found applications in the homology
theory of ample groupoids, topological full groups, and the structure
theory and $K$-theory of reduced groupoid $C^{*}$-algebras, and
deep connections to an increasing number of other important properties
have been established (see, for example, \cite{Matui2015Topological,Matui2017Topological,Nekrashevych2019Simple,Kerr2020Dimension,Suzuki2020Almost}). 

Focusing on the case of transformation groupoids (arising from \emph{topological
dynamical systems}, that is, groups acting on topological spaces),
Kerr \cite{Kerr2020Dimension} presented an insightful perspective
that links almost finiteness with \emph{regularity properties} in
the classification and structure theory of $C^{*}$-algebras, extending
the link between AF groupoids and AF $C^{*}$-algebras, as well as
paralleling the link between \emph{hyperfiniteness} for ergodic probability-measure-preserving
equivalence relations and hyperfiniteness for II$_{1}$ factors. Here,
``regularity properties'' refers to a handful of natural and intrinsic
properties of $C^{*}$-algebras pioneered by Winter that arose from
Elliott's classification program of simple separable nuclear $C^{*}$-algebras
and played pivotal roles in its eventual success: in short, they were
the missing piece needed to characterize the $C^{*}$-algebras classifiable
via the Elliott invariant \cite{GLN2020a, GLN2020b,ElliottGongLinNiu2015classification,TikuisisWhiteWinter2017Quasidiagonality,CGSTW,
Phillips2000classification}.
Moreover, as predicted by Toms and Winter, this handful of regularity
properties turn out to be (mostly) equivalent for simple separable
nuclear $C^{*}$-algebras (\cite{Rordam2004stable,MatuiSato2012Strict,Winter2012Nuclear,SatoWhiteWinter2015Nuclear,TikuisisWhiteWinter2017Quasidiagonality,CastillejosEvingtonTikuisisWhiteWinter2019Nuclear});
thus having one of them is good enough for classification. Among these
regularity properties, it was Hirshberg and Orovitz's \emph{tracial
$\CZ$-stability} (\cite{HirshbergOrovitz2013Tracially}) that
Kerr found to be the closest in spirit to Matui's almost finiteness.

Furthermore, Kerr did not merely translate Matui's almost finiteness
into the language of group actions \LyXbar \LyXbar{} this would just
mean that we have a group action on a totally disconnected space that
admits partitions into open (in fact, clopen) Rokhlin towers modelled
on F{\o}lner sets. He  installed an upgrade to its applicability: while
the elementary subgroupoids in Matui's almost finiteness are required
to cover the entire unit space, Kerr's almost finiteness allows a
``small'' \emph{remainder} to be left uncovered. This slight relaxation
has immense consequences: while for actions on a totally disconnected
(i.e., zero-dimensional) space, this ``small'' remainder can always
be absorbed into the Rokhlin towers and we thus recover Matui's definition,
yet Kerr's version now also applies to actions on higher-dimensional
spaces (though in this case we also need to explicitly require that
the levels of the Rokhlin towers all have tiny diameters). This idea
of approximation modulo a small remainder is well aligned with the
understanding of tracial $\CZ$-stability as an analogue of the McDuff
property that allows for a tracially small error \LyXbar \LyXbar{}
in a way similar to Lin's earlier definition of tracially AF $C^{*}$-algebras
(\cite{Lin2001Tracially,Lin2001introduction}). 

This pivotal upgrade, however, depends on the precise meaning of ``small''
sets. There are two natural ways to describe them: The first is measure-theoretic:
namely, a set is ``small'' if, with regard to every invariant measure,
its measure is smaller than a predetermined positive number $\eps$.
The other is topological: roughly speaking, a set is ``small'' if
it is \emph{dynamically subequivalent} to a (predetermined) nonempty
open set\footnote{Kerr's definition actually requires the remainder to be dynamically
subequivalent to a small portion of the unit space of the open elementary
subgroupoid, instead of a predetermined nonempty open set. This makes
it work better in the non-minimal setting. }, i.e., roughly speaking, the ``small'' set is able to be disassembled
and then translated, piece by piece via the group action, into ``non-touching''
positions inside the nonempty open set, where ``non-touching'' means
the closures of the translated pieces do not intersect. It is clear
that the second method yields a stricter sense of smallness, and it
turns out to be a desirable property of a topological dynamical system
for the two methods to agree. This is the essence of what Kerr dubbed
\emph{dynamical (strict) comparison}, after the analogous property
of \emph{strict comparison }of positive elements in a $C^{*}$-algebra,
which is also among the aforementioned handful of $C^{*}$-algebraic
regularity properties. 

Thus Kerr's almost finiteness comes in two flavors for higher-dimensional
spaces: the ordinary one uses dynamical subequivalence to express
smallness of the remainder and an auxiliary notion called \emph{almost
finiteness in measure}, which uses invariant measures instead. Kerr
and Szab\'{o} \cite{KerrSzabo2020Almost} showed that the former condition
is equivalent to the conjunction of the latter and dynamical strict
comparison, while, at least for free actions, the latter condition is equivalent to the \emph{small
boundary property} (which, in turn, is closely related to \emph{mean
dimension zero} as well as $C^*$-regularity properties; see \cite{GiolKerr2010Subshifts, ElliottNiu2017C, Niu2019Comparisona, Niu2019Comparisonb, Niu2020Z, Naryshkin2024URP}). To cement the link to $C^{*}$-algebraic regularity
properties, Kerr \cite{Kerr2020Dimension} proved that a free minimal almost
finite action on a compact metrizable space by a countable discrete amenable group gives rise to
a tracially $\CZ$-stable crossed product $C^{*}$-algebra. 
This was applied to show that free minimal actions on finite-dimensional metric spaces by a large class of amenable groups give rise to classifiable crossed product $C^{*}$-algebras \cite{KerrSzabo2020Almost,KerrNaryshkin2025Elementary,NaryshkinPetrakos2025Almost}.  
Based on these evidences, Kerr suggested, at least in the case of actions
by amenable groups, almost finiteness may be understood as a dynamical
regularity property. 

This great confluence of ideas from topological dynamics, topological
groupoid theory, and operator algebras opened the gate to a plethora
of new connections and applications. At the same time, it also left
us with a number of unanswered questions and problems. 
\begin{quesintro}
\label{question:Z-stability}Do almost finite groupoids in Matui's
sense always give rise to groupoid $C^{*}$-algebras satisfying regularity
properties such as tracial $\CZ$-stability? 
\end{quesintro}
Kerr's result above answered this in the affirmative for transformation
groupoids, but it remained largely unclear beyond that case. Kerr's
method appears difficult to generalize, for it makes use of the fact
that the levels in the Rokhlin towers witnessing almost finiteness
are labeled by group elements, which allows one to apply Ornstein
and Weiss' theory of quasi-tilings to carefully manipulate the towers.
 Ito, Whittaker and Zacharias \cite{ItoWhittakerZacharias-tiling}
managed to extend Kerr's method and result to the case of \'{e}tale groupoids
from aperiodic tilings, exploiting the fact that \'{e}tale groupoids from
aperiodic tilings are reductions of transformation groupoids associated
to $\BR^{n}$-actions. Nevertheless, this method appears unsuitable
for groupoids without obvious underlying group structures. 
\begin{quesintro}
\label{question:infinite}Is there a groupoid regularity property
that works both for groupoids with invariant measures on their unit
spaces and for those without? 
\end{quesintro}
Both Matui's and Kerr's notions of almost finiteness imply the existence
of invariant measures on the unit space of a groupoid and thus the
existence of traces in their groupoid $C^{*}$-algebras. While this
is often a useful feature, it does not line up with the fact that
$C^{*}$-algebraic regularity properties also applies to \emph{purely
infinite algebras}, which do not have any traces. Thus one would hope
a more relaxed notion could include groupoids without invariant measures
on their unit spaces. We point out that the correspondence between
dynamical strict comparison and pure infiniteness was established
by the first author \cite{Ma-purely}.
\begin{quesintro}
\label{question:amenability}Related to the previous item: Can we
isolate from Matui's definiton of almost finiteness an ``amenability
property'' that is responsible for the existence of invariant measures? 
\end{quesintro}
Both Kerr's and Matui's definitions made explicit use of F{\o}lner sets
or a F{\o}lner-type condition \LyXbar \LyXbar{} this is the immediate
reason for the existence of invariant measures. While Kerr's notion
forces the acting group to be amenable, it suggested an analogous
kind of amenability property for \'{e}tale groupoids might hide behind
Matui's notion. A possible candidate was topological amenability,
but this would be the wrong target, as it does not imply the existence
of invariant measures in general, and a number of examples have shown
almost finiteness does not imply topological amenability (\cite{AraBonickeBosaLi-almost,Elek-qualitative}).

Motivated by \ref{question:amenability}, we introduced a notion termed
\emph{(ubiquitous) fiberwise amenability} for \'{e}tale groupoids in \cite{MWfiberwise}. It
simply demands, for any $\eps$ and any compact subset $K$ in our
\'{e}tale groupoid $\CG$, there is a nonempty finite subset $F$ in $\CG$
such that 
\[
\frac{\left|KF\right|}{\left|F\right|}\leq1+\eps.
\]
Essentially, $F$ is what one may call a \emph{F{\o}lner set}. Moreover, a systematic coarse-geometric framework in \cite{MWfiberwise} was developed for \'{e}tale groupoids, within which we establish several foundational properties of the (ubiquitous) fiberwise amenability.

Notably, we established in \cite[Theorem D]{MWfiberwise} a F\o lner--paradoxical dichotomy for any \'{e}tale groupoid, where the alternative is determined by whether or not the groupoid satisfies fiberwise amenability. This dichotomy will serve as a key tool in the current paper to address \ref{question:Z-stability} and \ref{question:infinite}, which are the main goals of the current paper. The answers to these questions also enlarge the view of the dynamical/geometric study on the regularity properties of $C^*$-algebras from crossed products to the general groupoid $C^*$-algebras setting, which covers more examples. For instance, because of $K$-theoretical obstructions, unital AF algebras and the Jiang-Su algebra $\CZ$ cannot be written as $\Z$-crossed products while they all can be written as reduced groupoid $C^*$-algebras. On the other hand, it was proven in \cite{Li2020classifiability} that every classifiable $C^*$-algebra arises as a (twisted) groupoid $C^*$-algebra.

We summarize
our results, with the following standing assumption. 
\begin{description} \label{standing-assumption}
\item [{Assumption}] The groupoids below are $\sigma$-compact, locally
compact, Hausdorff, \'{e}tale topological groupoids. 
\end{description}

Motivated by \ref{question:infinite}, we introduce, for 
\'{e}tale groupoids with compact unit spaces, a new regularity
and approximation property called \textit{almost elementariness},
which extends both Matui's and Kerr's almost finiteness (c.f.,
Section~\ref{sec:Matui} and \ref{sec:Kerr}).   In fact, we establish the following: 
\begin{stateintro}[Theorem~\ref{thm:Kerr-AE}]
For transformation groupoids, Kerr's almost finiteness is equivalent
to the conjunction of almost elementariness and ubiquitous fiberwise
amenability (i.e., the amenability of the acting group). 
\end{stateintro}

\begin{stateintro}[Theorem~\ref{thm:Matui-AE}]
For ample $\sigma$-compact groupoids with a compact unit space, Matui's almost
finiteness is also equivalent to the conjunction of almost elementariness
and ubiquitous fiberwise amenability. 
\end{stateintro}

To formulate the notion of almost elementary groupoids, we first develop a systematic framework of  \textit{castles} in (possibly non-ample) \'{e}tale groupoids in Sections \ref{sec:castles} and \ref{sec:castles-2}, which may be understood as elementary subgroupoids plus a bit of extra data. Our definition of almost elementariness replaces the F{\o}lner-type conditions in Matui's and
Kerr's almost finiteness by a new condition that requires a castle in the approximation to be extendable to a larger
castle. This condition draws inspiration from
coarse geometry, e.g., from how, in the definition of the asymptotic
dimension of a metric space, we ask for a cover that is able to grow
by a predetermined large distance without increasing its multiplicity. 
Similarly, in our definition of almost elementary
groupoids, we would like our castle to be able to
grow by a predetermined distance (as measure by compact subsets in
the groupoid) while still being a castle (alternatively, it
should be able to shrink without jeopardizing the smallness of the
remainder). Replacing the F{\o}lner-type conditions by this new condition
allows our notion to break free from fiberwise amenability and apply
to groupoids both with and without invariant probability measures
on their unit spaces.  This inclusiveness with regard to invariant
measures underlies the fact that the groupoid $C^{*}$-algebras of
almost elementary groupoids include both stably finite algebras and
purely infinite ones\footnote{This is one reason why we decided not to use the word ``finite''
in the name of our new notion, after a helpful suggestion of George Elliott to the second author. }, a fact that makes almost elementariness a candidate for a closer
analogue of $C^{*}$-algebraic regularity properties. 

To support this analogy with $C^{*}$-algebraic regularity properties,
we establish several connections. First of all, we provide the following generalization of a result of \cite{KerrSzabo2020Almost}:  
\begin{stateintro}[Theorem~\ref{thm:AE-DSC}]
A minimal groupoid is almost elementary if and only if it has
groupoid strict comparison and is almost elementary in measure. 
\end{stateintro}

We also establish a direct link to tracial $\CZ$-stablility, providing
an affirmative answer to \ref{question:Z-stability}. Recall that
to show a $C^{*}$-algebra $A$ is tracially $\CZ$-stable, we need
to produce order zero maps from arbitrarily large matrix algebras
into $A$ with approximately central images and with only ``tracially
small'' defects from being unital. We also remind the reader of our standing assumption above. 
\begin{stateintro}[Theorem~\ref{thm: tracial Z stable}]
\label{stateintro:tracial-Z-stability}Let $\CG$ be a minimal \'{e}tale groupoid with a compact infinite unit space. Suppose $\CG$ is
almost elementary (e.g., almost finite). Then the reduced groupoid
$C^{*}$-algebra $C_{r}^{*}(\CG)$ is tracially $\CZ$-stable. 
\end{stateintro}
This answers \ref{question:Z-stability} above.  In addition, our theorem generalizes Kerr's result on almost finite actions of amenable
groups (\cite{Kerr2020Dimension}) in several ways: our theorem can
be applied to groupoids without obvious underlying group structures
and without topological amenability; even when restricted to transformation
groupoids, our result can now deal with possibly non-free actions by possibly nonamenable groups. For example, using a result of \cite{OrtegaScarparo2020almostfinite}, \ref{stateintro:tracial-Z-stability} can be immediately applied to all minimal actions on the Cantor set by the infinite dihedral group $\Z\rtimes \Z_2$, regardless of freeness (c.f., Example~\ref{8.15}). 

Moreover, as indicated in \ref{question:Z-stability}, our proof necessarily
takes an approach different from Kerr's. In place of the Ornstein-Weiss
tiling theory of amenable groups, we develop a ``nesting'' form
of almost elementariness, which is an approximation (modulo a small
remainder) of the groupoid $\CG$ by not one, but two open elementary
subgroupoids in a \emph{nested} position, a notion reminiscent of
how multi-matrix algebras embed into each other. Thus passing from
approximation by one elementary subgroupoid to approximation by a
nesting of two (or perhaps more) is akin to how, from the local definition
of an AF algebra, one can produce a tower of nested multi-matrix algebras
organized by a Bratteli diagram. Indeed, an open elementary subgroupoid
of $\CG$ will induce an order zero map from a multi-matrix algebra
into $C_{r}^{*}(\CG)$, and if another open elementary subgroupoid
is nested in the first one, then we have an embedding between two
multi-matrix algebras. By arranging the nesting to have a large \emph{multiplicity},
we can ensure this embedding has a large relative commutant. This
will essentially be the source of the desired large matrix algebra
together with an order zero map into $C_{r}^{*}(\CG)$ with an approximately
central image. The existence of the remainders in these approximations
unfortunately makes the proof appear technically complicated, but
their smallness will eventually guarantee that the resulting order
zero map only has a ``tracially small'' defect from being unital.

The following is a direct consequence of \ref{stateintro:tracial-Z-stability}
by combining results in \cite{CastillejosEvingtonTikuisisWhiteWinter2019Nuclear,ElliottGongLinNiu2015classification,GLN2020a,GLN2020b, CGSTW,HirshbergOrovitz2013Tracially,Phillips2000classification,TikuisisWhiteWinter2017Quasidiagonality}. 
\begin{stateintro}[Corollary~\ref{8.12}]
 Let $\CG$ be a second countable amenable minimal \'{e}tale groupoid
with a compact unit space. Suppose $\CG$ is almost elementary. Then
$C_{r}^{*}(\CG)$ is unital simple separable nuclear and $\CZ$-stable
and thus has nuclear dimension one. In addition, in this case $C_{r}^{*}(\CG)$
is classified by its Elliott invariant. Finally, if $M(\CG)\neq\emptyset$,
then $C_{r}^{*}(\CG)$ is quasidiagonal; if $M(\CG)=\emptyset$, then
$C_{r}^{*}(\CG)$ is a unital Kirchberg algebra. 
\end{stateintro}

The paper is organized as follows.  Section~\ref{sec:preliminaries} compiles the necessary background on \'{e}tale groupoids and their $C^*$-algebras, and establishes our notational conventions. Sections~\ref{sec:castles} and~\ref{sec:castles-2} lay out a comprehensive framework for castles, elementary subgroupoids, and their associated functions. In Sections~\ref{sec:in-measure}, \ref{sec:in-measure-extra}, \ref{sec:almost-elementary}, and~\ref{sec:comparison}, we formulate and analyze the notion of almost elementariness along with its variants---including the in-measure and relative versions---as well as the concept of groupoid strict comparison. Building on this foundation, Sections~\ref{sec:Matui} and~\ref{sec:Kerr} establish that the existing notion of almost finiteness in the literature is in fact equivalent to almost elementariness combined with the fiberwise amenability introduced in \cite{MWfiberwise}. Finally, Sections~\ref{sec:nesting} and~\ref{sec:Z-stable} present the nesting technique, through which we prove (tracial) $\mathcal{Z}$-stability for the $C^*$-algebras of minimal almost elementary groupoids.

\section{Preliminaries\label{sec:preliminaries}}
\subsection{Backgrounds on groupoids}
We recall some basic backgrounds on
\'{e}tale groupoids and $C^{*}$-algebras.
We refer to \cite{Renault1980groupoid} and \cite{Sims-groupoids}
as references for groupoids and we record several fundamental definitions
and results for locally compact Hausdorff \'{e}tale groupoids here. 
\begin{defn}
A \textit{groupoid} $\CG$ is a set equipped with a distinguished
subset $\CG^{(2)}\subseteq\CG\times\CG$, called the set of \textit{composable
pairs}, a product map $\CG^{(2)}\rightarrow\CG$, denoted by $(\gamma,\eta)\mapsto\gamma\eta$
and an inverse map $\CG\rightarrow\CG$, denoted by $\gamma\mapsto\gamma^{-1}$
such that the following hold 
\begin{enumerate}[label=(\roman*)]
\item If $(\alpha,\beta)\in\CG^{(2)}$ and $(\beta,\gamma)\in\CG^{(2)}$
then so are $(\alpha\beta,\gamma)$ and $(\alpha,\beta\gamma)$. In
addition, $(\alpha\beta)\gamma=\alpha(\beta\gamma)$ holds in $\CG$.
\item For all $\alpha\in\CG$ one has $(\gamma,\gamma^{-1})\in\CG^{(2)}$
and $(\gamma^{-1})^{-1}=\gamma$.
\item For any $(\alpha,\beta)\in\CG^{(2)}$ one has $\alpha^{-1}(\alpha\beta)=\beta$
and $(\alpha\beta)\beta^{-1}=\alpha$. 
\end{enumerate}
Every groupoid is equipped with a subset $\GU=\{\gamma\gamma^{-1}:\gamma\in\CG\}$
of $\CG$. We refer to elements of $\GU$ as \textit{units} and to
$\GU$ itself as the \textit{unit space}. We define two maps $s,r:\CG\rightarrow\GU$
by $s(\gamma)=\gamma^{-1}\gamma$ and $r(\gamma)=\gamma\gamma^{-1}$,
respectively, in which $s$ is called the \textit{source} map and
$r$ is called the \textit{range} map. 
\end{defn}
\label{paragraph:sections-in-groupoids}
When a groupoid $\CG$ is endowed with a locally compact Hausdorff
topology under which the product and inverse maps are continuous,
the groupoid $\CG$ is called a locally compact Hausdorff groupoid.
A locally compact Hausdorff groupoid $\CG$ is called \textit{\'{e}tale}
if the range map $r$ is a local homeomorphism from $\CG$ to itself,
which means for any $\gamma\in\CG$ there is an open neighborhood
$U$ of $\gamma$ such that $r(U)$ is open and $r|_{U}$ is a homeomorphism.
In this case, since the map of taking inverses is an involutive homeomorphism
on $\CG$ that intertwines $r$ and $s$, thus the source map $s$
is also a local homeomorphism. 
A subset $S$ in $\CG$ is called an \emph{$s$-section} (respectively, an \emph{$r$-section})
if there is an open neighborhood $U$ of $S$ in $\CG$ such that
the source map $s$ (respectively,
the range map $r$) restricts to a homeomorphism from $U$
onto an open subset of $\GU$. 
It is called a \emph{bisection} if it is both an $s$-section and an $r$-section. 
It is not hard to see a locally compact
Hausdorff groupoid is \'{e}tale if and only if its topology has a basis
consisting of open bisections. We say a locally compact Hausdorff
\'{e}tale groupoid $\CG$ is \textit{ample} if its topology has a basis
consisting of compact open bisections.

\begin{example}\label{exa:pair-groupoid}
	Let $X$ be a discrete space. The \emph{pair groupoid} over $X$ has the discrete space $X \times X$ as the underlying topological space, and the groupoid operations are defined so that for any $u,v,w \in X$, we have 
	\[
		r(u,v) = u, \quad s(u,v) = v, \quad (u,v)^{-1} = (v,u) \quad \text{ and } \quad (u,v)(v,w) = (u,w) \; .
	\]
\end{example}

\begin{example}
\label{exa:transformation-groupoid}
Let $X$ be a locally compact
Hausdorff space and $\Gamma$ be a discrete group. Then any action
$\Gamma\curvearrowright X$ by homeomorphisms induces a locally compact
Hausdorff \'{e}tale groupoid 
\[
X\rtimes\Gamma\coloneqq\{(\gamma x,\gamma,x):\gamma\in\Gamma,x\in X\}
\]
equipped with the relative topology as a subset of $X\times\Gamma\times X$.
In addition, $(\gamma x,\gamma,x)$ and $(\beta y,\beta,y)$ are composable
only if $\beta y=x$ and 
\[
(\gamma x,\gamma,x)(\beta y,\beta,y)=(\gamma\beta y,\gamma\beta y,y).
\]
One also defines $(\gamma x,\gamma,x)^{-1}=(x,\gamma^{-1},\gamma x)$
and announces that $\GU\coloneqq\{(x,e_{\Gamma},x):x\in X\}$. It
is not hard to verify that $s(\gamma x,\gamma,x)=x$ and $r(\gamma x,\gamma,x)=\gamma x$.
The groupoid $X\rtimes\Gamma$ is called a \textit{transformation
groupoid}. 
\end{example}
The following are several basic properties of locally compact Hausdorff
\'{e}tale groupoids whose proofs could be found in \cite{Sims-groupoids}.
\begin{prop}
Let $\CG$ be a locally compact Hausdorff \'{e}tale groupoid. Then $\GU$
is a clopen set in $\CG$. 
\end{prop}

For any set $D\subseteq\GU$, Denote by 
\[
\CG_{D}\coloneqq\{\gamma\in\CG:s(\gamma)\in D\},\ \CG^{D}\coloneqq\{\gamma\in\CG:r(\gamma)\in D\},\ \text{and}\ \ \CG_{D}^{D}\coloneqq\CG^{D}\cap\CG_{D}.
\]
For the singleton case $D=\{u\}$, we write $\CG_{u}$, $\CG^{u}$
and $\CG_{u}^{u}$ instead for simplicity. In this situation, we call
$\CG_{u}$ a \textit{source fiber} and $\CG^{u}$ a \textit{range
fiber}. In addition, each $\CG_{u}^{u}$ is a group, which is called
the \textit{isotropy} at $u$. We also denote by \[\opIso(\CG)=\bigcup_{u\in \GU}\CG^u_u=\{x\in \CG: s(x)=r(x)\}\] the isotropy of the groupoid $\CG$. We say a groupoid $\CG$ is \textit{principal}
if $\opIso(\CG)=\GU$. A groupoid $\CG$ is called \textit{topologically
principal} if the set $\{u\in\GU:\CG_{u}^{u}=\{u\}\}$ is dense in
$\GU$. The groupoid $\CG$ is also said to be \textit{effective} if $\opIso(\CG)^o=\GU$. Recall that effectiveness is equivalent to topological principalness  if $\CG$ is second countable (See \cite[Lemma 4.2.3]{Sims-groupoids}). Therefore, effectiveness is equivalent to the topological freeness of an action of a countable discrete group acting on a compact metrizable space by looking at the corresponding transformation groupoid.

A subset $D$ in $\GU$ is called $\CG$-\textit{invariant}
if $r(\CG D)=D$, which is equivalent to the condition $\CG^{D}=\CG_{D}$.
Note that $\CG|_{D}\coloneqq\CG_{D}^{D}$ is a subgroupoid of $\CG$
with the unit space $D$ if $D$ is a $\CG$-invariant set in $\GU$.
A groupoid $\CG$ is called \textit{minimal} if there are no proper
non-trivial closed $\CG$-invariant subsets in $\GU$.
The following lemma underlies the recurrence behavior of minimal groupoids
with compact unit spaces. We refer to, e.g., \cite[Lemma 5.13]{MWfiberwise} for a proof.
\begin{lem}\label{lem:minimal-recurrence}
Let $\CG$ be an
\'{e}tale groupoid. Let $K$ and $V$ be subsets
of $\CG^{(0)}$ such that $K$ is compact, $V$ is open, and $K \subseteq r \left( \CG \cdot V \right)$. Then there are compact bisections $K_{1},\dots,K_{n}$
such that $\bigcup_{i=1}^{n}r\left(K_{i}\right)\subseteq V$ and $K\subseteq\bigcup_{i=1}^{n}s\left(K_{i}\right)^{\operatorname{o}}$. 
\end{lem}

\begin{defn}
\label{def:invariant-measure}A measure on the unit space of an \'{e}tale groupoid $\CG$ is \emph{invariant} if $\mu(r(U))=\mu(s(U))$
for any Borel measuable bisection $U$. We write $M(\CG)$ for the collection
of all \emph{invariant regular Borel probability measures} on $\CG^{(0)}$. \end{defn}

\subsection{Fiberwise amenability of groupoids} The following concepts on the coarse geometric properties of groupoids were introduced in \cite{MWfiberwise}. We refer to  \cite{MWfiberwise} for more details.

\begin{defn}(\cite[Definition 5.1]{MWfiberwise})\label{def:groupoid-boundaries}
Let $\CG$ be a groupoid. For any subsets $A,K\subseteq\CG$, we define
the following boundary sets:
\begin{enumerate}[label=(\roman*)]
\item \emph{left outer $K$-boundary}: $\partial_{K}^{+}A=(KA)\setminus A=\{yx\in\CG\setminus A:y\in K,x\in A\}$;
\item \emph{left inner $K$-boundary}: $\partial_{K}^{-}A=A\cap(K^{-1}(\CG\setminus A))=\{x\in A:yx\in\CG\setminus A\mbox{ for some }y\in K\}$;
\item \emph{left $K$-boundary}: $\partial_{K}A=\partial_{K}^{+}A\cup\partial_{K}^{-}A$. 
\end{enumerate}
\end{defn}

\begin{rem}(\cite[Remark 5.2]{MWfiberwise})\label{rem: double control for boundries}
For any subsets $A,K\subseteq\CG$, it is straightforward
to see $\partial_{K}^{+}A\subseteq K\partial_{K}^{-}A$ and $\partial_{K}^{-}A\subseteq K^{-1}\partial_{K}^{+}A$. 
\end{rem}

\begin{defn}(\cite[Definition 5.3]{MWfiberwise})
\label{def:groupoid-Folner} Let $\CG$ be an \'{e}tale
groupoid. For any subset $K\subseteq\CG$ and $\varepsilon>0$, a finite
non-empty set $F\subseteq \CG$ is called $(K,\varepsilon)$-F{\o}lner if
it satisfies 
\[
\frac{|\partial_{K}F|}{|F|}\leq\varepsilon.
\]
We denote by $\opFol(K,\varepsilon)$ the collection of all $(K,\varepsilon)$-F{\o}lner
sets. 
\end{defn}
This leads to a natural definition of fiberwise amenability. 
\begin{defn}(\cite[Definition 5.4]{MWfiberwise})
\label{def:fiberwise-amenable} Let $\CG$ be a locally
compact \'{e}tale groupoid. 
\begin{enumerate}
\item We say $\CG$ is \textit{fiberwise amenable} if for any compact subset
$K$ of $\CG$ and any $\varepsilon>0$, there exists a $(K,\varepsilon)$-F{\o}lner
set.
\item We say $\CG$ is \textit{ubiquitously fiberwise amenable} if and only
if for any compact subset $K$ of $\CG$ and any $\varepsilon>0$, there
exists a compact subset $L$ of $\CG$ such that for any unit $u\in\CG^{(0)}$,
there is a $(K,\varepsilon)$-F{\o}lner set in $Lu\cup\{u\}$. 
\end{enumerate}
\end{defn}

These two notions coincide for minimal groupoids.

\begin{thm}(\cite[Theorem C]{MWfiberwise})\label{thm:minimal-fiberwise-amenability}
     Let $\CG$
be a $\sigma$-compact \'{e}tale groupoid. Suppose
$\CG$ is minimal. Then $\CG$ is fiberwise amenable if and only
if it is ubiquitously fiberwise amenable.
\end{thm}

We finally record the following dichotomy in this subsection. 

\begin{thm}(\cite[Theorem D]{MWfiberwise})\label{thm:fiberwise-amenable-dichotomy}
 Let $\CG$ be a 
$\sigma$-compact \'{e}tale groupoid with a compact unit space. 
Then we have the following dichotomy. 
\begin{enumerate}[label=(\roman*)]
\item\label{thm:fiberwise-amenable-dichotomy:amenable}
If $\CG$ is ubiquitous fiberwise amenable, then 
for any compact subset $K$ in $\CG$ and any $\varepsilon>0$, 
there is a compact set $L \subseteq \CG$ such that
for any finite set $M \subseteq\CG$, there
is a finite set $F$ satisfying 
\[
	M \subseteq F\subseteq LM \ \textrm{and}\ |KF|\leq(1+\varepsilon)|F|.
\]

\item\label{thm:fiberwise-amenable-dichotomy:non-amenable}
If $\CG$ is not fiberwise amenable, then 
for any compact set $K \subseteq\CG$
and any $n\in\mathbb{N}$, there is a compact set $L \subseteq\CG$ such that
for any finite set $M\subseteq\CG$, the set
$LM$ contains at least $n|M|$ many disjoint sets of the form $K x$ for $x \in \CG$. 
\end{enumerate}
\end{thm}
We remark that in the minimal case, Theorem \ref{thm:fiberwise-amenable-dichotomy}, combined with Theorem \ref{thm:minimal-fiberwise-amenability}, provides a genuine dichotomy.

\subsection{Groupoid $C^*$-algebras}
Let $\CG$ be a locally compact Hausdorff \'{e}tale groupoid. We define
a convolution product on $C_{c}(\CG)$ by 
\[
(f*g)(\gamma)=\sum_{\alpha\beta=\gamma}f(\alpha)g(\beta)
\]
and an involution by 
\[
f^{*}(\gamma)=\overline{f(\gamma^{-1})}.
\]
These two operations make $C_{c}(\CG)$ a $*$-algebra. Then the reduced
groupoid $C^{*}$-algebra $C_{r}^{*}(\CG)$ is defined to be the completion
of $C_{c}(\CG)$ with respect to the norm $\|\cdot\|_{r}$ induced
by all regular representation $\pi_{u}$ for $u\in\GU$, where $\pi_{u}:C_{c}(\CG)\rightarrow B(\ell^{2}(\CG_{u}))$
is defined by $\pi_{u}(f)\eta=f*\eta$ and $\|f\|_{r}=\sup_{u\in\GU}\|\pi_{u}(f)\|$.
It is well known that there is a $C^{*}$-algebraic embedding $\iota:C_{0}(\GU)\rightarrow C_{r}^{*}(\CG)$.
On the other hand, $E_{0}:C_{c}(\CG)\rightarrow C_{c}(\GU)$ defined
by $E_{0}(a)=a|_{\GU}$ extends to a faithful canonical conditional
expectation $E:C_{r}^{*}(\CG)\rightarrow C_{0}(\GU)$ satisfying $E(\iota(f))=f$
for any $f\in C_{0}(\GU)$ and $E(\iota(f)a\iota(g))=fE(a)g$ for
any $a\in C_{r}^{*}(\CG)$ and $f,g\in C_{0}(\GU)$.

As a typical example, it can be verified that for the transformation
groupoid in Example 2.2, the reduced groupoid $C^{*}$-algebra is
isomorphic to the reduced crossed product $C^{*}$-algebra of the
dynamical system. The following are some standard facts on reduced
groupoid $C^{*}$-algebras that could be found in \cite{Sims-groupoids}.
Throughout the paper, the notation $\spo(f)$ for a function $f$
on a topological space $X$ denotes the open support $\{x\in X:f(x)\neq0\}$
of $f$. In addition, we write $\supp(f)$ the usual support $\overline{\spo(f)}$
of $f$. We say an open set $O$ in a topological space $X$ is \textit{precompact}
if $\overline{O}$ is compact.
\begin{prop}
\label{2.6} Let $\CG$ be a locally compact Hausdorff \'{e}tale groupoid.
Any $f\in C_{c}(\CG)$ can be written as a sum $f=\sum_{i=0}^{n}f_{i}$
such that there are precompact open bisections $V_{0},\dots,V_{n}$
such that $V_{0}\subseteq\GU$ and $V_{i}\cap\GU=\emptyset$ for all
$0<i\leq n$ as well as $\supp(f_{i})\subseteq V_{i}$ for any $0\leq i\leq n$. 
\end{prop}

\begin{prop}\label{prop: product of functions supported on bisections}
Let $\CG$ be a locally compact Hausdorff \'{e}tale groupoid. Suppose
$U,V$ are open bisections and $f,g\in C_{c}(\CG)$ such that $\supp(f)\subseteq U$
and $\supp(g)\subseteq V$. Then $\supp(f*g)\subseteq U\cdot V$ and for
any $\gamma=\alpha\beta\in U\cdot V$ one has $(f*g)(\gamma)=f(\alpha)g(\beta)$. 
\end{prop}

\begin{prop}\label{prop: norms on bisections} Let $\CG$ be a locally compact Hausdorff \'{e}tale groupoid.
For $f\in C_{c}(\CG)$, one has $\|f\|_{\infty}\leq\|f\|_{r}.$ If
$f$ is supported on a bisection, then one has $\|f\|_{\infty}=\|f\|_{r}.$ 
\end{prop}
Let $\CG$ be a locally compact Hausdorff \'{e}tale groupoid. Suppose
$U$ is an open bisection and $f\in C_{c}(\CG)_{+}$ such that $\supp(f)\subseteq U$.
Define functions $s(f),r(f)\in C_{0}(\GU)$ by $s(f)(s(\gamma))=f(\gamma)$
and $r(f)(r(\gamma))=f(\gamma)$ for $\gamma\in\supp(f)$. Since $U$
is a bisection, so is $\supp(f)$. Then the functions $s(f)$ and
$r(f)$ are well-defined functions supported on $s(\supp(f))$ and
$r(\supp(f))$, respectively. Note that $s(f)=(f^{*}*f)^{1/2}$ and
$r(f)=(f*f^{*})^{1/2}$.

\textit{\emph{The }}\textit{Jiang-Su}\emph{ algebra} $\CZ$, introduced
in \cite{JiangSu1999simple} by Jiang and Su, is a infinite dimensional
unital nuclear simple separable $C^{*}$-algebra, but $KK$-equivalent
to $\C$ in the sense of Kasparov. We say a $C^{*}$-algebra $A$
is $\CZ$-\textit{stable} if $A\otimes\CZ\simeq A$.

\begin{notation}
	Throughout the paper, we write $B\sqcup C$ to indicate that
	the union of sets $B$ and $C$ is a disjoint union. In addition,
	we denote by $\bigsqcup_{i\in I}B_{i}$ for the disjoint union of
	the family $\{B_{i}:i\in I\}$. We also denote by $\ceilstar{\cdot}$
	the ceiling function and by $\floorstar{\cdot}$ the floor function
	from $[0,\infty)$ to $\N$. 
	For two collections $\CB$ and $\CC$ of subsets in a space $X$, we write $\CB \vee \CC := \left\{ B \cap C \colon B \in \CB, C \in \CC \right\}$. 
\end{notation}

Finally, we make the following assumption to simplify terminology.

\begin{notation} \label{standing-assumption-paper}
	Throughout the paper, we mean by ``\'{e}tale groupoids'' locally compact, Hausdorff, \'{e}tale topological groupoids. 
\end{notation}

\newcommand{\unitSpace}{^{(0)}}

\newcommand{\RedCD}[2]{\operatorname{Red}(#1,#2)}
\newcommand{\redCD}[2]{\operatorname{red}(#1,#2)}
\newcommand{\ResCD}[2]{\operatorname{Res}(#1,#2)}
\newcommand{\resCD}[2]{\operatorname{res}(#1,#2)}
\newcommand{\CommCD}[2]{\operatorname{Comm}(#1,#2)}
\newcommand{\commCD}[2]{\operatorname{comm}(#1,#2)}

\newcommand{\osupp}[1]{\operatorname{supp}^{\operatorname{o}}\left(#1\right)}

\section{Castles in \'{e}tale groupoids} \label{sec:castles}
In this section, we introduce the notion of castles on an \'{e}tale groupoid. They are akin to the notion with the same name defined by Kerr in the case of group actions (\cite{Kerr2020Dimension}) and can be thought of as elementary subgroupoids plus a bit more structure.  This coordinate-free description draws solely on the intrinsic algebraic structure of groupoids and thus it is not only of independent interest but also promises an applicability in subsequent research on groupoids. From a $C^*$-algebraic point of view, they give rise to c.p.c.\ order zero maps from multi-matrix algebras into groupoid $C^*$-algebras. This notion will be a centerpiece in one of the main regularity properties we study in this paper, almost-elementariness (Section~\ref{sec:almost-elementary}). Therefore, in this section, we provide an axiomatic approach to the concept of castles in a general groupoid setting as well as its relation to $C^*$-algebras.

Let us first introduce a precursory notion. 

\begin{defn} \label{defn:plurisection}
	Let $\CG$ be a groupoid. 
	A \emph{transverse system} 
	in $\CG$ is a disjoint collection $\CC$ of nonempty subsets in $\CG$ satisfying
	the following conditions: 
	\begin{enumerate}[label=(\roman*)]
		\item \label{defn:plurisection:unary} for any $C \in \CC$, we have $r(C), s(C), C^{-1} \in \CC$, and
		\item \label{defn:plurisection:binary} for any $C, D \in \CC$, we have $C D \in \CC \cup \{\varnothing\}$. 
	\end{enumerate}
	Given a transverse system $\CC$, we write $\CC^{(0)} := \left\{ C \in \CC \colon C \subseteq \CG^{(0)} \right\}$. 
	
	Given another transverse system $\CC'$, we write $\CC \preceq \CC'$ and say $\CC$ \emph{refines} $\CC'$ if for any $C \in \CC$, there is $C' \in \CC'$ such that $C\subseteq C'$. 
\end{defn}

We establish a few basic properties of transverse systems. 
	
\begin{lem} \label{lem:plurisection-basic}
	Let $\CC$ be a transverse system in $\CG$. 
	\begin{enumerate}
		\item \label{lem:plurisection-basic:levels-source} We have $\CC^{(0)} = \left\{ s(C) \colon C \in \CC \right\} = \left\{ r(C) \colon C \in \CC \right\}$. 
		\item \label{lem:plurisection-basic:inverse} For any $C \in \CC$, we have $s(C) = C^{-1} C$ and $r(C) = C C^{-1}$. 
		\item \label{lem:plurisection-basic:composability} For any $C, D \in \CC$, we have $C D \not= \varnothing$ if and only if $s(C) = r(D)$. 
		\item \label{lem:plurisection-basic:bisection} For any $C \in \CC$, the maps $s$ and $r$ are both injective when restricted to $C$. 
		\item \label{lem:plurisection-basic:subquotient} The union $\bigcup \CC$ is a subgroupoid of $\CG$, the collection $\CC$ also becomes a groupoid with operations induced from those of $\CG$, and the quotient map $\bigcup \CC \to \CC$ is a groupoid homomorphism, i.e., it preserves all the groupoid operations.  
		\item \label{lem:plurisection-basic:reduction} For any $\CD_0 \subseteq \CC\unitSpace$, the collection $\CD := \left\{ C \in \CC \colon r(C), s(C) \in \CD_0 \right\}$ defines a transverse system with $\CD\unitSpace = \CD_0$. 
		\item \label{lem:plurisection-basic:towers} Let $\sim_{\CC}$ be the relation on $\CC$ such that for any $C, D \in \CC$, we have $C \sim_{\CC} D$ if and only if there is $E \in \CC$ such that $s(C) = r(E)$ and $s(E) = r(D)$. Then $\sim_{\CC}$ is an equivalence relation on $\CC$, each of whose equivalence classes is itself a transverse system in $\CG$. 
	\end{enumerate}
\end{lem}

\begin{proof}
	\begin{enumerate}
		\item Since $s$ is the identity map on $\CG^{(0)}$, it is clear that $\left\{ C \in \CC \colon C \subseteq \CG^{(0)} \right\} \subseteq \left\{ s(C) \colon C \in \CC \right\}$. On the other hand, the opposite containment follows from Definition~\ref{defn:plurisection}\ref{defn:plurisection:unary}. The case of $r$ is similar. 
		
		\item Since $C \in \CC$, we have $C \not= \varnothing$. Fixing $x \in C$, we have $s(x) = x^{-1} x \in s(C) \cap C^{-1} C$. By Definition~\ref{defn:plurisection}, the sets $s(C)$ and $C^{-1} C$ are both contained in the disjoint family $\CC$ while having nonempty intersection, whence $s(C) = C^{-1} C$. Similarly we have $r(C) = C C^{-1}$. 
		
		\item Observe that $C D \not= \varnothing$ if and only if $s(C) \cap r(D) \not= \varnothing$. The latter condition is equivalent to $s(C) = r(D)$ since both sets are in the disjoint family $\CC$.

		\item 
		For any $x_1, x_2 \in C$ with $s(x_1) = s(x_2)$, we observe that 
		$x_1 x_2^{-1} \in C C^{-1} = r(C) \subseteq \CG^{(0)}$ by \eqref{lem:plurisection-basic:inverse}, whence $x_1 x_2^{-1} = r \left(x_1 x_2^{-1} \right) = r(x_1) = x_1 x_1^{-1}$ and thus $x_1 = x_2$ by cancellation, which proves $s$ is injective. The case of $r$ is similar. 
		
		\item It is straightforward to verify that $\bigcup \CC$ (respectively, $\CC$) are both groupoids with operations inherited (respectively, induced) from $\CG$, with $\bigcup \CC^{(0)}$ (respectively, $\CC^{(0)}$) as its unit space. The fact that the quotient map is a groupoid homomorphism is then tautological. 
		
		\item This is immediate from the definition. 
		
		\item It is straightforward to verify that $\sim_{\CC}$ is an equivalence relation, by using the facts that $\CC$ is closed under composition and taking inverses. Furthermore, it is clear that for any $C \in \CC$, we have $C \sim_{\CC} r(C)$, $C \sim_{\CC} s(C)$, and $C \sim_{\CC} C^{-1}$. For any $C , D \in \CC$ with $C D \not= \varnothing$, by \eqref{lem:plurisection-basic:composability}, we have $s(C) = r(D)$, which further implies $s \left( C^{-1} \right) = r(C) = r(CD)$, whence $C \sim_{\CC} D \sim_{\CC} CD$. Therefore we have verified the conditions in Definition~\ref{defn:plurisection} for any equivalence class of $\sim_{\CC}$. 
	\end{enumerate}
\end{proof}

Another important related concept is invariant functions on groupoids. 

\begin{defn} \label{defn:invariant-functions}
	Let $\CG$ be a groupoid and let $\CH$ be a subgroupoid of $\CG$. Define the equivalence relation $\sim_{\CH}$ on $\CG$ such that $x \sim_{\CH} y$ if and only if there are $z, w \in \CH \cup \CG^{(0)}$ such that $x = z y w$. 
	
	Let $X$ be a subset of $\CG$ and let $Y$ be a set. We say a map $f \colon X \to Y$ is \emph{$\CH$-invariant} if $f(x_1) = f(x_2)$ whenever $x_1 \sim_{\CH} x_2$. When $X$ is equipped with a locally compact Hausdorff topology, we may write $C_0(X)^{\CH}$ for the set of $\CH$-invariant $C_0$-functions on $X$. Similarly we define $C_{\operatorname{c}}(X)^{\CH}$, $C_{\operatorname{b}}(X)^{\CH}$, etc. 
\end{defn}

\begin{rem}
	Observe that the relation $\sim_{\CH}$ does not depend on the ambient groupoid $\CG$ in the sense that if $\CG'$ is an intermediate groupoid between $\CG$ and $\CH$, then the relation $\sim_{\CH}$ constructed on $\CG$, upon restriction to $\CG' \subseteq \CG$, agrees with the relation $\sim_{\CH}$ constructed on $\CG'$.  
\end{rem}

\begin{lem} \label{lem:invariant-functions-equivalence-relation}
	Let $\CG$ be a groupoid. Then the following relations agree: 
	\begin{enumerate}
		\item \label{lem:invariant-functions-equivalence-relation:original} The relation $\sim_{\CG}$ in Definition~\ref{defn:invariant-functions}, 
		\item \label{lem:invariant-functions-equivalence-relation:source-range} The smallest equivalence relation $\sim'_{\CG}$ such that for any $x \in \CG$, we have $r(x) \sim'_{\CG} x \sim'_{\CG} s(x)$. 
		\item \label{lem:invariant-functions-equivalence-relation:middleman} The relation $\sim''_{\CG}$ satisfying that for any $x, y \in \CG$, we have $x \sim''_{\CG} y$ if and only if there is $z \in \CG$ such that $s(x) = r(z)$ and $s(z) = r(y)$. 
	\end{enumerate}
\end{lem}

\begin{proof}
	To show $\sim_{\CG} {\Rightarrow} \sim'_{\CG}$, we fix arbitrary $x, y \in \CG$ with $x \sim_{\CG} y$, i.e., there are $z, w \in \CH \cup \CG^{(0)}$ such that $x = z y w$, and observe that $x \sim'_{\CG} s(x) = s(z y w) = s(w) \sim'_{\CG} w \sim'_{\CG} r(w) = s(y) \sim'_{\CG} y$. 
	
	To show $\sim'_{\CG} {\Rightarrow} \sim''_{\CG}$, we simply observe that the latter is an equivalence relation and, for any $x \in \CG$, we have $x \sim''_{\CG} s(x)$ since, writing $z = s(x)$, we have $s(x) = r(z)$ and $s(z) = r(y)$, and similarly we have $x \sim''_{\CG} r(x)$. 
	
	To show $\sim''_{\CG} {\Rightarrow} \sim_{\CG}$, we fix arbitrary $x, y \in \CG$ with $x \sim''_{\CG} y$, i.e., there is $z \in \CG$ such that $s(x) = r(z)$ and $s(z) = r(y)$ and observe that $x = (xz) y \left( y^{-1} z^{-1} \right)$. 
\end{proof}

\begin{rem}
	Lemma~\ref{lem:invariant-functions-equivalence-relation} implies that the equivalence relation $\sim_{\CC}$ in Lemma~\ref{lem:plurisection-basic}\eqref{lem:plurisection-basic:towers} agrees with the equivalence relation $\sim_{\CC}$ defined in Definition~\ref{defn:invariant-functions}, when we view $\CC$ as a groupoid as in Lemma~\ref{lem:plurisection-basic}\eqref{lem:plurisection-basic:subquotient}. 
\end{rem}

In the rest of the section, we assume $\CG$ is an \'{e}tale groupoid. Now we give the main definition of the section. 

\begin{defn} \label{defn:castle}
	A \emph{castle} in $\CG$ is a transverse system $\CC$ that is finite and \emph{principal} in the sense that the map $r \times s \colon \CC \to \CC^{(0)} \times \CC^{(0)}$, $C \mapsto (r(C) , s(C))$ is injective. 
	
	The members of $\CC^{(0)}$ are sometimes called the \emph{levels} of the castle $\CC$ or simply \emph{$\CC$-levels}, while
	the members of $\CC$ are sometimes called the \emph{ladders} of the castle $\CC$ or simply \emph{$\CC$-ladders}. 
	
	We say a castle is \emph{open} (respectively, \emph{closed}, \emph{clopen}, \emph{Borel}, \emph{compact} or \emph{precompact}) if each of its members is \emph{open} (respectively, \emph{closed}, \emph{clopen}, \emph{Borel}, \emph{compact} or \emph{precompact}) in $\CG$. 
			
	Let $\sim_{\CC}$ be the equivalence relation in Lemma~\ref{lem:plurisection-basic}\eqref{lem:plurisection-basic:towers}. Then the equivalence classes of $\sim_{\CC}$ are called the \emph{towers} in $\CC$ or simply \emph{$\CC$-towers}. 
	In particular, if $\sim_{\CC}$ has at most one equivalence class, then $\CC$ itself may also be called a \emph{tower}. 
	
	The \emph{height} of a tower $\CT$ in $\CC$ is the natural number $\left| \CT^{(0)} \right|$. The \emph{minimal height} of $\CC$ is the minimal height among all of its towers. 
\end{defn}

\begin{lem}\label{lem:castle-basic}
	Let $\CC$ be a castle. 
	\begin{enumerate}
		\item \label{lem:castle-basic:bisection} If $\CC$ is an open castle, then any $C \in \CC$ is a bisection. 
		\item \label{lem:castle-basic:refine-bisection} If $\CC$ refines an open castle, then any $C \in \CC$ is a bisection. 
		\item \label{lem:castle-basic:subquotient} The collection $\CC$ becomes a finite principal groupoid with operations induced from those of $\CG$, while the union $\bigcup \CC$ becomes a principal subgroupoid of $\CG$. 
		\item \label{lem:castle-basic:reduction} For any $\CD_0 \subseteq \CC\unitSpace$, the collection $\CD := \left\{ C \in \CC \colon r(C), s(C) \in \CD_0 \right\}$ defines a castle with $\CD\unitSpace = \CD_0$. 
		\item \label{lem:castle-basic:towers} Each tower in $\CC$ is also itself a castle. 
		\item \label{lem:castle-basic:height} For any tower $\CT$ in $\CC$ and any $u \in \bigcup \CT$, we have 
		\[
			\left| \CT\unitSpace \right| = \left| r \left( \left( \bigcup \CT \right) u \right) \right|  = \left| \left( \bigcup \CT \right) u \right|  = \left| \left( \bigcup \CC \right) u \right| \; .
		\]
	\end{enumerate}
\end{lem}

\begin{proof}
	The first five statements follow trivially from Lemma~\ref{lem:plurisection-basic}, where the principality condition in \eqref{lem:castle-basic:subquotient} follows from the principality requirement in Definition~\ref{defn:castle}. 
	For the last statement, we assume that $u \in C_0$ where $C_0 \in \CT$. 
	For any $C \in \CT\unitSpace$, by \eqref{lem:castle-basic:subquotient}, there is a unique $D \in \CT$ such that $C = r(D)$ and $s(D) = r(C_0)$, or equivalently, $C \cap r(D) \not= \varnothing$ and $s(D) \cap r(C_0) \not= \varnothing$, whence $\left| C \cap r \left( \left( \bigcup \CT \right) u \right) \right| = \left| r \left( D u \right) \right| = 1$ by Lemma~\ref{lem:plurisection-basic}\eqref{lem:plurisection-basic:bisection}. Since $r \left( \left( \bigcup \CT \right) u \right) \subseteq \bigcup \CT\unitSpace$, this shows the first equality $\left| \CT\unitSpace \right| = \left| r \left( \left( \bigcup \CT \right) u \right) \right|$. 
	The second equality follows from the principality of $\bigcup \CT$, a consequence of \eqref{lem:castle-basic:towers} and \eqref{lem:castle-basic:subquotient}. 
	The last equality follows from the fact that 
	for any $C \in \CC \setminus \CT$, we have $Cu = \varnothing$ by the definition of $\sim_\CG$.  
\end{proof}

\begin{rem}\label{rem: elementary groupoid}
	When $\CC$ is an open or closed castle, the subgroupoid $\bigcup \CC$ given in Lemma~\ref{lem:plurisection-basic}\eqref{lem:plurisection-basic:subquotient} is \emph{elementary} in the sense of \cite[III.1.1]{Renault1980groupoid}, and moreover, $\CC$ is clearly determined by the quotient map $\pi \colon \bigcup \CC \twoheadrightarrow \CC$ as the collection of fibers. However, it is not true that every finite principal quotient of an open (or closed) elementary subgroupoid of $\CG$ gives rise to a castle in this way. For example, consider the finite principle (and thus elementary) groupoid $\CG = \CG_0 \sqcup \CG_1$ where $\CG_0 = \{0\}$ is the trivial singleton groupoid and $\CG_1$ is the pair groupoid on $\{1,2\}$. Note that $\CG_1$ is a quotient of $\CG$ by mapping $\CG_0$ to $(1,1)$ and taking the identity on $\CG_1$, but the fibers of this quotient map $\pi$ do not constitute a castle, for $r \left( \pi^{-1} (1,2) \right) = r (\{1,2\}) = \{(1,1)\} \subsetneq \{ 0 , (1,1) \} = \pi^{-1} (1,1)$. 
\end{rem}

Now we present a more concrete model for castles and relate them to the following notion  introduced by Nekrashevych in \cite[Definition~3.1]{Nekrashevych2019Simple}: 
\begin{defn}
	\label{defn:multisection} A finite set of open
    bisections $\CT = \left\{C_{i,j} \colon i,j\in F \right\}$
	with a finite index set $F$ is called an \emph{open multisection}\footnote{We point out that in \cite[Definition~3.1]{Nekrashevych2019Simple},
		the author worked with ample groupoids and assumed multisections are
		clopen; we do not make such an assumption. } 
	if it
	satisfies 
	\begin{enumerate}
		\item $C_{i,j}C_{j,k}=C_{i,k}$ for $i,j,k\in F$; 
		\item $\left\{C_{i,i}:i\in F \right\}$ is a disjoint family of subsets of $\CG^{(0)}$. 
	\end{enumerate}
\end{defn}

\begin{rem}
	\label{rem:multisection-basics}
	Given an open multisection $\CT = \left\{C_{i,j} \colon i,j\in F \right\}$
	with a finite index set $F$, we have the following identities for any $i,j,k,l \in F$: 
	\begin{enumerate}
		\item $s(C_{i,j}) = C_{j,j}$, since we get $s(C_{i,j}) \subseteq C_{j,j}$ from $C_{i,j} = C_{i,j} C_{j,j}$, and $s(C_{i,j}) \supseteq C_{j,j}$ from $C_{j,j} = C_{j,i} C_{i,j}$; 
		\item $r(C_{i,j}) = C_{i,i}$, proved similarly; 
		\item $C_{i,j}^{-1} = C_{j,i}$, since we get $C_{j,i} = C_{j,i} C_{i,i} \subseteq C_{j,i} C_{i,j} C_{i,j}^{-1} = C_{jj} C_{i,j}^{-1} = C_{i,j}^{-1}$ and similarly $C_{i,j} \subseteq C_{j,i}^{-1}$; 
		\item $C_{i,j} C_{kl} = 
		\begin{cases}
		C_{i,l}  \in \CT, & \text{if } j = k \\
		C_{i,j} C_{j,j} C_{k,k}C_{k,l} = \varnothing , & \text{if } j \not= k
		\end{cases}$.  
	\end{enumerate}
\end{rem}

\begin{prop}
	\label{prop:tower-multisection}
	If $\CT$ is an open multisection with an index set $F$ and $\CT \not= \{\varnothing\}$, then $\CT$ is an open tower. Conversely, every open tower $\CC$ is an open multisection indexed by $\CC^{(0)}$. 
\end{prop}

\begin{proof}
	Assume $\CT=\{C_{i,j} \colon i,j\in F\}$ is an open multisection. 
	Observe that $\CT \not= \{\varnothing\}$ implies $\varnothing \notin \CT$, since if $C_{i_0,j_0}$ were empty for some $i_0, j_0 \in F$, then by Remark~\ref{rem:multisection-basics}, for any $i, j \in F$, we would have $s \left( C_{i, j_0} \right) = s \left( C_{i_0, j_0} \right) = \varnothing$, which would imply $C_{i, j_0} = \varnothing$ and similarly $C_{i, j} = \varnothing$, a contradiction. 
	It then follows from Remark~\ref{rem:multisection-basics} that $\CT$ is a finite open transverse system. 
	It is also clear that the map $r \times s \colon \CT \to \CT^{(0)} \times \CT^{(0)}$ is injective as $(r \times s) (C_{i,j}) = \left(r(C_{i,j}), s(C_{i,j})\right) = \left( C_{i,i} , C_{j,j} \right)$, which implies that $\CT$ is an open castle.  
	Moreover, observe that for any $i,j,k,l \in F$, we have $s(C_{i,j}) = C_{j,j} = r(C_{j,k})$ and $s(C_{j,k}) = C_{k,k} = r(C_{k,l})$, whence $C_{i,j} \sim_{\CT} C_{k,l}$. 
	Therefore $\CT$ is an open tower. 
	
	In the other direction, given an open tower $\CC$, we observe that the map $r \times s \colon \CC \to \CC^{(0)} \times \CC^{(0)}$ is not only injective, but also surjective, since for any $C, C' \in \CC^{(0)}$, it follows from $C \sim_{\CC} C'$ that there is $C'' \in \CT$ such that $(r \times s) (C'') = (C, C')$. For $i , j \in \CC^{(0)}$, define $C_{i,j} = (r \times s)^{-1} (i, j)$, which is a bisection by Lemma~\ref{lem:castle-basic}\eqref{lem:castle-basic:bisection}. It is straightforward to verify that $\CT = \left\{C_{i,j} \colon i,j\in \CC^{(0)} \right\}$ is an open multisection. 
\end{proof}

Hence an open tower may be viewed as a ``coordinate-free'' presentation of an open multisection. Such an approach, albeit more abstract initially, will help us greatly simplify notations later on. We get a greater simplification when we deal with castles consisting of multiple towers, as the ``coordinate-dependent'' presentation will need another layer of indices to enumerate the towers. 

\begin{cor} \label{cor:castle-multisections}
	A collection $\CC$ of nonempty open bisections in $\CG$ is an open castle if and only if there are finite index sets $I$ and $F_l$ for $l \in I$ and an enumeration $\CC= \left\{C_{i,j}^{l} \colon i,j\in F_{l}, l\in I \right\}$,
	satisfying the following conditions: 
	\begin{enumerate}[label=(\roman*)]
		\item for each $l \in I$, the set  $\CC_l := \left\{C_{i,j}^{l} \colon i,j\in F_{l} \right\}$ is a multisection indexed by $F_{l}$; 
		\item $C_{i,i}^{l} \cap C_{i',i'}^{l'}=\emptyset$ for any distinct $l, l' \in I$, any $i \in F_l$ and any $i' \in F_{l'}$. 
	\end{enumerate}
	Moreover, in this case, the collections $\{C_{i,j}^{l}:i,j\in F_{l}\}$, for $l \in I$, are exactly the towers in $\CC$. 
\end{cor}

\begin{proof}
	The ``only if'' direction follows from applying Proposition~\ref{prop:tower-multisection} to the towers in $\CC$. The ``if'' direction is also straightforward by Proposition~\ref{prop:tower-multisection} upon noticing that $C_{i,j}^{l} C_{i',j'}^{l'}=\emptyset$ for any distinct $l, l' \in I$, any $i,j \in F_l$ and any $j',j' \in F_{l'}$. 
\end{proof}

\begin{rem} \label{rem:castle-matrices}
	Let $\CC$ be an open castle in $\CG$ and let $\CC = \left\{C_{i,j}^{l} \colon i,j\in F_{l},l\in I \right\}$ be an enumeration as in Corollary~\ref{cor:castle-multisections}. Then, viewing $\CC$ as a finite groupoid as in Lemma~\ref{lem:castle-basic} and identifying its groupoid $C^*$-algebra with its groupoid algebra $\mathbb{C} \CC = \left\{ \displaystyle \sum_{l\in I} \sum_{i,j\in F_{l}} \lambda_{i,j}^{l} C_{i,j}^{l} \colon  \lambda_{i,j}^{l}  \in \mathbb{C} \text{ for } i,j \in F_{l} \text{ and } l\in I \right\}$, we have, thanks to Definition~\ref{defn:multisection} and Remark~\ref{rem:multisection-basics}, a $*$-isomorphism
	\[
		C^*(\CC) \xrightarrow{\cong} \bigoplus_{l \in I} M_{|F_l|} (\mathbb{C}) \; , \quad C_{i,j}^{l} \mapsto e_{i,j}^{l} \; ,
	\]
	where $e_{i,j}^{l}$ stands for the matrix unit at the $(i,j)$-th position in the $l$-th matrix block. 
\end{rem}


\section{More on castles}\label{sec:castles-2}

In this section, we discuss a few constructions and operations involving castles. They will be of great use to us in the rest of the paper. 

We start with a procedure that allows us to split a castle into a ``finer'' one. 

\begin{lem} \label{lem:castle-split}
	Let $\CC$ be a castle in $\CG$ and let $\CP := \left\{ P_1, \ldots , P_n \right\}$ be a partition of $\bigcup \CC$. Consider the collection 
	\begin{align*}
		\CC_\CP &:= \left\{ Q_{C, \beta} \colon C \in \CC, \, \beta \in \{1,\ldots,n\}^{\Xi(C)} \right\} \setminus \{\varnothing\} \; 
	\end{align*}
	of nonempty subsets of $\CG$, where 
	\[
		\Xi(C) := \left\{ (D,E) \in \CC \colon DCE \not=\varnothing \right\} 
		\quad \text{ and } \quad 
		Q_{C, \beta} := 
		\bigcap_{ (D,E) \in \Xi(C) } 
		D^{-1} \, P_{\beta(D,E)} \, E^{-1} \; .
	\]
	Then the following hold:
	\begin{enumerate}
		\item \label{lem:castle-split:exists} 
		The family $\CC_\CP$ is a castle that refines both $\CC$ and $\CP$ and satisfies $\bigcup \CC_\CP = \bigcup \CC$. 
		
		\item \label{lem:castle-split:refine-basic}
		If $\CD$ is a castle that refines $\CC$ and satisfies $\bigcup \CD = \bigcup \CC$, then for any $D \in \CD$ and $C, E \in \CC$, we have $CDE \in \CD \cup \{\varnothing\}$.  
		
		\item \label{lem:castle-split:universal} 
		Any castle $\CD$ refining both $\CC$ and $\CP$ and satisfying $\bigcup \CD = \bigcup \CC$ refines $\CC_\CP$. 
		
		\item \label{lem:castle-split:extra} 
		If $\CC$ is open (respectively, compact) and $P_1, \ldots, P_n$ are open (respectively, compact), then so is $\CC_\CP$. 
	\end{enumerate}
	
\end{lem}

\begin{proof}
	To prove \eqref{lem:castle-split:exists}, 
	we first observe that for any $C \in \CC$ and $\beta \in \{1,\ldots,n\}^{\Xi(C)}$, we have \[Q_{C,\beta} \subseteq r(C) P_{\beta(r(C),s(C))} s(C) \subseteq P_{\beta(r(C),s(C))} \cap \left(r(C) \left( \bigcup \CC \right) s(C)  \right) = P_{\beta(r(C),s(C))} \cap C.\]
    This shows $\CC_\CP$ refines both $\CC$ and $\CP$, which in turn also implies $\bigcup \CC_\CP \subseteq \bigcup \CC$. 
	
	To see $\bigcup \CC_\CP \supseteq \bigcup \CC$, we fix an arbitrary $x \in \bigcup \CC$ and let $C \in \CC$ be such that $x \in C$. We define $\beta_x \in \{1,\ldots,n\}^{\Xi(C)}$ as follows: for any $(D,E) \in \Xi(C)$, there are $y \in D$ and $z \in E$ such that $s(y) = r(x)$ and $r(x) = s(z)$, which allows us to define $\beta_x(D,E)$ so that $yxz \in P_{\beta(D,E)}$, which implies $x \in D^{-1} \, P_{\beta_x(D,E)} \, E^{-1}$. 
	This implies $x \in Q_{C, \beta_x}$, as desired. 
	
	To see $\CC_\CP$ is a castle, we first prove it is a disjoint family. Indeed, for any $C, C' \in \CC$, any $\beta \in \{1,\ldots,n\}^{\Xi(C)}$ and any $\beta' \in \{1,\ldots,n\}^{\Xi(C')}$, since we have shown $Q_{C, \beta} \subseteq C$ and $Q_{C', \beta'} \subseteq C'$, we have $Q_{C, \beta} \cap Q_{C', \beta'} = \varnothing$ if $C \not= C'$; on the other hand, if $C = C'$ but $\beta \not= \beta'$, then there is $(D,E) \in \Xi(C)$ such that $\beta(D,E) \not= \beta'(D,E)$, which implies $Q_{C, \beta} \cap Q_{C', \beta'} \subseteq \left( D^{-1} \, P_{\beta(D,E)} \, E^{-1} \right) \cap \left( D^{-1} \, P_{\beta'(D,E)} \, E^{-1} \right) = \varnothing$. This shows $\CC_\CP$ is disjoint. 
	
	Now observe that for any $C, C' \in \CC$ with $s(C) = r(C')$, there are bijections 
	\begin{align*}
		\rho_{C, C'} &\colon \Xi(C) \to \Xi(CC') \, ,  & (D,E) \mapsto \left( D, (C')^{-1} E \right) \; , \\
		\lambda_{C, C'} &\colon \Xi(C') \to \Xi(CC') \, ,  & (D,E) \mapsto \left( DC^{-1} , E \right) \; ,	
	\end{align*}
	satisfying $\rho_{C, s(C)} = \operatorname{id}_{\Xi(C)} = \lambda_{r(C), C}$, $\rho_{C, C'} \circ \rho_{CC', C''} = \rho_{C, C'C''}$ and $\rho_{C', C''} \circ \rho_{C, C'C''} = \rho_{CC', C''}$ for any $C'' \in \CC$ with $s(C') = r(C'')$,  
	and it is easy to verify that $Q_{C,\beta \circ \rho_{C,C'}} C' = Q_{CC', \beta} = C Q_{C',\beta \circ \lambda_{C,C'}}$ for any $\beta \in \{1,\ldots,n\}^{\Xi(CC')}$. 
	For any $C \in \CC$ and any $\beta \in \{1,\ldots,n\}^{\Xi(C)}$, since $Q_{C, \beta} \subseteq C$ and $r(C) = CC^{-1}$, we have $r \left( Q_{C, \beta} \right) = Q_{C, \beta} C^{-1} = Q_{r(C), \beta \circ \rho_{r(C), C}}$ and similarly $s \left( Q_{C, \beta} \right) = Q_{s(C), \beta \circ \lambda_{C, s(C)}}$, 
	while since $C^{-1} = C^{-1} C C^{-1}$, we have $Q_{C, \beta}^{-1} = C^{-1} Q_{C, \beta} C^{-1} = C^{-1} Q_{r(C), \beta \circ \rho_{r(C), C}} = Q_{C^{-1} , \beta \circ \rho_{r(C), C} \circ \lambda_{C, C^{-1}}}$. 
	For any $C, C' \in \CC$, any $\beta \in \{1,\ldots,n\}^{\Xi(C)}$ and any $\beta' \in \{1,\ldots,n\}^{\Xi(C')}$, we have 
	\begin{align*}
		& Q_{C, \beta} Q_{C', \beta'} \\
		=&\ C Q_{s(C), \beta \circ \lambda_{C, s(C)}} Q_{r(C'), \beta' \circ \rho_{r(C'), C'}} C' \\
		=&\ 
		\begin{cases}
			Q_{C, \beta} C' = Q_{CC', \beta \circ \lambda_{CC', (C')^{-1}}} \, , & \text{if } s(C) = r(C') \text{ and } \beta \circ \lambda_{C, s(C)} = \beta' \circ \rho_{r(C'), C'} \\
			\varnothing \, , & \text{otherwise} 
		\end{cases}
		\; .
	\end{align*}
	This shows $\CC_\CP$ is a transverse system. It is clearly finite, while its principality follows from that of $\CC$. Therefore $\CC_\CP$ is a castle. 
	
	To prove \eqref{lem:castle-split:refine-basic}, we fix $\CD$, $D$, $C$, and $E$ as in the statement and first show that $CD \in \CD$. 
	To this end, let $F \in \CC$ be such that $D \subseteq F$. 
	It suffices to consider the case where $C D \not= \varnothing$. 
	Now we have $\varnothing \not= C D \subseteq C F \subseteq \bigcup \CC = \bigcup \CD$, whence there is $D' \in \CD$ such that $C D \cap D' \not= \varnothing$, which implies $C F \cap D' \not= \varnothing$ and thus $D' \subseteq C F$ as $\CD$ refines the disjoint family $\CC$. But on the other hand, since $s(D') \cap s(D) \supseteq s(D') \cap s(C D) \not= \varnothing$, we have $s(D') = s(D) \in \CD$ and thus $s(CD) = s(D')$, which implies $C D = D' \in \CD$ by applying Lemma~\ref{lem:plurisection-basic}\eqref{lem:plurisection-basic:bisection} to $C F$. Similarly we may further deduce $C D F \in \CD$, as desired. 
	
	To prove \eqref{lem:castle-split:universal}, we fix a castle $\CD$ that refines both $\CC$ and $\CP$ and satisfies $\bigcup \CD = \bigcup \CC$. 
	Also fix an arbitrary $D \in \CD$ and let $C \in \CC$ be such that $D \subseteq C$. 
	Now we define $\beta_{D} \in \{1,\ldots,n\}^{\Xi(C)}$ as follows: 
	for any $(E,F) \in \Xi(C)$, since $EDF \in \CD$ by \eqref{lem:castle-split:refine-basic} and $\CD$ refines $\CP$, we may define $\beta_{D} (E,F) \in \{1, \ldots, n\}$ so that $EDF \subseteq P_{\beta_{D} (E,F)}$, which implies that $D = E^{-1} E D F F^{-1} \subseteq E^{-1} P_{\beta_{D} (E,F)} F^{-1}$. Consequently we have $D \subseteq Q_{C, \beta_D} \in \CC_\CP$, as desired. 
	
	Statement~\eqref{lem:castle-split:extra} is obvious from the definition of $\CC_\CP$. 
\end{proof}


We proceed to discuss how castles interact with measures and functions on $\CG$. 

\begin{rem} \label{rem:castle-measure}
	Let $\mu \in M(\CG)$ be an invariant regular Borel probability measure. Let $\CC$ be a Borel castle in $\CG$. Then the following are immediate from the definitions: 
	\begin{enumerate}
		\item \label{rem:castle-measure:tower} For any $C$ and $C'$ in the same tower in $\CC$, we have $\mu(C) = \mu(C')$. 
		\item \label{rem:castle-measure:castle} Writing $\CC / \sim_{\CC}$ for the collection of towers in $\CC$ and fixing an arbitrary set $\left\{ C_{\CT} \in \CT\unitSpace \colon \CT \in \CC / \sim_{\CC} \right\}$ of representatives for the relation $\sim_{\CC}$ restricted to $\CC\unitSpace$, we have 
		\[
			\mu \left( \bigcup \CC\unitSpace \right) = \sum_{ \CT \in \CC / \sim_{\CC} }  \left| \CT\unitSpace \right| \mu \left( C_{\CT} \right) \; .
		\]
	\end{enumerate}
\end{rem}

We then study invariant functions on castles, in the sense of Definition~\ref{defn:invariant-functions}. They will be used for several purposes. Let us first discuss the existence of these functions.

\begin{rem}\label{rem:castle-fundamental-domain}
	Let $\CC$ be an open castle and view $\bigcup \CC$ as an open subgroupoid of $\CG$. 
	Then we have the following isomorphisms: 
	\begin{enumerate}
		\item \label{rem:castle-fundamental-domain:castle-tower}
		Writing $\CC / \sim_{\CC}$ for the collection of towers in $\CC$, we have
		\[
			C_0 \left(\bigcup \CC \right)^{\bigcup \CC} 
			\cong \bigoplus_{ \CT \in \CC / \sim_{\CC} }  C_0\left(\bigcup \CT \right)^{\bigcup \CT} \; ,
		\]
		which is obtained from the restriction of the decomposition $\displaystyle C_0 \left(\bigcup \CC \right)	\cong \bigoplus_{ \CT \in \CC / \sim_{\CC} }  C_0\left(\bigcup \CT \right)$ that is induced from the relatively clopen partition $\displaystyle \bigcup \CC = \bigsqcup_{ \CT \in \CC / \sim_{\CC} }  \left(\bigcup \CT\right)$. Indeed, this follows from the observation that 
		the equivalence classes of $\sim_{\bigcup \CC}$ are precisely the equivalence classes of $\sim_{\bigcup \CT}$, as $\CT$ ranges over all towers in $\CC$. 
		
		\item \label{rem:castle-fundamental-domain:tower-ladder}  
		For any tower $\CT$ in $\CC$ and any $C \in \CT$, we have an isomorphism
		\[
			C_0(C) \cong C_0\left(\bigcup \CT \right)^{\bigcup \CT}  \; ,
		\]
		which is induced from the restriction map $C_0(\bigcup \CT) \twoheadrightarrow C_0(C)$, since each member of the relatively clopen partition $\bigcup \CT = \bigsqcup_{C' \in \CT} C'$ constitutes a set of representatives for $\sim_{\bigcup \CT}$, as is evident from Proposition~\ref{prop:tower-multisection}. 
		
		\item \label{rem:castle-fundamental-domain:castle-ladder} 
		Consequently, for any set $\left\{ C_{\CT} \in \CT \colon \CT \in \CC / \sim_{\CC} \right\}$ of representatives of $\sim_{\CC}$, we have an isomorphism
		\[
			C_0 \left(\bigcup \CC \right)^{\bigcup \CC} 
			\cong \bigoplus_{ \CT \in \CC / \sim_{\CC} } C_0 \left( C_{\CT} \right) \; .
		\]
	\end{enumerate}
\end{rem}

Let us present a useful method to produce invariant functions. 

\begin{lem} \label{lem:castle-invarint-function-smallest}
	Let $\CC$ be an open castle and fix $f \in C_0 \left( \bigcup \CC , \mathbb{R} \right)$. Consider the function 
	$f_{\uparrow\CC} \colon \bigcup \CC \to \mathbb{R}$ such that 
	\[
		f_{\uparrow\CC}(x) := \sup \left\{ f(y) \colon y \in \left( \bigcup \CC \right) x \left( \bigcup \CC \right) \right\} \quad \text{ for any } x \in \bigcup \CC \; . 
	\]
	Then $f_{\uparrow\CC}$ is the smallest element in $C_0 \left( \bigcup \CC  , \mathbb{R} \right)^{\bigcup \CC}$ such that $f \leq f_{\uparrow\CC}$. 
	It also satisfies $\displaystyle \sup_{x \in \bigcup \CC} f_{\uparrow\CC} (x) = \sup_{x \in \bigcup \CC} f (x)$.

	Moreover, 
	the assignment $f \mapsto f_{\uparrow\CC}$ is order-preserving, $1$-Lipschitz in the uniform norm, and equal to the identity map when restricted to $C_0 \left( \bigcup \CC  , \mathbb{R} \right)^{\bigcup \CC}$. 
\end{lem}

\begin{proof}
	It is immediate from the definition that 
	$f_{\uparrow\CC} \geq f$, 
	$\displaystyle \sup_{x \in \bigcup \CC} f_{\uparrow\CC} (x) = \sup_{x \in \bigcup \CC} f (x)$, 
	and
	$f_{\uparrow\CC}$ is $\bigcup \CC$-invariant, 
	while for any $g \in C_0 \left( \bigcup \CC , \mathbb{R}  \right)$, we also have  
	$f_{\uparrow\CC} \leq g_{\uparrow\CC}$ if $f \leq g$, 
	$g_{\uparrow\CC} = g $ if $g$ is $\bigcup \CC$-invariant, 
	and 
	$\left\| f_{\uparrow\CC} - g_{\uparrow\CC} \right\| \leq \|f - g\|$. 
	By Lemma~\ref{lem:castle-basic}\eqref{lem:castle-basic:bisection}, for any $C, C' \in \CC$ and any $x \in \CG$, we have $|CxC'| \leq 1$. This allows us to write 
	\[
		f_{\uparrow\CC}(x) = \max_{C,C' \in \CC} \max  \left\{  f(y) \colon y \in C x C' \right\}  \quad \text{ for any } x \in \bigcup \CC \; , 
	\]
	which shows $f_{\uparrow\CC}$ is continuous. To see $f_{\uparrow\CC}$ is $C_0$, we fix an arbitrary $\varepsilon > 0$ and apply a simple functional calculus (e.g., using the function $t \mapsto t + \frac{|t- \varepsilon| - |t + \varepsilon|}{2}$) to obtain  $g  \in C_{\operatorname{c}} \left( \bigcup \CC  , \mathbb{R} \right)$ such that $\|f - g\| \leq \varepsilon$. Let $K$ be a compact subset of $\bigcup \CC$ that contains the support of $g$. Then it is evident from the definition that $g_{\uparrow\CC}$ is supported in $\left( \bigcup \CC \right) K \left( \bigcup \CC \right)$, or equivalently, in the union 
	\[
		\bigcup \left\{ C (K \cap C') C'' \colon C,C',C'' \in \CC \text{ satisfying } s(C) = r(C') \text{ and } s(C') = r(C'') \right\} \; ,
	\]
	which is compact since for any $C, C',C'' \in \CC$ as above, $K \cap C'$ is relatively clopen in $K$ and thus compact, and $C$ and $C'$ are bisections. Since $\left\| f_{\uparrow\CC} - g_{\uparrow\CC} \right\| \leq \|f - g\| \leq \varepsilon$ and $\varepsilon$ was chosen arbitrarily, we see that $f_{\uparrow\CC} \in C_0 \left( \bigcup \CC  , \mathbb{R} \right)$. Finally, for any $g \in C_0 \left( \bigcup \CC , \mathbb{R}  \right)^{\bigcup \CC}$ satisfying $f \leq g$, we have $f_{\uparrow\CC} \leq  g_{\uparrow\CC} = g$ by what we have proved, which shows $f_{\uparrow\CC}$ is the smallest element in $C_0 \left( \bigcup \CC  , \mathbb{R} \right)^{\bigcup \CC}$ such that $f \leq f_{\uparrow\CC}$. 
\end{proof}

We record one useful feature of invariant functions on a castle. Later we will use this to ``trim'' castles in the horizontal direction. 

\begin{lem} \label{lem:castle-trim}
	Let $\CC$ be an open castle. Then for any $f \in C_0 \left(\bigcup \CC \right)^{\bigcup \CC}$ and any $Y \subseteq \mathbb{C}$, the collection 
	$\left\{ f^{-1} (Y) \right\} \vee \CC := \left\{ f^{-1} (Y) \cap C \colon C \in \CC \right\} \setminus \{ \varnothing \}$ is again a castle. 
\end{lem}

\begin{proof}
	Applying Lemma~\ref{lem:invariant-functions-equivalence-relation}\eqref{lem:invariant-functions-equivalence-relation:source-range} to $\bigcup \CC$, we see that $f(x) = f(s(x)) = f(r(x)) = f \left( x^{-1} \right) = f(xy)$ for any composable $x, y \in \bigcup \CC$. 
	It follows that for any $C, D \in \CC$, we have 
	$s \left( f^{-1} (Y) \cap C  \right) = f^{-1} (Y) \cap s(C)$, 
	$r \left( f^{-1} (Y) \cap C  \right) = f^{-1} (Y) \cap r(C)$, 
	$\left( f^{-1} (Y) \cap C  \right)^{-1} = f^{-1} (Y) \cap C^{-1}$, and 
	$\left( f^{-1} (Y) \cap C  \right) \left( f^{-1} (Y) \cap D  \right) = f^{-1} (Y) \cap CD$; hence these sets all fall back into $\left\{ f^{-1} (Y) \right\} \vee \CC$. Hence $\left\{ f^{-1} (Y) \right\} \vee \CC$ is a finite transverse system. Its principality follows directly from that of $\CC$. 
\end{proof}

In preparation for our discussion of regulartiy properties of groupoid $C^*$-algebras, we discuss how an open castle in $\CG$ gives rise to order zero maps from the above multi-matrix algebra into $C^*(\CG)$. 
Recall from \cite[Definition~2.3]{WinterZacharias2009Completely} that a c.p.c.\ map $\psi:A\to B$ between two $C^*$-algebras is
said to be \textit{order zero} if for any positive elements $a,b\in A$ with $ab=0$,
we have $\psi(a)\psi(b)=0$ as well. 
Since it is easy to see that a c.p.c.\ order zero map from a finite abelian $C^*$-algebra $\C^n$ to an abelian $C^*$-algebra $C_0(X)$ amounts to a collection of $n$ ``bump functions'' with disjoint (open) supports, 
the idea of our construction is to make use of ``invariant bump functions'' with disjoint (open) supports. 
In the following result, we again view $\CC$ as a finite groupoid as in Lemma~\ref{lem:castle-basic} and identify its groupoid $C^*$-algebra with its groupoid algebra, with $\CC$ serving as a linear basis.  

\begin{prop} \label{prop:castle-order-zero}
	Let $\CC$ be an open castle in $\CG$ and let $f$ be a positive contraction in $C_0 \left(\bigcup \CC \right) \cap C_{\operatorname{c}}(\CG)$ that is $\left(\bigcup \CC\right)$-invariant. We write $f_{|C} = f \cdot \chi_{C} \in C_0(C) \cap C_{\operatorname{c}}(\CG)$ for any $C \in \CC$. Then the linear map 
	\[
		\psi_f \colon C^*(\CC) \to  C^*(\CG) \quad \text{ determined by } \quad  \CC \ni C \mapsto f_{|C} \in C_{\operatorname{c}} (\CG) \subseteq C^*(\CG)  
	\]
	is a c.p.c.\ order zero map that takes $C\left( \CC^{(0)} \right) \subseteq C^*(\CC)$ into $C_0 \left( \CG^{(0)}\right) \subseteq C^*(\CG)$. 
\end{prop}

\begin{proof}
	As the last assertion is evident, we focus on showing $\psi_f$ is a c.p.c.\ order zero map. 
	Since $C^*(\CC)$ is finite-dimensional with $\CC$ serving as a basis, every element in $C_0((0,1]) \otimes C^*(\CC)$ can be written uniquely as $\sum_{C \in \CC} g_C \otimes C$, where $g_C \in C_0((0,1])$ for $C \in \CC$. Consider the map 
	\[
		\widetilde{\psi}_f \colon C_0((0,1]) \otimes C^*(\CC) \to  C^*(\CG) \; , \quad  \sum_{C \in \CC} g_C \otimes C \mapsto \sum_{C \in \CC} (g_C \circ f)_{|C} \in C_{\operatorname{c}} (\CG) \subseteq C^*(\CG)  \; ,
	\]
	where $(g_C \circ f)_{|C} := (g_C \circ f) \cdot \chi_{C} \in C_0(C) \cap C_{\operatorname{c}}(\CG)$. Since clearly $\psi_f = \widetilde{\psi}_f \circ  \iota $, where $\iota \colon C^*(\CC) \to C_0((0,1]) \otimes C^*(\CC)$ maps any $a \in C^*(\CC)$ to $\operatorname{id}_{(0,1]} \otimes a$, by \cite[Corollary~4.1]{WinterZacharias2009Completely}, it suffices to show that $\widetilde{\psi}_f$ is a $*$-homomorphism.

	To do this, we first observe that $\widetilde{\psi}_f$ is linear. 
	It remains to verify that 
	for any $C, D \in \CC$ and any $g, h \in C_0((0,1])$, 
	we have 
	\begin{equation} \tag{\ensuremath{*}} \label{prop:castle-order-zero:proof:eq}
		\widetilde{\psi}_f \left( g \otimes C \right)  \widetilde{\psi}_f \left( h \otimes D \right) 
		= 
		\widetilde{\psi}_f \left( (g \otimes C)  (h \otimes D) \right)
		\text{ and } 
		\left(\widetilde{\psi}_f \left( g \otimes C \right)\right)^* 
		= 
		\widetilde{\psi}_f \left( (g \otimes C) ^* \right).
	\end{equation}	
	To this end, we fix an arbitrary $x \in \CG$ and compute, using Lemma~\ref{lem:castle-basic}\eqref{lem:castle-basic:bisection}, that
	\begin{align*}	
		& \left( (g \circ f)_{|C} \ast  (h \circ f)_{|D} \right) (x) \\
		= & \sum_{\left\{ (y,z) \in \CG^{(2)} \colon x = yz \right\}} (g \circ f)_{|C} (y) \cdot  (h \circ f)_{|D} (z) \\
		= & \begin{cases}
			0 \, , & \text{if } x \notin CD \\
			(g \circ f) \left( \left(r|_C\right)^{-1}(r(x)) \right) \cdot  (h \circ f)  \left( \left(s|_D\right)^{-1}(s(x)) \right) \, , & \text{if } x \in CD 
		\end{cases}
		\; .
	\end{align*}
	Since $f$ is $\left( \bigcup \CC \right)$-invariant, by Lemma~\ref{lem:invariant-functions-equivalence-relation}\eqref{lem:invariant-functions-equivalence-relation:source-range}, the expression in the second case can be simplified to $(g \circ f) \left( x \right) \cdot  (h \circ f)  \left( x \right) = ((g h ) \circ f) \left( x \right)$. On the other hand, since 
	\[
		(g \otimes C) \cdot (h \otimes D) = 
		\begin{cases}
			0 \, ,  & \text{if } CD = \varnothing \\
			(gh) \otimes (CD) \, , & \text{if } CD \in \CC
		\end{cases}
		\; ,
	\]
	we have 
	\begin{align*}
		\widetilde{\psi}_f \left( (g \otimes C) \cdot (h \otimes D) \right) (x) 
		= & 
		\begin{cases}
		0 \, ,  & \text{if } CD = \varnothing \\
		((gh) \circ f)_{|CD}  (x) \, , & \text{if } CD \in \CC
		\end{cases} \\
		= & 
		\begin{cases}
		0 \, ,  & \text{if } x \notin CD \\
		((gh) \circ f) (x) \, , & \text{if } x \in CD 
		\end{cases} 
	\end{align*}
	This verifies the first part of \eqref{prop:castle-order-zero:proof:eq}. For the second part, since $f(x) = f(s(x)) = f\left( r\left(x^{-1}\right) \right) = f\left(x^{-1}\right)$ by Lemma~\ref{lem:invariant-functions-equivalence-relation}\eqref{lem:invariant-functions-equivalence-relation:source-range}, we thus have
	\begin{align*}
		\widetilde{\psi}_f \left( (g \otimes C) ^* \right) (x) & = \widetilde{\psi}_f \left( \overline{g} \otimes C^{-1} \right) (x) = (\overline{g} \circ f)_{|C^{-1}} (x) = \overline{(g \circ f)} (x) \cdot  \chi_{C^{-1}} (x)   \\
		& = \overline{(g \circ f) \left(x^{-1}\right) \cdot \chi_{C} \left(x^{-1}\right) } =  \overline{ ({g} \circ f)_{|C}  \left(x^{-1}\right) } = \left(\widetilde{\psi}_f \left( g \otimes C \right)\right)^*  (x) \; ,
	\end{align*}
	as desired. 
\end{proof}

In the following sections, it is often necessary to ``horizontally shrink'' the levels and ladders of an open castle, by exploiting the normality of $\CG\unitSpace$. The following definition is convenient for this purpose. 

\begin{defn}
	\label{defn:fortified-castle}
	A \emph{fortified castle} is a collection $\CF$ of pairs $C := (C_+, C_-)$ of nonempty open subsets in $\CG$, such that $\CF_+ := \left\{ C_+ \colon C \in \CF \right\}$ and $\CF_- := \left\{ C_- \colon C \in \CF \right\}$ are both open castles in $\CG$, and $\overline{ C_- } \subseteq C_+$ for any $C \in \CF$. 
	Here $\CF_+$ and $\CF_-$ are, respectively, called the \emph{outer castle} and the \emph{inner castle} of $\CF$. 
	We also write $\CF\unitSpace$ for $\left\{ C \in \CF \colon C_+ \in \CF_+\unitSpace \right\}$. 
	
	Given any $C, C' \in \CF$, we write $C \sim_{\CF} C'$ if $C_+ \sim_{\CF_+} C'_+$. This defines an equivalence relation on $\CF$, whose equivalence classes are the \emph{fortified towers} in $\CF$. 
	
	We say an open castle $\CC$ is \emph{sandwiched} in $\CF$ if there is a bijection $\varphi \colon \CF \to \CC$ such that $\overline{ C_- } \subseteq \varphi(C) \subseteq \overline{ \varphi(C) } \subseteq C_+$ for each $C \in \CF$. Note that such a bijection $\varphi$ is unique if it exists. 
\end{defn} 

\begin{rem}
	\label{rem:fortified-castle-basics}
	Given a fortified castle $\CF$ in $\CG$, the following hold: 
	\begin{enumerate}
		\item \label{rem:fortified-castle-basics:operations} For any $C, C' \in \CF$, we have $(r(C_+), r(C_-)), (s(C_+), s(C_-)), \left( C_+^{-1}, C_-^{-1} \right) \in \CF$ and $\left( C_+ C'_+ , C_- C'_- \right) \in \CF \cup \{ (\varnothing, \varnothing) \}$, and thus $C_+ \sim_{\CF_+} C'_+$ if and only if $C_- \sim_{\CF_-} C'_-$. 
		
		\item \label{rem:fortified-castle-basics:subcollection} For any $\CD \subseteq \CF_+$ that forms a castle, the collection $\left\{ C \in \CF \colon C_+ \in \CD \right\}$ is also a fortified castle. An analogous statement also holds with $\CF_+$ and $C_+$ replaced by $\CF_-$ and $C_-$. 
		
		\item \label{rem:fortified-castle-basics:tower} In particular, any fortified tower in $\CF$ is again a fortified castle. 
		
		\item \label{rem:fortified-castle-basics:sandwich} If an open castle $\CC$ is sandwiched in $\CF$, then both $\left\{ (C_+, \varphi(C)) \colon C \in \CF \right\}$ and $\left\{ ( \varphi(C) , C_-) \colon C \in \CF \right\}$ are fortified castles, where $\varphi$ is the bijection in Definition~\ref{defn:fortified-castle}. 
		
		\item \label{rem:fortified-castle-basics:precompact} It follows from Lemma~\ref{lem:castle-basic}\eqref{lem:castle-basic:bisection} that 
		if $\CF_-\unitSpace$ consists of precompact sets, then $\CF_-$ is a precompact castle\footnote{It may happen that an open castle has non-precompact ladders while all of its levels are precompact. Consider the transformation groupoid of an irrational rotation $\alpha \colon \mathbb{Z} \curvearrowright \mathbb{T}$. 
		We use the density of orbits to 
		recursively pick a sequence of pairs $\{(U_k,n_k)\}_{k=1}^\infty$ where $U_k$ is a nonempty open interval on $\mathbb{T}$ and $n_k$ is a positive integer, such that 
		$0 < k_1 < k_2 < \ldots$ and $\{ \overline{U_k}, \alpha_{n_k} (\overline{U_k}) \colon k = 1, 2, \ldots \}$ is a disjoint family.  
		Define an open castle $\CC = \left\{ C_{i,j} \colon i,j \in \{1,2\} \right\}$ where 
			$C_{1,2} = \bigsqcup_{k=1}^\infty \{ (x, n_k , \alpha_{n_k} (x)) \colon x \in U_k \}$ and the other sets are defined accordingly, i.e., $C_{1,1} = r(C_{1,2})$, $C_{2,2} = s(C_{1,2})$ and $C_{2,1} = C_{1,2}^{-1}$. }. 
		
		\item \label{rem:fortified-castle-basics:Urysohn} There is $f \in C_0 \left(\bigcup \CF_+ \right)^{\bigcup \CF_+}$ with its range in $[0,1]$ such that $f(x) = 1$ for any $x \in \bigcup \CF_-$.  
		To see this, we simply choose a set $\left\{ C_{\CT} \in \CT \colon \CT \in \CF / \sim_{\CF} \right\}$ of representatives of $\sim_{\CF}$, apply Urysohn's lemma to obtain $f_{C_{\CT}} \in C_0\left( \left( C_{\CT} \right)_+ \right)$ with $f_{C_{\CT}} (x) = 1$ for any $x \in \overline{\left( C_{\CT} \right)_-}$, and then apply Remark~\ref{rem:castle-fundamental-domain}\eqref{rem:castle-fundamental-domain:castle-ladder} to obtain the desired function $f$. 
		
		\item \label{rem:fortified-castle-basics:construction} 
		On the other hand, thanks to Lemma~\ref{lem:castle-trim}, we also have: given any open castle $\CC$, any $f \in C_0 \left(\bigcup \CC \right)^{\bigcup \CC}$ and any nonempty relatively open subset $Y \subseteq f \left(\bigcup \CC \right) \subseteq \mathbb{C}$ with $0 \notin \overline{Y}$, the collection 
		\[
		\left\{ \left( C, f^{-1} (Y) \cap C \right) \colon C \in \CC  \right\}
		\]
		forms a fortified castle. 
		Moreover, for any relatively open subset $Y' \subseteq f \left(\bigcup \CC \right) \subseteq \mathbb{C}$ with $\overline{Y} \subseteq Y' \subseteq \overline{Y'} \subseteq \mathbb{C} \setminus \{0\}$, the collection 
		\[
		\left\{ f^{-1} (Y') \cap C \colon C \in \CC  \right\}
		\]
		forms an open castle sandwiched in the above fortified castle. 
		
		\item \label{rem:fortified-castle-basics:sandwich-exist} 
		It follows from \eqref{rem:fortified-castle-basics:Urysohn} and \eqref{rem:fortified-castle-basics:construction} that there always exists an open castle sandwiched in $\CF$. 
		Repeating this process using \eqref{rem:fortified-castle-basics:sandwich}, we obtain infinitely many open castles sandwiched in $\CF$. 
	\end{enumerate}
\end{rem}

\section{Almost elementary \'{e}tale groupoids in measure} \label{sec:in-measure}
In this section, we will introduce the notion of \emph{almost-elementariness in measure} for \'{e}tale groupoids.  
It will be a precursor to one of the main regularity properties we study in this paper, \emph{almost-elementariness}. 

With the understanding that open castles are desirable structures to find in an \'{e}tale groupoid $\CG$, 
the intuition behind almost-elementariness (in measure) is to approximate $\CG$ by open castles. 
This idea of approximation breaks down into three kinds of parameters to control: how ``small'' the \emph{remainder} $\CG^{(0)} \setminus \left( \bigcup \CC^{(0)} \right)$ is, how ``fine'' the levels of $\CC$ are, and how ``tall'' the towers in $\CC$ are.
The two notions (i.e., with and without ``in measure'') differ in how to measure the ``smallness'' of the remainder. In this section, where we deal with the ``in measure'' version, we will essentially use the supremum of $\mu \left( \CG^{(0)} \setminus \left( \bigcup \CC^{(0)} \right) \right)$ as $\mu$ ranges over all invariant Borel probability measures on $\CG^{(0)}$. 
On the other hand, the two notions agree in how they measure the ``fineness'' of the levels and the ``tallness'' of the towers. 

In Definition~\ref{defn:castle}, we introduced the notion of heights for towers. However, such a numerical quantity is too weak to control the ``tallness'' of the towers in $\CC$. We shall need a more elaborate and indirect approach that should remind the reader of our discussion of fiberwise amenability in \cite{MWfiberwise}.

\begin{defn}\label{defn:castle-shrinking}
	Let $\CC$ be a castle in an \'{e}tale groupoid $\CG$ and let $K$ be a subset of $\CG$. 
	We define the \textit{left $K$-shrinking} of $\CC$ by
	\[L(\CC,-K) := \{C\in \CC: KC\subseteq\bigcup \CC\} \; ,\] 
	the \textit{right $K$-shrinking} of $\CC$ by
	\[R(\CC,-K) := \{C\in \CC: CK\subseteq\bigcup \CC\} \; , \]
	and the \emph{$K$-shrinking} of $\CC$ by 
	\begin{align*}
		\CC^{-K} :=&\ L(\CC, -K)\cap L(\CC, -K^{-1})\cap R(\CC, -K)\cap R(\CC, -K^{-1}) \\
		=&\ \left\{ C \in \CC \colon \left( K \cup K^{-1} \right) C \subseteq \bigcup \CC \text{ and }  C \left( K \cup K^{-1} \right) \subseteq \bigcup \CC \right\} \; .
	\end{align*}

	Similarly, if $\CF$ is a fortified castle, the \emph{$K$-shrinking} $\CF^{-K}$ is the fortified castle $\left\{ C \in \CF \colon C_+ \in (\CF_+)^{-K} \right\}$. Note that $(\CF_+)^{-K} = \left( \CF^{-K} \right)_+$ and $(\CF_-)^{-K} \supseteq \left( \CF^{-K} \right)_-$ but equality may not hold in the latter. 
\end{defn}

\begin{lem}
	\label{lem:castle-shrinking}
	Let $\CC$ be a castle (respectively, a fortified castle) in an \'{e}tale groupoid $\CG$ and let $K$ be a subset of $\CG$. Then the following hold: 
	\begin{enumerate}
		\item \label{lem:castle-shrinking:castle} The set $\CC^{-K}$ is again a castle (respectively, a fortified castle). 
		\item \label{lem:castle-shrinking:tower} For any tower (respectively, fortified tower) $\CT$ in $\CC$, we have $\CT^{-K} = \CT \cap \CC^{-K}$, and if $\CT^{-K} \not= \varnothing$, then it is a tower (respectively, fortified tower) in $\CC^{-K}$. 
	\end{enumerate}
\end{lem}

\begin{proof}
	We only prove the case of castles, from which the case of fortified castles follows, in view of Remark~\ref{rem:fortified-castle-basics}\eqref{rem:fortified-castle-basics:subcollection}. 
	
	Statement \eqref{lem:castle-shrinking:castle} follows the observation that $\CC^{-K}$ is equal to the castle constructed from $\CD_0 := \CC\unitSpace \cap \CC^{-K}$ as in Lemma~\ref{lem:castle-basic}\eqref{lem:castle-basic:reduction}.

For \eqref{lem:castle-shrinking:tower}, 
	we first observe that $\CT^{-K} \subseteq \CT \cap \CC^{-K}$ by definition. 
	On the other hand, for any $C \in \CT \cap \CC^{-K}$, we have $\left( K \cup K^{-1} \right) C \subseteq \bigcup \CT$ since for any $D \in \CC \setminus \CT$, we have $s(D) \cap s \left( \left( K \cup K^{-1} \right) C \right) \subseteq s(D) \cap s(C) = \varnothing$, and similarly we have $C \left( K \cup K^{-1} \right) \subseteq \bigcup \CT$, whence $C \in \CT^{-K}$. This shows $\CT^{-K} = \CT \cap \CC^{-K}$. 

	We then observe that for any $C, D \in \CT^{-K}$, since $C,D \in \CT$, there is $E \in \CT$ such that $s(C) = r(E)$ and $s(E) = r(D)$, but we also have $E \in \CT^{-K}$ as $\left( K \cup K^{-1} \right) E \subseteq \left( K \cup K^{-1} \right) s(C) E \subseteq \bigcup \CT$ and similarly $E \left( K \cup K^{-1} \right) \subseteq \bigcup \CT$, whence $C \sim_{\CT^{-K}} D$. 
	On the other hand, for any $D' \in \CC^{-K} \setminus \CT^{-K} = \CC^{-K} \setminus \CT \subseteq \CC \setminus \CT$, we have $C \not\sim_{\CT} D$ and thus $C \not\sim_{\CT^{-K}} D$. 
	This shows $\CT^{-K}$, if nonempty, is a tower in $\CC^{-K}$. 
\end{proof}

We then define the dual notion of shrinking. 

\begin{defn}\label{defn:castle-extendable}
	We say a castle $\CC$ (respectively, the fortified castle $\CF$) in an \'{e}tale groupoid $\CG$ is \emph{$K$-extendable} if there is a castle $\CD$ (respectively, a fortified castle $\CE$) in $\CG$ such that $\CC \subseteq \CD^{-K}$ (respectively, a fortified castle $\CF \subseteq \CE^{-K}$). 
	
	When we need to refer to the castle  $\CD$ (respectively, the fortified castle $\CE$) in the above, we say $\CC$ (respectively, $\CF$) is \emph{$K$-extendable to} $\CD$ (respectively, $\CE$). 
\end{defn}

\begin{rem}\label{rem: left extendbility enough}
	In Definition~\ref{defn:castle-extendable}, if $K$ is symmetric in the sense that $K^{-1} = K$, then we have $\CC \subseteq \CD^{-K}$ if and only if $\CC\subseteq L(\CD, -K)$ if and only if $\CC\subseteq R(\CD, -K)$. This is due to the fact that $\CC$ is closed under the (set) inverse operation. 
\end{rem}

\begin{rem}\label{rem: multisection-castle-extendability}
	When using the language of open multisections to describe open castles as in Corollary~\ref{cor:castle-multisections}, we may characterize extendability as follows:  
	Suppose $K$ a set in $\CG$ and $\CC=\{C_{i,j}^{l}:i,j\in F_{l},l\in I\}$ satisfies the condition in Corollary~\ref{cor:castle-multisections}. 
	Then $\CC$ is $K$-extendable to a castle $\CD$ if and only if there is an enumeration $\CD=\{D_{i',j'}^{l'}:i',j'\in F'_{l'},l'\in I'\}$ that satisfies the condition in Corollary~\ref{cor:castle-multisections} and the requirements that 
	$I \subseteq I'$, $F'_{l}\subseteq F_{l}$ for any $l \in I$, $C_{i,j}^{l}=D_{i.j}^{l}$ for any $l \in I$ and any $i,j\in F_{l}$, and 
	\[
		K\cdot\bigsqcup_{i,j\in F_{l}}C_{i,j}^{l}\subseteq\bigsqcup_{i,j\in E_{l}}D_{i,j}^{l} \quad \text{ for any } l \in I \; .
	\]
\end{rem}

We point out that the lack of principality in an \'{e}tale groupoid $\CG$ can be an obstruction to the extendability of castles in $\CG$. 

\begin{lem} \label{lem:castle-extendable-obstruction}
	Let $\CG$ be an \'{e}tale groupoid and let $K$ be a subset of $\CG$ containing two distinct elements $x$ and $y$ with $s(x) = s(y)$ and $r(x) = r(y)$. Then for any $K$-extendable castle $\CC$ in $\CG$, we have $(\bigcup \CC) \cap \{ s(x), r(x), x, x^{-1}, y , y^{-1} \} = \varnothing$. 
\end{lem}

\begin{proof}
	We prove that $s(x) \notin \bigcup \CC$, the case for $r(x)$ being similar, while the cases for $x$, $x^{-1}$, $y$, $y^{-1}$ will then follow, as $\bigcup \CC$ is closed under $s$ and $r$. 
	This amounts to proving that for any castle $\CD$ in $\CG$, we have $s(x) \notin \bigcup \CD^{-K}$. 
	Suppose the contrary, i.e., $s(x) \in D$ for some $D \in \CD^{-K}$. By Definition~\ref{defn:castle-shrinking}, we have $\{x,y\} = \{x,y\} s(x) \subseteq K D \subseteq \bigcup \CD$ and thus there are $D', D'' \in \CD$ such that $x \in D'$ and $y \in D''$. 
	Since $D'$ and $D''$ are bisections while $s$ is not injective on $\{x, y\}$, we have $D' \not= D''$. 
	On the other hand, since $s(D') \cap s(D'') \not= \varnothing$ and thus $s(D) = s(D')$, and similarly $r(D') = r(D'')$, by the injectivity requirement of $r \times s \colon \CD \to \CD\unitSpace \times \CD\unitSpace$ in Definition~\ref{defn:castle}, we also have $D' = D''$, a contradiction. 	
\end{proof}

Now we arrive at the following definition as our first groupoid regularity property in this paper. 

\begin{defn} \label{defn:AE-in-measure} 
	We say that an \'{e}tale groupoid $\CG$ with a compact unit space is \emph{almost elementary in measure}
	if for any compact set $K$ in $\CG$, 
	any $\varepsilon > 0$, 
	and any finite open cover $\CV$ of $\CG^{(0)}$, 
	there is an open castle $\CC$
	satisfying 
	\begin{enumerate}[label=(\roman*)]
		\item \label{defn:AE-in-measure:extendable} $\CC$ is $K$-extendable; 
		\item \label{defn:AE-in-measure:fine} $\CC\unitSpace$ refines $\CV$, i.e., any level $C \in \CC\unitSpace$ is contained in an open set $V\in\CV$; 
		\item \label{defn:AE-in-measure:full} $\bigcup \CC^{(0)}$ is \emph{full} in $\CG$ in the sense that $r \left( \CG \left(\bigcup \CC^{(0)}\right) \right) = \CG\unitSpace$; 
		\item \label{defn:AE-in-measure:remainder} $\mu \left( \CG^{(0)}\setminus\ \left( \bigcup \CC^{(0)} \right) \right) \leq \varepsilon$ for any invariant regular Borel probability measure $\mu$ on $\CG^{(0)}$.  
	\end{enumerate}
	In this case, the castle $\CC$ is sometimes referred to as an \emph{almost elementary approximation of $\CG$ in measure} and the set $\CG^{(0)}\setminus\ \left( \bigcup \CC^{(0)} \right)$ is the \emph{remainder} of this approximation. 
\end{defn}

\begin{rem}
	Note that we do not rule out the possibility that the set $M(\CG)$ of invariant regular Borel probability measures on $\CG\unitSpace$ is empty, in which case  condition~\ref{defn:AE-in-measure:remainder} in Definition~\ref{defn:AE-in-measure} becomes vacuous. 
\end{rem}

\begin{rem} \label{rem:AE-in-measure-fullness}
	Condition~\ref{defn:AE-in-measure:full} in Definition~\ref{defn:AE-in-measure} may take various equivalent forms: 
	\begin{enumerate}
		\item \label{rem:AE-in-measure-fullness:compact} It follows from the compactness of $\CG\unitSpace$ and Lemma~\ref{lem:minimal-recurrence} that condition~\ref{defn:AE-in-measure:full} is equivalent to the a priori stronger requirement that there is a compact subset $L \in \CG$ such that $r \left( L \left(\bigcup \CC^{(0)}\right) \right) = \CG\unitSpace$. 
		\item \label{rem:AE-in-measure-fullness:minimal} When $\CG$ is minimal and nonempty, condition~\ref{defn:AE-in-measure:full} is equivalent to the a priori weaker requirement that $\CC$ is nonempty. 
		\item \label{rem:AE-in-measure-fullness:finite} Moreover, when $\CG$ is minimal and $M(\CG) \not= \varnothing$, condition~\ref{defn:AE-in-measure:full} may be dropped since it is implied by the last condition (as soon as $\varepsilon < 1$).   	\end{enumerate}
\end{rem}

To help establish intuition, we first investigate how this notion behaves in the case of \'{e}tale groupoids with finite unit spaces, even though our main interest lies in the case of infinite unit spaces. 

\begin{prop} \label{prop:AEm-finite-pair}
	An \'{e}tale groupoid $\CG$ with a finite unit space is almost elementary in measure if and only if it is principal if and only if it is a disjoint union of pair groupoids. 
\end{prop}

\begin{proof}
	It is straightforward to see that an \'{e}tale groupoid $\CG$ with a finite unit space is principal if and only if $|u \CG v| \leq 1$ for any $u, v \in \CG\unitSpace$ if and only if it is a (finite) disjoint union of (finite) pair groupoids. 
	
	We also observe that if $\CG$ is a finite disjoint union of finite pair groupoids, then it is almost elementary in measure, since Definition~\ref{defn:AE-in-measure} can be witnessed by the castle consisting of all singletons in $\CG$, no matter the choices of $K$, $\varepsilon$ and $\CV$. 
	
	It remains to show that if $\CG$ is almost elementary in measure, then $|u \CG v| \leq 1$ for any $u, v \in \CG\unitSpace$. We prove this by contradiction. Suppose there are distinct elements $x, y \in \CG$ with $s(x) = s(y)$ and $r(x) = r(y)$. Then taking $K := \{x,y\}$, we deduce from Lemma~\ref{lem:castle-extendable-obstruction} that any $K$-extendable open castle does not cover $s(x)$. On the other hand, notice by Definition~\ref{def:invariant-measure} that $M(\CG)$ includes the uniform measure on $\CG\unitSpace$, i.e., the one assigning the weight $\frac{1}{|\CG\unitSpace|}$ to every singleton in $\CG\unitSpace$. 
	Hence, taking $\varepsilon = \frac{1}{2 |\CG\unitSpace|}$, we observe that any open castle $\CC$ satisfying condition~\ref{defn:AE-in-measure:remainder} in Definition~\ref{defn:AE-in-measure} must satisfy $\bigcup \CC\unitSpace = \CG\unitSpace$, a contradiction. 
\end{proof}

In practice, it is more convenient to use the notion of $K$-shrinking directly instead of that of $K$-extendability, which leads us to the equivalent characterization in Lemma~\ref{lem:AE-in-measure-shrinking}\eqref{lem:AE-in-measure-shrinking:shrinking} below. We also present a stronger version in  Lemma~\ref{lem:AE-in-measure-shrinking}\eqref{lem:AE-in-measure-shrinking:fine-entourage} that will be handy later. 
We shall first need an auxiliary result.

\begin{lem}
	\label{lem:entourage-fineness-control}
	Let $K$ be a compact subset in an \'{e}tale groupoid $\CG$ and let $\CV$ be a collection of open subset in $\CG$ that covers $K$. 
	Then there is a finite open cover $\CW$ of $\CG\unitSpace$ such that 
	for any $x \in K$ and any $W \in \CW$ containing $s(x)$, there is $V \in \CV$ containing $x$ such that $W \subseteq s(V)$. 
\end{lem}

\begin{proof}
	By compactness, there exists a finite subcollection $\left\{ V_{1}, \ldots, V_{n} \right\}$ of $\CV$ that covers $K$. 
	Define a collection of open subsets of $\CG\unitSpace$: 
	\begin{align*}
	\CW &:= \left\{ W_{F} \colon F \subseteq \{1, \ldots, n\} \right\} 
	&\text{where } W_{F} := \CG\unitSpace \cap \left( \bigcap_{i \in F} s \left( V_{i} \right) \right) \setminus s \left( K \setminus \left( \bigcup_{i \in F} V_{i}  \right) \right) \; .
	\end{align*}

	We claim that $\CW$ covers $\CG\unitSpace$. Indeed, for any $u \in \CG\unitSpace$, if we let $F_u := \{ i \in \{1, \ldots, n\} \colon u \in s \left( V_{i} \right) \}$, then $K u = K \cap s^{-1} (u) \subseteq K \cap \left( \bigcup_{i \in F_u} V_{i}  \right)$ and thus $u \notin s \left( K \setminus \left( \bigcup_{i \in F_u} V_{i}  \right) \right)$, whence $u \in W_{F_u}$, proving the claim. 
	
	Now, fixing arbitrary $x \in K$ and $F \subseteq \{1, \ldots, n\}$ with $s(x) \in W_{F}$, we deduce that $x \notin K \setminus \left( \bigcup_{i \in F} V_{i}  \right)$ and thus $x \in \bigcup_{i \in F} V_{i}$. 
	Choosing $i_0 \in F$ such that $x \in V_{i_0}$, we have $W_{F} \subseteq s(V_{i_0})$ by construction. This completes the proof. 
\end{proof}

\begin{lem} \label{lem:AE-in-measure-shrinking} 
	Let $\CG$ be an \'{e}tale groupoid with a compact unit space. 
	The following are equivalent: 
	\begin{enumerate}
		\item \label{lem:AE-in-measure-shrinking:original} The groupoid $\CG$ is {almost elementary in measure}. 
		
		\item \label{lem:AE-in-measure-shrinking:shrinking} For any compact set $K$ in $\CG$, 
		any $\varepsilon > 0$, 
		and any finite open cover $\CV$ of $\CG^{(0)}$, 
		there is an open castle $\CD$
		satisfying 
		\begin{enumerate}[label=(\roman*)]
			\item \label{lem:AE-in-measure-shrinking:refine} $\CD \unitSpace$ refines $\CV$; 
			\item \label{lem:AE-in-measure-shrinking:full} $\bigcup \left( \CD^{-K} \right)^{(0)}$ is {full} in $\CG$; 
			\item \label{lem:AE-in-measure-shrinking:remainder} 
			$\mu \left( \bigcup \left( \CD^{-K} \right)^{(0)} \right) \geq 1 - \varepsilon$ for any $\mu \in M(\CG)$.  
		\end{enumerate}
	
		\item \label{lem:AE-in-measure-shrinking:fine-entourage} For any compact sets $K, L \subseteq \CG$, 
		any $\varepsilon > 0$, 
		and any collection $\CV$ of open subsets in $\CG$ that covers $L$, 
		there is an open castle $\CD$
		satisfying 
		\begin{enumerate}[label=(\roman*)]
			\item \label{lem:AE-in-measure-shrinking:fine-entourage:refine} for any $D \in \CD$ with $D \cap L \not= \varnothing$, there is $V \in \CV$ such that $D \subseteq V$; 
			\item \label{lem:AE-in-measure-shrinking:fine-entourage:full} $\bigcup \left( \CD^{-K} \right)^{(0)}$ is {full} in $\CG$; 
			\item \label{lem:AE-in-measure-shrinking:fine-entourage:remainder} 
			$\mu \left( \bigcup \left( \CD^{-K} \right)^{(0)} \right) \geq 1 - \varepsilon$ for any $\mu \in M(\CG)$.  
		\end{enumerate}	
	\end{enumerate}
\end{lem}

Before starting the proof, let us point out where the difficulty lies: 
in the context of \eqref{lem:AE-in-measure-shrinking:fine-entourage} (assuming $\CG\unitSpace \subseteq K$), the condition that $\left( \CD^{-K} \right) \unitSpace$ refines $\CV$ is easily guaranteed by \eqref{lem:AE-in-measure-shrinking:original}, 
while \eqref{lem:AE-in-measure-shrinking:shrinking} amounts to the requirement that $\CD \unitSpace$ refines $\CV$, 
but \eqref{lem:AE-in-measure-shrinking:fine-entourage} demands a typically larger collection in $\CD$ also refines $\CV$. 

\begin{proof}
	The implication \eqref{lem:AE-in-measure-shrinking:shrinking} $\Rightarrow$ \eqref{lem:AE-in-measure-shrinking:original} 
	is clear, 
	since if a castle $\CD$ satisfies the conditions in \eqref{lem:AE-in-measure-shrinking:shrinking}, then the $K$-shrinking $\CC :=  \CD^{-K}$ will satisfy the conditions in Definition~\ref{defn:AE-in-measure}. 
	
	The implication \eqref{lem:AE-in-measure-shrinking:fine-entourage} $\Rightarrow$ \eqref{lem:AE-in-measure-shrinking:shrinking} 
	is also clear by taking $L := \CG\unitSpace$. 
    
	To show \eqref{lem:AE-in-measure-shrinking:original} $\Rightarrow$ \eqref{lem:AE-in-measure-shrinking:fine-entourage}, 
	we fix $K$, $\varepsilon$ and $\CV$ as in \eqref{lem:AE-in-measure-shrinking:fine-entourage}. 
	By replacing $\CV$ with $\CV \cup \{ \CG \setminus L \}$ and then enlarging both $K$ and $L$ to the same compact set $K \cup L \cup K^{-1} \cup L^{-1} \cup \CG\unitSpace$, we may assume without loss of generality that $K = K^{-1} = L = L^{-1}\supseteq \CG\unitSpace$. 
	Since $\CG$ is locally compact, Hausdorff and \'{e}tale, by refining $\CV$ if necessary, we may further assume without loss of generality that $\CV$ consists of precompact open bisections. 
	By compactness, there exists a finite subcollection $\left\{ V_{1}, \ldots, V_{n} \right\}$ of $\CV$ that covers $K$. 
	Applying Lemma~\ref{lem:entourage-fineness-control}, we obtain a finite open cover $\CU$ of $\CG\unitSpace$ such that for any $x \in K$ and any $U \in \CU$ containing $s(x)$, there is $i \in \{1, \ldots, n\}$ such that $x \in V_i$ and $U \subseteq s(V_i)$. 
	Applying Lemma~\ref{lem:entourage-fineness-control} again with $\CV$ replaced by $\left\{ V_i \cap r^{-1}(U) \colon i \in \{1, \ldots, n\}, U \in \CU \right\}$, we obtain a finite open cover $\CW$ of $\CG\unitSpace$ such that for any $x \in K$ and any $W \in \CW$ containing $s(x)$, there are $i \in \{1, \ldots, n\}$ and $U \in \CU$ such that $x \in V_i$, $r(x) \in U$ and $W \subseteq s \left( V_i \cap r^{-1}(U) \right)$. 
	Without loss of generality, we may assume $\CW$ refines $\CU$. 
	Define $K' := \bigcup_{i, j =1}^n \overline{V_{i}} \cdot \overline{V_{j}}$, which is a compact set in $\CG$ that contains $K^2$ as well as $\bigcup_{i =1}^n \overline{V_{i}}$ and $K$, since $\CG\unitSpace \subseteq K$.

	By \eqref{lem:AE-in-measure-shrinking:original}, we may apply Definition~\ref{defn:castle-shrinking} to find an open castle $\CC$ 
	satisfying 
	\begin{enumerate}[label=(\roman*)]
		\item $\CC$ is $K'$-extendable; 
		\item $\CC\unitSpace$ refines $\CW$; 
		\item $\bigcup \CC\unitSpace$ is full in $\CG$; 
		\item $\mu \left( \CG^{(0)}\setminus\ \left( \bigcup \CC^{(0)} \right) \right) \leq \varepsilon$ for any $\mu \in M(\CG)$.  
	\end{enumerate}
	By Definition~\ref{defn:castle-shrinking}, there exists an open castle $\CE$ such that $\CC \subseteq \CE^{-K'}$. 
	
	For any $C \in \CC\unitSpace$ and any $j \in \{1,\ldots, n\}$ such that $C \subseteq s\left( V_{j} \right)$, we claim that the open bisection $B := r ( V_{j} C )$ is a relatively clopen subset in $\bigcup \CE$. 
	Indeed, since $V_{j} C \subseteq K' C \subseteq \bigcup \CE$, we have $ r ( V_{j} C ) \subseteq \bigcup \CE$. 
	Moreover, it follows from the containment $V_{j} C \subseteq \bigcup \CE$ that there are $E_1, \ldots, E_m \in \CE$ such that $V_{j} C \subseteq \bigsqcup_{k=1}^m E_k$ and $(V_{j} C) \cap E_k \not= \varnothing$ for any $k \in \{1, \ldots, m\}$, whence $\left\{ (V_{j} C) \cap E_k \colon k \in \{1, \ldots, m\} \right\}$ forms a relatively clopen partition of $V_{j} C$. 
	Since $V_{j} C$ is a bisection with $s (V_{j} C) = C$, we see that $\left\{ s ((V_{j} C) \cap E_k)  \colon k \in \{1, \ldots, m\} \right\}$ forms a relatively clopen partition of $C$. 
	On the other hand, for any $k \in \{1, \ldots, m\}$, 
	since $C$ and $E_k$ are both in the castle $\CE$ and $C \cap s\left(E_k\right) = s\left( V_{j} C \right) \cap s\left(E_k\right) \supseteq s \left( (V_{j} C) \cap E_k \right) \not=\varnothing$, it follows that $C = s(E_k)$ and thus $s ((V_{j} C) \cap E_k)$ is relatively clopen in $s(E_k)$, 
	and furthermore, 
	since $E_k$ is also a bisection, 
	we have $(V_{j} C) \cap E_k = E_k \cdot s \left( (V_{j} C) \cap E_k \right)$,
	which is thus a relatively clopen subset in $E_k$ and thus also a relatively clopen subset in $\bigcup \CE$. It follows that $V_{j} C$ is relatively clopen in $\bigcup \CE$ and thus so is $r(V_{j} C)$, proving the claim.
	
	Next, 
	for any $C \in \CC\unitSpace$ and any $i, j \in \{1,\ldots, n\}$ such that $C \subseteq s\left( V_{i} V_{j} \right)$ and thus $r ( V_{j} C ) \subseteq s(V_i)$, we claim that the open bisection $B := V_{i} \cdot r ( V_{j} C )$ is a relatively clopen subset in $\bigcup \CE$. 
	Indeed, since $V_{j} C \subseteq K' C \subseteq \bigcup \CE$ and $V_{i} V_{j} C \subseteq K' C \subseteq \bigcup \CE$, it follows that $B = V_{i} ( V_{j} C ) ( V_{j} C )^{-1} \subseteq (\bigcup \CE) \cdot (\bigcup \CE)^{-1} = \bigcup \CE$. 
	Similar to the above, it follows from the containment $B \subseteq \bigcup \CE$ that there are $F_1, \ldots, F_n \in \CE$ such that $\left\{ B \cap F_k \colon k \in \{1, \ldots, l\} \right\}$ forms a relatively clopen partition of $V_{j} C$. Since $B$ is a bisection with $s(B) = r ( V_{j} C )$, we see that $\left\{ s (B \cap F_k)  \colon k \in \{1, \ldots, l\} \right\}$ forms a relatively clopen partition of $r ( V_{j} C )$. 
	Since we proved above that $r ( V_{j} C )$ is relatively clopen in $\bigcup \CE$, it follows that for any $k \in \{1, \ldots, l\}$, the set  $s (B \cap F_k)$ is also relatively clopen in $\bigcup \CE$ and thus in $s(F_k)$, 
	and furthermore,  
	since $F_k$ is also a bisection, 
	we have $B \cap F_k = F_k \cdot s \left( B \cap F_k \right)$, 
	which is a relatively clopen subset in $F_k$ and thus also a relatively clopen subset in $\bigcup \CE$. 
	It follows that $B$ is relatively clopen in $\bigcup \CE$, proving the claim.

	Hence for any $C \in \CC\unitSpace$ and any $i, j \in \{1,\ldots, n\}$ such that $C \subseteq s\left( V_{i} V_{j} \right)$, we obtain a relatively clopen (and thus open in $\CG$) partition 
	\[
		\CP_{C,i,j} := \left\{ r \left( V_{j} C \right), \left(\bigcup \CE\right) \setminus r \left( V_{j} C \right) \right\} \vee \left\{ V_{i} \cdot r \left( V_{j} C \right), \left(\bigcup \CE\right) \setminus \left( V_{i} \cdot r \left( V_{j} C \right) \right)\right\}
	\]
	 of $\bigcup \CE$. 
	Let $\CP$ be the common refinement of all $\CP_{C,i,j}$, where $(C, i, j)$ ranges over all pairs where $C \in \CC$ and $i, j \in \{1,\ldots, n\}$ satisfying $C \subseteq s\left( V_{i} V_{j} \right)$. By Lemma~\ref{lem:castle-split}\eqref{lem:castle-split:exists} and~\eqref{lem:castle-split:extra}, there is an open castle $\CE_\CP$ that refines both $\CE$ and $\CP$ and satisfies $\bigcup \CE_\CP = \bigcup \CE$. 
	
	Define 
	\[
		\CD_0 := \left\{ D \in \CE_\CP\unitSpace \colon D K \left( \bigcup \CC\unitSpace \right) \not= \varnothing  \right\}
	\]
	and let $\CD \subseteq \CE_\CP$ be the open castle constructed from $\CD_0$ as in Lemma~\ref{lem:castle-basic}\eqref{lem:castle-basic:reduction}, so that $\CD\unitSpace = \CD_0$. By our assumption that $K = K^{-1} \subseteq K'$, we deduce from $\CC \subseteq \CE^{-K'}$ that $\bigcup \CC \subseteq \bigcup \CD^{-K}$. 
	In particular, $\bigcup \left( \CD^{-K} \right)^{(0)}$ is full and satisfies 
	\begin{equation}\label{lem:AE-in-measure-shrinking:eq:remainder-control}
	\mu \left( \bigcup \left( \CD^{-K} \right)^{(0)} \right) \geq \mu \left( \bigcup \CC^{(0)} \right) = 1 - \mu \left( \CG^{(0)}\setminus\ \left( \bigcup \CC^{(0)} \right) \right) \geq 1 - \varepsilon \; . 
	\end{equation}
	It remains to verify condition~\ref{lem:AE-in-measure-shrinking:fine-entourage:refine} in \eqref{lem:AE-in-measure-shrinking:fine-entourage}. To this end, we fix an arbitrary $D \in \CD$ with $D \cap K \not= \varnothing$. 
	By our assumption above, there is $C \in \CC\unitSpace$ and $x \in K$ such that $D x C \not= \varnothing$, i.e., $s(x) \in C$ and $r(x) \in s(D)$. 
	Since $\CC\unitSpace$ refines $\CW$, there is $W \in \CW$ such that $C \subseteq W$. 
	By our construction of $\CW$, there is $j \in \{1, \ldots, n\}$ and $U \in \CU$ such that $x \in V_j$, $r(x) \in U$ and $W \subseteq s \left( V_j \cap r^{-1}(U) \right)$. 
	By our construction of $\CP$, we see that $\CE_\CP$ refines $\left\{ r \left( V_{j} C \right), \left( \bigcup \CE \right) \setminus r \left( V_{j} C \right) \right\}$ and thus so does $\CD$. Since $r(x) \in s(D) \cap r \left( V_{j} C \right)$, it follows that $s(D) \subseteq r \left( V_{j} C \right) \subseteq r \left( V_{j} \cdot s \left( V_j \cap r^{-1}(U) \right) \right) \subseteq r \left( V_{j} \right) \cap U$. 
	Now choose $y \in D \cap K$. 
	Since $s(y) \in s(D) \subseteq U$, 
	by our construction of $\CU$, there is $i \in \{1, \ldots, n\}$ such that $y \in V_i$ and $U \subseteq s \left( V_i \right)$. 
	By our construction of $\CP$, we see that $\CE_\CP$ refines $\left\{ V_i \cdot r \left( V_{j} C \right), \left( \bigcup \CE \right) \setminus \left( V_i \cdot r \left( V_{j} C \right) \right) \right\}$ and thus so does $\CD$. Since $y \in D \cap \left( V_i \cdot s(D) \right) \subseteq D \cap \left( V_i \cdot r \left( V_{j} C \right) \right)$, it follows that $D \subseteq V_i \cdot r \left( V_{j} C \right) \subseteq V_i$. This verifies condition~\ref{lem:AE-in-measure-shrinking:fine-entourage:refine} in \eqref{lem:AE-in-measure-shrinking:fine-entourage}. 
\end{proof}

\section{Relations with effectiveness and the groupoid small boundary property}\label{sec:in-measure-extra}

In this section, we discuss how almost elementariness in measure interacts with a number of other structural properties of \'{e}tale groupoids, including effectiveness, groupoid small boundary property, and fiberwise amenability. 

We begin by reformulating effectiveness in terms of open castles. 

\begin{prop} \label{prop:effective-castle}
	The following are equivalent for an \'{e}tale groupoid $\CG$: 
	\begin{enumerate}
		\item \label{prop:effective-castle:effective} The groupoid $\CG$ is effective. 
		
		\item \label{prop:effective-castle:castle} For any compact subset $K \subseteq \CG$ and any nonempty open subset $O \subseteq \CG\unitSpace$, there exists a $K$-extendable open castle $\CC$ in $\CG$ such that $O \cap (\bigcup \CC\unitSpace) \not= \varnothing$. 
		
		\item \label{prop:effective-castle:tower-fine} 
		For any compact subset $K \subseteq \CG$ and any collection $\CV$ of open subsets of $\CG\unitSpace$ with $\bigcup \CV \not= \varnothing$, there exists a nonempty $K$-extendable open tower $\CC$ in $\CG$ such that any $C \in \CC\unitSpace$ is contained in some $V \in \CV$.  
	\end{enumerate}
\end{prop}
\begin{proof}
	The implication \eqref{prop:effective-castle:tower-fine} $\Rightarrow$ \eqref{prop:effective-castle:castle} is trivial by taking $\CV = \left\{ O \right\}$.

	We prove \eqref{prop:effective-castle:castle} $\Rightarrow$ \eqref{prop:effective-castle:effective} 
	by contradiction. 
	Suppose 	
	$\CG$ is not effective, i.e., the open set $\opIso(\CG)^{\operatorname{o}}\setminus\CG^{(0)}$ is nonempty.  
	Then there is a nonempty precompact open
	bisection $V$ such that $V\subseteq\overline{V}\subseteq\opIso(\CG)^{o}\setminus\CG^{(0)}$. 
	Then define a compact set $K := \overline{V}$ and a nonempty open set $O=s(V)=r(V)\subseteq\CG^{(0)}$. 
	We claim that any $K$-extendable open castle $\CC$ in $\CG$ satisfies $O \cap (\bigcup \CC\unitSpace) = \varnothing$. 
	Indeed, for any $u \in O$, since $V$ is a bisection in $\opIso(\CG)^{o}\setminus\CG^{(0)}$, we have $V u = \{ x \}$ for some $x \in \CG \setminus \CG\unitSpace$ such that $s(x) = r(x) = u$ but $x \not= u$, 
	whence by Lemma~\ref{lem:castle-extendable-obstruction}, $u \notin \bigcup \CC\unitSpace$. This contradicts \eqref{prop:effective-castle:castle}.

	To prove \eqref{prop:effective-castle:effective} $\Rightarrow$ \eqref{prop:effective-castle:tower-fine}, we assume $\CG$ is effective and fix $K$ and $\CV$ as in \eqref{prop:effective-castle:tower-fine}. Without loss of generality, we may assume $K = K^{-1}$. Fix $V_0 \in \CV$ and $x_0 \in V_0$, and enumerate the finite subset $K x_0 \setminus \{ x_0 \}$ of the source fiber $\CG_{x_0}$ as $\left\{ x_1, \ldots, x_n \right\}$. 
	Choose disjoint precompact open bisections $W_{i} \subseteq \CG \setminus \CG\unitSpace$ with $x_i \in W_{i}$ and $s \left( \overline{ W_{i} } \right) \subseteq V_0$ for $i \in \{1,\ldots,n\}$. 
	Applying Lemma~\ref{lem:entourage-fineness-control} with $\CV$ replaced by the collection $\CW := \left\{ W_1, \ldots, W_n, W_{0} := \bigcap_{i=1}^n s \left( W_{i} \right), \CG \setminus (\CG x_0) \right\}$ of open sets covering $K$, we obtain a finite open cover $\CU$ of $\CG\unitSpace$ such that for any $x \in K$ and $U \in \CU$ containing $s(x)$, there is $W \in \CW$ containing $x$ such that $U \subseteq s(W)$. 
	Pick $U_0 \in \CU$ that contains $x_0$. 
	Since $U_0 \not\subseteq s(\CG \setminus (\CG x_0))$, for any $x \in K u$ with $u=s(x) \in U_0$, it follows from the construction of $\CU$ that $x \in \bigcup_{i =0}^n W_{i}$, whence
	\[
		K U_0 \subseteq \bigcup_{i =0}^n W_{i} \; . 
	\]
	Define nonempty precompact open bisections
	\[
	U_{i} :=  W_{i} U_{0} = \left( s \middle|_{W_{i}} \right)^{-1} \left( U_{0} \right) \subseteq W_{i} \quad \text{ for } i \in \{1, \ldots, n\}  \; , 
	\]
	which satisfy $s \left( U_{i} \right) = U_{0}$ for any $i \in \{1, \ldots, n\} $. Observe that for any $i , j \in \{0, 1, \ldots, n\}$, the set $U_{i} U_{j}^{-1}$ is a nonempty precompact open bisection homeomorphic to $U_0$ via the homeomorphism 
	\[
	\varphi_{i,j} := s \circ \left( r \middle|_{U_{i}} \right)^{-1} \circ r = s \circ \left( r \middle|_{U_{j}} \right)^{-1} \circ s \colon U_{i} U_{j}^{-1} \to U_{0} \; ,
	\]
	which satisfies the identites $\varphi_{i,i} \circ r = \varphi_{i,j} = \varphi_{j,j} \circ s$ and $\varphi_{i,j} \left( (-)^{-1} \right) = \varphi_{j,i}(-)$ for any $i , j \in \{0, 1, \ldots, n\}$. 
	
	Since $\CG$ is effective, the set $\left(\CG\setminus \opIso(\CG)\right) \cup \CG\unitSpace$ is dense and open in $\CG$. It follows that for any $i , j \in \{0, 1, \ldots, n\}$, the set 
	\[
	Z_{i,j} := \left( U_{i} U_{j}^{-1} \right) \cap \left( \left(\CG\setminus \opIso(\CG)\right) \cup \CG^{(0)} \right)
	\]
	is dense and open in $U_{i} U_{j}^{-1}$, whence $\varphi_{i,j} \left( Z_{i,j} \right)$ is dense and open in $U_{0}$. Hence the finite intersection 
	\[
	Z := \bigcap_{i,j \in \{ 0, 1, \ldots, n \}} \varphi_{i,j} \left( Z_{i,j} \right)
	\]
	is also dense and open in $U_{0}$, and in particular nonempty. 
	Fix $u_0$ in $Z$. 
	For any distinct $i, j \in \{ 0, 1, \ldots, n\} $, since $U_i \cap U_j \subseteq W_i \cap W_j = \varnothing$, we have $\left(U_{i} U_{j}^{-1}\right) \cap \CG\unitSpace = \varnothing$ and thus 
	$\varphi_{i,j}^{-1} \left( u_0 \right) \in Z_{i,j} = Z_{i,j} \setminus \CG\unitSpace \subseteq \CG\setminus \opIso(\CG)$, whence $\varphi_{i,i}^{-1} \left( u_0 \right) = r \left( \varphi_{i,j}^{-1} \left( u_0 \right) \right) \not= s \left( \varphi_{i,j}^{-1} \left( u_0 \right) \right) = \varphi_{j,j}^{-1} \left( u_0 \right)$. 
	Choose disjoint open sets $Y_i$ satisfying $\varphi_{i,i}^{-1} \left( u_0 \right) \in Y_i \subseteq \varphi_{i,i}^{-1} \left( Z \right)$ for $i \in \{ 0, 1, \ldots, n\}$ and define 
	\[
	C_{0,0} := \bigcap_{i \in \{ 0, 1, \ldots, n \}} \varphi_{i,i} \left( Y_i \right) \subseteq Z \subseteq U_0 \subseteq W_0 \subseteq V_0
	\]
	which is open and nonempty since it contains $u_0$. 
	Define nonempty precompact open bisections
	\[
	C_{i,j} := \varphi_{i,j}^{-1} \left( C_{0,0} \right)  = U_{i} C_{0,0} U_{j}^{-1} \subseteq Z_{i,j} \subseteq U_{i} U_{j}^{-1} \quad \text{ for } i, j \in \{ 0, 1, \ldots, n\}  \; , 
	\]
	which satisfies 
	$C_{i,i} \subseteq Y_i$, 
	$r \left( C_{i,j} \right) = C_{i,i}$, $s \left( C_{i,j} \right) = C_{j,j}$, 
	$C_{i,j} C_{j,k} = C_{i,k}$ and $C_{i,j}^{-1} = C_{j,i}$ for any $i, j, k \in \{ 0, 1, \ldots, n\} $. 
	Note there is no conflict in the definition of $C_{0,0}$. 
	It follows from the disjointness of $\left\{ Y_i \right\}_{i \in \{0, 1, \ldots, n\}}$ that the family $\left\{ C_{i,i} \right\}_{i \in \{0, 1, \ldots, n\}}$ is also disjoint, whence the family $\CT := \left\{ C_{i,j} \colon i,j \in \{0, 1, \ldots, n\} \right\}$ is an open multisection consisting of nonempty sets and thus an open tower by Proposition~\ref{prop:tower-multisection}. 
	Finally, for any $u \in C_{0,0} \subseteq U_0$, since we have proved $K u \subseteq \bigcup_{i =0}^n W_{i}$, we also have 
	\[
	K u \subseteq \bigcup_{i =0}^n W_{i} u = \bigcup_{i =0}^n U_{i} u = \left\{ \varphi_{i,0}^{-1} \left( u_0 \right) \colon i \in \{0, 1, \ldots, n\} \right\} \subseteq \bigcup_{i =0}^n C_{i,0} \subseteq \bigcup \CT
	\]
	and thus also $u K \subseteq \bigcup_{i =0}^n C_{0, i}\subseteq \bigcup \CT$. It follows that $C_{0,0} \in \CT^{-K}$. Therefore the tower $\CC = \left\{ C_{0,0} \right\}$ satisfies all the requirements. 
\end{proof}

Proposition~\ref{prop:effective-castle}\eqref{prop:effective-castle:castle} takes a simpler form for minimal groupoids with compact unit spaces. 

\begin{cor} \label{cor:effective-castle-minimal}
	A nonempty minimal \'{e}tale groupoid $\CG$ with a compact unit space is effective if and only if for any compact subset $K \subseteq \CG$, there exists 
	a nonempty $K$-extendable open castle. 
\end{cor}

\begin{proof}
	Since $\CG$ is nonempty, the condition above immediately follows from Proposition~\ref{prop:effective-castle}\eqref{prop:effective-castle:castle} upon fixing $O := \CG\unitSpace$. This proves the ``only if'' direction. 
	
	For the ``if'' direction, 
	we fix $K$ and $O$ as in Proposition~\ref{prop:effective-castle}\eqref{prop:effective-castle:castle}. 
	Since $\CG\unitSpace$ is compact and $\CG\unitSpace = r \left( \CG \cdot O \right)$ by minimality, 
	we may apply Lemma~\ref{lem:minimal-recurrence} to find compact bisections $K_{1},\dots,K_{n}$
	such that $\bigcup_{i=1}^{n}r\left(K_{i}\right)\subseteq O$ and $\CG\unitSpace \subseteq\bigcup_{i=1}^{n}s\left(K_{i}\right)^{\operatorname{o}}$. 
	Now we define a compact set $L = \bigcup_{i=1}^{n} (K \cup K^{-1} \cup \CG\unitSpace ){K_i}$. 
	By our assumption, there is a nonempty $L$-extendable open castle, or equivalently, an open castle $\CD$ such that $\CD^{-L}$ is nonempty. 
	Fix $D_0 \in \left( \CD^{-L} \right)\unitSpace$ and $u \in D_0$. 
	Since $\CG\unitSpace \subseteq\bigcup_{i=1}^{n}s\left(K_{i}\right)^{\operatorname{o}}$, there is $i \in \{1, \ldots, n\}$ such that $u \in s(K_i)^{\operatorname{o}}$. 
	Since $K_i$ is a bisection with $r(K_i) \subseteq O$, there is $x \in K_i$ such that $s(x) = u$ and $r(x) \in O$. Since $x \in L$, by Definition~\ref{defn:castle-shrinking}, there is $D_1 \in \CD$ such that $x \in D_1$. Note that $s(D_1) = D_0$. 
	Let $E_0 := D_1 \cap K_i^{\operatorname{o}}$, which is again an open bisection containing $x$. 
	Let $C := r(E_0) \subseteq r(D_1) \cap r(K_i^{\operatorname{o}})$ and $\CE := \left\{ DCD' \colon D, D' \in \CD \right\}$. It is routine to verify that $\CE$ is an open castle (in fact, an open tower) containing $C$. 
	
	We claim that $C \in \CE^{-K}$. Indeed, by our construction and the assumption that $D_0 \in \left( \CD^{-L} \right)\unitSpace$, we have 
	\begin{align*}
		(K \cup K^{-1}) C &= (K \cup K^{-1}) E_0 s(E_0) E_0^{-1} C \\
		&\subseteq (K \cup K^{-1}) K_i s(D_1) D_1^{-1} C \\
		&\subseteq L D_0 D_1^{-1} C \\
		&\subseteq (\bigcup \CD) D_1^{-1} C = D_1^{-1} C \in \CE
	\end{align*}
	Hence $(K \cup K^{-1}) C \subseteq \bigcup \CE$ and similarly $C (K \cup K^{-1}) \subseteq \bigcup \CE$, which proves the claim. Hence the nonempty open castle $\{C\}$ is $K$-extendable. 
\end{proof}

\begin{thm} 
	\label{thm:effective}
	If a minimal \'{e}tale groupoid $\CG$ with a compact unit space is almost elementary in measure, then
	$\CG$ is effective. 
	The converse holds when $M(\CG) = \varnothing$. 
\end{thm}

\begin{proof}
	In view of Lemma~\ref{lem:AE-in-measure-shrinking}\eqref{lem:AE-in-measure-shrinking:shrinking}, the first statement follows directly from the implication \eqref{prop:effective-castle:castle} $\Rightarrow$ \eqref{prop:effective-castle:effective} in Lemma~\ref{prop:effective-castle}. 
	The second statement follows directly from the implication \eqref{prop:effective-castle:effective} $\Rightarrow$ \eqref{prop:effective-castle:tower-fine} in Lemma~\ref{prop:effective-castle}, since when $M(\CG) = \varnothing$, condition~\ref{defn:AE-in-measure:remainder} in Definition~\ref{defn:AE-in-measure} is vacuous. 
\end{proof}

\begin{rem}\label{OSnonfree}
	We remark that we cannot expect almost elementariness in measure to imply the stronger property of principality instead of effectiveness. Recall that for transformation groupoids, principality corresponds to the freeness of the action while effectiveness corresponds to topological freeness. For example, it was shown in \cite{OrtegaScarparo2020almostfinite} that all minimal actions on the Cantor set by the infinite dihedral group $\Z\rtimes \Z_2$ give rise to transformation groupoids that are almost finite in Matui's sense (and thus almost elementary in measure; see Theorem~\ref{thm:Matui-AE} and \ref{thm:AE-DSC} ), but some of these actions \LyXbar \LyXbar{} in particular, some odometer actions \LyXbar \LyXbar{} are not free. 
\end{rem}

Another consequence of almost elementariness in measure is a generalization of Lindenstrauss and Weiss' small boundary property (\cite{LindenstraussWeiss2000Mean}) to the groupoid setting. 
Let us start with a topological concept that is extracted from the dynamical context. 

\begin{defn}
	\label{7.7} 
	\label{defn:SBP-topological}
	Let $X$ be a compact Hausdorff space.  Let $\Omega$ be a weak{*}- closed subset of the set $M(X)$ of regular Borel probability measures on $X$.
	We say $X$ has the $\Omega$\emph{-small boundary property} if 
	for   
	any $x \in X$ and any open neighborhood $U$ of $x$, 
	there is an	open set $V$ such that $x\in V\subseteq\overline{V}\subseteq U$
	and $\mu(\del V) =0$ for any $\mu\in\Omega$. 
\end{defn}

\begin{defn}
	\label{7.6} 
	\label{defn:GSBP}
	An \'{e}tale groupoid $\CG$ with a compact unit space has 
	the \emph{groupoid small boundary property} (GSBP for short) 
	or simply the \emph{small boundary property}
	if $\CG\unitSpace$ has the $M(\CG)$-small boundary property. 
\end{defn}

\begin{rem}
	If $M(\CG) = \varnothing$, then $\CG$ has the groupoid small boundary property vacuously. 
\end{rem}

\begin{example}
	\label{exa:GSBP-transformation-groupoid}
	Let $\alpha \colon \Gamma \curvearrowright X$ be an action of a discrete group on a compact Hausdorff space. Then the transformation groupoid $X \rtimes \Gamma$ has the groupoid small boundary property if and only if $\alpha$ has the small boundary property, since we have $M(X \rtimes \Gamma) = M_\Gamma (X)$. 
\end{example}

Before investigating this notion further, let us first establish a few lemmas that are more or less known to experts. 

\begin{lem}
	\label{lem:measure-semicontinuous}
	Let $X$ be a compact Hausdorff space. 
	Then 
	\begin{enumerate}
		\item For any open subset $U$ in $X$, the function $M(X) \to [0,1]$, $\mu \mapsto \mu(U)$ is lower semicontinuous, i.e., $\left\{ \mu \in M(X) \colon \mu(U) > t \right\}$ is open for any $t \in \mathbb{R}$. 
		\item For any compact subset $K$ in $X$, the function $M(X) \to [0,1]$, $\mu \mapsto \mu(K)$ is upper semicontinuous, i.e., $\left\{ \mu \in M(X) \colon \mu(K) < t \right\}$ is open for any $t \in \mathbb{R}$.  
	\end{enumerate}
\end{lem}

\begin{proof}
	This follows from the facts that $\mu(U) = \sup \left\{ \int_X f \mathrm{d} \mu \colon f \in C(X)_+ , \,  f \leq \chi_U \right\}$ and $\mu(K) = \inf \left\{ \int_X f \mathrm{d} \mu \colon f \in C(X)_+ , \,  f \geq \chi_K \right\}$
\end{proof}

The following is a generalization 
to \cite[Lemma 9.1]{Kerr2020Dimension} and \cite[Lemma 3.2]{Ma2019Invariant}.

\begin{lem}
	\label{5.10} 
	\label{lem:small-neighborhood-in-measure}
	Let $X$ and $\Omega$ be as in Definition~\ref{defn:SBP-topological}. 
	\begin{enumerate}
		\item For any compact subset $K$ in $X$
		and any lower semicontinuous function $f \colon \Omega \to \mathbb{R}$ with $\mu(K) < f(\mu)$ for any $\mu \in \Omega$, 
		there is an open neighborhood $V$ of $K$ such that $\mu(\overline{V}) < f(\mu)$
		for any $\mu\in\Omega$. 
		
		\item For any open subset $U$ in $X$
		and any upper semicontinuous function $g \colon \Omega \to \mathbb{R}$ with $\mu(U) > g(\mu)$ for any $\mu \in \Omega$, 
		there is an open subset $W$ with $\overline{W} \subseteq U$ such that $\mu({W}) > g(\mu)$
		for any $\mu\in\Omega$. 
	\end{enumerate}
	
\end{lem}

\begin{proof}
	It suffices to prove the first statement, as it will then imply the second by taking $K := X \setminus U$ and $W := X \setminus \overline{V}$. 
	
	To this end, for any $\mu \in \Omega$, by the outer regularity of $\mu$, there is an open neighborhood $U_\mu$ of $K$ such that $\mu(U_\mu) < f(\mu)$. Since $X$ is normal, there is an open set $V_\mu$ such that $K \subseteq V_\mu$ and $\overline{V_\mu} \subseteq U_\mu$. 
	Since $f(\mu) - \mu \left( \overline{V_\mu} \right) \geq f(\mu) - \mu \left( {U_\mu} \right) > 0$, by Lemma~\ref{lem:measure-semicontinuous}, there is an open neighborhood $W_\mu$ of $\mu$ in $\Omega$ such that $f(\nu) - \nu \left( \overline{V_\mu} \right) > 0$ for any $\nu \in W_\mu$. Hence we have an open cover $\left\{ W_\mu \colon \mu \in \Omega \right\}$ of $\Omega$. By compactness, there is a finite subcover $\left\{ W_{\mu_i} \colon i \in \{1, \ldots, n\} \right\}$. It is clear that $V := \bigcap_{i=1}^n V_{\mu_i}$ satisfies $\mu(\overline{V}) \leq \min_{i \in \{1, \ldots, n\}} \mu \left( \overline{V_{\mu_i}} \right) < f(\mu)$
	for any $\mu\in\Omega$, as desired. 
\end{proof}

\begin{lem}\label{lem:SBP-equivalent}
	Let $X$ and $\Omega$ be as in Definition~\ref{defn:SBP-topological}. 
	The following are equivalent: 
	\begin{enumerate}
		\item\label{lem:SBP-equivalent:local-zero} The space $X$ has the $\Omega$-small boundary property. 
		
		\item\label{lem:SBP-equivalent:global-zero} For any compact set $K$ and any open set $U$ with $K \subseteq U \subseteq X$, 
		there is an	open set $V$ such that $K \subseteq  V\subseteq \overline{V} \subseteq U$
		and $\mu(\del V) = 0$ for any $\mu\in\Omega$. 
		
		\item\label{lem:SBP-equivalent:global-epsilon} For any $\varepsilon > 0$, any compact set $K$ and any open set $U$ with $K \subseteq U \subseteq X$, 
		there is an	open set $V$ such that $K \subseteq  V\subseteq \overline{V} \subseteq U$
		and $\mu(\del V) < \varepsilon$ for any $\mu\in\Omega$. 
		
		\item\label{lem:SBP-equivalent:global-thick} For any $\varepsilon > 0$, any compact set $K$ and any open set $U$ with $K \subseteq U \subseteq X$, 
		there are open sets $W$ and $W'$ such that $K \subseteq  W$, $\overline{W} \subseteq W'$, $\overline{W'} \subseteq U$
		and $\mu(\overline{W'} \setminus W) < \varepsilon$ for any $\mu\in\Omega$. 
		
		\item\label{lem:SBP-equivalent:cover-zero} For any open cover $\CU$ of $X$, there is a finite disjoint collection $\CW$ of open subsets of $X$ that refines $\CU$ and satisfies $\mu\left(X \setminus \bigcup \CW\right) = 0$  for any $\mu\in\Omega$. 
		
		\item\label{lem:SBP-equivalent:cover-epsilon} For any $\varepsilon > 0$ and any open cover $\CU$ of $X$, there is a disjoint collection $\CW$ of open subsets of $X$ that refines $\CU$ and satisfies $\mu\left(X \setminus \bigcup \CW \right) < \varepsilon$  for any $\mu\in\Omega$. 
	\end{enumerate}
\end{lem}

\begin{proof}
	We prove \eqref{lem:SBP-equivalent:local-zero} $\Rightarrow$ \eqref{lem:SBP-equivalent:global-zero} $\Rightarrow$ \eqref{lem:SBP-equivalent:global-epsilon} $\Rightarrow$ \eqref{lem:SBP-equivalent:global-thick} $\Rightarrow$ \eqref{lem:SBP-equivalent:local-zero} and \eqref{lem:SBP-equivalent:local-zero} $\Rightarrow$ \eqref{lem:SBP-equivalent:cover-zero} $\Rightarrow$ \eqref{lem:SBP-equivalent:cover-epsilon} $\Rightarrow$ \eqref{lem:SBP-equivalent:global-epsilon}. 
	Note that the implications \eqref{lem:SBP-equivalent:global-zero} $\Rightarrow$ \eqref{lem:SBP-equivalent:global-epsilon} and \eqref{lem:SBP-equivalent:cover-zero} $\Rightarrow$ \eqref{lem:SBP-equivalent:cover-epsilon} are trivial. 
	
	To prove \eqref{lem:SBP-equivalent:local-zero} $\Rightarrow$ \eqref{lem:SBP-equivalent:global-zero}, we fix $K$ and $U$ as in \eqref{lem:SBP-equivalent:global-zero}. For any $x \in K$, applying \eqref{lem:SBP-equivalent:local-zero} gives us an open neighborhood $V_x$ of $x$ satisfying $\overline{V_x} \subseteq U$ and $\mu\left( \del V_x \right) = 0$ for any $\mu \in \Omega$. 
	By compactness, we obtain a finite open cover $\left\{ V_{x_i} \colon i \in \{1, \ldots, n\} \right\}$. It is clear that $V := \bigcup_{i=1}^n V_{x_i}$ satisfies the requirement. 
	
	The implication \eqref{lem:SBP-equivalent:global-epsilon} $\Rightarrow$ \eqref{lem:SBP-equivalent:global-thick} follows from applying Lemma~\ref{lem:small-neighborhood-in-measure} to the boundary $\del V$ obtained in \eqref{lem:SBP-equivalent:global-epsilon} and the constant function $\varepsilon$ to obtain an open neighborhood $Y$ of $\del V$, where we may assume without loss of generality that $\overline{Y} \subseteq U \setminus K$, and then defining $W := V \setminus \overline{Y}$ and $W' := V \cup Y$. 
	
	To prove \eqref{lem:SBP-equivalent:global-thick} $\Rightarrow$ \eqref{lem:SBP-equivalent:local-zero}, we fix $x$ and $U$ as in Definition~\ref{defn:SBP-topological} and recursively apply \eqref{lem:SBP-equivalent:global-thick} to $\varepsilon_n$, $K_n$, $U_n$ to obtain $W_n$ and $W'_n$, where $\varepsilon_n := \frac{1}{2^n}$ and 
	\[
		K_n := 
		\begin{cases}
			\{x\} \, , & n = 0 \\
			\overline{W_{n-1}} \, , & n > 0
		\end{cases}
		\; , 
		U_n := 
		\begin{cases}
			U \, , & n = 0 \\
			{W'_{n-1}} \, , & n > 0
		\end{cases}
		\; .
	\]
	Hence we obtain an increasing sequence $W_0 \subseteq W_1 \subseteq \ldots$ and a decreasing sequence $W'_0 \supseteq W'_1 \supseteq \ldots$ such that $\mu \left( \overline{W'_n} \setminus W_n \right) < \frac{1}{2^n}$ for any $n \in \mathbb{N}$. The open set $V := \bigcup_{n = 0}^{\infty} W_n$ thus satisfies the requirements. 
	
	To prove \eqref{lem:SBP-equivalent:local-zero} $\Rightarrow$ \eqref{lem:SBP-equivalent:cover-zero}, we fix an open cover $\CU$ as in \eqref{lem:SBP-equivalent:cover-zero}. For any $x \in X$, there is an open neighborhood $U_x \in \CU$ of $x$, and thus, by \eqref{lem:SBP-equivalent:local-zero}, another open neighborhood $V_x$ of $x$ such that $\overline{V_x} \subseteq U_x$ and $\mu \left(\del V_x \right) = 0$ for any $\mu \in \Omega$. By the compactness of $X$, there is a finite open cover $\left\{ V_{x_i} \colon i = 1, \ldots, n \right\}$. It is routine to verify that the collection $\CW := \left\{ W_i \colon i = 1, \ldots, n \right\}$ satisfies the requirement, where $W_i := V_{x_i} \setminus \left( \bigcup_{k=1}^{i-1} \overline{V_{x_{k}}} \right)$, 
	since $\bigcup \CV \supseteq X \setminus \left( \bigcup_{i=1}^{n} \del V_{x_i} \right)$. 
	
	To prove \eqref{lem:SBP-equivalent:cover-epsilon} $\Rightarrow$ \eqref{lem:SBP-equivalent:global-epsilon}, 
	we fix $\varepsilon$, $K$ and $U$ as in \eqref{lem:SBP-equivalent:global-epsilon} and apply normality of $X$ to obtain open sets $U_0, V_0$ such that $K \subseteq V_0 \subseteq \overline{V_0} \subseteq U_0 \subseteq \overline{U_0} \subseteq U$.  
	Define $\CU := \left\{ U_0, X \setminus \overline{V_0} \right\}$, and obtain $\CW$ as in \eqref{lem:SBP-equivalent:cover-epsilon}. Without loss of generality, we may assume $\CW$ is finite, since we have $\displaystyle \mu \left(X \setminus \bigcup \CW \right) = \inf_{\CF \subseteq \CW \text{ finite}} \mu\left(X \setminus \bigcup \CF \right)$. We observe that $V := V_0 \cup \left( \bigcup \left\{ W \in \CW \colon W \cap \overline{V_0} \not= \varnothing \right\} \right)$ satisfies the requirement, since we have $\del V \subseteq X \setminus \left( \bigcup \CW \right)$. 
\end{proof}

The following is a partial generalization of \cite[Theorem~5.6]{KerrSzabo2020Almost}. 

\begin{thm}
	\label{7.9} 
	\label{thm:GSBP}
	If an \'{e}tale groupoid $\CG$ with a compact unit space is almost elementary in measure, 
	then $\CG$ has the groupoid small boundary property. 
\end{thm}

\begin{proof}
	It suffices to verify Lemma~\ref{lem:SBP-equivalent}\eqref{lem:SBP-equivalent:cover-epsilon}. To this end, we fix $\varepsilon >0$ and an open cover $\CU$ of $\CG\unitSpace$. By Definition~\ref{defn:AE-in-measure}, there is a nonempty open castle $\CC$ satisfying 
	\begin{enumerate}[label=(\roman*)]
		\item $\CC$ is $\CG\unitSpace$-extendable; 
		\item $\CC\unitSpace$ refines $\CU$;
		\item $\mu \left( \CG^{(0)}\setminus\ \left( \bigcup \CC^{(0)} \right) \right) \leq \varepsilon$ for any $\mu \in M(\CG)$.
	\end{enumerate}
	Then $\CW := \CC\unitSpace$ satisfies the requirement in Lemma~\ref{lem:SBP-equivalent}\eqref{lem:SBP-equivalent:cover-epsilon}. 
\end{proof}

\begin{rem}
	It was in fact shown in \cite[Theorem~5.6]{KerrSzabo2020Almost} that the converse to Theorem~\ref{thm:GSBP} also holds for principal transformation groupoids, 
	while it follows from Theorem~\ref{thm:effective} that this converse also holds for effective groupoids with $M(\CG) = \varnothing$. We do not know to what extent this converse holds among effective \'{e}tale groupoids with compact unit spaces.  
\end{rem}

Lastly, recall that it was shown in \cite[Theorem B]{MWfiberwise} that fiberwise amenability implies that $M(\CG)$ is not empty. Now we show that among almost elementary groupoids in measure, the converse also holds.

\begin{thm}
	Suppose a nonempty $\sigma$-compact \'{e}tale groupoid $\CG$ with a compact unit space is almost elementary in measure. Then $\CG$ is fiberwise amenable if and only if $M(\CG) \not= \varnothing$. 
\end{thm}

\begin{proof}
	The ``only if'' direction is precisely \cite[Theorem B]{MWfiberwise}. For the ``if'' direction, we assume $\CG$ is not fiberwise amenable and show $M(\CG) = \varnothing$. To this end, we apply Theorem~\ref{thm:fiberwise-amenable-dichotomy}\ref{thm:fiberwise-amenable-dichotomy:non-amenable} (with $L = \CG\unitSpace$ and $n = 2$) to obtain a compact set $K \subseteq \CG$ such that for any compact set $M \subseteq \CG$ and any $u \in \CG\unitSpace$, the set $KMu$ contains at least $2|Mu|$ many elements. 
	Applying Lemma~\ref{lem:AE-in-measure-shrinking}\eqref{lem:AE-in-measure-shrinking:shrinking} to the nonempty compact set $K$, $\varepsilon = \frac{1}{3}$ and the open cover $\CV = \left\{ \CG\unitSpace \right\}$, 
	we obtain 
	an open castle $\CD$ satisfying $\mu \left( \bigcup \left( \CD^{-K} \right)^{(0)} \right) \geq \frac{2}{3}$ for any $\mu \in M(\CG)$.

	We claim that for any $\mu \in M(\CG)$, we have 
	\[
		\mu \left( \bigcup \CD^{(0)} \right) \geq 2 \mu \left( \bigcup \left( \CD^{-K} \right)^{(0)} \right) \quad \text{ for any } \mu \in M(\CG) \; ,
	\]
	which will imply $M(\CG) = \varnothing$ since it implies $\mu \left( \bigcup \CD^{(0)} \right) 
	\geq \frac{4}{3}$, 
	which would be absurd. 
	To prove this claim, we first observe that 
	by our choice of $K$, for any $u \in \bigcup \left(\CD^{-K}\right)\unitSpace$, since $\left( \bigcup \CD^{-K} \right) u = \bigcup_{D \in \CD^{-K}} D u$ is a finite set, we have, 
	\[
	\left| \left( \bigcup \CD^{-K} \right) u \right| 
	\leq \frac{1}{2} \left| K {\left( \bigcup \CD^{-K} \right)} u \right| 
	\leq \frac{1}{2} \left| \left( \bigcup \CD \right) u \right| \; .
	\]
	By Remark~\ref{rem:fortified-castle-basics}\eqref{rem:fortified-castle-basics:tower}, Lemma~\ref{lem:castle-basic}\eqref{lem:castle-basic:height} and Lemma~\ref{lem:castle-shrinking}\eqref{lem:castle-shrinking:tower}, for any tower $\CT$ in $\CD$ and any $u \in \bigcup \left(\CT^{-K}\right)\unitSpace $, we have 
	\[
	\left|   \CT  \unitSpace \right|
	= \left| \left( \bigcup \CD \right) u \right|
	\geq 2 \left| \left( \bigcup \CD^{-K} \right) u \right| 
	= 2 \left|  \left(\CT^{-K}\right)  \unitSpace \right| 
	\; .
	\]
	Choose a set $\left\{ D_{\CT} \in  \CT \unitSpace \colon \CT \in \CD / \sim_{\CD} \right\}$ of representatives for the relation $\sim_{\CD}$ restricted to $\CD\unitSpace$. 	
	By Lemma~\ref{lem:castle-shrinking}\eqref{lem:castle-shrinking:tower}, we may assume without loss of generality that $D_{\CT} \in \left( \CT^{-K} \right)\unitSpace$ for any $\CT \in \CD / \sim_{\CD}$ with $\CT^{-K} \not= \varnothing$. 
	By Remark~\ref{rem:castle-measure}\eqref{rem:castle-measure:castle}, 
	we have 
	\begin{align*}
	\mu \left( \bigcup \CD^{(0)} \right) 
	&= \sum_{ \CT \in \CD / \sim_{\CD} }  \left| \CT\unitSpace \right| \mu  \left(D_{\CT} \right)  \\
	&\geq 2 \sum_{ \CT \in \CD / \sim_{\CD} }  \left| \left(\CT^{-K}\right)  \unitSpace \right| \mu  \left(D_{\CT} \right) 
	= 2 \mu \left( \bigcup \left( \CD^{-K} \right)^{(0)} \right)
	\end{align*}
	as desired. 
\end{proof}

To end this section, we remark that the almost elementary in measure has been applied by the first author in \cite{Ma-fiberwise} to establish the uniform property $\Gamma$ for reduced ample groupoid $C^*$-algebras

\section{Almost elementary \'{e}tale groupoids} \label{sec:almost-elementary}

In this section, we introduce one of the main concepts in the paper, namely almost elementariness for \'{e}tale groupoids, considered as a regularity property for \'{e}tale groupoids (with compact unit spaces). It differs from \textemdash\ and in fact strengthens (see Theorem~\ref{thm:AE-DSC}) \textemdash\ almost elementariness in measure in the following way: 
In the latter notion, the remainder control (see condition~\ref{defn:AE-in-measure:remainder} in Definition~\ref{defn:AE-in-measure}) is ``measure-theoretic'' \textemdash\ we control the remainder of a castle by weighing it under all invariant measures. 
In contrast, in almost elementariness, this is replaced by a finer, ``topological'' control, which is called for in a number of applications. 
This is implemented by the following notions of subequivalence in an \'{e}tale groupoid.

\begin{defn}
	\label{defn:subequivalence}
	Let $\CG$ be an \'{e}tale groupoid. Let $K$ be a compact subset of $\CG^{(0)}$ and let $U,V \subseteq \CG^{(0)}$. 
	\begin{enumerate}[label=(\roman*)]
		\item We write $K\prec_{\CG}V$ if there is an open $s$-section $S$ such that $K \subseteq r(S)$ and $s(S) \subseteq V$. 
		\item We write $U\precsim_{\CG}V$
		if $K\prec_{\CG}V$ for every compact subset $K\subseteq U$. 
	\end{enumerate}
\end{defn}

We establish a number of basic properties. 

\begin{rem}
	\label{rem:subequivalence-basic}
	Let $K$, $U$, and $V$ be as in Definition~\ref{defn:subequivalence}. 
	\begin{enumerate}
		\item \label{rem:subequivalence-basic:subset} If $K \subseteq V$, then $K\prec_{\CG}V$. Hence if $U \subseteq V$, then $U\precsim_{\CG}V$. 
		
		\item \label{rem:subequivalence-basic:clopen} 
		When $U$ is compact, we have $U \precsim_{\CG} V$ if and only if $U \prec_{\CG} V$. 

        \item \label{rem:subequivalence-basic: compact to open} From the definition, if $K\prec_\CG V$. Then there exists an open set $U\supseteq K$ such that $U \precsim_\CG V$. Indeed, one obtains this property by defining the open set $U$ by $U=r(S)$ where $S$ is the open $s$-section that witnesses $K\prec_\CG V$.
		
		\item \label{rem:subequivalence-basic:interior} It is clear that $K\prec_{\CG}V$ if and only if $K\prec_{\CG} V^{\operatorname{o}}$, and thus $U\precsim_{\CG}V$ if and only if $U\precsim_{\CG} V^{\operatorname{o}}$. 

        \item \label{rem:subequivalence-basic:precompact-section} We also have $K\prec_{\CG}V$ if and only if there is a compact $s$-section $T$ such that $K = r(T)$ and $s(T) \subseteq V^{\operatorname{o}}$. 
Indeed, for the ``if'' direction, it follows from the definition of an $s$-section on page~\pageref{paragraph:sections-in-groupoids} that there is an open $s$-section $T' \supseteq T$. Then the open  $s$-section $S := T' V^{\operatorname{o}} = (s|_{T'})^{-1}(s(T') \cap V^{\operatorname{o}}) \subseteq T'$ satisfies $s(S) = s(T') \cap V ^{\operatorname{o}}\supseteq s(T)$, whence $S \supseteq T$ and thus $r(S) \supseteq r(T) = K$. 
        For the ``only if'' direction, assuming, by \eqref{rem:subequivalence-basic:interior}, there is an open $s$-section $S$ such that $K \subseteq r(S)$ and $s(S) \subseteq V^{\operatorname{o}}$, we use the normality of $\CG\unitSpace$ and the fact that $s|_S \colon S \to s(S)$ is a homeomorphism to write $S$ as the union of all the precompact open subsets whose closures are in $S$, and then use the compactness of $K$ and the fact that $r$ is an open map to find some precompact open subset $S'$ such that $\overline{S'} \subseteq S$ and $K \subseteq r(S')$, whence $T := (r_{\overline{S'}})^{-1}(K)$ is a compact subset in $S$ with $r(T) = K$ and $s(T) \subseteq s(S) \subseteq V^{\operatorname{o}}$.

		\item \label{rem:subequivalence-basic:compact} It follows from \eqref{rem:subequivalence-basic:precompact-section} that $K\prec_{\CG}V$ if and only if $K\prec_{\CG} K'$ for some compact subset $K' \subseteq V^{\operatorname{o}}$. 
        Indeed, for the nontrivial ``only if'' direction, we apply \eqref{rem:subequivalence-basic:interior} and \eqref{rem:subequivalence-basic:precompact-section} to obtain a compact $s$-section $T$ such that $K = r(T)$ and $s(T) \subseteq V^{\operatorname{o}}$, and then use normality and local compactness to obtain some $K' \subseteq V^{\operatorname{o}}$ with $s(T) \subseteq (K')^{\operatorname{o}}$. Applying \eqref{rem:subequivalence-basic:precompact-section} again, we see $K\prec_{\CG} K'$.

		\item \label{rem:subequivalence-basic:auxiliary} It follows from \eqref{rem:subequivalence-basic:compact} that if $K\prec_{\CG}U$ and $U\precsim_{\CG}V$, then $K\prec_{\CG}V$. 
		
		\item \label{rem:subequivalence-basic:preorder}  
		It follows from \eqref{rem:subequivalence-basic:subset} and \eqref{rem:subequivalence-basic:auxiliary} that $\precsim_{\CG}$ is reflexive and transitive, i.e., a preorder. 
		
		\item \label{rem:subequivalence-basic:additive-prec}
		It is immediate from the definition that for any compact subset $L$ of $\CG\unitSpace$, if $K\prec_{\CG}U$, $L\prec_{\CG}V$, and $U \cap V = \varnothing$, then we have $K \cup L \prec_{\CG} U \sqcup V$.

		\item \label{rem:subequivalence-basic:additive-precsim}
		For any open subsets $U', V' \subseteq \CG\unitSpace$, if $U'\precsim_{\CG}U$, $V'\precsim_{\CG}V$, and $U \cap V = \varnothing$, then we have $U' \cup V' \precsim_{\CG} U \sqcup V$. Indeed, for any compact subset $L \subseteq U' \cup V'$, we may exploit the local compactness of $\CG\unitSpace$ to obtain precompact open subsets $U'', V'' \subseteq \CG\unitSpace$ such that $\overline{U''} \subseteq U'$, $\overline{V''} \subseteq V'$, and $L \subseteq U'' \cup V''$, whence we have $\overline{U''} \prec_{\CG} U$ and $\overline{V''} \prec_{\CG} V$ by \eqref{rem:subequivalence-basic:subset} and \eqref{rem:subequivalence-basic:auxiliary}, 
		$\overline{U''} \cup \overline{V''} \prec_{\CG} U \sqcup V$ by \eqref{rem:subequivalence-basic:additive-prec}, and thus  $L \prec_{\CG} U \sqcup V$ by \eqref{rem:subequivalence-basic:subset}, \eqref{rem:subequivalence-basic:clopen} and \eqref{rem:subequivalence-basic:auxiliary} again, as desired. 

		\item \label{rem:subequivalence-basic:additive-prec-precsim}
        For any open subset $W \subseteq \CG\unitSpace$, if $K \setminus W \prec_\CG U$, $W \precsim_\CG V$, and $U \cap V = \varnothing$, then we have $K \prec_\CG U \sqcup V$. Indeed, by \eqref{rem:subequivalence-basic: compact to open}, there is an open subset $U' \supseteq K$ such that $U' \precsim_\CG U$. It then follows from \eqref{rem:subequivalence-basic:additive-prec} that $K = (K \setminus W) \sqcup W \subseteq U' \cup W \precsim_\CG U \sqcup V$, whence $K \prec_\CG U \sqcup V$ by  \eqref{rem:subequivalence-basic:subset}, \eqref{rem:subequivalence-basic:clopen} and \eqref{rem:subequivalence-basic:auxiliary}. 
		
		\item \label{rem:subequivalence-basic:subgroupoid}
		Let $\CH$ be an open subgroupoid of $\CG$. 
		It follows from the definition that if $K, V \subseteq \CH\unitSpace$ and $K \prec_{\CH} V$, then $K \prec_{\CG} V$. 
		Hence, furthermore, if $U, V \subseteq \CH\unitSpace$ and $U \precsim_{\CH} V$, then $U \prec_{\CG} V$.

		\item \label{rem:subequivalence-basic:castle}
		Let $\CC$ be an open castle in $\CG$ and let $\CX , \CY \subseteq \CC\unitSpace$. 
Then we have $ \CX \precsim_{ \CC}  \CY$ in the finite groupoid $\CC$ (as in \Cref{lem:castle-basic}\eqref{lem:castle-basic:subquotient}) if and only if $|\CX \cap \CT\unitSpace| \leq |\CY \cap \CT\unitSpace|$ for each tower $\CT$ in $\CC$. Note that  $ \CX \precsim_{ \CC}  \CY$  implies $\bigcup \CX \precsim_\CG \bigcup \CY$.  
		
		\item \label{rem:subequivalence-basic:measure-auxiliary}
		If $K\prec_{\CG}V$, then we have $\mu(K)\leq\mu(V)$ for any $\mu\in M(\CG)$. Indeed, assuming there is an open $s$-section $S$ such that $K \subseteq r(S)$ and $s(S) \subseteq V$, we may again use the compactness of $K$ and the fact that $\CG$ is \'{e}tale to find open bisections $S_1, \ldots, S_n \subseteq S$ such that $K \subseteq \bigcup_{i=1}^n r \left( S_i \right)$, 
		which give rise to the mutually disjoint Borel measurable bisections $T_i := S_i \setminus \left( \bigcup_{j = 1}^i S_j \right)$ for $i \in \{1,\ldots, n\}$, 
		whence for any $\mu\in M(\CG)$, we have $\displaystyle \mu(K) \leq \mu \left( \bigcup_{i=1}^n r \left( S_i \right) \right) = \mu \left( \bigcup_{i=1}^n r \left( T_i \right) \right) \leq \sum_{i=1}^n \mu \left( r \left( T_i \right) \right) = \sum_{i=1}^n \mu \left( s \left( T_i \right) \right) = \mu \left( \bigcup_{i=1}^n s \left( T_i \right) \right) \leq \mu(s(S)) \leq \mu(V)$. 
		
		\item \label{rem:subequivalence-basic:measure-subequivalence}
		It follows from \eqref{rem:subequivalence-basic:measure-auxiliary} and the regularity of measures that if $U \precsim_{\CG}V$, then we have $\mu(U)\leq\mu(V)$ for any $\mu\in M(\CG)$. 

	\end{enumerate}
\end{rem}

As \eqref{rem:subequivalence-basic:measure-auxiliary} and \eqref{rem:subequivalence-basic:measure-subequivalence} above indicate, the subequivalences $\prec_{\CG}$ and $\precsim_{\CG}$ provide ``topological'' reasons for one set to weigh less than another under every invariant measure. 

\begin{example}
	\label{exa:subequivalence-actions}
	Let $X \rtimes_\alpha \Gamma$ be the transformation groupoid associated to an action $\alpha$ of a discrete group $\Gamma$ on a compact space $X$ (see Example~\ref{exa:transformation-groupoid} for notations). 
	Let $K$ be a compact subset of $X$ and let $U,V \subseteq X$. 
	We write $\prec_{\alpha}$ instead of $\prec_{X\rtimes_{\alpha}\Gamma}$, and $\precsim_{\alpha}$ instead of $\precsim_{X\rtimes_{\alpha}\Gamma}$, for the sake of brevity. 
	
	Let us work out what the relations $\prec_{\alpha}$ and $\precsim_{\alpha}$ mean explicitly in terms of the action. 
	\begin{enumerate}
		\item \label{exa:subequivalence-actions:auxiliary} We have $K\prec_{\alpha}V$ if and only if there exist mutually disjoint open subsets $W_1, \ldots, W_n$ of $V^{\operatorname{o}}$ and $\gamma_1 , \ldots, \gamma_n \in \Gamma$ such that $K \subseteq \bigcup_{i=1}^n \alpha_{\gamma_i} \left( W_i \right)$. 
		Indeed, for the ``if'' direction, the set 
		\[
		S := \bigcup_{i=1}^n  \left\{ \left(\alpha_{\gamma_i} (x) , \gamma_i, x \right) \in X \rtimes_\alpha \Gamma \colon x \in W_i \right\}
		\]
		will fulfill the role in Definition~\ref{defn:subequivalence}. For the ``only if'' direction, we may argue as in Remark~\ref{rem:subequivalence-basic}\eqref{rem:subequivalence-basic:compact} that the $S$ in Definition~\ref{defn:subequivalence} may be taken to be precompact, and thus have a finite image under the projection $X \rtimes_\alpha \Gamma \to \Gamma$, which allows us to write $S$ as a finite union as above, where $W_i := \left\{ x \in X \colon  \left(\alpha_{\gamma_i} (x) , \gamma_i, x \right) \in S \right\}$, and verify the desired equivalent characterization with these $W_1, \ldots, W_n$ and $\gamma_1 , \ldots, \gamma_n$. 
		
		\item \label{exa:subequivalence-actions:subequivalence} It follows from \eqref{exa:subequivalence-actions:auxiliary} that the preorder $\precsim_{\alpha}$ coincides with the relation $\prec_0$ in \cite[Definition~3.1]{Kerr2020Dimension}. 
	\end{enumerate}
\end{example}

\begin{rem}\label{rem:subequivalence-change}
    Our definitions of $\prec_{\CG}$ and $\precsim_{\CG}$ in \Cref{defn:subequivalence} are somewhat weaker than those given in \cite[Definition~3.2]{Ma-purely}, since the latter additionally requires, as phrased in the context of \Cref{defn:subequivalence}, that the open $s$-section $S$ is a finite disjoint union of open bisections. We record the following observations regarding this discrepancy: 
    \begin{enumerate}
        \item \label{rem:subequivalence-change:strong_enough} Our weaker relations are still strong enough for all of our purposes. Most crucially, for any open subsets $U,V \subseteq \GU$, $U \precsim_{\CG} V$ still implies $\mu(U) \leq \mu(V)$ for any $\mu \in M(\CG)$ (see \Cref{rem:subequivalence-basic}\eqref{rem:subequivalence-basic:measure-subequivalence}), and for any $f,g \in C_0(\GU)$, $\supp^o(f) \precsim_{\CG} \supp^o (g)$ still implies $f \precsim g$ in $C^*(\CG)$ (see \Cref{lem:groupoid subequi to Cuntz subequi}). 
        \item \label{rem:subequivalence-change:comparison} Consequently, groupoid strict comparison, when phrased in these weaker relations (see \Cref{defn:comparison}), still makes sense and is more likely to hold. 
        \item \label{rem:subequivalence-change:purely_infinite} Thanks to \eqref{rem:subequivalence-change:strong_enough} and the fact that these weaker relations satisfy the same basic permanence properties as the original ones (see \Cref{rem:subequivalence-basic}), results in \cite{Ma-purely} still hold when we use these weaker relations in place of the original. 
        \item \label{rem:subequivalence-change:same} Whenever an \'{e}tale groupoid $\CG$ admits a partition into clopen bisections (this is the case for transformation groupoids and ample groupoids), our relations agree with the original, since restricting said partition to any precompact open $s$-section $S$ identifies the latter as a finite disjoint union of open bisections. 
        \item \label{rem:subequivalence-change:example} On the other hand, we present an example below in \Cref{example:subequivalence-change} where our new relations are strictly weaker than the original. We do not know if such examples exist among \emph{minimal} groupoids. 
    \end{enumerate}
\end{rem}

\begin{example}\label{example:subequivalence-change}
We fulfill the claim in \Cref{rem:subequivalence-change}\eqref{rem:subequivalence-change:example}. 
    Let $X := S^2$ and let $\alpha \colon \mathbb{Z}/2 \curvearrowright X$ be the antipodal action. Let $Y$ be the quotient space of this action (so $Y \cong \mathbb{R}P^2$), and write $Z := X \sqcup Y$. Consider the compact Hausdorff space $\mathcal{G} := (X \times \{0,1, \uparrow, \downarrow\} ) \sqcup Y$ and identify $Z$ canonically with the subspace $(X \times \{0\} ) \sqcup Y$ of $\mathcal{G}$. Equip $\mathcal{G}$ with the structure of a principal groupoid, with $Z$ being the unit space, by embedding $\CG$ into the (non-\'{e}tale) pair groupoid $Z \times Z$ via mapping each $[x] \in Y$ to $([x], [x]) \in Y \times Y \subseteq Z \times Z$, and mapping 
    \[
    X \times \{0,1, \uparrow, \downarrow\} \ni (x , j)  \longmapsto
        \begin{cases}
            (\alpha_i(x),x) \in X\times X \subseteq Z \times Z \, ,  & j \in \{0,1\} \\
            (x,[x]) \in X\times Y \subseteq Z \times Z \, ,  & j = \uparrow \\
            ([x], x) \in Y\times X \subseteq Z \times Z \, , & j \mathrel{=} \downarrow
        \end{cases}
        \; .
    \]
    One readily verifies that $\CG$ is a principal \'{e}tale compact Hausdorff groupoid. 

    Viewing $X$ and $Y$ as clopen subsets of the unit space $Z$, we immediately have $Y \prec_{\CG} X$ and $Y \precsim_{\CG} X$, as witnessed by the clopen $s$-section $X \times \{\downarrow\}$. 

    On the other hand, we can show that  $Y \prec_{\CG} X$ (and thus also  $Y \precsim_{\CG} X$ ) does not hold in the sense of \cite[Definition~3.2]{Ma-purely}, that is, one cannot find a finite collection $\{ B_i \colon i \in I\}$ of open bisections such that $Y \subseteq \bigcup_{i \in I} r(B_i)$ and $\bigsqcup_{i \in I} s(B_i) \subseteq X$. Indeed, suppose such a finite collection $\{ B_i \colon i \in I\}$ exists. Without loss of generality, we may assume $r(B_i) \subseteq Y$ for each $i \in I$ (otherwise we may replace $B_i$ by the open subset $Y B_i$). For each $i \in I$, since $s(B_i) \subseteq X$, a look at the groupoid structure of $\CG$ tells us that $B_i \subseteq X \times \{\downarrow\}$. 
    Note that restricting $s$ to $X \times \{\downarrow\}$ yields a homeomorphism $s \colon X \times \{\downarrow\} \to X$, under which the restriction of $r$ to $X \times \{\downarrow\}$ can be identified with the quotient map $X \to Y$; in particular, for each $[x] \in Y$, we have $(r|_{X \times \{\downarrow\}})^{-1} ([x]) = \{x, \alpha_1(x)\} \times \{\downarrow\}$. For each $i \in I$, since $r$ restricts to an injection on $B_i$, it follows that $B_i \cap \alpha_1 (B_i) = \varnothing$. Now a form of the Borsuk-Ulam theorem implies that there exists $[x] \in Y$ and distinct $i_1, i_2, i_3 \in I$ such that 
    $[x] \in r(B_{i_1}) \cap r(B_{i_2}) \cap r(B_{i_3})$
    (for a direct proof, we argue that if this were not the case, then after choosing a partition of unity $\{f_i \colon i \in I\}$ of $Y$ subordinate to the open cover $\{r(B_i) \colon i \in I \}$, we would be able to construct, in contradiction to the Borsuk-Ulam theorem, an odd continuous map $X = S^2 \to \mathbb{R}^2 \setminus\{0\}$ by mapping each $x \in X$ to 
    \[
        \big(\mathrm{sgn}_1(x) \, (2 f_{\max}(x) - 1), \;  \mathrm{sgn}_2(x) \, (1- f_{\max}(x)) \big)
    \]
    where, noting that our hypothesis implies there is at most two distinct $i,j \in I$ with nonzero $f_i([x])$ and $f_j([x])$, and thus at most one $i \in I$ with $f_i([x]) > \frac{1}{2}$, we fix a total order of $I$ and define 
    \begin{align*}
        \mathrm{sgn}_1(x) &:= \begin{cases}
            1 \, , & \text{ if there is $i \in I$ such that $f_i([x]) > \frac{1}{2}$ and $x \in s(B_i)$,} \\  
            -1 \, , & \text{ if there is $i \in I$ such that $f_i([x]) > \frac{1}{2}$ and $x \in \alpha_1 (s(B_i))$,} \\
            0 \, , & \text{ otherwise,}
        \end{cases}  \\
        \mathrm{sgn}_2(x) &:= 
        \begin{cases}
            1 \, , & \text{ if there are distinct $i,j \in I$ such that $\min \{ f_i([x]), f_j([x])\} > 0$ } \\  
            & \quad \text{ with $i < j$ and $x \in s(B_i) $,} \\
            -1 \, , & \text{ if there are distinct $i,j \in I$ such that $\min \{ f_i([x]), f_j([x])\} > 0$} \\  
            & \quad \text{ with $i < j$ and $x \in \alpha_1(s(B_i)) $,} \\
            0 \, , & \text{ otherwise,}
        \end{cases} 
    \end{align*}
   and $f_{\max}(x) := \max\{f_i([x]) \colon i \in I \}$). 
    However, since $(r|_{X \times \{\downarrow\}})^{-1} ([x])$ consists of two elements, the pigeonhole principle tells us that the collection $\{B_{i_1}, B_{i_2}, B_{i_3}\}$ (and thus also $\{s(B_{i_1}), s(B_{i_2}), s(B_{i_3})\}$) cannot be mutually disjoint, contradicting our assumption. 

    With the same proof method, this construction can be generalized by taking $\alpha \colon G \curvearrowright X$ to be any topological action of a finite group on a compact Hausdorff space such that the $G$-index $\mathrm{ind}_{G}(X) \geq |G|$ (see, e.g., \cite[Definition~5.3.1]{matousek_using_2008} for the definition of the $G$-index). 
\end{example}

Now we are ready to introduce almost elementariness. 
Since our main focus is on minimal \'{e}tale groupoids, we will restrict ourselves to this setting in this section. 
This has the advantage that our first definition below takes a simpler form. 
In the next section, we will give an equivalent definition that also works in the non-minimal case. 

\begin{defn} \label{defn:AE} 
	We say that a minimal \'{e}tale groupoid  $\CG$ with a compact unit space is \emph{almost elementary}
	if for any compact set $K$ in $\CG$,  
	any nonempty open subset $O$ in $\CG\unitSpace$, 
	and any finite open cover $\CV$ of $\CG^{(0)}$, 
	there is an open castle $\CC$
	satisfying 
	\begin{enumerate}[label=(\roman*)]
		\item \label{defn:AE:extendable} $\CC$ is $K$-extendable; 
		\item \label{defn:AE:fine} $\CC\unitSpace$ refines $\CV$; 
		\item \label{defn:AE:nonempty} $\CC$ is nonempty;
		\item \label{defn:AE:remainder} $\CG^{(0)}\setminus\ \left( \bigcup \CC^{(0)} \right) \prec_{\CG} O$.  
	\end{enumerate}
	In this case, the castle $\CC$ is sometimes referred to as an \emph{almost elementary approximation of $\CG$} and the set $\CG^{(0)}\setminus\ \left( \bigcup \CC^{(0)} \right)$ is the \emph{remainder} of this approximation. 
\end{defn}

As pointed out in Remark~\ref{rem:AE-in-measure-fullness}, since we restrict ourselves to the minimal case, the nonemptiness requirement in condition~\ref{defn:AE:nonempty} above is equivalent to the fullness requirement as in condition~\ref{defn:AE-in-measure:full} in Definition~\ref{defn:AE-in-measure}. 

As in the case of almost elementariness in measure, it is often convenient to use the notion of $K$-shrinking directly instead of that of $K$-extendability.

\begin{lem} \label{lem:AE-shrinking} 
	Let $\CG$ be a minimal \'{e}tale groupoid with a compact unit space. 
	The following are equivalent: 
	\begin{enumerate}
		\item \label{lem:AE-shrinking:original} The groupoid $\CG$ is {almost elementary}. 
		
		\item \label{lem:AE-shrinking:shrinking} For any compact set $K$ in $\CG$, 
		any non-empty open set $O$ in $\GU$, 
		and any finite open cover $\CV$ of $\CG^{(0)}$, 
		there is an open castle $\CD$
		satisfying 
		\begin{enumerate}[label=(\roman*)]
			\item \label{lem:AE-shrinking:refine} $\CD \unitSpace$ refines $\CV$; 
			\item \label{lem:AE-shrinking:full} $\CD^{-K}$ is nonempty; 
			\item \label{lem:AE-shrinking:remainder} $\CG^{(0)}\setminus\ \left( \bigcup \left( \CD^{-K} \right)^{(0)} \right) \prec_{\CG} O$. 
		\end{enumerate}
		
		\item \label{lem:AE-shrinking:fine-entourage} For any compact sets $K, L \subseteq \CG$, 
		any non-empty open set $O$ in $\GU$, 
		and any collection $\CV$ of open subsets in $\CG$ that covers $L$, 
		there is an open castle $\CD$
		satisfying 
		\begin{enumerate}[label=(\roman*)]
			\item \label{lem:AE-shrinking:fine-entourage:refine} for any $D \in \CD$ with $D \cap L \not= \varnothing$, there is $V \in \CV$ such that $D \subseteq V$; 
			\item \label{lem:AE-shrinking:fine-entourage:full} $\CD^{-K}$ is nonempty; 
			\item \label{lem:AE-shrinking:fine-entourage:remainder} $\CG^{(0)}\setminus\ \left( \bigcup \left( \CD^{-K} \right)^{(0)} \right) \prec_{\CG} O$. 
		\end{enumerate}	
	\end{enumerate}
\end{lem}

\begin{proof}
	We follow the proof of Lemma~\ref{lem:AE-in-measure-shrinking} almost verbatim, with the only significant modification being to replace equation~\eqref{lem:AE-in-measure-shrinking:eq:remainder-control} with 
	\[
		\GU \setminus \left( \bigcup \left( \CD^{-K} \right)^{(0)} \right) \subseteq \GU \setminus \left( \bigcup \CC^{(0)} \right) \prec_\CG O \; ,
	\]
	which implies $\GU \setminus \left( \bigcup \left( \CD^{-K} \right)^{(0)} \right) \prec_\CG O$ by Remark~\ref{rem:subequivalence-basic}\eqref{rem:subequivalence-basic:subset} and~\eqref{rem:subequivalence-basic:preorder}. 
\end{proof}

As a consequence, we have an analog of Theorem~\ref{thm:effective}. 

\begin{thm} 
	\label{thm:effective-AE}
	If a minimal \'{e}tale groupoid with a compact unit space is almost elementary, then
	it is effective. 
\end{thm}

\begin{proof}
	This follows directly from the ``if'' direction of Corollary~\ref{cor:effective-castle-minimal}. 
\end{proof}

We again investigate the case of finite groupoids, this time using the above result about effectiveness. 

\begin{cor} \label{cor:AE-finite-pair}
	A minimal \'{e}tale groupoid with a finite unit space is almost elementary if and only if it is a finite {pair groupoid}. 
\end{cor}

\begin{proof}
	The ``if'' direction follows from the observation that if $\CG$ is a finite pair groupoid, then we can verify Definition~\ref{defn:AE} by taking $\CC = \left\{ \{x\} \colon x \in \CG \right\}$ regardless of $K$, $\varepsilon$ and $\CV$, since $\CC^{-K} = \CC$ for any $K \subseteq \CG$. 
		
	In view of Theorem~\ref{thm:effective-AE}, the ``only if'' direction follows from the observation that a minimal effective \'{e}tale groupoid with a discrete unit space must be a pair groupoid. Indeed, by the minimality of $\CG$ and the discreteness of $\CG\unitSpace$, for any $u,v \in \CG\unitSpace$, there is $x_{(u,v)} \in \CG$ such that $s \left( x_{(u,v)} \right) = v$ and $r \left( x_{(u,v)} \right) = u$. Now, for any $y \in \CG$ with $s(y) = v$ and $r(y) = u$ we must have $y = x_{(u,v)}$, since otherwise $y^{-1} x_{(u, v)} \in \opIso(\CG)^{\operatorname{o}}\setminus\CG^{(0)}$, which is a contradiction to the effectivity of $\CG$. 
	It is then straightforward to verify that $(u,v) \mapsto x_{(u,v)}$ yields a groupoid isomorphism from the pair groupoid on the set $\CG\unitSpace$ to $\CG$. 
\end{proof}

In preparation for the following sections, 
we shall need another equivalent formulation of almost elementariness in terms of fortified castles (see Definition~\ref{defn:fortified-castle}). 

\begin{lem} \label{lem:AE-fortified} 
	A minimal \'{e}tale groupoid $\CG$ with a compact unit space is {almost elementary}
	if and only if 
	for any compact sets $K, L \subseteq \CG$, 
	any nonempty open set $O$ in $\GU$ 
	and any finite open cover $\CV$ of $L$, 
	there is a fortified castle $\CF$
	satisfying 
	\begin{enumerate}[label=(\roman*)]
		\item \label{lem:AE-fortified:refine} for any $F_+ \in \CF_+$ with $F_+ \cap L \not= \varnothing$, there is $V \in \CV$ such that $F_+ \subseteq V$;   
		\item \label{lem:AE-fortified:precompact} $\bigcup \CF_+$ is precompact in $\CG$;  
		\item \label{lem:AE-fortified:nonempty} $\CF^{-K}$ is nonempty;  
		\item \label{lem:AE-fortified:remainder} $\GU \setminus \left( \bigcup \left( \CF^{-K} \right)_-^{(0)} \right) \prec_\CG O$.  
	\end{enumerate}
\end{lem}

\begin{proof}
	In view of Lemma~\ref{lem:AE-shrinking}, the ``if'' direction is straightforward, since the castle $\CF_-$ clearly fulfills the role of $\CD$ in Lemma~\ref{lem:AE-shrinking}\eqref{lem:AE-shrinking:fine-entourage}. 
	
	To prove the ``only if'' direction, we fix $K$, $\varepsilon$ and $\CV$ as in the lemma and apply Lemma~\ref{lem:AE-shrinking}\eqref{lem:AE-shrinking:fine-entourage} to obtain an open castle $\CD$
	satisfying 
	\begin{enumerate}[label=(\roman*)]
		\item for any $D \in \CD$ with $D \cap L \not= \varnothing$, there is $V \in \CV$ such that $D \subseteq V$; 
		\item $\CD^{-K}$ is nonempty; 
		\item $\CG^{(0)}\setminus\ \left( \bigcup \left( \CD^{-K} \right)^{(0)} \right) \prec_{\CG} O$.   
	\end{enumerate}
	By Definition~\ref{defn:subequivalence}, 
	there is an open $s$-section $S$ such that $\CG^{(0)}\setminus\ \left( \bigcup \left( \CD^{-K} \right)^{(0)} \right)  \subseteq r(S)$ and $s(S) \subseteq O$. 
	Observe that $r(S)$ is open since $r$ is an open map. 
	Since $\CG\unitSpace \setminus \left( \bigcup \left( \CD^{-K} \right)^{(0)} \right) \not= \CG\unitSpace$, we may assume without loss of generality that $r(S) \not= \CG\unitSpace$. 
	Thus we may apply Urysohn's lemma to obtain a continuous function $f \colon \CG\unitSpace \to [0,1]$ that sends $\CG^{(0)}\setminus\ \left( \bigcup \left( \CD^{-K} \right)^{(0)} \right)$ to $0$ and $\CG\unitSpace \setminus r(S)$ to $1$.  
	Viewing $f$ as an element in $C_0 \left( \bigcup \left( \CD^{-K} \right)^{(0)} , \mathbb{R} \right) \subseteq C_0 \left( \bigcup \CD , \mathbb{R} \right)$, 
	we apply Lemma~\ref{lem:castle-invarint-function-smallest} to obtain $f_{\uparrow\CC} \in C_0 \left( \bigcup \CD , \mathbb{R} \right)^{\bigcup \CD}$ satisfying $0 \leq f \leq f_{\uparrow\CC} \leq \sup_{x \in \bigcup \CD} f(x) = 1$. 
	Hence applying Remark~\ref{rem:fortified-castle-basics}\eqref{rem:fortified-castle-basics:construction} to $f_{\uparrow\CC}$, 
	we obtain a fortified castle 
	\[
		\CF := \left\{ \left( f_{\uparrow\CC}^{-1} \left( \left(\frac{1}{3},1 \right] \right) \cap D, f_{\uparrow\CC}^{-1} \left( \left(\frac{1}{2},1 \right] \right) \cap D \right) \colon D \in \CD  \right\} 
	\]
	with $\CF_+$ refining $\CD$. 
	It follows that $\CF_+$ refines $\CV$ and $\CF^{-K}$ is nonempty. 
	It follows that conditions~\ref{lem:AE-fortified:refine} and~\ref{lem:AE-fortified:nonempty} follow from the corresponding conditions for $\CD$, while the precompactness of $\bigcup \CF_+$ follows from the fact that $f_{\uparrow\CC}$ is a $C_0$-function, which proves condition~\ref{lem:AE-fortified:precompact}.  
	It remains to show condition~\ref{lem:AE-fortified:remainder}.   
	To this end, it follows from $f \leq f_{\uparrow\CC}$ that $\bigcup \CF_- = f_{\uparrow\CC}^{-1} \left( \left(\frac{1}{2},1 \right] \right) \supseteq f_{\uparrow\CC}^{-1} (1) \supseteq \CG\unitSpace \setminus r(S)$. Since we also have $\CG\unitSpace \setminus r(S) \subseteq \bigcup \left( \CD^{-K} \right)^{(0)}$, 
 we have $\CG\unitSpace \setminus r(S) \subseteq \left( \bigcup \CF_- \right) \cap  \left( \CD^{-K} \right)^{(0)} = \bigcup \left( \CF^{-K} \right)_-^{(0)}$, whence 
	$\GU \setminus \bigcup \left( \CF^{-K} \right)_-^{(0)} \subseteq r(S)$ and thus $\GU \setminus \bigcup \left( \CF^{-K} \right)_-^{(0)} \prec_\CG O$ by the definition of $S$. 
\end{proof}

\section{Groupoid strict comparison} \label{sec:comparison}

In this section, we study another regularity property of \'{e}tale groupoids, namely, groupoid strict comparison. This is the key link to connect almost elementariness to its in-measure version. 

Comparing Definition~\ref{defn:AE-in-measure} against Definition~\ref{defn:AE}, we see that the difference between the two notions lies in the remainder conditions \textemdash\ the former uses a ``measure-theoretical'' subequivalence while the latter uses the ``topological'' subequivalence $\prec_{\CG}$ (or, equivalently in this context by Remark~\ref{rem:subequivalence-basic}\eqref{rem:subequivalence-basic:clopen}, $\precsim_{\CG}$). Note that Remark~\ref{rem:subequivalence-basic}\eqref{rem:subequivalence-basic:measure-auxiliary} and~\eqref{rem:subequivalence-basic:measure-subequivalence} provide a one-way relationship between the two kinds of subequivalences. The property of groupoid strict comparison below is nothing but a partial converse
of Remark~\ref{rem:subequivalence-basic}\eqref{rem:subequivalence-basic:measure-subequivalence}.
For the ease of our exposition, we restrict our attention to minimal groupoids, leaving the study of non-minimal groupoids to future work.

\begin{defn}
	\label{5.2} 
	\label{defn:comparison}
	We say a minimal \'{e}tale groupoid $\CG$ with a compact unit space has \emph{groupoid strict comparison} (or simply \emph{groupoid
		comparison} or\emph{ comparison}) if, for any nonempty open sets $U,V\subseteq\CG^{(0)}$,
	we have $U\precsim_{\CG}V$ whenever $\mu(U)<\mu(V)$ for all $\mu\in M(\CG)$.
\end{defn}

\begin{rem}\label{rem:comparison-infinite}
	Note that we do not rule out the case where $M(\CG) = \varnothing$, in which case having groupoid strict comparison amounts to saying $\GU \precsim_\CG V$ for any nonempty open set $V \subseteq \GU$, by virtue of Remark~\ref{rem:subequivalence-basic}\eqref{rem:subequivalence-basic:subset} and ~\eqref{rem:subequivalence-basic:preorder}. 
\end{rem}

\begin{example} \label{exa:GSBP-pair-groupoid}
	It is clear that the pair groupoid $X \times X$ on a finite set $X$ has groupoid strict comparison, since in this case the only invariant measure on the unit space is the uniform measure on $X$, while for any $U, V \subseteq X$, $U \precsim_{X \times X} V$ means there is a surjection from a subset of $V$ onto $U$. Therefore, such a groupoid satisfies a stronger comparison property with the strict inequality $<$ in Definition~\ref{defn:comparison} replaced by $\leq$. 
\end{example}

\begin{example}
	\label{exa:comparison-actions}
	Let $X \rtimes_\alpha \Gamma$ be a minimal transformation groupoid as in Example~\ref{exa:subequivalence-actions}. Note that $X \rtimes_{\alpha}\Gamma$ is minimal if and only if $\alpha$ is a minimal action. 
	\begin{enumerate}		
		\item \label{exa:comparison-actions:comparison} It follows from Example~\ref{exa:subequivalence-actions}\eqref{exa:subequivalence-actions:subequivalence} and the canonical identification of $M(X \rtimes_\alpha \Gamma)$ with the set $M_\Gamma(X)$ of $\Gamma$-invariant regular Borel probability measures on $X$ that in the case of minimal actions, $X \rtimes_\alpha \Gamma$ has groupoid strict comparison if and only if the action $\alpha$ has dynamical strict	comparison defined in \cite[Definition 3.2]{Kerr2020Dimension}. 
		
		\item \label{exa:comparison-actions:non-strict} It has been known in the study of Cantor $\mathbb{Z}$-systems that the stronger comparison property with $<$ replaced by $\leq$ (see Example~\ref{exa:GSBP-pair-groupoid}) has no chance to hold for general \'{e}tale transformation groupoids (see .............). 
	\end{enumerate}
\end{example}

The existence of open castles in an \'{e}tale groupoid greatly helps us with questions such as groupoid strict comparison, by virtue of the following observation, which provides concrete criteria for groupoid subequivalence. 

\begin{lem} \label{lem:castle-comparison}
	Let $\CC$ be an open castle in an \'{e}tale groupoid $\CG$. Recall from Lemma~\ref{lem:castle-basic}\eqref{lem:castle-basic:subquotient} that $\bigcup \CC$ is an open subgroupoid of $\CG$ while $\CC$ itself forms a finite groupoid. Let $\CX, \CY \subseteq \CC\unitSpace$. 
	Then the following are equivalent: 
	\begin{enumerate}
		\item \label{lem:castle-comparison:subequivalence} $\bigcup \CX \precsim_{\bigcup \CC} \bigcup \CY$; 
		\item \label{lem:castle-comparison:cardinality} for any tower $\CT$ in $\CC$, we have $|\CX \cap \CT\unitSpace| \leq |\CY \cap \CT\unitSpace|$; 
		\item \label{lem:castle-comparison:bisection} there are $\CZ \subseteq \CY$ and a bisection $\CB$ in the finite groupoid $\CC$ such that $s(\CB) = \CX$ and $r(\CB) = \CZ$. 
	\end{enumerate}
	Moreover, these conditions imply $\bigcup \CX \precsim_{\CG} \bigcup \CY$. 
\end{lem}

\begin{proof}
The fact that $\bigcup \CX \precsim_{\bigcup \CC} \bigcup \CY$ implies $\bigcup \CX \precsim_{\CG} \bigcup \CY$ follows from Remark~\ref{rem:subequivalence-basic}\eqref{rem:subequivalence-basic:subgroupoid}.	 We now prove 
	\eqref{lem:castle-comparison:subequivalence} $\Rightarrow$ \eqref{lem:castle-comparison:cardinality} $\Rightarrow$ \eqref{lem:castle-comparison:bisection} $\Rightarrow$ \eqref{lem:castle-comparison:subequivalence}.
 
\eqref{lem:castle-comparison:subequivalence} $\Rightarrow$ \eqref{lem:castle-comparison:cardinality}:
Let $\CT$ be a tower in $\CC$, which is an equivalence class of the relation $\sim_{\CC}$ in the sense of Definition \ref{defn:castle}. Then for each $C\in \CX\cap \CT\unitSpace$, choose $u_C\in C$ such that for any $C, C'\in \CX\cap \CT\unitSpace$ there exists a $D\in \CT$ with $r(Du_C)=u_{C'}$. Note that this implies $u_C\sim_{\bigcup \CC} u_{C'}$ in the sense of Definition \ref{defn:invariant-functions}.

Denote by $K=\{u_C: C\in \CX\cap \CT\unitSpace\}$, which is finite. Note that $|K|=|\CX\cap \CT\unitSpace|$. Now the condition \eqref{lem:castle-comparison:subequivalence} implies that there exists an open $s$-section $S\subseteq \bigcup \CC$ such that $K\subseteq r(S)$ and $s(S)\subseteq \bigcup \CY$. Then for each $u_C\in K$, choose $\gamma_C\in S$ with $r(\gamma_C)=u_C$. Denote by $L=\{\gamma_C\in S: C\in \CX\cap \CT\unitSpace\}$. By definition, observe that $|K|=|L|=|s(L)|$ because $S$ is an $s$-section. Now, for any $\gamma_C\in L$, there exists a $D\in \CC$ such that $\gamma_C\in D$. This implies that $CD\neq \emptyset$ and one has $r(D)=C$ by Lemma \ref{lem:plurisection-basic}\eqref{lem:plurisection-basic:composability}. This further implies  that $s(D)\sim_{\CC} C$ and therefore $s(D)\in \CT\unitSpace$. On the other hand, $\gamma_C\in S$ implies that $s(\gamma_C)\in C'$ for some $C'\in \CY$. Then, similarly by Lemma \ref{lem:plurisection-basic}\eqref{lem:plurisection-basic:composability}, $DC'\neq \emptyset$ implies that $C'=s(D)\in \CT\unitSpace\cap \CY$. 

Note that for different $\gamma_1, \gamma_2\in L$, one necessarily has $s(\gamma_1)\neq s(\gamma_2)$ because $L\subseteq S$, which is a $s$-section. Now, let $C_1, C_2\in \CT\unitSpace\cap \CY$ such that $s(\gamma_1)\in C_1$ and $s(\gamma_2)\in C_2$. Now, fix a $u_C$ for some $C\in \CX \cap \CT\unitSpace$. Then by our construction, one has $s(\gamma_i)\sim_{\bigcup\CC} u_C$ for all $i=1,2$. Then, because any $C\in \CC$ is a bisection by Lemma \ref{lem:castle-basic}\eqref{lem:castle-basic:bisection}, there is at most one $x\in C$ such that $x\sim_{\bigcup\CC} u_C$. This implies that $C_1\neq C_2$ and thus $C_1\cap C_2=\emptyset$ by Lemma \ref{lem:plurisection-basic}\eqref{lem:plurisection-basic:composability}. Then, one has 
\[|\CX\cap \CT\unitSpace|=|s(L)|\leq|\{C\in \CY\cap \CT\unitSpace: s(\gamma)\in C, \gamma\in L\}|\leq  |\CT\unitSpace\cap \CY|.\]

\eqref{lem:castle-comparison:cardinality} $\Rightarrow$ \eqref{lem:castle-comparison:bisection}: For each tower $\CT$ in $\CC$, choose an injective map $f_\CT: \CX\cap \CT\unitSpace\to \CY\cap \CT\unitSpace$. Then define $\CB_{\CT}=\{C\in \CT: s(C)\in \CX\text{ and }r(C)=f_\CT(s(C))\}$ and $\CB=\bigsqcup\{\CB_\CT: \CT\text{ is a tower in }\CC\}$ is the desired bisection in  the groupoid $\CC$.

\eqref{lem:castle-comparison:bisection} $\Rightarrow$ \eqref{lem:castle-comparison:subequivalence}:
Let $\CZ\subseteq \CY$ and $\CB$ be the bisection in the groupoid $\CC$ such that $s(\CB)=\CX$ and $r(\CB)=\CZ$. Define $S=\bigcup \CB$, which is an open bisection in $\bigcup\CC$. In particular, $S^{-1}$ is an $s$-section in $\bigcup \CC$. Observe $K\subseteq \bigcup \CX=r(S^{-1})$ for any compact set $K\subseteq \bigcup \CX$. Furthermore, $s(S^{-1})\subseteq \bigcup \CZ\subseteq \bigcup \CY$. This establishes \eqref{lem:castle-comparison:subequivalence}.
\end{proof}

Our main task for the rest of this section is to relate almost elementariness, its in-measure version, and groupoid strict comparison. 
The most intricate step is to show almost elementariness implies groupoid strict comparison. To do this, we investigate the fiberwise amenable case and the non-fiberwise-amenable case separately. 

We start with the former case, where we produce castles with ``F\o{}lner shapes'' in order to facilitate comparison of sets. 

\begin{lem}\label{lem: double shrinking}
 Let $\CG$ be an \'{e}tale  groupoid with a compact unit space and let $K, L\subseteq \CG$ be compact set. Let $\CC$ be a castle in $\CG$. Suppose $K$ contains $\GU$ and suppose $L'\supset L$ is another compact set in $\CG$ containing $\GU$ and satisfies that for any $C\in \CC$, if $C\cap L\neq \emptyset$ then $C\subseteq L'$. Then one has
 \[\CC^{-KL'}\subseteq (\CC^{-K})^{-L}\subseteq \CC^{-KL}.\]
\end{lem}
\begin{proof}
First, let $C\in (\CC^{-K})^{-L}$. for any $x\in C$ and $y\in L$, then by definition, one has $y\cdot x\in D$ for some $D\in \CC^{-K}$. Then for any $z\in K$, one obtains that $z\cdot x\cdot y\in \bigcup \CC$. This implies that $C\in \CC^{-KL}$ and thus we have $(\CC^{-K})^{-L}\subseteq \CC^{-KL}$.

Then, let $C\in \CC^{-KL'}$, $x\in C$ and $y\in L$. Since $\GU\subseteq K$ and $L\subseteq L'$, one obtains 
\[y\cdot r(x)=s(y)\cdot y\cdot r(x)\in \bigcup \CC,\]
which implies that there exits a $C_1\in \CC$ such that $y\in C_1$ and $s(C_1)=r(C)$. This further entails $C_1\cap L\neq \emptyset$ and therefore $C_1\subseteq L'$ by our assumption. Now note that 
\[K\cdot C_1\subseteq KL'\cdot r(C)\subseteq \bigcup \CC\]
as $r(C)\in \CC^{-KL'}$ by Lemma \ref{lem:castle-shrinking}\eqref{lem:castle-shrinking:castle}.  This shows that actually $C_1\in \CC^{-K}$. In addition, because $L'$ contains $\GU$, one obtains $\CC^{-KL'}\subseteq \CC^{-K}$
thus one obtains $y\cdot x\in C_1\cdot C$, which is an element in $\CC^{-K}$. This thus shows $L\cdot C\subseteq \bigcup \CC^{-K}$, which means $C\in (\CC^{-K})^{-L}$. As a consequence, $\CC^{-KL'}\subseteq (\CC^{-K})^{-L}$ holds.
\end{proof}

\begin{lem} \label{lem:AE-fiberwise-fortified} 
	Let $\CG$ be a minimal \'{e}tale groupoid with a compact unit space. 
	Suppose $\CG$ is $\sigma$-compact and fiberwise amenable. 
	If $\CG$ is almost elementary, then 
	for any nonempty compact set $K$ in $\CG$, 
	any $\varepsilon > 0$, 
	any non-empty open set $O$ in $\GU$ 
	and any finite open cover $\CV$ of $\CG^{(0)}$, 
	there is a nonempty fortified castle $\CF$
	satisfying 
	\begin{enumerate}[label=(\roman*)]
		\item \label{lem:AE-fiberwise-fortified:refine} 
		$\CF_+ \unitSpace$ refines $\CV$;  
		\item\label{lem:AE-fiberwise-fortified:precompact} $\bigcup \CF_+$ is precompact in $\CG$; 
		\item \label{lem:AE-fiberwise-fortified:remainder} 
		$\GU \setminus \left( \bigcup \CF_-^{(0)} \right) \prec_\CG O$; 
		\item \label{lem:AE-fiberwise-fortified:Folner} 
		for any $u \in \CG\unitSpace$, we have $|K (\bigcup \CF_+) u| \leq (1+\varepsilon) |(\bigcup \CF_+) u|$. 
	\end{enumerate}
\end{lem}

As will be clear later, the above condition actually gives yet another equivalent characterization of almost elementariness. 

\begin{proof}
	Assume $\CG$ is almost elementary and let $K$, $\varepsilon$, $O$ and $\CV$ be given as in the statement above. 
	Since $\CG$ is locally compact and Hausdorff, we may find a finite cover $\CW$ of $K \cup \CG\unitSpace$ by precompact open bisections. 
	Define another compact set $K' := \bigcup_{W \in \CW} \overline{W}$ in $\CG$.
	By Theorems~\ref{thm:minimal-fiberwise-amenability} and~\ref{thm:fiberwise-amenable-dichotomy}\ref{thm:fiberwise-amenable-dichotomy:amenable}, there is a compact set $L$ in $\CG$ such that for any finite set $M \subseteq\CG$, there
	is a finite set $N$ satisfying 
	\[
		M \subseteq N\subseteq LM \ \textrm{and}\ |K'N|\leq(1+\varepsilon)|N|.
	\]
Choose an open covering $\CO$ of $L\cup \GU$ consisting of  precompact open bisections and define $L'=\bigcup_{O\in \CO}\overline{O}$.  Define another open cover 
\[\CU := (\CW\cup\{\CG\setminus K \}) \vee (\CO\cup \{\CG\setminus L\}) \vee (\CV \cup \{ \CG \setminus \CG\unitSpace \})\] of $K \cup  L\cup\CG\unitSpace$.     
	Now by Lemma~\ref{lem:AE-fortified}, there is a fortified castle $\CE$
	satisfying 
	\begin{enumerate}[label=(\alph*)]
		\item \label{lem:AE-fiberwise-fortified:proof:refine}  for any $E_+ \in \CE_+$ with $E_+ \cap (K\cup L\cup \GU) \not= \varnothing$, there is $U \in \CU$ such that $E_+ \subseteq U$;   
		\item \label{lem:AE-fiberwise-fortified:proof:precompact}  $\bigcup \CE_+$ is precompact in $\CG$;  
		\item \label{lem:AE-fiberwise-fortified:proof:remainder}  $\GU \setminus \left( \bigcup \left( \CE^{-K'L'} \right)_-^{(0)} \right) \prec_\CG O$.  
	\end{enumerate}
	Decompose $\CE_+$ into its towers as $\bigsqcup_{i \in I} \CT_i$. For any $i \in I$, choose an $E_{i,+} \in \CT_i \unitSpace$ and a $u_i \in E_{i,+}$, and define $M_i := \left( \bigcup \CT_i^{-K'L'} \right) u_i$, which is a finite set since $\CT_i^{-K'L'}$ is a finite collection of bisections by Lemma~\ref{lem:castle-basic}. Moreover, for each $i\in I$ note that 
    \[L\cdot (\bigcup T_i^{-K'L'})\subseteq \bigcup T_i^{-K'}\]
    by Lemma \ref{lem: double shrinking} . By our assumption on $L$ above, for any $i \in I$, there is a finite set $N_i$ satisfying $M_i \subseteq N_i\subseteq LM_i\subseteq (\bigcup T_i^{-K'})u_i$ and $|K'N_i|\leq(1+\varepsilon)|N_i|$. 
	Writing $N :=  \bigsqcup_{i \in I} N_i$, we define $\CF$ to be the subcollection
	\[
		\left\{ (F_+, F_-) \in \CE \colon s(F_+) \cap r(N) \not= \varnothing \not=  r(F_+) \cap r(N) \right\} \; .
	\]
	From the construction, observe that $\CE^{-K'L}_+ \subseteq \CF_+\subseteq \CE^{-K'}_+$. Hence it follows directly from conditions~\ref{lem:AE-fiberwise-fortified:proof:refine}, \ref{lem:AE-fiberwise-fortified:proof:precompact} and~\ref{lem:AE-fiberwise-fortified:proof:remainder} above that $\CF$ satisfies conditions~\ref{lem:AE-fiberwise-fortified:refine}, \ref{lem:AE-fiberwise-fortified:precompact} and~\ref{lem:AE-fiberwise-fortified:remainder} in the statement. 
	
	It remains to prove condition~\ref{lem:AE-fiberwise-fortified:Folner}. 
	To this end, let $u \in \CG\unitSpace$ be given. Observe that it suffices to restrict to the case where $u \in \bigcup \CF_+\unitSpace$ because other with the inequality in ~\ref{lem:AE-fiberwise-fortified:Folner} holds trivially as $|K (\bigcup \CF_+) u| = |(\bigcup \CF_+) u|=0$. Now, since $\CF_+\subseteq \CE^{-K'}_+$, one has $K\cdot\bigcup \CF_+\subseteq \bigcup\CE_+$. Let $u\in \bigcup \CF_+\unitSpace$. Then for any $x\in \bigcup\CF_+u$ there exists an $i\in I$ and a level $F_{i, +}\in \CF_+\unitSpace\cap \CT_i\unitSpace$ such that $r(x)\in F_{i, +}$. Then for the base level $E_{i, +}\in \CT_i\unitSpace$ containing $u_i$ defined above, there exists a unique $x'\in \bigcup\CF_+$ such that $s(x')=u_i$ and $r(x')\in F_{i, +}$. We write $f(x)$ for the $x'$ determined by $x$ above.  Now, suppose $Kx\neq \emptyset$ and let $y\in Kr(x)$. Then because $\CF_+\subseteq \CE_+^{-K'}$ and $K\subseteq K'$, there exists a ladder $E_{y,+}\in \CT_i$ such that $y\in E_{y, +}$, which implies that $E_{y, +}\cap K\neq \emptyset$. Then the condition (a) above for $\CE_+$ implies that $E_{y, +}\subseteq U$ for some $U\in \CU$ and therefore $E_{y, +}\subseteq K'$.  Observe that for all $x\in \bigcup \CF_+u$ with $Kx\neq  \emptyset$ this process yields a map $g_x$ from $Kx$ to $K'f(x)$, defined by $g_x(yx)=zf(x)$ such that $y\in Kr(x)$ and $z\in E_{y,+}$ with $s(z)=r(f(x))$. Now define a map $g: K(\bigcup\CF_+)u\to K'f((\bigcup\CF_+)u)$ by
 \[g(yx)=g_x(yx)\]
 for $y\in K$ and $x\in \CF_+u$. 
 
 We claim that the map $g$ is well-defined. Indeed, let $y_1x_1=y_2x_2\in K(\bigcup\CF_+)u$ with $y_1, y_2\in K$, $x_1, x_2\in (\bigcup\CF_+)u$. Then let $x_1\in F_{1, +}$ and $x_2\in F_{2, +}$, where $F_{1, +}$ and $F_{2, +}$ are $\CF_+$-ladders. This implies that $E_{y_1,+}\cdot F_{1,+}=E_{y_2, +}\cdot F_{2, +}$ because they are $\CE$-ladders consisting of a common element $y_1x_1=y_2x_2$. This implies that 
 \[g(y_1x_1)=g_{x_1}(y_1x_1)=g_{x_2}(y_2x_2)=g(y_2x_2)\]
 because both $g_{x_1}(y_1x_1)$ and $g_{x_2}(y_2x_2)$ are the unique element $\gamma$ in $\bigcup \CE$ such that $s(\gamma)=u_i$ and $r(\gamma)\in r(E_{y_1,+}\cdot F_{1,+})=r(E_{y_2,+}\cdot F_{2,+})$.

 Moreover, the function $g$ is injective. Indeed, if $y_1x_1\neq y_2x_2$ in $K (\bigcup \CF_+) u$, then $r(y_1)\neq r(y_2)$. This implies that  $r(E_{y_1, +})\neq r(E_{y_2, +})$. Then $g(y_1x_1)\neq g(y_2x_2)$ by definition because each $g(y_jx_j)$ is the element $\gamma\in \bigcup\CE$ with $s(\gamma)=u_i$ and $r(\gamma)\in r(E_{y_j, +})$ for $j=1,2$. On the other hand, observe that $K'f((\bigcup \CF_+)u)\subseteq K'(\bigcup \CF_+)u_i$ and therefore one has
 \[|K(\bigcup \CF_+)u|\leq |K'f((\bigcup \CF_+)u)|\leq |K'(\bigcup \CF_+)u_i|.\]
 Next, note that $|(\bigcup \CF_+)u|=|(\bigcup \CF_+)u_i|$ because $\CF_+$ is a tower. Then using the fact that  $(\bigcup \CF_+)u_i=N_i$, one has
 \[|K'(\bigcup \CF_+)u_i|=|K'N_i|\leq (1+\varepsilon)|N_i|=(1+\varepsilon)|(\bigcup \CF_+)u_i|=(1+\varepsilon)|(\bigcup \CF_+)u|.\]
 Finally, by combining these two inequalities, one has 
 \[|K(\bigcup \CF_+)u|\leq (1+\varepsilon)|(\bigcup \CF_+)u|.\]
 This has established \ref{lem:AE-fiberwise-fortified:Folner}. \end{proof}

\begin{lem} \label{lem:O-meas-sup}
	Let $\CG$ be a minimal \'{e}tale groupoid with a compact unit space. 
	Suppose $\GU$ is an infinite set. 
	Then for any nonempty open subset $V$ of $\GU$ and any natural number $N$, 
	there is an open tower $\CT$ such that $\bigcup \CT^{(0)} \subseteq V$ and $\left| \CT^{(0)} \right| \geq N$. 
	In particular, we have $\displaystyle \sup_{\mu \in M(\CG)} \mu (C) \leq \frac{1}{N}$ for any $C \in \CT^{(0)}$. 
\end{lem}

\begin{proof}
Let $u\in V$ and $k\geq N$ be an integer. The minimality of $\CG$ implies that there are pairwise different $x_1, \dots, x_k\in 
\CG$ such that $s(x_i)=u$ and $r(x_i)\in V$ for all $1\leq i\leq k$. Then choose open bisections $U_i\ni x_i$ for $i\leq k$ such that $u\in \bigcap_{i=1}^ks(U_i)$ and  $\{r(U_i): 1\leq i\leq k\}$ is a disjoint family of open sets contained in $V$. Now, choose an open set $O$ with $u\in O\subseteq \bigcap_{i=1}^ks(U_i)$. Then define $C_{i, j}=U_iOU^{-1}_j$ for all $1\leq i, j\leq k$. Now, it is direct to verify that $\CT=\{C_{i, j}: 1\leq i, j\leq k\}$ is a tower in Definition \ref{defn:castle}. Moreover, by our construction, one has $|\CT\unitSpace|=k\geq N$ and $\bigcup\CT\unitSpace\subseteq\bigsqcup_{i=1}^kr(U_i)\subseteq V$.
\end{proof}

\begin{prop}
	 \label{prop:AE-DSC-fiberwise} 
	Let $\CG$ be a minimal $\sigma$-compact \'{e}tale groupoid
	with a compact unit space. 
	If $\CG$ is fiberwise amenable and almost elementary, then $\CG$ has groupoid strict comparison. 
\end{prop}
\begin{proof}
	In view of Corollary~\ref{cor:AE-finite-pair} and Example~\ref{exa:GSBP-pair-groupoid}, it suffices to focus on the case where $\CG\unitSpace$ if infinite. 
	Now let two nonempty open sets $U,V\subseteq\CG^{(0)}$ be given and suppose $\mu(U)<\mu(V)$ for any $\mu\in M(\CG)$. To show $U\precsim_{\CG}V$, it suffices to show, after fixing an arbitrary compact subset $L \subseteq U$, that $L\prec_{\CG}V$. 
	
	To this end, we first exploit the strictness of the inequality $\mu(U)<\mu(V)$ to create some ``wiggle room'' around $L$ and inside $V$: 
	Since $\CG\unitSpace$ is locally compact and Hausdorff, there is an open neighborhood $U'$ of $L$ with $\overline{U'} \subseteq U$. 
	By Lemma~\ref{lem:measure-semicontinuous}, the function $M(\CG\unitSpace) \to [-1,1]$, defined by $\mu \mapsto \mu(V) - \mu(\overline{U'})$, is lower semicontinuous, and thus it has a minimum $\eta$ on the compact set $M(\CG) \subseteq M(\CG\unitSpace)$. Note that $\eta > 0$. 
	By Lemma~\ref{lem:O-meas-sup}, there is a nonempty open set $O_0 \subseteq V$ such that $\mu(O_0) \leq \eta/2$. Again since $\CG\unitSpace$ is compact and Hausdorff,  we may choose a nonempty open set $O$ with $\overline{O} \subseteq O_0$. Let $W := V \setminus \overline{O}$. Since 
    \[\mu(W) = \mu(V) - \mu(\overline{O}) \geq \mu(V) - \mu({O_0}) \geq \mu(V) - \eta/2 > \mu(V) - \eta \geq \mu(\overline{U'})\]  holds for any $\mu \in M(\CG)$, we may apply Lemma~\ref{lem:small-neighborhood-in-measure} to the open set $W$ and the upper continuous function $\mu \mapsto \mu(\overline{U'})$ to obtain an open subset $W'$ with $\overline{W'} \subseteq W$ such that $\mu({W'}) > \mu(\overline{U'})$ for any $\mu \in M(\CG)$. Apply Lemma~\ref{lem:measure-semicontinuous} again to obtain $\eta' := \min_{\mu \in M(\CG)} (\mu({W'}) - \mu(\overline{U'})) > 0$. 
	
	Now we claim that there are a compact set $K$ and $\varepsilon > 0$ such that for any finite subset $F$ in $\CG$ with $|KF| < (1+\varepsilon) |F|$, we have 
 \[\frac{1}{|F|}\sum_{x\in F}1_{\overline{U'}}(r(x))+\frac{\eta'}{2}<\frac{1}{|F|}\sum_{x\in F}1_{W'}(r(x)).\]
 Suppose the contrary, for any increasing sequence of compact sets $\{K_n: n\in \N\}$ be in $\CG$ with $\bigcup_{n\in \N}K^o_n=\CG$ and any decreasing sequence of positive numbers$\{\varepsilon_n: n\in \N\}$ converging to $0$, there exists a finite $F_n$ in $\CG$ such that $|K_nF_n|<(1+\varepsilon_n)|F_n|$ and 
 \[\mu_n(\overline{U'})+\frac{\eta'}{2}\geq\mu_n(W'),\]
 where $\mu_{n}=\frac{1}{|F_{n}|}\sum_{x\in F_{n}}\delta_{r(x)}$. Then the proof of \cite[Theorem B]{MWfiberwise} implies that any cluster point $\mu$ of the sequence $\{\mu_n: n\in \N\}$ belongs to $M(\CG)$. This entails that
 \[\mu(W')+\frac{\eta'}{2}\leq \liminf_{n\to \infty}\mu_n(W')\leq \limsup_{n\to \infty}\mu_n(\overline{U'})+\frac{\eta'}{2}\leq \mu(\overline{U'})+\eta'.\]
 This is a contradiction to the choice of $\eta'$ above and therefore the claim above holds.
	
 Let $\CV := \left\{ U' , \CG\unitSpace \setminus L \right\} \vee \left\{ W , \CG\unitSpace \setminus \overline{W'} \right\}$. For the compact set $K$, the number $\varepsilon>0$, the open set $O$ and the cover $\CV$, Lemma~\ref{lem:AE-fiberwise-fortified} yields a fortified castle $\CF$ satisfying
 \begin{enumerate}[label=(\roman*)]
		\item  $\CF_+ \unitSpace$ refines $\CV$;  
		\item $\bigcup \CF_+$ is precompact in $\CG$; 
		\item $\GU \setminus \left( \bigcup \CF_-^{(0)} \right) \prec_\CG O$; 
		\item  \label{prop:AE-DSC-fiberwise:Folner}  for each fortified tower $\CT$ in $\CF$, one may choose a $u_\CT\in \bigcup\CT\unitSpace$ such that $|K (\bigcup \CT_+) u_\CT| \leq (1+\varepsilon) |(\bigcup \CT_+) u_\CT|$. 
  \end{enumerate}
Define $L_1=L\cap (\CG^0\setminus \bigcup \CF\unitSpace_-)$ and $L_2=L\cap \left(\bigcup\{\overline{C_-}: C_-\in \CF\unitSpace_-\}\right)$. First, observe $L_1\prec_\CG O$ because $L_1$ is a subset of $\GU \setminus \left( \bigcup \CF_-^{(0)} \right) \prec_\CG O$.  Define $\CX=\{C_+\in \CF\unitSpace_+: C_+\cap L\neq \emptyset\}$ and $\CY=\{C_+\in \CF\unitSpace_+: C_+\cap \overline{W'}\neq \emptyset\}$. In addition, define $\CX'=\{C_+\in \CF\unitSpace_+: C_+\subseteq U'\}$ and $\CY'=\{C_+\in \CF\unitSpace_+: C_+\subseteq W\}$. Because $\CF\unitSpace_+$ refines $\CV$, one has $C_+\subseteq U'$ holds for any $C_+\in \CX$, and $C_+\subseteq W$ holds for any $C_+\in \CY$. Therefore, one actually has $\CX\subseteq \CX'$ and $\CY\subseteq \CY'$. 

Now, let $\CT$ be a fortified tower in $\CF$. Let $u_\CT$ be a unit satisfying \ref{prop:AE-DSC-fiberwise:Folner} above and define $F_\CT=(\bigcup\CT_+)u_\CT$. Then the choice of $K$ and $\varepsilon$ implies that $|r(F_\CT)\cap U'|<|r(F_\CT)\cap \overline{W'}|$ because $F_\CT$ satisfies
\[\sum_{x\in F_\CT}1_{\overline{U'}}(r(x))<\sum_{x\in F_\CT}1_{W'}(r(x))\leq \sum_{x\in F_\CT}1_{\overline{W'}}(r(x)).\]
This further implies that 
\[|\CX\cap \CT\unitSpace|\leq |\CX'\cap \CT\unitSpace|\leq |\CY'\cap \CT\unitSpace|.\]
Lemma~\ref{lem:castle-comparison} then entails that $\bigcup \CX\precsim_\CG \bigcup \CY'\subseteq W$. Moreover, since $L_2\subseteq\bigcup\CX$ and $O, W$ are disjoint open sets in $V$, Remark~\ref{rem:subequivalence-basic}\eqref{rem:subequivalence-basic:additive-prec} entails that $L=L_1\cup L_2\prec_\CG V$.
\end{proof}

\begin{rem}\label{rem: minimal-perfect}
  If $\CG$ is a minimal \'{e}tale groupoid such that the unit space $\GU$ is infinite compact, then $\GU$ is \textit{perfect} in the sense that it has no isolated points. Suppose the contrary, let $\{u\}$ be an open set in $\GU$. Since $\CG$ is minimal, there are open bisections $O_1,\dots, O_n$ such that $\GU=\bigcup_{i=1}^ns(O_i)$ and $\bigcup_{i=1}^nr(O_i)\subseteq \{u\}$, which implies that $\GU$ is finite. This is a contradiction to the assumption that $\GU$ is infinite.
  \end{rem}

\begin{prop} \label{prop:AE-DSC-non-fiberwise} 
	Let $\CG$ be a minimal $\sigma$-compact \'{e}tale groupoid
	with a compact unit space. 
	If	$\CG$ is non-fiberwise-amenable and almost elementary, then $M(\CG)=\emptyset$ and $\CG$ has groupoid strict comparison. 
\end{prop}

\begin{proof}
	It still suffices to work on the case that $\GU$ is infinite by  Corollary~\ref{cor:AE-finite-pair} and Example~\ref{exa:GSBP-pair-groupoid}. To this end, we will show $\GU\prec_\CG O_0$ for any non-empty open set $O_0$. Remark \ref{rem: minimal-perfect} implies that $\GU$ is perfect and therefore, there are two disjoint non-empty open sets $U, V$ in $O$. For the $V$, choose a non-empty open set $W$ with $W\subseteq \overline{W}\subseteq V$.  Since $\CG$ is minimal, Lemma \ref{lem:minimal-recurrence} implies that there are compact bisections $K_1, \dots, K_m$ in $\CG$ such that $\bigcup_{i=1}^ms(K_i)^o=\GU$ and $\bigcup_{i=1}^mr(K_i)\subseteq W$. Define $K=\bigcup_{i=1}^mK_i\cup \GU$, which is a compact set in $\CG$. Then Theorem \ref{thm:fiberwise-amenable-dichotomy}\eqref{thm:fiberwise-amenable-dichotomy:non-amenable} implies that there exists a compact set $L$ containing $\GU$ such that for any $n\in \N$ and any finite set $M\subseteq \CG$, the set $LM$ contains at least $n|M|$ many disjoint sets of the form $Kx$ for $x\in \CG$. Define an open cover $\CU=\{V, \GU\setminus \overline{W}\}\vee \{s(K_1)^o,\dots, s(K_n)^o\}$ of $\GU$ and an open cover $\CV=\CU\cup\{\CG\setminus \GU\}$ of $\CG$. Then apply Lemma~\ref{lem:AE-fortified} to compact sets $L, \GU$ in $\CG$, the cover $\CV$ and the open set $U$ to obtain a fortified castle $\CF$
satisfying
\begin{enumerate}[label=(\roman*)]
		\item\label{prop:AE-DSC-non-fiberwise-condition1} 
  for any $F_+ \in \CF_+$ with $F_+ \cap \GU \not= \varnothing$, there is $V \in \CV$ such that $F_+ \subseteq V$;   
		\item  $\bigcup \CF_+$ is precompact in $\CG$;  
		\item  $\CF^{-L}$ is nonempty;  
		\item\label{prop:AE-DSC-non-fiberwise-condition4} 
  $\GU \setminus \left( \bigcup \left( \CF^{-L} \right)_-^{(0)} \right) \prec_\CG U$.  
	\end{enumerate}
For any fortified tower $\CT$ in $\CF$, choose a unit $u_\CT\in (\CT^{-L})\unitSpace_-$ and define a finite set $M_\CT=(\bigcup\CT^{-L}_-)u_\CT$. Then by our choice of $L$, the set $LM_\CT$ contains at least $|M_\CT|$ many disjoint non-empty finite subsets of the form $Kx$ for $x\in \CG$ with $s(x)=u_\CT$. Denote by $m_\CT=|M_\CT|$ for simplicity and write these $Kx$ explicitly by $\{Kx_1, \dots, Kx_{m_\CT}\}$. Moreover, one has
\[\bigsqcup_{i=1}^{m_\CT}Kx_i\subseteq LM_\CT\subseteq (\bigcup\CT_-)u_\CT.\]
Choose a bijective map $\varphi_\CT: (\CT^{-L})\unitSpace_+\to \{x_1,\dots, x_{m_\CT}\}$, which is possible because $m_\CT=|(\bigcup\CT^{-L}_-)u_\CT|=|(\CT^{-L})\unitSpace_+|.$ Then for each $x_i$, choose one $K_{j_i}\subseteq K$ such that $r(x_i)\in s(K_{j_i})^o$. Note that $K_{j_i}x_i$ is a singleton and we denote by $\{\gamma_i\}=K_{j_i}x_i$ for simplicity. On the other hand, since $r(K_{j_i})\subseteq W$ and $\{\gamma_i\}=K_{j_i}x_i\subseteq (\bigcup\CT_-)u_\CT$ hold, there exists a unique $F_+\in \CT\unitSpace_+$ such that $r(\gamma_i)\in F_+\cap W$. This determines an injective map $\psi_\CT: \{x_1, \dots, x_{m_\CT}\} \to \{F_+\in \CT\unitSpace_+: F_+\cap W\neq \emptyset\}$. Note that 
condition \ref{prop:AE-DSC-non-fiberwise-condition1} above for $\CF$ implies for all $F_+\in \CF\unitSpace$, there exists an $O\in \CU$ such that $F_+\subseteq O$. Thus, by definition of $\CU$, if $F_+\cap W\neq \emptyset$ then $F_+\subseteq V$ holds necessarily. Define $\CY_\CT=\{F_+\in \CT\unitSpace_+: F_+\subseteq V\}$. Then for each fortified tower $\CT$ in $\CF$, the process above yields an injective map $\Phi_\CT=\psi_\CT\circ\varphi_\CT: (\CT^{-L})\unitSpace_+\to \CY_\CT$, which implies that $|(\CT^{-L})\unitSpace_+|\leq |\CY_\CT|$.

Now, define $\CY=\{F_+\in \CF\unitSpace_+: F_+\subseteq V\}$. Observe that for any fortified tower $\CT$ in $\CF$, one has $(\CF^{-L})\unitSpace_+\cap \CT_+\unitSpace=(\CT^{-L})\unitSpace_+$ and $\CY\cap \CT_+\unitSpace=\CY_\CT$. Then Lemma \ref{lem:castle-comparison} entails that $\bigcup (\CF^{-L})\unitSpace_+ \precsim_\CG \bigcup \CY\subseteq V$. Then, together with the condition \ref{prop:AE-DSC-non-fiberwise-condition4} above, Remark~\ref{rem:subequivalence-basic}\eqref{rem:subequivalence-basic:additive-prec} entails that $\GU\prec_\CG U\sqcup V\subseteq O_0$.

Finally, suppose the contrary that $M(\CG)\neq\emptyset$. Let $\mu\in M(\CG)$ and $O_1, O_2$ are disjoint non-empty open set in $\GU$.  Then Remark~\ref{rem:subequivalence-basic}\eqref{rem:subequivalence-basic:measure-auxiliary} implies that $1=\mu(\GU)\leq \mu(O_i)\leq 1$ for $i=1, 2$ because we have established $\GU\prec_\CG O_i$ for $i=1, 2$. But this implies that $\mu(O_1\sqcup O_2)=2$, which is a contradiction.
\end{proof}

We need one more lemma before the final conclusions of the section. 

\begin{lem} \label{lem:O-meas-inf}
	Let $\CG$ be a minimal \'{e}tale groupoid with a compact unit space and let $O$ be a nonempty open subset of $\GU$. Then $\displaystyle \inf_{\mu \in M(\CG)} \mu(O) > 0$. 
\end{lem}

\begin{proof}
	Since $O$ is compact Hausdorff, there is a nonempty open set $O'$ with $\overline{O'} \subseteq O$. It follows from Lemma~\ref{lem:measure-semicontinuous} that the function $M(\CG) \to [0,\infty)$ given by $\mu \mapsto \mu(\overline{O'})$ is upper semicontinuous and thus has a minimal, which must be nonzero since $\mu(\overline{O'}) \geq \mu({O'}) > 0$ by the minimality of $\CG$ (which can be seen, e.g., by applying Lemma~\ref{lem:minimal-recurrence}). 
\end{proof}

\begin{thm}\label{thm:AE-DSC} 
	Let $\CG$ be a minimal $\sigma$-compact \'{e}tale groupoid
	with a compact unit space. 
	Then $\CG$ is almost elementary if and only if it is almost elementary in measure and has groupoid strict comparison. 
\end{thm}

\begin{proof}
	When $\GU$ is finite, the desired equivalence follows from Corollary~\ref{cor:AE-finite-pair} and Proposition~\ref{prop:AEm-finite-pair} as well as Example~\ref{exa:GSBP-pair-groupoid}. 
	
	Now suppose $\GU$ is infinite. 
	In view of Lemmas~\ref{lem:O-meas-inf} and~\ref{lem:O-meas-sup}, by comparing Definitions~\ref{defn:AE-in-measure} and~\ref{defn:AE}, we see that almost elementariness implies almost elementariness in measure and the converse holds in the presence of groupoid strict comparison. 
	Finally, combining Proposition \ref{prop:AE-DSC-fiberwise} and \ref{prop:AE-DSC-non-fiberwise} gives us that almost elementariness implies groupoid strict comparison. This completes the proof. 
\end{proof}

\begin{cor}
	Let $\CG$ be as in Theorem~\ref{thm:AE-DSC}. Suppose $M(\CG) = \varnothing$. Then $\CG$ is almost elementary if and only if it is effective and has groupoid strict comparison. 
\end{cor}

\begin{proof}
	This follows immediately from Theorems~\ref{thm:AE-DSC} and~\ref{thm:effective}. 
\end{proof}

\begin{rem}\label{rem:purely-infinite}
    By \cite[Theorem~5.1]{Ma-purely}, a groupoid $\CG$ as in Theorem~\ref{thm:AE-DSC} satisfies both $M(\CG) = \varnothing$ and groupoid strict comparison if and only if it has \emph{paradoxical decomposition} in the sense that for any non-empty open set $O$ in $\GU$ and any compact subset $K \subseteq O$, there are disjoint open subset $O_1, O_2 \subseteq O$ such that $K \prec_{\CG} O_1$ and  $K \prec_{\CG} O_2$ (see also \Cref{rem:subequivalence-change}\eqref{rem:subequivalence-change:purely_infinite}). This latter notion generalizes Matui's \emph{pure infiniteness}  \cite{Matui2015Topological} for ample groupoids (see \cite[Corollary~5.5]{Ma-purely} and also \Cref{rem:subequivalence-change}\eqref{rem:subequivalence-change:same}). In particular, a purely infinite minimal ample \'{e}tale groupoid with a compact unit space is almost elementary. 
\end{rem}

\section{A relative remainder control} \label{sec:relative}

We introduce an equivalent characterization of almost elementary groupoids, using a ``relative'' remainder control, i.e., by comparing the remainder with (a small portion of) the castle in the almost elementary approximation. This turns out to be useful in a number of the results. It also has the advantage that it is a meaningful concept without assuming $\CG$ is minimal, unlike the original Definition~\ref{defn:AE}, where in a non-minimal groupoid, the remainder may be forced to be empty. 

For the moment, let us distinguish this characterization from the original one by calling it almost elementariness with relative remainder control. 

\begin{defn} \label{defn:AE-relative} 
	Let $\CG$ be an \'{e}tale groupoid with a compact unit space. 
	We say that $\CG$ is \emph{almost elementary with relative remainder control}
	if for any nonempty compact set $K$ in $\CG$, 
	any $\varepsilon > 0$, 
	and any finite open cover $\CV$ of $\CG^{(0)}$, 
	there is a nonempty open castle $\CC$
	satisfying 
	\begin{enumerate}[label=(\roman*)]
		\item \label{defn:AE-relative:extendable} $\CC$ is $K$-extendable; 
		\item \label{defn:AE-relative:fine} $\CC\unitSpace$ refines $\CV$; 
		\item \label{defn:AE-relative:remainder} there is a subcollection $\CX \subseteq \CC^{(0)}$ such that $|\CX \cap \CT^{(0)}| \leq \varepsilon |\CT^{(0)}|$ for every tower $\CT$ in $\CC$ and
		$\CG^{(0)}\setminus\ \left( \bigcup \CC^{(0)} \right) \prec_{\CG} \bigcup \CX$.  
	\end{enumerate}
\end{defn}

\begin{rem}\label{rem: AE-relative-property-full}
Observe that condition~\ref{defn:AE-relative:remainder} in Definition~\ref{defn:AE-relative} implies that $\bigcup \CC\unitSpace$ is full in $\CG$. Indeed, by Definition~\ref{defn:subequivalence}, it follows from $\CG^{(0)}\setminus\ \left( \bigcup \CC^{(0)} \right) \prec_{\CG} \bigcup \CX$ that $\CG^{(0)}\setminus\ \left( \bigcup \CC^{(0)} \right) \subseteq r(\CG \cdot (\bigcup \CX)) \subseteq r(\CG \cdot (\bigcup \CC\unitSpace))$, whence $\CG^{(0)} \subseteq r(\CG \cdot (\bigcup \CC\unitSpace))$. 
	\end{rem}

One advantage of Definition~\ref{defn:AE-relative} is that it relates to almost elementariness in measure in a more direct way than Definition~\ref{defn:AE} does.

\begin{prop} \label{prop:AE-relative-in-measure} 
	If an \'{e}tale groupoid $\CG$ with a compact unit space is almost elementary with relative remainder control, then it is almost elementary in measure. 
\end{prop}
\begin{proof}
    With the help of Remark \ref{rem: AE-relative-property-full}, it is left to show $\mu(\GU\setminus \bigcup \CC^{(0)})<\varepsilon$ for any $\mu\in M(\CG)$. To this end, let $\mu\in M(\CG)$ and $\CX$ be in Definition \ref{defn:AE-relative}\ref{defn:AE-relative:remainder}. Note that 
    \begin{align*}
        \mu(\bigcup \CX)&=\sum_{\CT\text{ is a tower in }\CC} \mu(\bigcup (\CT^{(0)}\cap \CX))\\
        &\leq \sum_{\CT\text{ is a tower in }\CC}\varepsilon\cdot \mu(\bigcup \CT^{(0)})=\varepsilon\cdot \mu(\bigcup \CC^{(0)})\leq \varepsilon.
        \end{align*}
   Then $\GU\setminus \bigcup \CC^{(0)}\prec_\CG \bigcup\CX$ entails that $\mu(\GU\setminus \bigcup \CC^{(0)})\leq \varepsilon$ by Remark \ref{rem:subequivalence-basic}\eqref{rem:subequivalence-basic:measure-subequivalence}.
   \end{proof}

We shall go through a number of results about this variant that are parallel to those in Section~\ref{sec:almost-elementary}, before showing that it is an equivalent characterization in the minimal setting. 

The first result we need aims at using the notion of $K$-shrinking directly instead of that of $K$-extendability. We first record several general properties of splitting a castle in the sense of Lemma \ref{lem:castle-split}.

\begin{prop}\label{prop:castle-split-more}
    Let $\CC, \CD$ be castles in $\CG$ such that $\CD$ refines $\CC$ and $\bigcup \CD=\bigcup \CC$.  Then the following hold.
    \begin{enumerate}[label=(\roman*)]
        \item\label{prop:castel-split-more-general}
        For any $\CX\subseteq \CC$,  if we denote by 
        \[\CX^{\CD}=\{D\in \CD: D\subseteq C
        \text{ for some } C\in \CX\},\]
         then $\bigcup \CX=\bigcup \CX^{\CD}$ holds.
         \item\label{prop:castel-split-more-towers}         Let $\CT$ be a tower in $\CC$ and $\CT'$ a tower in $\CD$ refining $\CT$. Then one has $|\CT^{(0)}|=|\CT'^{(0)}|$. Moreover, for any       
         $\CX\subseteq \CC^{(0)}$, one has $|\CX\cap \CT^{(0)}|=|\CX^{\CD}\cap \CT'^{(0)}|$.
             \end{enumerate}
\end{prop}
\begin{proof}
   First, let $C\in \CC$ and $D\in \CD$ such that $C\cap D\neq \emptyset$. Then $C$ is the unique member such that $D\subseteq C$. Indeed, since $\CD$ refines $\CC$, for the $D\in \CD$, there exists a $C_0\in \CC$ such that $D\subseteq C_0$. Now because $\CC$ is a disjoint family, one has that $D$ is disjoint from all other bisections in $\CC$. This thus implies that $C_0=C$.

   Now let $C\in \CC$. First, $\bigcup\{D\in \CD, D\subseteq C\}\subseteq C$ holds trivially. For the converse direction, let $x\in C$. Because of $\bigcup \CC=\bigcup \CD$, there exists a $D\in \CD$ such that $x\in D$, which entails $x\in C\cap D$. Thus, $x\in D\subseteq C$ holds necessarily and thus one has $C=\bigcup\{D\in \CD, D\subseteq C\}$.  
    Therefore, let $\CX$ be a subfamily of $\CC$ and denote by $\CX^{\CD}=\{D\in \CD: D\subseteq C \text{ for some }C\in \CX\}$, one has $\bigcup \CX=\bigcup \CX^{\CD}$. This has established \ref{prop:castel-split-more-general}.

    Let $\CT$ be a tower in $\CC$ and  $C_0\in \CT$.  Then for any $\CD$-member $D_0\subseteq C_0$, we claim that the collection $\CT'=\{C_1D_0C_2: C_1, C_2\in \CC\}\setminus \{\emptyset\}$ is a tower in $\CD$ refining $\CT$. Indeed, Lemma \ref{lem:castle-split}\eqref{lem:castle-split:refine-basic} implies that $\CT'\subseteq \CD$ and actually $\CT'$ is contained in a $\CD$-tower by the definition of $\sim_{\CD}$ recorded in Lemma \ref{lem:plurisection-basic}\eqref{lem:plurisection-basic:towers}. 
    We now show that such $\CT'$ is already a $\sim_{\CD}$-class and thus a $\CD$-tower. To this end, let $D\sim_\CD D_0$, which means that there exists an $E\in \CD$ such that $s(D)=r(E)$ and $s(E)=r(D_0)$. This further entails that $D=(DE)\cdot D_0\cdot (ED_0)^{-1}$. Because $\CD$ refines $\CC$, there exists $C_1, C_2\in \CC$ such that $DE\subseteq C_1$ and $(ED_0)^{-1}\subseteq C_2$. Since $C_1$ and $C_2$ are bisections, one has that $D=C_1\cdot D_0\cdot C_2$ and therefore $D\in \CT'$. This shows that $\CT'$ is already an equivalence class of the relation $\sim_\CD$ and thus a $\CD$-tower. Finally, because $C\in \CT$ and $D_0\subseteq C_0$, if either $C_1$ or $C_2$ is not in $\CT$, then $C_1D_0C_2=\emptyset$. Therefore, $\CT'$ refines $\CT$ and actually $\CT'=\{C_1D_0C_2: C_1, C_2\in \CT\}\setminus \{\emptyset\}$. 
    
    We then show that all $\CD$-tower $\CT'$ refining a $\CC$-tower $\CT$ is of the form $\{C_1D_0C_2: C_1, C_2\in \CT
    \}\setminus \{\emptyset\}$ for some fixed $D_0\in \CT'$. Indeed, note that any $D\in \CT'$ can be written as $D=D_1D_0D_2$ for some $D_1, D_2\in \CT'$. Choose $C_1, C_2\in \CT$ such that $D_1\subseteq C_1$ and $D_2\subseteq C_2$. Then one has $D=C_1D_0C_2$ since $C_1, C_2, D_1, D_2$ are all bisections.    
    
    Now, let $\CT$ be a $\CC$-tower and $C_0\in \CT$. Let  $\CT'=\{C_1D_0C_2: C_1,C_2\in \CC\}\setminus \{\emptyset\}$ be an arbitrary $\CD$-tower refining $\CT$ for a $D_0\in \CD$ with $D_0\subseteq C_0$. We define a map $f: \CT'\to \CT$ by $f(C_1D_0C_2)=C_1C_0C_2$. This map is well-defined because $C=C_1C_0C_2$ is the unique bisection in $\CT$ containing $D=C_1D_0C_2$. Moreover, the definition of $f$ shows that $f$ is surjective. For injectivity, let $D_1, D_2\in \CT'$ with $D_1\neq D_2$. This means that $s(D_1)\cap s(D_2)=\emptyset$ or $r(D_1)\cap r(D_2)=\emptyset$.  Then, write $D_1=C_1D_0C_2$ and $D_2=C'_1D_0C'_2$ for some $C_1, C_2, C'_1, C'_2\in \CT$. Then necessarily one has $r(C_1)\neq r(C'_1)$ or $s(C_2)\neq s(C'_2)$. But this implies that 
    \[f(D_1)=C_1C_0C_2\neq C'_1C_0C'_2=f(D_2)\]
     and thus $f$ is injective. Moreover, note that by definition, $f(\CT'^{(0)})=\CT^{(0)}$ and therefore one has $|\CT'^{(0)}|=|\CT^{(0)}|$.

Let $\CX\subseteq \CC^{(0)}$ and write $\CX^\CD=\{D\in \CD: D\subseteq C \text{ for some }C\in \CX\}$. Let $\CT$ be a $\CC$-tower and $\CT'$ a $\CD$-tower refining $\CT$.  By definition of $f: \CT'\to \CT$, one has $f(\CX^\CD\cap \CT'^{(0)})=\CX\cap \CT^{(0)}$. This shows that $|\CX^\CD\cap \CT'^{(0)}|=|\CX\cap \CT^{(0)}|$. We thus have established \ref{prop:castel-split-more-towers}.
 \end{proof}

\begin{lem} \label{lem:AE-relative-shrinking}  
	A groupoid $\CG$ with a compact unit space is {almost elementary with relative remainder control}
	if and only if for any nonempty compact set $K, L$ in $\CG$, 
	any $\varepsilon > 0$, 
	and any finite open cover $\CV$ of $L$, 
	there is an open castle $\CD$
	satisfying 
	\begin{enumerate}[label=(\roman*)]
		\item\label{lem:AE-relative-shrinking-refine}   for any $D\in \CD$ with 
  $D\cap L\neq \emptyset$, there exists a $V\in \CV$ such that $D\subseteq V$.

		\item \label{lem:AE-relative-shrinking-small reminder}   
		there is a subcollection $\CX \subseteq \left( \CD^{-K} \right)^{(0)}$ such that $|\CX \cap \CT^{(0)}| \leq \varepsilon |\CT^{(0)}|$ for every tower $\CT$ in $\CD^{-K}$ and
		$\CG^{(0)}\setminus\ \left( \bigcup \left( \CD^{-K} \right)^{(0)} \right) \prec_{\CG}  \bigcup \CX$.  
	\end{enumerate}
\end{lem}

\begin{proof}
 The ``if'' part is trivial by choosing $\CC=\CD^{-K}$ and $L=\GU$. For the reverse direction, we follow the proof of \eqref{lem:AE-in-measure-shrinking:original} $\Rightarrow$ \eqref{lem:AE-in-measure-shrinking:fine-entourage} of Lemma~\ref{lem:AE-in-measure-shrinking} almost verbatim with the help of Proposition \ref{prop:castle-split-more}.  The only significant modifications are the following. First, let $K, \CV, \varepsilon$ in the statement be given. Then the same construction in the proof of  Lemma~\ref{lem:AE-in-measure-shrinking} yields a compact set $K'$ containing $K^2$,  an open cover $\CW$ of $\GU$. Then almost elementariness with relative remainder control for $\CG$ implies that there exists an open castle $\CC$ satisfying
 \begin{enumerate}
     \item $\CC$ is $K'$-extendable;
     \item $\CU$ refines $\CW$;
     \item there is a subcollection $\CY\subseteq \CC^{(0)}$ such that $|\CY\cap \CS^{(0)}|\leq \varepsilon|\CS^{(0)}|$ for every tower $\CS$ in $\CC$ and 
     $\GU \setminus \left( \bigcup \CC^{(0)} \right) \prec_\CG \bigcup \CY$.
 \end{enumerate}
By Definition \ref{defn:castle-shrinking}, there exists an open castle $\CE$ such that  $\CC\subseteq \CE^{-K'}$. Then for the exactly same partition $\CP$ constructed in Lemma~\ref{lem:AE-in-measure-shrinking}  for $\bigcup \CE$, Lemma \ref{lem:castle-split}\eqref{lem:castle-split:exists} and \eqref{lem:castle-split:extra} imply that there exists an open castle $\CE_\CP$ that refines both $\CE$ and $\CP$ and satisfies $\bigcup \CE_\CP=\bigcup \CE$.
 
 Then still following the same constructions in the proof of  Lemma~\ref{lem:AE-in-measure-shrinking}, one obtains a castle $\CD$ as a subcastle of $\CE_\CP$ such that $\bigcup \CC\subseteq \bigcup \CD^{-K}$ so that \ref{lem:AE-relative-shrinking-refine} holds.

 For \ref{lem:AE-relative-shrinking-small reminder}, Recall that $\CC\subseteq \CE^{-K'}\subseteq \CE$ and $\bigcup \CE=\bigcup \CE_\CP$ as well as $\CE_\CP$ refines $\CE$. Then Proposition \ref{prop:castle-split-more}\ref{prop:castel-split-more-general} implies that the family
 \[\CD'=\{D\in \CE_\CP: D\subseteq C \text{ for some }C\in \CC\}\]
 satisfies $\bigcup \CD'=\bigcup \CC\subseteq \bigcup \CD^{-K}$. Therefore, $\CD'$ is a subcastle of $\CD^{-K}$. Let $\CT$ be a tower in $\CD^{-K}$, we denote by $\CT_{\CD'}$ its natural restriction to $\CD'$, which is a thus a $\CD'$-tower. Now for the $\CY\subseteq \CC$ above, we denote by 
 \[\CX=\CY^{\CD'}=\{D\in \CD': D\subseteq C \text{ for some }C\in \CY\},\]
 which is a subfamily of $D^{-K}$ satisfying $\bigcup \CX=\bigcup \CY$ by Proposition \ref{prop:castle-split-more}\ref{prop:castel-split-more-general}. 
 Moreover, Proposition \ref{prop:castle-split-more}\ref{prop:castel-split-more-towers} shows that $|\CX\cap \CT^{(0)}_{\CD'}|\leq \varepsilon|\CT^{(0)}_{\CD'}|$.
 Also, observe that $\CX\cap \CT^{(0)}=\CX\cap \CT^{(0)}_{\CD'}$ and $\CT_{\CD'}\subseteq \CT$. These further entail
 \[|\CX\cap \CT^{(0)}|=|\CX\cap \CT^{(0)}_{\CD'}|\leq \varepsilon|\CT^{(0)}_{\CD'}|\leq \varepsilon|\CT^{0}|.\]
 Finally, one has
	\[
	\GU \setminus \left( \bigcup \left( \CD^{-K} \right)^{(0)} \right) \subseteq \GU \setminus \left( \bigcup \CC^{(0)} \right) \prec_\CG \bigcup \CY=\bigcup \CX \; ,
	\]
	which implies $\GU \setminus \left( \bigcup \left( \CD^{-K} \right)^{(0)} \right) \prec_\CG \bigcup\CX$ by Remark~\ref{rem:subequivalence-basic}\eqref{rem:subequivalence-basic:subset} and~\eqref{rem:subequivalence-basic:preorder}. This has established \ref{lem:AE-relative-shrinking-small reminder}. 
\end{proof}

\begin{cor} \label{cor:AE-relative-finite-pair}
	Suppose $\CG^{(0)}$ is finite. Then $\CG$ is almost elementary with relative remainder control if and only if it is a finite union of finite {pair groupoids}. 
\end{cor}

\begin{proof}
	First, if $\CG$ is a finite union of pair groupoid, then for any $K\subseteq \CG, \varepsilon>0$ and open cover $\CV$ of $\GU$,  one may verify Definition \ref{defn:AE-relative} by taking $\CC=\{\{x\}: x\in\CG\}$ since $\CC^{-K}=\CC$ and $\GU\setminus \bigcup \CC^{(0)}=\emptyset$. Thus $\CG$ is almost elementary with relative remainder control. For the converse, by Proposition \ref{prop:AE-relative-in-measure}, the groupoid $\CG$ with a finite unit space, is actually almost elementary in measure. Then Proposition \ref{prop:AEm-finite-pair} shows that $\CG$ is a finite union of finite pair groupoids.
\end{proof}

We also need a ``fortified'' version. 

\begin{lem} \label{lem:AE-relative-fortified} 
	An \'{e}tale groupoid $\CG$ with a compact unit space  is {almost elementary with relative remainder control}
	if and only if 
	for any compact sets $K$ and $L$ in $\CG$, 
	any $\varepsilon > 0$, 
	and any finite open cover $\CV$ of $\CG^{(0)}$, 
	there is a fortified castle $\CF$
	satisfying 
	\begin{enumerate}[label=(\roman*)]
		\item \label{lem:AE-relative-fortified:refine} 
		for any $F_+ \in \CF_+$ with $F_+ \cap L \not= \varnothing$, there is $V \in \CV$ such that $F_+ \subseteq V$; 
		\item\label{lem:AE-relative-fortified:precompact} 
		$\bigcup \CF_+$ is precompact in $\CG$;  
		\item \label{lem:AE-relative-fortified:remainder} 
		there is a subcollection $\CW \subseteq \left( \CF^{-K} \right)_-^{(0)}$ such that $|\CW \cap \CT^{(0)}| \leq \varepsilon |\CT^{(0)}|$ for every tower $\CT$ in $\left( \CF^{-K} \right)_-$ and
		$\CG^{(0)}\setminus\ \left( \bigcup \left( \CF^{-K} \right)_-^{(0)} \right) \prec_{\CG}  \bigcup \CW$.   \end{enumerate}
\end{lem}
\begin{proof}
	We follow the proof of Lemma~\ref{lem:AE-fortified} almost verbatim, with the only significant modification being constructing $f_{\uparrow\CC}$ with additional properties.
 ``If'' part is straightforward by Lemma \ref{lem:AE-relative-shrinking} since the castle $\CF_-$ plays a role of $\CD$ in Lemma \ref{lem:AE-relative-shrinking}. 
 For the converse, suppose $\CG$ is almost elementary with relative remainder control, Lemma \ref{lem:AE-relative-shrinking} implies that for compact sets $K, L$ in $\CG$, number $\varepsilon>0$ and an open cover $\CV$ for $L$, there exists a castle  $\CD$ satisfying
\begin{enumerate}
		\item for any $D\in \CD$ with 
  $D\cap L\neq \emptyset$, there exists a $V\in \CV$ such that $D\subseteq V$.
\item There is a subcollection $\CX \subseteq \left( \CD^{-K} \right)^{(0)}$ such that $|\CX \cap \CT^{(0)}| \leq \varepsilon |\CT^{(0)}|$ for every tower $\CT$ in $\CD^{-K}$ and
		$\CG^{(0)}\setminus\ \left( \bigcup \left( \CD^{-K} \right)^{(0)} \right) \prec_{\CG}  \bigcup \CX$.  
	\end{enumerate} 
 By Definition \ref{defn:subequivalence} and Remark \ref{rem:subequivalence-basic}\eqref{rem:subequivalence-basic:compact}, there exists an open $s$-section $S$ and a compact set $M\subseteq \bigcup \CX$ such that $\CG^{(0)}\setminus\ \left( \bigcup \left( \CD^{-K} \right)^{(0)} \right)\subseteq r(S)$ and $s(S)\subseteq M$. By shrinking $S$ if necessary, one may assume that $r(S)\cap M=\emptyset$. Indeed, to do this, note that the compact set $\CG^{(0)}\setminus\ \left( \bigcup \left( \CD^{-K} \right)^{(0)}\right)$ is disjoint from $M$. Then choose an open set $O$ in $\GU$ such that $\CG^{(0)}\setminus\ \left( \bigcup \left( \CD^{-K} \right)^{(0)}\right)\subseteq O$ and $O\cap M=\emptyset$  and define $S'=S\cap r^{-1}(O)$. Note that $S'$ is still an open $s$-section and one has $\CG^{(0)}\setminus\ \left( \bigcup \left( \CD^{-K} \right)^{(0)}\right)\subseteq r(S')$ and $r(S')\cap M=\emptyset$. In addition, note that  $s(S')\subseteq s(S)\subseteq M$. Then $S'$ is the $s$-section that we need.

 Then choose a continuous function $f:\GU\to [0, 1]$ that sends $\CG^{(0)}\setminus\ \left( \bigcup \left( \CD^{-K} \right)^{(0)} \right)$ to $0$ and $\GU\setminus r(S)$ to $1$, respectively. We apply Lemma~\ref{lem:castle-invarint-function-smallest} to obtain $f_{\uparrow\CC} \in C_0 \left( \bigcup \CD , \mathbb{R} \right)^{\bigcup \CD}$ satisfying $0 \leq f \leq f_{\uparrow\CD} \leq \sup_{x \in \bigcup \CD} f(x) = 1$.  
	Hence applying Remark~\ref{rem:fortified-castle-basics}\eqref{rem:fortified-castle-basics:construction} to $f_{\uparrow\CD}$, 
	we obtain a fortified castle 
	\[
		\CF := \left\{ \left( f_{\uparrow\CD}^{-1} \left( \left(\frac{1}{3},1 \right] \right) \cap D, f_{\uparrow\CD}^{-1} \left( \left(\frac{1}{2},1 \right] \right) \cap D \right) \colon D \in \CD  \right\} 
	\]
	with $\CF_+$ refining $\CD$. 
	It follows that $\CF_+$ refines $\CV$ and $\CF^{-K}$ is nonempty. Define $\CW=\{f_{\uparrow\CC}^{-1} \left( \left(\frac{1}{2},1 \right] \right) \cap D: D\in \CX\}$, which is a subset of $(\CF^{-K})^{(0)}_-$. Moreover, by the construction of $\CF$, it is straightforward to obtain $|\CW\cap \CT^{(0)}|\leq \varepsilon|\CT^{(0)}|$ for every tower $\CT$ in $(\CF^{-K})_-$.

On the other hand, observe that $f$ is constant $1$ on $M$ and thus so is $f_{\uparrow\CC}$ by definition.  This implies that $M\subseteq \bigcup \CW$ because
$M\subseteq \bigcup \CX$. Then,  using the same argument in the final part of Lemma~\ref{lem:AE-fortified}, one observes that $\CG^{(0)}\setminus\ \left( \bigcup \left( \CF^{-K} \right)_-^{(0)} \right)\subseteq r(S)$. This implies that 
 $\CG^{(0)}\setminus\ \left(\bigcup \left( \CF^{-K} \right)_-^{(0)}\right)\prec_\CG \bigcup \CW$ because $s(S)\subseteq M\subseteq \bigcup \CW$
 \end{proof}

Another additional advantage of Definition~\ref{defn:AE-relative} is that it combines well with ubiquitous fiberwise amenability to give an equivalent characterization using castles with ``F\o{}lner shapes''. The equivalent conditions below may be understood as a generalized definition of \emph{almost finiteness} (see Sections~\ref{sec:Matui} and~\ref{sec:Kerr}).

\begin{prop} \label{prop:AE-relative-fiberwise-fortified} 
	Let $\CG$ be a $\sigma$-compact \'{e}tale groupoid with a compact unit space. 
	Then the following are equivalent: 
	\begin{enumerate}
		\item\label{prop:AE-relative-fiberwise-fortified:conjunction}  The groupoid $\CG$ is ubiquitously fiberwise amenable and almost elementary with relative remainder control.

        \item \label{prop:AE-relative-fiberwise-fortified:Foelner-0} 
        for any nonempty compact set $K$ in $\CG$, 
		any $\varepsilon > 0$, 
		and any finite open cover $\CV$ of $K$, 
		there is a nonempty open castle $\CC$
		satisfying 
		\begin{enumerate}[label=(\roman*)]
			\item \label{prop:AE-relative-fiberwise-fortified:Foelner-0:refine} For any $C\in\CC$ with $C\cap K\neq \emptyset$, there exists a $V\in \CV$ such that $C\subseteq V$;  
			\item \label{prop:AE-relative-fiberwise-fortified:Foelner-0:remainder} there is a subcollection $\CX \subseteq \CC^{(0)}$ such that $|\CX \cap \CT^{(0)}| \leq \varepsilon |\CT^{(0)}|$ for every tower $\CT$ in $\CC$ and
			$\GU \setminus \left( \bigcup \CC^{(0)} \right) \prec_{\CG}  \bigcup \CX$; 
			\item \label{prop:AE-relative-fiberwise-fortified:Foelner-0:Foelner} for any $u \in \CG\unitSpace$, we have $|K (\bigcup \CC) u| \leq (1+\varepsilon) |(\bigcup \CC) u|$. 
		\end{enumerate}

		\item\label{prop:AE-relative-fiberwise-fortified:Foelner-fortified} 
            for any nonempty compact sets $K, L$ in $\CG$, 
		any $\varepsilon > 0$, 
		and any finite collection $\CV$ of open sets in $\CG$ that covers $L$, 
		there is a nonempty fortified castle $\CF$
		satisfying 
		\begin{enumerate}[label=(\roman*)]
			\item\label{prop:AE-relative-fiberwise-fortified:Foelner-fortified:refine} for any $C \in \CF_+$ with $C \cap L \not= \varnothing$, there is $V \in \CV$ such that $C \subseteq V$;  
			\item\label{prop:AE-relative-fiberwise-fortified:Foelner-fortified:precompact} $\bigcup \CF_+$ is precompact in $\CG$;  
			\item\label{prop:AE-relative-fiberwise-fortified:Foelner-fortified:remainder} there is a subcollection $\CX \subseteq \CF_-^{(0)}$ such that $|\CX \cap \CT^{(0)}| \leq \varepsilon |\CT^{(0)}|$ for every tower $\CT$ in $\CF_-$ and
			$\GU \setminus \left( \bigcup \CF_-^{(0)} \right) \prec_{\CG}  \bigcup \CX$; 
			\item\label{prop:AE-relative-fiberwise-fortified:Foelner-fortified:Foelner} for any $u \in \CG\unitSpace$, we have $|\partial_K ((\bigcup \CF_+) u)| \leq \varepsilon |(\bigcup \CF_+) u|$, where the notation $\partial_K$ is defined in \Cref{def:groupoid-boundaries}.  
		\end{enumerate}
	\end{enumerate}
\end{prop}

We remark that the main difficulty in proving the above result lies in upgrading ``fiberwise F\o{}lner conditions'' such as \eqref{prop:AE-relative-fiberwise-fortified:Foelner-0}\ref{prop:AE-relative-fiberwise-fortified:Foelner-0:Foelner} and \eqref{prop:AE-relative-fiberwise-fortified:Foelner-fortified}\ref{prop:AE-relative-fiberwise-fortified:Foelner-fortified:Foelner} to a type of ``ladder-wise F\o{}lner condition'' (see \eqref{prop:AE-relative-fiberwise-fortified::proof:Foelner-fortified:shrunken-size} below). 
To overcome this difficulty, we shall need the following simple lemma in point-set topology. 

\begin{lem}[{\cite[Lemma~4.11]{HirshbergWu2023Long}}]\label{lem:cover-clumping}
	Let $ X $ be a locally compact Hausdorff space and let $K \subseteq X $ be a compact subset. Let $\CV$ be a finite collection of open sets in $X$ such that $K \subseteq \bigcup \CV$. 
	Then there exists a finite collection $\CW$  of open sets in $X$ such 
	that $K \subseteq \bigcup \CW$ and, for any subcollection $\CW' 
	\subseteq \CW$, if $\bigcap \CW' \not= \varnothing$, then there 
	exists $V \in \CV$ such that $\bigcup \CW' \subseteq V$. 
\end{lem}

\begin{proof}[Proof of \Cref{prop:AE-relative-fiberwise-fortified}]
	The implication \eqref{prop:AE-relative-fiberwise-fortified:Foelner-fortified} $\Rightarrow$ \eqref{prop:AE-relative-fiberwise-fortified:Foelner-0} is obvious when we set $L = \CG^{(0)}$ and observe that $K (\bigcup \CF_+) u \subseteq (\bigcup \CF_+) u \cup \partial_K ((\bigcup \CF_+) u)$ for any $u \in \CG^{(0)}$.

	Next we show the implication \eqref{prop:AE-relative-fiberwise-fortified:Foelner-0} $\Rightarrow$ \eqref{prop:AE-relative-fiberwise-fortified:Foelner-fortified}. Let $K, L, \CV$ and $\varepsilon>0$ be given. Define $K_0=K\cup L$ and $\delta=\varepsilon/(1+\sup\{|K^{-1}_0u|: u\in \GU\})$. Then \eqref{prop:AE-relative-fiberwise-fortified:Foelner-0} implies that there exists an open castle $\CC$ such that 
\begin{enumerate}[label=(\roman*)]
			\item  For any $C\in\CC$ with $C\cap K_0\neq \emptyset$, there exists a $V\in \CV$ such that $C\subseteq V$;  
			\item  there is a subcollection $\CX \subseteq \CC^{(0)}$ such that $|\CX \cap \CT^{(0)}| \leq \varepsilon |\CT^{(0)}|$ for every tower $\CT$ in $\CC$ and
			$\GU \setminus \left( \bigcup \CC^{(0)} \right) \prec_{\CG}  \bigcup \CX$; 
			\item for any $u \in \CG\unitSpace$, we have $|K_0 (\bigcup \CC) u| \leq (1+\delta) |(\bigcup \CC) u|$. 
		\end{enumerate}
    Following the proof of Lemma \ref{lem:AE-relative-fortified}, we obtain a fortified castle
    \[
		\CF := \left\{ \left( f_{\uparrow\CC}^{-1} \left( \left(\frac{1}{3},1 \right] \right) \cap C, f_{\uparrow\CC}^{-1} \left( \left(\frac{1}{2},1 \right] \right) \cap C \right) \colon C \in \CC  \right\} 
	\]
for some properly chosen function $f_{\uparrow\CC}$ such that the first three conditions in \eqref{prop:AE-relative-fiberwise-fortified:Foelner-fortified} hold for $\CF$. It remains to show \ref{prop:AE-relative-fiberwise-fortified:Foelner-fortified:Foelner} in \eqref{prop:AE-relative-fiberwise-fortified:Foelner-fortified}, but this follows from the inequalities
\begin{align*}
    |\partial_{K}(\bigcup\CF_+)u|&\leq|\partial_{K_0}(\bigcup\CF_+)u|\leq (1+\sup_{u\in \GU}|K_0^{-1}u|)\cdot|\partial^+_{K_0}(\bigcup\CF_+)u|\\
    &=(1+\sup_{u\in \GU}|K_0^{-1}u|)\cdot|K_0(\bigcup\CF_+)u\setminus(\bigcup\CF_+)u|\leq \varepsilon|\bigcup(\CF_+)u|
\end{align*}
for any $u\in \bigcup\CF_+\unitSpace$ by Remark \ref{rem: double control for boundries} and the condition \ref{prop:AE-relative-fiberwise-fortified:Foelner-0:Foelner} in \eqref{prop:AE-relative-fiberwise-fortified:Foelner-0}. Note that if $u\notin \bigcup\CF_+\unitSpace$, the inequality is trivial.

The implication \eqref{prop:AE-relative-fiberwise-fortified:conjunction}$\Rightarrow$\eqref{prop:AE-relative-fiberwise-fortified:Foelner-fortified} follows from a similar argument to the proof of Lemma \ref{lem:AE-fiberwise-fortified}. Suppose \eqref{prop:AE-relative-fiberwise-fortified:conjunction} holds. Let the compact set $K$, $\varepsilon>0$, and the open cover $\CV$ of $K$ be given. Without loss of generality, one may assume $\GU\subseteq K$.  Then choose an cover $\CW$ of $K$, consisting of precompact open bisections and define  a new compact set $K'=\bigcup_{W\in \CW}\overline{W}$. By Theorem~\ref{thm:fiberwise-amenable-dichotomy}\ref{thm:fiberwise-amenable-dichotomy:amenable}, there is a compact set $L$ in $\CG$ such that for any finite set $M \subseteq\CG$, there
	is a finite set $N$ satisfying 
	\[
		M \subseteq N\subseteq LM \ \textrm{and}\ |K'N|\leq(1+\varepsilon)|N|.
	\]
Choose an open covering $\CO$ of $L\cup \GU$ consisting of  precompact open bisections and define $L'=\bigcup_{O\in \CO}\overline{O}$.  Define another open cover 
\[\CU := (\CW\cup \{\CG\setminus K\}) \vee (\CO\cup\{\CG\setminus (L\cup \GU)\})\] of $K \cup  L\supset \GU$. Then by Lemma \ref{lem:AE-relative-shrinking} implies that there exists an open castle $\CD$ such that
\begin{enumerate}[label=(\roman*)]
		\item for any $D\in \CD$ with 
  $D\cap (K\cup L)\neq \emptyset$, there exists a $U\in \CU$ such that $D\subseteq U$.
		\item 
		there is a subcollection $\CX \subseteq \left( \CD^{-K} \right)^{(0)}$ such that $|\CX \cap \CT^{(0)}| \leq \varepsilon |\CT^{(0)}|$ for every tower $\CT$ in $\CD^{-K'L'}$ and
		$\CG^{(0)}\setminus\ \left( \bigcup \left( \CD^{-K'L'} \right)^{(0)} \right) \prec_{\CG}  \bigcup \CX$.  
	\end{enumerate}
    Then the same argument as in Lemma \ref{lem:AE-fiberwise-fortified} allows us to choose a castle $\CC$ with  $\CD^{-K'L'}\subseteq\CC\subseteq \CD^{-K'}$ satisfying the requirements. 

    We finally show \eqref{prop:AE-relative-fiberwise-fortified:Foelner-fortified}  $\Rightarrow$\eqref{prop:AE-relative-fiberwise-fortified:conjunction}. To see   \eqref{prop:AE-relative-fiberwise-fortified:Foelner-fortified} implies $\CG$ is ubiquitously fiberwise amenable, we verify the equivalent characterization recorded in \cite[Proposition 5.5]{MWfiberwise}. Let $K$ be a compact set in $\CG$, $\varepsilon>0$ and a finite  open cover $\CV$ of $\GU$. Denote by $\CF$ the fortified castle and by $\CX\subseteq\CF_-\unitSpace$ satisfying \eqref{prop:AE-relative-fiberwise-fortified:Foelner-fortified}.  Then Remark \ref{rem:subequivalence-basic}\eqref{rem:subequivalence-basic:precompact-section} implies that there is a compact $s$-section $T$ such that $\GU\setminus (\bigcup \CF_-\unitSpace)= r(T)$ and $s(T)\subseteq \bigcup \CX$. Now define 
    \[L=(\bigcup\{\overline{F_-}: F_-\in \CF_-\})\cup (\bigcup\{\overline{F_-}: F_-\in \CF_-\})\cdot{T^{-1}},\]
    which is a compact set in $\CG$.
Then if $u\in \bigcup \CF_-\unitSpace$,  defining $F=(\bigcup\CF)u\subseteq Lu\bigcup \{u\}$, one has 
    \[|K(\bigcup\CF_-)u|\leq (1+\varepsilon)|(\bigcup \CF_-)u|.\]
Otherwise, $u\in \GU\setminus (\bigcup \CF\unitSpace)\subseteq r(T)$. Then choose one $x\in T^{-1}$ such that $s(x)=u$ and $r(x)\in \bigcup \CX$. In this case, define $F=(\bigcup \CF)x\subseteq Lu\cup \{u\}$. Then one has
\[|K(\bigcup \CF_-)x|=|K(\bigcup \CF_-)r(x)|\leq (1+\varepsilon)|(\bigcup \CF_-)r(x)|=(1+\varepsilon)|(\bigcup \CF_-)x|\]
as $r(x)\in \bigcup\CX\subseteq \bigcup\CF_-\unitSpace$. Thus, $\CG$ is ubiquitous fiberwise amenable by \cite[Proposition 5.5]{MWfiberwise}.

It is left to show  \eqref{prop:AE-relative-fiberwise-fortified:Foelner-fortified} implies $\CG$ is almost elementary with relative remainder control. 
Following \Cref{defn:AE-relative}, we let $\varepsilon \in (0,1)$, a nonempty compact set $K$ in $\CG$ and a finite open cover $\CV$ of $\CG^{(0)}$ be given. 
Without loss of generality, we assume $K = K^{-1}$. 
Choose $\delta \in (0,1) $ such that $\frac{2 \delta}{1 - \delta} \leq \varepsilon$, and 
choose a finite family $\CU$ of precompact open bisections such that $K\subseteq \bigcup \CU$. 
By \Cref{lem:cover-clumping}, there is a finite family $\CB$ of open sets such that $K\subseteq \bigcup \CB$ and for any $B_1, B_2 
	\in \CB$, if $B_1 \cap B_2 \not= \varnothing$, then there exists $U \in \CU$ such that $B_1 \cup B_2 \subseteq U$. 
    In particular, $\CB$ refines $\CU$, and thus each $B \in \CB$ is also a precompact open bisection. 
    By a standard argument using a partition of unity, we may choose an open set $Y_B$ for each $B \in \CB$ such that $\overline{Y_B} \subseteq B$ for each $B \in \CB$ and $K \subseteq \bigcup _{B \in \CB} Y_B$. 
Define a compact set  
$L:= \bigcup _{B \in \CB} \overline{ Y_B}$, which satisfies $K \subseteq \bigcup _{B \in \CB} Y_B \subseteq L \subseteq \bigcup \CB$. 
By \Cref{lem:entourage-fineness-control}, there is a finite open cover $\CW$ of $\CG\unitSpace$ such that for any $x \in K$ and any $W \in \CW$ containing $s(x)$, there is $B \in \CB$ such that $x \in Y_B$ and $W \subseteq s(Y_B)$. 
Write $\CW \vee \CV := \{ W \cap V \colon W \in \CW, V \in \CV \}$ for the common refinement of $\CW$ and $\CV$. 
Note that $(\CW \vee \CV ) \cup \{ \CG \setminus \CG\unitSpace\}$ forms a finite open cover of $\CG$. Define $\CA$ to be the common refinement of $(\CW \vee \CV ) \cup \{ \CG \setminus \CG\unitSpace\}$ and $\CB$, whence $L\subseteq \bigcup \CA$. 
Then the assumption in  \eqref{prop:AE-relative-fiberwise-fortified:Foelner-fortified} implies that there exists an open castle $\CD$ such that 
\begin{enumerate}[label=(\roman*$^{\circ}$)]
    \item\label{prop:AE-relative-fiberwise-fortified::proof:Foelner-fortified:refine}  for any $D\in \CD$ with $D\cap L\neq \emptyset$, there exists a $A \in \CA$ such that $D\subseteq A$ (in particular, there is $B \in \CB$ such that $D \subseteq B$ and, if $D \in \CD\unitSpace$, then there are also $W \in \CW$ and $V \in \CV$ such that $D \in W \cap V$); 
    \item \label{prop:AE-relative-fiberwise-fortified::proof:Foelner-fortified:precompact} there is a subcollection $\CX\subseteq \CD\unitSpace$ such that $|\CX\cap \CT\unitSpace|\leq \delta \, |\CT\unitSpace|$ for any tower $\CT$ in $\CD$ and $\GU\setminus \bigcup \CD\unitSpace\prec_\CG \bigcup \CX$; 
    \item\label{prop:AE-relative-fiberwise-fortified::proof:Foelner-fortified:remainder} the condition $|\partial_L (\bigcup \CD)u|\leq \delta \,|(\bigcup \CD)u|$ holds for any $u\in \GU$. 
\end{enumerate}
Define the shrunken castle $\CC := \CD^{-K}$, which is then $K$-extendable by definition, verifying \Cref{defn:AE-relative}\ref{defn:AE-relative:extendable}. It follows from \ref{prop:AE-relative-fiberwise-fortified::proof:Foelner-fortified:refine} that $\CC$ also satisfies \Cref{defn:AE-relative}\ref{defn:AE-relative:fine}. 
To verify \Cref{defn:AE-relative}\ref{defn:AE-relative:remainder}, we claim that 
\begin{equation}\label{prop:AE-relative-fiberwise-fortified::proof:Foelner-fortified:shrunken-size} 
    |\CC\unitSpace \cap \CT\unitSpace| \geq (1-\delta) |\CT\unitSpace| \quad  \text{ for any tower } \CT \text{ in } \CD \; ,
\end{equation}
which we shall verify at the end of the proof. Now, observe that this claim implies,  together with \ref{prop:AE-relative-fiberwise-fortified::proof:Foelner-fortified:precompact}, that for any tower $\CT$ in $\CD$,  
\[
    \max \{ |(\CD\unitSpace \setminus \CC\unitSpace) \cap \CT\unitSpace| , |\CX\cap \CT\unitSpace| \} \leq \delta \, |\CT\unitSpace| \leq \frac{\delta}{1-\delta} \, |\CC\unitSpace \cap \CT\unitSpace| \leq \frac{\varepsilon}{2} \, |\CC\unitSpace \cap \CT\unitSpace| \leq \frac{1}{2} \, |\CC\unitSpace \cap \CT\unitSpace| \; ,
\]
whence we can choose disjoint subsets $\CY, \CZ \subseteq \CC\unitSpace$ such that, for any tower $\CT$ in $\CD$,  
\[
    |\CY \cap \CT\unitSpace| = |\CX \cap \CT\unitSpace| \quad \text{ and } \quad   |\CZ \cap \CT\unitSpace| = |(\CD\unitSpace \setminus \CC\unitSpace) \cap \CT\unitSpace| \; ,
\]
which immediately implies 
\[
    |(\CY \sqcup \CZ) \cap \CT\unitSpace| \leq \varepsilon  \, |\CC\unitSpace \cap \CT\unitSpace| \; .
\]
On the other hand, it follows from \eqref{prop:AE-relative-fiberwise-fortified::proof:Foelner-fortified:shrunken-size}, \ref{prop:AE-relative-fiberwise-fortified::proof:Foelner-fortified:precompact} , and \Cref{rem:subequivalence-basic}\eqref{rem:subequivalence-basic:castle} and \eqref{rem:subequivalence-basic:additive-prec-precsim} that 
\[
    \CG\unitSpace \setminus \bigcup \CC\unitSpace = (\CG\unitSpace \setminus \bigcup \CD\unitSpace) \sqcup \bigcup (\CD\unitSpace \setminus \CC\unitSpace) \prec_\CG \bigcup (\CY \sqcup \CZ) \; . 
\]
Since by \Cref{lem:castle-shrinking}\eqref{lem:castle-shrinking:tower}, any $\CC$-tower is given by intersecting a $\CD$-tower with $\CC$, the last two formulas combine to verify \Cref{defn:AE-relative}\ref{defn:AE-relative:remainder}. Hence we have showed $\CG$ is almost elementary with relative remainder control, modulo the claim in  \eqref{prop:AE-relative-fiberwise-fortified::proof:Foelner-fortified:shrunken-size}.  

It remains to prove \eqref{prop:AE-relative-fiberwise-fortified::proof:Foelner-fortified:shrunken-size}. In view of \Cref{lem:castle-basic}\eqref{lem:castle-basic:height}, it suffices to show that for any $u \in \bigcup \CD\unitSpace$, we have $|(\bigcup \CC\unitSpace) (\bigcup \CD) u| \geq (1-\delta) |(\bigcup \CD) u|$, which, by  \ref{prop:AE-relative-fiberwise-fortified::proof:Foelner-fortified:remainder}, can be further reduced to showing 
\[
    (\bigcup \CC\unitSpace) (\bigcup \CD) u \supseteq  \big((\bigcup \CD) u \big) \setminus \partial^-_L \big((\bigcup \CD) u \big)\; .
\]
To prove this, it suffices to show that for any $D \in \CD$ and $u \in s(D)$, if $L D u \subseteq (\bigcup \CD) u$, then $r(D) \in \CC\unitSpace$. Further observe that by replacing $D$ with $r(D)$, it suffices to consider in the above statement only $D \in \CD\unitSpace$. So, we fix arbitrary $D \in \CD\unitSpace$ and $u \in D$ satisfying $L D u \subseteq (\bigcup \CD) u$ and aim to show $D \in \CC$, that is, $KD \subseteq \bigcup \CD$ (which is also equivalent to $DK \subseteq \bigcup \CD$ since $K$ is symmetric). To this end, we fix an arbitrary $x \in K$ with $v := s(x) \in D$. 
By \ref{prop:AE-relative-fiberwise-fortified::proof:Foelner-fortified:refine} and our choice of the collection $\CW$, there is an open bisection $B \in \CB$ such that $x \in Y_B \subseteq B$ and $D \subseteq s(Y_B)$. Since $Y_B$ is a bisection, we may write $Y_B u = \{y\}$. Observe that since $Y_B \subseteq L$, we have $y \in L u \subseteq (\bigcup \CD) u$ and thus $y \in L \cap (\bigcup \CD)$, whence by \ref{prop:AE-relative-fiberwise-fortified::proof:Foelner-fortified:refine}, there is $D' \in \CD$ and $B' \in \CB$ such that $y \in D' \subseteq B'$. Note that $s(D') = D$, for $s(y) = u \in s(D') \cap D$. Now since $y \in B \cap B'$, it follows from our choice of $\CB$ that there is an open bisection $U \in \CU$ such that $B \cup B' \subseteq U$, whence $D' \subseteq U$ and $\{x\} = Bv = Uv = D'v \subseteq \bigcup \CD$. Since $x$ was chosen arbitrarily, this shows $KD \subseteq \bigcup \CD$. 

This completes the proof. 
\end{proof}

We are now ready to restrict to the minimal case and show that Definition~\ref{defn:AE-relative} is equivalent to Definition~\ref{defn:AE} in this setting. 

\begin{prop}
	\label{prop:AE-relative-DSC-fiberwise} 
	Let $\CG$ be a minimal $\sigma$-compact \'{e}tale groupoid with a compact unit space. 
	Suppose $\CG$ is fiberwise amenable. 
	If $\CG$ is almost elementary with relative remainder control, then $\CG$ has groupoid strict comparison. 
\end{prop}

\begin{proof}
	We follow the proof of Proposition~\ref{prop:AE-DSC-fiberwise}. Suppose $\GU$ is finite. Corollary \ref{cor:AE-relative-finite-pair} shows that $\CG$ has to be a union of pair groupoids.  Since $\CG$ is also minimal,  the groupoid $\CG$ has to be a pair groupoid and thus $\CG$ has groupoid strict comparison.

    Suppose $\GU$ is infinite. Let $U, V$ be non-empty open set in $\GU$ satisfying $\mu(U)<\mu(V)$ for any $\mu\in M(\CG)$. It suffices to show $L\prec_\CG V$ for an arbitrarily fixed compact subset $L\subseteq U$. 

    To this end, using a similar argument to the one in Proposition~\ref{prop:AE-DSC-fiberwise}, we find an open set $U'$ with $L\subseteq U'\subseteq \overline{U'}\subseteq U$, an open set $V'$ with $V'\subseteq \overline{V'}\subseteq V$ and an $\eta'>0$ so that one can find a compact set $K\subseteq \CG$ and $\varepsilon>0$ such that for any finite $F$ in $\CG$ with $|KF|<(1+\varepsilon)|F|$, one always has
    \begin{equation}\label{equi: folner}
       \frac{1}{|F|}\sum_{x\in F}1_{\overline{U'}}(r(x))+\frac{\eta'}{2}<\frac{1}{|F|}\sum_{x\in F}1_{V'}(r(x)). 
    \end{equation}
   Choose $n>\max\{\frac{3}{\eta'}, \frac{1}{\varepsilon}\}$ and define an open cover $\CV := \left\{ U' , \CG\unitSpace \setminus L \right\} \vee \left\{ V , \CG\unitSpace \setminus \overline{V'} \right\}$. 
   Since a minimal fiberwise amenable groupoid is also ubiquitously fiberwise amenable by \Cref{thm:minimal-fiberwise-amenability}, we see that $\CG$ satisfies   \Cref{prop:AE-relative-fiberwise-fortified} \eqref{prop:AE-relative-fiberwise-fortified:Foelner-fortified}, i.e., there exists an open castle $\CC$ such that 
    \begin{enumerate}[label=(\roman*)]
			\item\label{refine} $\CF_+ \unitSpace$ refines $\CV$;  
			\item $\bigcup \CF_+$ is precompact in $\CG$;  
			\item\label{remainder} there is a subcollection $\CX \subseteq \CF_-^{(0)}$ such that $|\CX \cap \CT^{(0)}| \leq \frac{1}{n} |\CT^{(0)}|$ for every tower $\CT$ in $\CF_-$ and
			$\GU \setminus \left( \bigcup \CF_-^{(0)} \right) \prec_{\CG}  \bigcup \CX$; 
			\item\label{folner} for any $u \in \CG\unitSpace$, we have $|K (\bigcup \CF_+) u| \leq (1+\frac{1}{n}) |(\bigcup \CF_+) u|\leq (1+\varepsilon)|\bigcup \CF_+u|$.  
			   \end{enumerate}
 Let $\CT$ be a tower in $\CF_-$, i.e, $\CT\in \CF_-/\sim_{\CF_-}$. Define $S_{\CT, 1}=\{F_-\in \CT: F_-\cap L\neq \emptyset\}$ and $S_{\CT, 2}=\{F_-\in \CT: F_-\cap \overline{V'}\neq \emptyset\}$. Now, for any level $F_+\in \CF_+\unitSpace$, it follows from the condition \ref{refine} that $F_+\subseteq U'$ if its $F_-\in S_{\CT,1}$ and $F_+\subseteq V$  whenever its $F_-\in S_{\CT, 2}$. Choose $u_\CT\in\bigcup \CT\unitSpace$. Then inequality \eqref{equi: folner} and  the condition \ref{folner} for $\bigcup \CT u_\CT$ imply that
 \[\frac{|S_{\CT, 1}|}{|\CT\unitSpace|}+\frac{\eta'}{2}<\frac{|S_{\CT, 2}|}{|\CT\unitSpace|}.\]
 Furthermore, by the choice of $n$, one has 
 \[\frac{|S_{\CT, 1}|}{|\CT\unitSpace|}+\frac{|\CX\cap \CT\unitSpace|}{|\CT\unitSpace|}<\frac{|S_{\CT, 2}|}{|\CT\unitSpace|}.\]
This allows us to choose two injective maps $\varphi_\CT: S_{\CT, 1}\to S_{\CT, 2}$ and $\psi_\CT: \CX\cap \CT\unitSpace\to S_{\CT, 2}$ with disjoint images, i.e., $\varphi_{\CT}(S_{\CT, 1})\cap \psi_\CT(\CX\cap \CT\unitSpace)=\emptyset$.
For the map $\varphi_\CT$, Lemma \ref{lem:castle-comparison} implies that 
\[\bigcup \{\CS_{\CT, 1}: \CT\in \CF_-/\sim_{\CF_-}\}\precsim_{\bigcup \CF_-} \bigcup \{\varphi_{\CT}(\CS_{\CT, 1}): \CT\in \CF_-/\sim_{\CF_-}\},\]
and by \ref{lem:castle-comparison}\eqref{lem:castle-comparison:bisection}  there exists a bisection $\CB_-\subseteq \CF_-$ such that 
\[\bigcup \{\CS_{\CT, 1}: \CT\in \CF_-/\sim_{\CF_-}\}=s(\bigcup \CB_-)\] and 
\[\bigcup \{\varphi_{\CT}(\CS_{\CT, 1}): \CT\in \CF_-/\sim_{\CF_-}\}\subseteq r(\bigcup \CB_-).\]
Note that $\bigcup\CB_-$ is an open bisection in $\CG$ and the induced $B=\bigcup\{F_+\in \CF_+: \CF_-\in \CB_-\}$ is also an open bisection in $\CG$.
Define a compact set $A=\bigcup\{\overline{F_-}: F_-\in \CS_{\CT, 1}, \CT\in \CF_-/\sim_{\CF_-}\}$ and an open set $O_1=\bigcup\{F_+: F_-\in \varphi_\CT(\CS_{\CT, 1}), \CT\in \CF_-/\sim_{\CF_-}\}\subseteq V$. Using $B$,  one actually has $A\prec_\CG O_1$ and thus $L\cap A\prec_\CG O_1$.
In addition, for the $\psi_\CT$, Lemma \ref{lem:castle-comparison} entails that
\[\bigcup\CX= \bigcup \{\CX\cap \CT\unitSpace: \CT \text{ is a tower of } \CF_-\}\precsim_{\bigcup \CF_-} \bigcup \{\psi_{\CT}(\CS_{\CT, 1}): \CT \text{ is a tower of } \CF_-\}.\]    
We denote by $O_2=\bigcup \{\psi_{\CT}(\CS_{\CT, 1}): \CT \text{ is a tower of } \CF_-\}$, which is an open set in $V$ but disjoint from $O_1$. Then because of \ref{remainder}, Remark \ref{rem:subequivalence-basic}\eqref{rem:subequivalence-basic:auxiliary} implies that $\GU\setminus \bigcup (\CF_-\unitSpace)\prec_\CG O_2
$ and thus $L\cap (\GU\setminus \bigcup (\CF_-\unitSpace))\prec_\CG O_2$. Then apply Remark \ref{rem:subequivalence-basic}\eqref{rem:subequivalence-basic:additive-prec} to $L\cap A$ and $L\cap (\GU\setminus \bigcup (\CF_-\unitSpace))$, one has $L\prec_\CG O_1\sqcup O_2\subseteq V$ as desired.
\end{proof}

\begin{prop}
	\label{prop:AE-relative-DSC-non-fiberwise} 
	Let $\CG$ be a minimal $\sigma$-compact \'{e}tale groupoid with a compact unit space. 
	Suppose $\CG$ is non-fiberwise-amenable. 
	If	$\CG$ is almost elementary with relative remainder control, then $M(\CG)=\emptyset$ and $\CG$ has groupoid strict comparison. 
\end{prop}
\begin{proof}
 It still suffices to study the case that $\GU$ is infinite as otherwise $\CG$ is pair groupoid and thus fiberwise amenable. This is a contradiction to the assumption. First, the minimality of $\CG$ implies that $\GU$ is perfect by Remark \ref{rem: minimal-perfect}. Let $O$ be non-empty open set in $\CG$. We will show $\GU\prec_\CG O$. To this end, We mainly follow the argument for Proposition \ref{prop:AE-DSC-non-fiberwise}. By the minimality, one chooses an open bisection $U$ such that $s(U), r(U)$ are disjoint open sets in $O$. Choose another open set $W\subseteq s(U)$ such that $\overline{W}\subseteq s(U)$ holds as well. Then like in Proposition \ref{prop:AE-DSC-non-fiberwise}, the minimality of $\CG$ also implies that there are compact bisections $K_1, \dots, K_m$ in $\CG$ such that $\bigcup_{i=1}^ms(K_i)^o=\GU$ and $\bigcup_{i=1}^mr(K_i)\subseteq W$. Define $K=\bigcup_{i=1}^mK_i\cup \GU$, which is a compact set in $\CG$. Then Theorem \ref{thm:fiberwise-amenable-dichotomy}\eqref{thm:fiberwise-amenable-dichotomy:non-amenable} implies that there exists a compact set $L$ containing $\GU$ such that for any finite set $M\subseteq \CG$, the set $LM$ contains at least $2|M|$ many disjoint sets of the form $Kx$ for $x\in \CG$. Define an open cover $\CU=\{V, \GU\setminus \overline{W}\}\vee \{s(K_1)^o,\dots, s(K_n)^o\}$ of $\GU$ and an open cover $\CV=\CU\cup\{\CG\setminus \GU\}$ of $\CG$. Then apply Lemma~\ref{lem:AE-relative-fortified} to compact sets $L, \GU$ in $\CG$ and the cover $\CV$ to obtain a fortified castle $\CF$
satisfying 
\begin{enumerate}[label=(\roman*)]
		\item $\CF_+\unitSpace$ refines $\CV$; 
		\item
		$\bigcup \CF_+$ is precompact in $\CG$;  
		\item  
		there is a subcollection $\CW \subseteq \left( \CF^{-L} \right)_-^{(0)}$ such that $|\CW \cap \CT^{(0)}| \leq \varepsilon |\CT^{(0)}|$ for every tower $\CT$ in $\left( \CF^{-L} \right)_-$ and
		$\CG^{(0)}\setminus\ \left( \bigcup \left( \CF^{-L} \right)_-^{(0)} \right) \prec_{\CG}  \bigcup \CW$.   
	\end{enumerate}

   For any fortified tower $\CT$ in $\CF$, choose a unit $u_\CT\in (\CT^{-L})^{(0)}_-$ and define a finite set $M_\CT=(\bigcup \CT^{-L}_-)u_\CT$. Then the choice of $L$ implies that the set $LM_\CT$ contains at least $|M_\CT|$ many disjoint non-empty finite subsets of the form $Kx$ for $x\in \CG$ with $s(x)=u_\CT$. Denote by $m_\CT=|M_\CT|$ for simplicity. Then enumerate these $Kx$ explicitly by 
   $\{Kx_1,\dots, Kx_{m_\CT}\}$.
   Then the same argument in Proposition \ref{prop:AE-DSC-non-fiberwise} shows that $\bigcup (\CF^{-L}_+)^{(0)}\prec_\CG s(U)$. 
In addition, using the bisection $U$, one has 
\[\GU\setminus \bigcup (\CF^{-L}_-)^{(0)}\prec_\CG \bigcup \CW\subseteq \bigcup (\CF^{-L}_-)^{(0)}\prec_\CG s(U)\prec_\CG r(U).\]
Then, Remark \ref{rem:subequivalence-basic}\eqref{rem:subequivalence-basic:additive-prec} shows that \[\GU\prec_\CG s(U)\sqcup r(U)\subseteq O.\]  Finally, the same argument in Proposition \ref{prop:AE-DSC-non-fiberwise} shows that $M(\CG)=\emptyset$.
\end{proof}

\begin{thm}
	\label{thm:AE-relative-AE} 
	Let $\CG$ be a minimal $\sigma$-compact \'{e}tale groupoid
	with a compact unit space. 
	Then $\CG$ is almost elementary if and only if it is almost elementary with relative remainder control. 
\end{thm}
\begin{proof}
Let $\CG$ be minimal. In view of Theorem~\ref{thm:AE-DSC}, it suffices to show $\CG$ is almost elementary with relative remainder control if and only if it is almost elementary in measure and has groupoid strict comparison. 
    
Suppose $\CG$ is almost elementary with relative remainder control. Then Propositions~\ref{prop:AE-relative-in-measure}, \ref{prop:AE-relative-DSC-fiberwise} and~\ref{prop:AE-relative-DSC-non-fiberwise} implies that $\CG$ is indeed almost elementary in measure and has groupoid strict comparison. Thus, $\CG$ is almost elementary.

For the converse, suppose $\CG$ is almost elementary. Let compact $K\subseteq \CG$, a number $1>\varepsilon>0$, and an open cover $\CV$ of $\GU$ be given. In the case that $\CG$ is fiberwise amenable, first choose an integer $n> 1/\varepsilon>1$. By compactness of $K$ and $\GU$, one may choose an open cover $\CW$ of $K\cup \GU$ consisting of precompact open bisections such that $|\{Wu: W\in \CW\}|\geq 2n(n+1)$ for any $u\in \GU$. Define a compact set $K'=\bigcup\{\overline{W}: W\in \CW\}$. Choose an open set $O$ in $\GU$ with $\max_{\mu\in M(\CG)}\mu(O)<\frac{1}{n+1}$ by Lemma \ref{lem:O-meas-sup}.

Now, for the compact set $K'$ and the open cover $\CV$, Lemma \ref{lem:AE-fiberwise-fortified} implies that there exists an open castle $\CC$ (i.e., the $\CF^+$ in  Lemma \ref{lem:AE-fiberwise-fortified}) such that 
\begin{enumerate}[label=(\roman*)]
    \item $\CC$ is $K'$-extendable and thus also $K$-extendable;
    \item $\CC\unitSpace$ refines $\CV$;
    \item $\mu(\GU\setminus(\bigcup \CC\unitSpace))<\mu(O)<\frac{1}{n+1}$ for any $\mu\in M(\CG)$;
    \item\label{long tower} $|K'\bigcup\CC u|\leq (1+\frac{1}{n})|\bigcup \CC u|$ for any $u\in \GU$.
\end{enumerate}
 Let $u$ be an arbitrary unit in $\bigcup\CC\unitSpace$. Then the condition \ref{long tower} above entails that
 \[|\bigcup\CC u|\leq |K'\bigcup \CC u|/(1+\frac{1}{n})\geq |K'u|/(1+\frac{1}{n})> n(n+1),\]
 which further implies that 
 \[\frac{|\bigcup \CC u|}{n}-\frac{|\bigcup \CC u|}{n+1}>1\]
 and thus there exists an integer in the interval $[\frac{|\bigcup \CC u|}{n+1}, \frac{|\bigcup \CC u|}{n}]$. This allows us to choose a $\CX\subseteq \CC\unitSpace$ such that 
 \[\frac{1}{n+1}|\CT\unitSpace|\leq |\CX\cap \CT\unitSpace|\leq \frac{1}{n}|\CT\unitSpace|<\varepsilon|\CT\unitSpace|\]
 for any tower $\CT$ in $\CC$. This further implies
\[\mu(\bigcup \CX)\geq \frac{1}{n+1}\]
and thus
\[\mu(\GU\setminus \bigcup \CC\unitSpace)<\frac{1}{n+1}\leq\mu(\bigcup \CX)\]
for any $\mu\in M(\CG)$.
 Recall that $\CG$ has groupoid strict comparison by Theorem \ref{thm:AE-DSC} and thus $\GU\setminus \bigcup \CC\unitSpace\prec_\CG \CX$. Thus $\CG$ is almost elementary with relative remainder control.

 On the other hand, suppose $\CG$ is not fiberwise amenable. Let the compact set $K$, the number $\varepsilon>0$, and the open cover $\CV$ be given. Still choose an integer $n>1/\varepsilon$. Now, since $\CG$ is almost elementary, there exists an open castle $\CC$ such that 
 \begin{enumerate}[label=(\roman*)]
    \item $\CC$ is $K$-extendable;
    \item $\CC\unitSpace$ refines $\CV$.
\end{enumerate} 
This means there exists an open castle $\CD$ such that $\CC\subseteq \CD^{-K}$.
Fix one level $C\in \CC\unitSpace$. Then Lemma \ref{lem:O-meas-sup} implies that one can construct a tower $\CT$ such that $\bigcup\CT\unitSpace\subseteq C$ with $|\CT\unitSpace|>n$. Define a new tower
\[\CS=\{C_1UC_2: U\in \CT, C_1, C_2\in \CD\}\]
and observe that $\CS^{-K}$ is not empty since it contains the tower 
\[\CR=\{C_1UC_2: U\in \CT, C_1, C_2\in \CC\}.\]
Therefore $\CR$ is $K$-extendable and note also $\CR\unitSpace$ refines $\CC\unitSpace$ and thus $\CR\unitSpace$ refines $\CV$. Finally, choose a non-empty subfamily $\CX\subseteq \CT\unitSpace$ with $|\CX|\leq\frac{1}{n}|\CT\unitSpace|<\varepsilon|\CR\unitSpace|$. Recall that $\CG$ is not fiberwise amenable and almost elementary. Then Proposition \ref{prop:AE-DSC-non-fiberwise} implies that $\GU\setminus(\bigcup\CR\unitSpace)\prec_\CG \bigcup\CX$. Combining these two cases, one has $\CG$ is almost elementary with relative remainder control.
\end{proof}

In view of Theorem~\ref{thm:AE-relative-AE}, we will henceforth drop the phrase ``with relative remainder control'' and treat Definition~\ref{defn:AE-relative} as a generalization of Definition~\ref{defn:AE} beyond the minimal setting. 

\section{Almost finite ample groupoids in the sense of Matui} \label{sec:Matui}

In this section, we investigate the relationship between our almost elementariness and Matui's original notion of almost finiteness for
ample groupoids (see \cite{Matui2012Homology}).

\begin{defn}[{\cite[Definition 6.2]{Matui2012Homology}}]
	\label{6.1} \label{defn:Matui}
	An ample \'{e}tale
	groupoid $\CG$ with a compact unit space is called \textit{almost
		finite} if, for any compact set $K$ in $\CG$ and $\varepsilon>0$, there
	is a compact open elementary subgroupoid $\CH$ of $\CG$ with $\HU=\GU$
	such that for any $u\in\GU$, one has 
	\[
	|K\CH u|<(1+\varepsilon)|\CH u|
	\; .
	\]
\end{defn}

\begin{rem}
	The original formulation in \cite[Definition 6.2]{Matui2012Homology} asks for $|K\CH u\setminus\CH u|<\varepsilon|\CH u|$ instead of $|K\CH u|<(1+\varepsilon)|\CH u|$, but these formulations are equivalent since, by enlarging $K$ if necessary, we may assume without generality that $K \supseteq \CG\unitSpace$. 
\end{rem}

Definition~\ref{defn:Matui} already bears some similarity with Proposition~\ref{prop:AE-relative-fiberwise-fortified}\eqref{prop:AE-relative-fiberwise-fortified:Foelner-0}, which is equivalent to Proposition~\ref{prop:AE-relative-fiberwise-fortified}\eqref{prop:AE-relative-fiberwise-fortified:conjunction}, i.e., the conjunction of almost elementariness and ubiquitous fiberwise amenability. There are, however, three main differences: 
\begin{enumerate}
	\item In Proposition~\ref{prop:AE-relative-fiberwise-fortified}\eqref{prop:AE-relative-fiberwise-fortified:Foelner-0}\ref{prop:AE-relative-fiberwise-fortified:Foelner-0:refine}, we use open castles, while Matui uses, essentially, compact open (or equivalently, clopen) castles (elementary subgroupoids are closely related to castles; see Lemma~\ref{lem:castle-basic}\eqref{lem:castle-basic:subquotient}). 
	
	\item We require in Proposition~\ref{prop:AE-relative-fiberwise-fortified}\eqref{prop:AE-relative-fiberwise-fortified:Foelner-0}\ref{prop:AE-relative-fiberwise-fortified:Foelner-0:refine} that the castles are ``fine'', while Matui makes no such requirements. 
	
	\item We allow in Proposition~\ref{prop:AE-relative-fiberwise-fortified}\eqref{prop:AE-relative-fiberwise-fortified:Foelner-0}\ref{prop:AE-relative-fiberwise-fortified:Foelner-0:remainder} a ``small'' remainder, while Matui allows no remainder.  
\end{enumerate}
It turns out that these a priori differences may be neutralized by the assumption that $\CG$ is ample. 
Let us start by addressing the first difference. 

\begin{prop} \label{prop:AE-relative-fiberwise-ample} 
	Let $\CG$ be an \emph{ample} $\sigma$-compact \'{e}tale groupoid with a compact unit space. 
	Then the following are equivalent: 
	\begin{enumerate}
		\item\label{prop:AE-relative-fiberwise-ample:conjunction}  The groupoid $\CG$ is ubiquitously fiberwise amenable and almost elementary (with relative remainder control). 
		
		\item\label{prop:AE-relative-fiberwise-ample:Foelner} For any non-empty compact open set $K, L$ in $\CG$, 
		any $\varepsilon > 0$,  
		and any finite compact open cover $\CV$ of $L$, 
		there is a nonempty compact open castle $\CC$
		satisfying 
		\begin{enumerate}[label=(\roman*)]
			\item\label{prop:AE-relative-fiberwise-ample:Foelner:refine} If $C\cap L\neq \emptyset$ then $C\subseteq V$ for some $V\in \CV$;  
			\item\label{prop:AE-relative-fiberwise-ample:Foelner:remainder} there is a subcollection $\CW \subseteq \CC^{(0)}$ such that $|\CW \cap \CT^{(0)}| \leq \varepsilon |\CT^{(0)}|$ for every tower $\CT$ in $\CC$ and
			$\GU \setminus \left( \bigcup \CC^{(0)} \right) \prec_{\CG}  \bigcup \CW$; 
			\item\label{prop:AE-relative-fiberwise-ample:Foelner:Foelner} 
			for any $u \in \CG\unitSpace$, we have $|K (\bigcup \CC) u| \leq (1+\varepsilon) |(\bigcup \CC) u|$. 
		\end{enumerate}
	\end{enumerate}
\end{prop}

\begin{proof}
For \eqref{prop:AE-relative-fiberwise-ample:Foelner}$\Rightarrow$\eqref{prop:AE-relative-fiberwise-ample:conjunction}, it suffices to verify Proposition~\ref{prop:AE-relative-fiberwise-fortified}\eqref{prop:AE-relative-fiberwise-fortified:Foelner-0}. Let $K$ be a compact set in $\CG$, $\varepsilon>0$ and an open cover $\CU$ of $K$. Since $\CG$ is ample, one may enlarge $K$ to a compact open $K'$ and define an open cover $\CU'=\CU\cup \{K'\setminus K\}$ of $K'$. Now choose a compact open cover $\CV$ of $K$ refining $\CU'$. For the $K'$, $\varepsilon$, and $\CV$,  the condition \eqref{prop:AE-relative-fiberwise-ample:Foelner} yields a compact open castle $\CC$ satisfying the conditions \ref{prop:AE-relative-fiberwise-ample:Foelner:refine}, \ref{prop:AE-relative-fiberwise-ample:Foelner:remainder}, and 
\ref{prop:AE-relative-fiberwise-ample:Foelner:Foelner}, which implies Proposition~\ref{prop:AE-relative-fiberwise-fortified}\eqref{prop:AE-relative-fiberwise-fortified:Foelner-0}.

For the converse direction, Recall that \eqref{prop:AE-relative-fiberwise-ample:conjunction} is equivalent to Proposition~\ref{prop:AE-relative-fiberwise-fortified}\eqref{prop:AE-relative-fiberwise-fortified:Foelner-fortified}, which will be shown to imply the statement \eqref{prop:AE-relative-fiberwise-ample:Foelner}.
Indeed, Proposition~\ref{prop:AE-relative-fiberwise-fortified}\eqref{prop:AE-relative-fiberwise-fortified:Foelner-fortified} entails that for the compact open $K, L$, the number $\varepsilon>0$ and the compact open cover $\CV$ of $L$, there is a nonempty fortified castle $\CF$
		satisfying 
		\begin{enumerate}[label=(\roman*)]
			\item\label{prop:AE-relative-fiberwise-fortified:Foelner-fortified:refine} for any $F_+ \in \CF_+$ with $F_+ \cap L \not= \varnothing$, there is $V \in \CV$ such that $F_+ \subseteq V$ 
			\item\label{prop:AE-relative-fiberwise-fortified:Foelner-fortified:precompact} $\bigcup \CF_+$ is precompact in $\CG$;  
			\item\label{prop:AE-relative-fiberwise-fortified:Foelner-fortified:remainder} there is a subcollection $\CX \subseteq \CF_-^{(0)}$ such that $|\CX \cap \CT^{(0)}| \leq \varepsilon |\CT^{(0)}|$ for every tower $\CT$ in $\CF_-$ and
			$\GU \setminus \left( \bigcup \CF_-^{(0)} \right) \prec_{\CG}  \bigcup \CX$; 
			\item for any $u \in \CG\unitSpace$, we have $|K (\bigcup \CF_+) u| \leq (1+\varepsilon) |(\bigcup \CF_+) u|$.  
		\end{enumerate}
Now, for any $\overline{F_-}\subseteq F_+$ with $F_-\in \CF_-\unitSpace$ and $F_+\in \CF_+\unitSpace$, since $\CG$ is ample, there exists a compact open set $C$ such that $\overline{F_-}\subseteq C\subseteq F_+$. Denote by $f$ the sum of all indicator functions $1_C$ of all such $C$. Then $f\in C_0(\bigcup \CF_+\unitSpace, \IR)\subseteq C_0(\bigcup \CF, \IR)$ and thus we obtain an invariant function $g=f_{\uparrow\CF_+}\in C_0(\bigcup \CF_+, \IR)^{\bigcup \CF_+}$ by Lemma \ref{lem:castle-invarint-function-smallest} and the support of $g$ is compact open and satisfies $\bigcup \CF_-\subseteq \operatorname{supp}(g)\subseteq \bigcup \CF_+$.
Then Lemma~\ref{lem:castle-trim} allows to define a clopen castle 
\[\CC=\{g^{-1}(\{1\})\} \vee \CF_+ := \left\{ g^{-1}(\{1\})  \cap F_+ \colon F_+ \in \CF_+ \right\} \setminus \{ \varnothing \}.\]
First, by the construction, $\CC\unitSpace$ refines $\CF_+\unitSpace$ and thus  $\CC\unitSpace$ refines $\CV$.
Then, it is direct to see
\[|K(\bigcup \CC)u|=|K (\bigcup \CF_+) u| \leq (1+\varepsilon) |(\bigcup \CF_+) u|=(1+\varepsilon) |(\bigcup \CC) u|\]
for any $u\in \bigcup\CC\unitSpace$. On the other hand, if $u\notin \bigcup\CC\unitSpace$, the inequality in \ref{prop:AE-relative-fiberwise-ample:Foelner:Foelner} holds trivially. Moreover, define $\CW=\bigcup\{C\in \CC\unitSpace: F_-\subseteq C, F_-\in \CX\}$. Then by the construction, one has $|\CW\cap \CT\unitSpace|<\varepsilon|\CT\unitSpace|$ for any tower $\CT$ in $\CC$ and thus 
\[\GU\setminus \bigcup \CC\unitSpace\subseteq \GU\setminus \bigcup \CF_-\unitSpace\prec_\CG \bigcup \CX\subseteq \bigcup \CW.\]
This implies  $\GU\setminus \bigcup \CC\unitSpace\prec_{\CG}\bigcup \CW$.
\end{proof}

The second and the third differences are resolved by partitioning castles with the help of Lemma~\ref{lem:castle-split}. We start with a simple fact about totally disconnected spaces, whose proof we leave to the reader. 

\begin{lem} \label{lem:ample-tighten}
	Let $X$ be a totally disconnected space. Let $O$ be a compact open set in $X$ and let $\CV$ be a collection of open sets in $X$ covering $O$. Then there is a finite clopen partition $O = \bigsqcup_{i \in I} O_i$ such that $\left\{ O_i \colon i \in I \right\}$ refines $\CV$. 
	\qed
\end{lem}

\begin{rem}\label{rem: subequivalence-ample}
	We make an observation about the subequivalence relation introduced in Definition~\ref{defn:subequivalence}. In an ample \'{e}tale groupoid $\CG$, if $U, V$ are compact open sets in $\GU$ and satisfying $U\precsim_{\CG}V$, 
	applying Lemma~\ref{lem:ample-tighten}, one actually can find a collection $\{A_1,\dots, A_n\}$ of compact open bisections such that $U=\bigsqcup_{i=1}^ns(A_i)$ and $\bigsqcup_{i=1}^nr(A_i)\subseteq V$.
\end{rem}

\begin{thm} \label{thm:Matui-AE} 
	Let $\CG$ be an {ample} $\sigma$-compact \'{e}tale groupoid with a compact unit space. 
	Then $\CG$ is almost finite in the sense of Matui if and only if it is almost elementary (with relative remainder control) and ubiquitously fiberwise amenable. 
\end{thm}

\begin{proof}
	For the ``only if'' direction, we verify Proposition~\ref{prop:AE-relative-fiberwise-ample}\eqref{prop:AE-relative-fiberwise-ample:Foelner}. Let compact open sets $K, L$, $\varepsilon>0$,  and a finite compact open cover $\CV$ of $L$ be given. Almost finiteness of $\CG$ implies that there exists a compact open elementary subgroupoid $\CH$ of $\CG$ with $\HU=\GU$ and $|K\CH u|\leq (1+\varepsilon)|\CH u|$ holds for any $u\in \GU$. Recall that the elementary groupoid $\CH$ admits a decomposition as compact open castles by Remark \ref{rem: elementary groupoid}, i.e., there exists a compact open castle $\CC$ such that $\CH=\bigcup \CC$. To establish \ref{prop:AE-relative-fiberwise-ample:Foelner:refine} in Proposition~\ref{prop:AE-relative-fiberwise-ample}\eqref{prop:AE-relative-fiberwise-ample:Foelner}, Let $V\in \CV$, define a clopen partition of $\CH$ by $\CP_V=\{\CH\setminus V, V\}$  and denote by $\CP$ the common refinement for all of these $\CP_V$. Then Lemma~\ref{lem:castle-split} shows that there exists another compact open castle $\CD$ satisfying that $\CD$ refines $\CC$ and $\CP$, and  $\bigcup\CD=\bigcup \CC$. This implies that for any $V\in \CV$ and $D\in \CD$, either $D\subseteq V$ or $D\cap V=\emptyset$. Then if $D\cap L\neq \emptyset$, then because $\CV$ is an cover of $L$, there has to be a $V\in \CV$ such that $C\cap V\neq\emptyset$ and therefore, $C\subseteq V$. This shows \ref{prop:AE-relative-fiberwise-ample:Foelner:refine} in Proposition~\ref{prop:AE-relative-fiberwise-ample}\eqref{prop:AE-relative-fiberwise-ample:Foelner}. Then, the condition \ref{prop:AE-relative-fiberwise-ample:Foelner:remainder}  in Proposition~\ref{prop:AE-relative-fiberwise-ample}\eqref{prop:AE-relative-fiberwise-ample:Foelner} holds trivially because $\GU\setminus \bigcup\CD\unitSpace=\emptyset$.
    Finally, let $u\in \GU$. The fact $\bigcup \CD=\bigcup \CC$ implies $\bigcup \CD u=\bigcup \CC u=\CH u$ and therefore one has
    \[|K(\bigcup \CD)u|=|K \CH u| \leq (1+\varepsilon) |\CH u|=(1+\varepsilon) |(\bigcup \CD) u|,\]
    which establishes the condition
    \ref{prop:AE-relative-fiberwise-ample:Foelner:Foelner} in Proposition~\ref{prop:AE-relative-fiberwise-ample}\eqref{prop:AE-relative-fiberwise-ample:Foelner}.
    
  For the ``if'' direction, it suffices to verify Definition~\ref{6.1} under the assumption of Proposition~\ref{prop:AE-relative-fiberwise-ample}\eqref{prop:AE-relative-fiberwise-ample:Foelner}. Let $K$ be a compact set in $\CG$ and $\varepsilon>0$. By compactness of $K$, the number $N=\sup_{x\in \CG}|Kx|$ is finite.
  Choose $\delta<\varepsilon/2N$. Proposition~\ref{prop:AE-relative-fiberwise-ample}\eqref{prop:AE-relative-fiberwise-ample:Foelner} implies that there exists a non-empty clopen castle $\CC$ such that
\begin{enumerate}
			\item there is a subcollection $\CW \subseteq \CC^{(0)}$ such that $|\CW \cap \CT^{(0)}| \leq \delta |\CT^{(0)}|$ for every tower $\CT$ in $\CC$ and
			$\GU \setminus \left( \bigcup \CC^{(0)} \right) \prec_{\CG}  \bigcup \CW$; 
			\item for any $u \in \CG\unitSpace$, we have $|K (\bigcup \CC) u| \leq (1+\delta) |(\bigcup \CC) u|$. 
		\end{enumerate}  
    Denote by clopen sets $R=\GU \setminus \left( \bigcup \CC^{(0)} \right)$ and $W=\bigcup \CW$ for simplicity. Then, Remark \ref{rem: subequivalence-ample} implies that there exists a collection $\CA$ of compact open bisections such that $R=\bigsqcup\{s(A): A\in \CA\}$ and by $\bigsqcup \{r(A): A\in \CA\}\subseteq W$. Then define a clopen partition $\CP_A=\{r(A), (\bigcup \CC)\setminus r(A)\}$ for every $A\in \CA$ and denote by $\CP$ the common refinement fo all of these $\CP_A$ for $A\in \CA$. Then for $\CP$, Lemma~\ref{lem:castle-split} shows that there exists another clopen castle $\CD$ such that $\bigcup\CD=\bigcup \CC$, and $\CD$ refines both $\CC$ and $\CP$.

    Define $\CW^\CD=\{D\in \CD: D\subseteq C \text{ for some }C\in \CW\}$ and $\CY=\{D\in \CW^{\CD}: D\subseteq r(A) \text{ for some }A\in \CA\}$. Proposition \ref{prop:castle-split-more} shows that $\bigcup \CW=\bigcup \CW^{\CD}$ and for any $\CD$-tower $\CT'$, one has 
    \[|\CY\cap \CT'^{(0)}|\leq|\CW^{\CD}\cap \CT'^{(0)}|=|\CW\cap \CT\unitSpace|\leq \delta|\CT\unitSpace|=\delta|\CT'^{(0)}|\]
    where $\CT$ is the unique $\CC$-tower that the $\CD$-tower $\CT'$ refines. Moreover, by definition, one has $\bigsqcup \{r(A): A\in \CA\}=\bigcup \CY$.

    Now, let $\CT'$ be a $\CD$-tower. For any $D\in \CY\cap \CT'^{(0)}$, since the range $\{r(A): A\in \CA\}$ is a disjoint family, there exists a unique $A_D\in \CA$ such that $D\subseteq r(A_D)$.  Then, $\{r(A^{-1}_DD): D\in \CY\cap \CT'^{(0)}\}$ is disjoint family of clopen sets in $R$ because $\{s(A): A\in \CA\}$ is also a disjoint family. Moreover this family defines a clopen tower $\CS_{\CT'}$ such that $\CS^{(0)}_{\CT'}=\{r(A^{-1}_DD): D\in \CY\cap \CT'^{(0)}\}$ by adding all bisections of the form $A^{-1}_{D_2}DA_{D_1}$ in, where $D$ is the unique $\CD$-ladder such that $s(D)=D_1$ and $r(D)=D_2$ for all $D_1, D_2\in \CY\cap \CT'^{(0)}$. Our construction implies $|\CS_{\CT'}^{(0)}|=|\CY\cap \CT'^{(0)}|\leq \delta|\CT'^{(0)}|$. In addition, observe $\CS^{(0)}_{\CT'}$ is disjoint from $\CT'^{(0)}$ and it is straightforward to find a clopen tower $\CF_{\CT'}$ such that $\CF_{\CT'}\unitSpace=\CT'^{(0)}\sqcup \CS^{(0)}_{\CT'}$. 
    
   By construction, note that $\bigcup\{\bigcup\CS_{\CT'}:\CT'\text{ is a }\CD\text{-tower}\}=R=\GU\setminus (\bigcup \CC\unitSpace)$. This also implies that $\CE=\{\CF_{\CT'}: \CT'\text{ is a }\CD\text{-tower}\}$ is a disjoint family and thus defines a clopen castle with $\bigcup\CE\unitSpace=\GU$.
    Finally, let $u\in \GU$. By our construction, one has $u\in \CF^{(0)}_{\CT'}$ for a unique $\CD$-tower $\CT'$. Then one has
    \begin{align*}
   |K\bigcup\CE u|&=|K\bigcup(\CT'\sqcup \CS_{\CT'})u|\leq |K(\bigcup \CD)u|+|K\bigcup\CS_{\CT'} u|\leq |K(\bigcup\CC)u|+N|\CS_{\CT'}|\\
   &\leq (1+\delta)|(\bigcup \CD) u|+N\delta|\CT'^{(0)}|\leq (1+\varepsilon)|\bigcup\CE u|.
    \end{align*}
Thus, $\CG$ is almost finite. 
\end{proof}

\section{Almost elementary group actions and Kerr's almost finiteness} \label{sec:Kerr}

In this section, we focus on the case of transformation groupoids, i.e., those arising from group actions. We will establish a connection to Kerr's almost finiteness and also derive what almost elementariness means for a group action. 

Throughout the section, we write $\alpha:\Gamma\curvearrowright X$ for an action of a countable
discrete group $\Gamma$ on a compact Hausdorff space $X$. Following Example~\ref{exa:transformation-groupoid}, we write $X\rtimes_{\alpha}\Gamma$ for the resulting transformation groupoid. We also identify $X$ with the unit space of $X\rtimes_{\alpha}\Gamma$. 

Let us first recall the definitions of towers and castles for group actions. To distinguish them from the notions we introduced in Section~\ref{sec:castles}, we will add the adjective ``dynamical'' to them when necessary. 

\begin{defn}[{\cite[Definition 4.1]{Kerr2020Dimension}}]
	Let $S\subseteq\Gamma$
	be a finite subset and $V$ be a set in $X$. 
	We say $(S,V)$ is a
	\emph{(dynamical) tower} if $\{sV:s\in S\}$ is a disjoint family. A (dynamical) tower is
	called open if $V$ is open. 
	
	A finite
	family $\{(S_{i},V_{i}):i\in I\}$ of (dynamical) towers is called a \emph{(dynamical) castle} if
	$sV_{i}\cap tV_{j}=\emptyset$ for any $s\in S_{i}$ and $t\in S_{j}$
	and different $i,j\in I$. A (dynamical) castle is called open if all (dynamical) towers inside
	are open.
	
	We shall assume henceforth, without loss of generality, that any castle $\{(S_{i},V_{i}):i\in I\}$ satisfies $S_i \not= \varnothing$ and $V_i \not= \varnothing$ for any $i \in I$. 
\end{defn}

\begin{rem} \label{rem:dynamical-vs-groupoid-castles}
	Let us specify the relationship between dynamical castles and castles in \'{e}tale groupoids. 
	A key ingredient here is a partition (by clopen bisections), namely 
	\[
		\CP_{\alpha} := \left\{ B_\gamma \colon \gamma \in \Gamma \right\} 
		\; ,\quad \text{where }  
		B_\gamma := \left\{ \left( \gamma x, \gamma, x \right) \colon x \in X  \right\}
		\; .
	\]
	Note that $B_{\gamma'} B_{\gamma} = B_{\gamma' \gamma}$ and $B_{\gamma}^{-1} = B_{\gamma^{-1}}$ for any $\gamma', \gamma \in \Gamma$. Moreover, Let $F\subseteq \Gamma$ be a finite set. Then $|F|=|\bigcup_{\gamma\in F}B_\gamma u|$ for any $u\in X$.
	\begin{enumerate}
		\item \label{rem:dynamical-vs-groupoid-castles:assignment} 
		We define a map from the collection of all dynamical castles in the system $(X, \Gamma, \alpha)$ to the collection of all castles in $X \rtimes_\alpha \Gamma$ that refines $\CP_\alpha$: 
		Given any dynamical castle $\{(S_{i},V_{i}) \colon i\in I\}$, it is straightforward to check that the collection $\CC$ given by
		\[
			\left\{ B_{\gamma'} V_i B_{\gamma}^{-1} \colon \gamma', \gamma \in S_{i}, i\in I \right\} 
			,\text{ where }  
			B_{\gamma'} V_i B_{\gamma}^{-1} = \left\{ \left(\gamma' x, \gamma' \gamma^{-1} ,\gamma x \right) \colon x \in V_i \right\} \; ,
		\]
		is a castle in $X \rtimes_\alpha \Gamma$ that refines $\CP_\alpha$, with each dynamical tower $(S_{i},V_{i})$ giving rise to a tower in $\CC$. 
		
		\item \label{rem:dynamical-vs-groupoid-castles:injective} It is straightforward to verify that two dynamical castles $\{(S_{i},V_{i}) \colon i\in I\}$ and $\{(S'_{j},V'_{j}) \colon j\in J\}$ give rise to the same castle in $X \rtimes_\alpha \Gamma$ using the map in \eqref{rem:dynamical-vs-groupoid-castles:assignment} if and only if there are a bijection $\varphi \colon I \to J$ and $\gamma_i \in S_i$ with $S'_{\varphi({i})} = S_{i} \gamma_i^{-1}$ and $V'_{\varphi({i})} = \gamma_i \cdot V_i$ for any $i \in I$. 
		
		\item \label{rem:dynamical-vs-groupoid-castles:surjective} We claim that the map in \eqref{rem:dynamical-vs-groupoid-castles:assignment} is surjective.  
		Indeed
		, given any (open) castle $\CC$ in $X \rtimes_\alpha \Gamma$ that refines $\CP_\alpha$, 
		partitioning $\CC$ into its towers as $\bigsqcup_{i \in I} \CT_i$ and picking $V_i \in \CT_i\unitSpace$ for each $i \in I$, 
		we see that since by Lemma~\ref{lem:plurisection-basic}\eqref{lem:plurisection-basic:towers}, each $C \in \CT_i$ with $s(C)=V_i$ can be written as $B_\gamma V_i$ for some $
        \gamma\in \Gamma$ as $\CC$ refines $\CP_\alpha$. Define\[S_i=\{\gamma\in \Gamma: \text{there exists }  C\in \CT_i \text{ such that } C\subseteq B_\gamma,\text{ and }s(C)=V_i\},\] which is a finite set as $\CT_i$ is finite. Then it is straightforward to see $\{(S_i, V_i): i\in I\}$ is an open dynamical castle that the image is $\CC$ by the map described in \eqref{rem:dynamical-vs-groupoid-castles:assignment}. By construction, one has $S_iV_i=\bigcup \CT_i\unitSpace$.
	\end{enumerate}
\end{rem}

The following is a slight generalization of Kerr's original definition. 

\begin{defn}[{\cite[Definition 8.2]{Kerr2020Dimension}}]
	\label{6.5} \label{defn:Kerr}
	We say $\alpha:\Gamma\curvearrowright X$
	is \emph{almost finite} if for any finite set $F \subseteq\Gamma$, 
	any $\delta>0$, and any finite open cover $\CV$ of $X$,  
	there is  an open castle $\{(S_{i},V_{i}):i\in I\}$ such that 
	\begin{enumerate}[label=(\roman*)]
		\item\label{defn:Kerr:refine} the collection $\left\{ sV_{i} \colon s\in S_{i}, i \in I \right\}$ of levels refines $\CV$,  
		\item\label{defn:Kerr:Foelner} $|F S_{i}|\leq(1+\delta)|S_{i}|$ for any $i \in I$, and
		\item\label{defn:Kerr:remainder} there are sets $S_{i}'\subseteq S_{i}$ with $|S_{i}'|< \delta |S_{i}|$ for $i \in I$ such that 
		\[
		X\setminus\bigsqcup_{i\in I}S_{i}V_{i}\prec_{\alpha}\bigsqcup_{i\in I}S'_{i}V_{i}
		\]
		in the sense of Example~\ref{exa:comparison-actions}.
	\end{enumerate}
\end{defn}

\begin{rem}
	We point out a few superficial differences between our definition and Kerr's original definition: 
	\begin{enumerate}
		\item We drop the requirement that $X$ is metrizable, which results in our using a finite cover $\CV$ to control the ``fineness'' of the castle, instead of using a diameter control. 
		\item We do not assume the action is free. 
		\item We have combined some of the universal quantifiers in Kerr's definition. 
	\end{enumerate}
\end{rem}

\begin{thm}
	\label{6.6} \label{thm:Kerr-AE}
	Let $\alpha:\Gamma\curvearrowright X$ be an action of a countable
	discrete group $\Gamma$ on a compact Hausdorff space $X$. 
	Then $\alpha$ is almost finite in the sense of Kerr
	if and only if the transformation
	groupoid $X\rtimes_{\alpha}\Gamma$ is 
	almost elementary with relative reminder control and (ubiquitously) fiberwise amenable.
\end{thm}

Note that by \cite[Remark 5.6]{MWfiberwise}, the transformation groupoid $X\rtimes_{\alpha}\Gamma$ is ubiquitously fiberwise amenable if and only if it is fiberwise amenable if and only if $\Gamma$ is amenable, which is why we put the word ``ubiquitously'' in parentheses. 

\begin{proof}
``$\Leftarrow$'': Suppose $\CG=X\rtimes_\alpha \Gamma$ is almost elementary and (ubiquitous) fiberwise amenable. Let $F\subseteq \Gamma$ be a finite set, $\delta>0$, and $\CV$ an open cover of $X$. Define $K=\bigsqcup_{\gamma\in F} B_{\gamma}$ and $L=X$, which are compact sets in $X\rtimes_\alpha \Gamma$.
Then Proposition~\ref{prop:AE-relative-fiberwise-fortified}\eqref{prop:AE-relative-fiberwise-fortified:Foelner-fortified} shows that there is a non-empty fortified castle $\CF$
		satisfying 
		\begin{enumerate}[label=(\roman*)]
			\item\label{refine2} for any $C \in \CF_+$ with $C \cap L \not= \varnothing$, there is $V \in \CV$ such that $C \subseteq V$;  ;  
			\item $\bigcup \CF_+$ is precompact in $\CG$;  
			\item\label{reminder2} there is a subcollection $\CX \subseteq \CF_-^{(0)}$ such that $|\CX \cap \CT^{(0)}| \leq \varepsilon |\CT^{(0)}|$ for every tower $\CT$ in $\CF_-$ and
			$\GU \setminus \left( \bigcup \CF_-^{(0)} \right) \prec_{\CG}  \bigcup \CX$; 
			\item\label{folner2} for any $u \in \CG\unitSpace$, we have $|K (\bigcup \CF_+) u| \leq (1+\varepsilon) |(\bigcup \CF_+) u|$.  
		\end{enumerate}
        Now, note that $\CG=\bigsqcup_{\gamma\in \Gamma}B_\gamma$, which is a disjoint union of compact open sets $B_\gamma$. Since $\bigcup \CF_+$ is precompact, so is $\CF_-$. Then there is a finite set $M\subseteq \Gamma$ such that $\bigcup \CF_-\subseteq \bigsqcup_{\gamma\in M}B_\gamma$ and therefore by defining $P_\gamma= (\bigcup \CF_-)\cap B_\gamma$, the family $\CP=\{P_\gamma: \gamma\in M\}$ forms an open partition of $\bigcup \CF_-$. Then Lemma~\ref{lem:castle-split} shows that there exists an open castle $\CE$ satisfying $\bigcup \CE=\bigcup \CF_-$ and refining $\CF_-$ and $\CP$. Now, note that $\CE$ refines $\CP_\alpha$ and we decompose $\CE=\bigsqcup_{i\in I} \CT_i$ into towers. The construction in Remark~\ref{rem:dynamical-vs-groupoid-castles}\eqref{rem:dynamical-vs-groupoid-castles:surjective} shows that there exists a dynamical open castle $\{(S_i, V_i): i\in I\}$ such that $V_i\in \CT_i\unitSpace$ is a chosen level and
        \[S_i=\{\gamma\in \Gamma: \text{there exists }  E\in \CT_i \text{ such that } E\subseteq B_\gamma,\text{ and }s(E)=V_i\}.\] Moreover, one has $S_iV_i=\bigcup \CT_i\unitSpace$. For the subcollection $\CX\subseteq \CF^{(0)}_-$ in condition \ref{reminder2}, define $\CX^\CE=\{E\in \CE: E\subseteq C \text{ for some }C\in \CX\}$. Then Proposition \ref{prop:castle-split-more} shows that $\bigcup \CX=\bigcup \CX^\CE$ and $|\CX^\CE\cap \CT_i\unitSpace|\leq \varepsilon|\CT_i\unitSpace|$ for any tower $\CT_i$ in $\CE$. Define 
        \[S'_i=\{\gamma\in \Gamma: \text{there is }  E\in \CT_i \text{ such that } E\subseteq B_\gamma, s(E)=V_i, \text{ and }r(E)\in \CX^\CE\}.\]
        Then, one has $\bigsqcup_{i\in I}S'_iV_i=\bigcup \CX^{\CE}=\bigcup \CX$.

        Now, the condition \ref{refine2} implies that $\CF_+\unitSpace$ refines $\CV$. Therefore, the collection $\{sV_i: s\in S_i, i\in I\}=\CE\unitSpace$ refines $\CV$ because $\CE\unitSpace$ refines $\CF_-\unitSpace$ and thus also refines $\CF_+\unitSpace$. 

        Then for each $i\in I$, choose a $u\in V_i$. By Remark~\ref{rem:dynamical-vs-groupoid-castles} and condition \ref{folner2}, one has
        \begin{align*}
            |FS_i|&=|(\bigcup_{\gamma\in FS_i}B_\gamma)u|=|(\bigcup_{g\in F}B_g)\cdot(\bigcup_{h\in S_i} B_h)u|=|K\cdot (\bigcup\CT_i) u|=|K\cdot (\bigcup\CF_+) u|\\  
             &\leq (1+\varepsilon)|\bigcup \CF_+ u|=(1+\varepsilon)|(\bigcup\CT_i) u|=(1+\varepsilon)|(\bigcup_{h\in S_i}B_h) u|=(1+\varepsilon)|S_i|.
             \end{align*}

        Finally, the condition \ref{reminder2} implies that 
        \[X\setminus \bigsqcup_{i\in I} S_iV_i=\GU\setminus (\bigcup \CF_-\unitSpace)\prec_\CG  \bigcup \CX= \bigcup \CX^\CE=\bigsqcup_{i\in I}S'_iV_i.\]
        This actually shows that $X\setminus \bigsqcup_{i\in I}S_iV_i\prec_\alpha \bigsqcup_{i\in I}S'_iV_i$ as relations $\prec_\CG$ and $\prec_\alpha$ coincide by Example \ref{exa:subequivalence-actions}. Thus, $\alpha$ is almost finite in the sense of Definition \ref{defn:Kerr}.

``$\Rightarrow$'':  Now suppose $\alpha$ is almost finite. Still denote by $\CG=X\rtimes_\alpha \Gamma$ for simplicity. Let $K$ be a compact set in $\CG$, and $\varepsilon>0$, as well as $\CV$ an open cover of $\GU=X$.  The compactness of $K$ implies that there exists a finite $F\subseteq \Gamma$ such that $K\subseteq \bigcup_{\gamma\in F}B_\gamma$. Without loss of generality, one may assume $e_\Gamma\in F$ and $F$ is symmetric, i.e. $F^{-1}=F$. Choose a $\delta>0$ such that $\delta<\min\{\varepsilon/|F^{-1}F|, \varepsilon/(1-\varepsilon)\}$. The almost finiteness of $\alpha$ implies that there exists an open dynamical castle $\{(S_i, V_i): i\in I\}$ such that 
\begin{enumerate}[label=(\roman*)]
		\item the collection $\left\{ sV_{i} \colon s\in S_{i}, i \in I \right\}$ of levels refines $\CV$,  
		\item\label{folner3} $|F S_{i}|\leq(1+\delta)|S_{i}|$ for any $i \in I$, and
		\item there are sets $S_{i}'\subseteq S_{i}$ with $|S_{i}'|< \delta |S_{i}|$ for $i \in I$ such that 
		\[
		X\setminus\bigsqcup_{i\in I}S_{i}V_{i}\prec_{\alpha}\bigsqcup_{i\in I}S'_{i}V_{i}
		\]
     \end{enumerate}   
         First, since $K$ is an arbitrary compact set and $\delta$ could be arbitrary small, the condition \ref{folner3} implies that $\Gamma$ is amenable and therefore $\CG$ is (ubiquitous) fiberwise amenable by \cite[Remark 5.6]{MWfiberwise}.
 Note that from the castle $\{(S_i, V_i): i\in I\}$, one obtains a (groupoid) castle $\CD$ by  Remark~\ref{rem:dynamical-vs-groupoid-castles}\eqref{rem:dynamical-vs-groupoid-castles:assignment}.   Now, for each $i\in I$ define $T_i=\bigcap_{\gamma\in F}\gamma^{-1}S_i$, which is a subset of $S_i$ as $e_\Gamma\in F$. Note that it is a standard fact in the amenable group $\Gamma$ that $|T_i|\geq (1-\varepsilon)|S_i|$ by the condition \ref{folner3} and the choice of $\delta$.
 This process provides us another dynamical open castle $\{(T_i, V_i):i\in I\}$, which also produces an open (groupoid) castle $\CC$ in $\CG$ by Remark~\ref{rem:dynamical-vs-groupoid-castles}\eqref{rem:dynamical-vs-groupoid-castles:assignment}. 
 
 Denote by $K'=\bigcup_{\gamma\in F} B_{\gamma}$ for simplicity, which is a symmetric compact set in $\CG$. We claim $\CC\subseteq \CD^{-K'}$. Indeed, let $y\in \bigcup \CC$. Then, Remark~\ref{rem:dynamical-vs-groupoid-castles}\eqref{rem:dynamical-vs-groupoid-castles:assignment} entails that  $y=(gx, gh^{-1}, hx)$ for some $x\in V_i, g, h\in T_i$ and $i\in I$. Then one has
 \[K'y=\{(\gamma gx, \gamma gh^{-1}, hx): g, h\in T_i,\gamma\in F, x\in V_i\}\subseteq \bigcup \CD\] because $FT_i\subseteq S_i$. Remark \ref{rem: left extendbility enough} then shows that $\CC\subseteq \CD^{-K'}$. This implies $\CC\subseteq \CD^{-K}$ by Definition \ref{defn:castle-shrinking} and the fact $K\cup K^{-1}\subseteq K'$. Therefore, $\CC$ is $K$-extendable.   

 Finally, for each $i\in I$ choose a $T'_i\subseteq T_i$ with $|T'_i|=|S'_i|$. Then $\{tV_i: t\in T'_i, i\in I\}$ determines a  family $\CW\subseteq \CC\unitSpace$ such that $\bigcup \CW=\bigsqcup_{i\in I}T'_iV_i$ and
 \[|\CW\cap \CT_i\unitSpace|=|T'_i|=|S'_i|\leq \delta|S_i|\leq \varepsilon|T_i|=\varepsilon |\CT_i\unitSpace| \]
 for any tower $\CT_i$ in $\CC$. On the other hand, similarly, the collection $\{sV_i: t\in S'_i, i\in I\}$ also determines a family $\CY\subseteq \CD\unitSpace$ such that $\bigcup \CY=\bigsqcup_{i\in I}S'_iV_i$ and $|\CY\cap \CS|=|\CW\cap \CS|$ for any tower $\CS$ in $\CD$. This shows $\bigcup \CY\prec_\CG \bigcup \CW$ by Lemma \ref{lem:castle-comparison}. Thus, one has 
 \[\GU\setminus (\bigcup \CC\unitSpace)=X\setminus \bigsqcup_{i\in I}S_iV_i\prec_\alpha \bigsqcup S'_iV_i=\bigcup \CY\prec_\CG \bigcup \CW,\]
 which implies $\GU\setminus (\bigcup \CC\unitSpace)\prec_\CG \bigcup \CW$ as relations $\prec_\CG$ and $\prec_\alpha$ coincide in this case. Thus, $\CG$ is almost elementary with relative reminder control.
 \end{proof}

It is also helpful to remove ubiquitous fiberwise amenability from the picture and derive a dynamical condition that corresponds to almost elementariness alone. 

\begin{rem}
	In light of Remark~\ref{rem:dynamical-vs-groupoid-castles}, let us clarify what shrinking (see Definition~\ref{defn:castle-shrinking}) means for dynamical castles. Let $\CT=\{(T_i, V_i): i\in I\}$ be a dynamical castle and $F\subseteq\Gamma$ be a finite set. We say $\CT$ is \textit{$F$-shrinkable}, if there exists a dynamical castle $\CS=\{(S_i, V_i): i\in I\}$ such that $FS_i\subseteq T_i$ for each $i\in I$. Denote by $K_F=\bigcup_{\gamma\in F}B_\gamma$ the compact set in the transformation groupoid $X\rtimes_\alpha\Gamma$ and we still denote by $\CT$ the corresponding castle in $X\rtimes_\alpha \Gamma$ generated by $\CT=\{(T_i, V_i): i\in I\}$ as in Remark~\ref{rem:dynamical-vs-groupoid-castles}\eqref{rem:dynamical-vs-groupoid-castles:assignment}. Then it is straightforward to see that $\CT$ is $F$-shrinkable as a dynamical castle if and only if $\CT$, as a castle in $X\rtimes_\alpha\Gamma$ is $K_F$-shrinkable, i.e, $\CT^{-K_F}\neq \emptyset$.
\end{rem}

\begin{thm}
	\label{thm:AE-actions}
	Let $\alpha:\Gamma\curvearrowright X$ be an action of a countable
	discrete group $\Gamma$ on a compact Hausdorff space $X$. 
	Then the transformation
	groupoid $X\rtimes_{\alpha}\Gamma$ is 
	almost elementary with relative reminder control if and only if 
	for any finite symmetric set $F \subseteq \Gamma$, 
	any $\delta>0$, and any finite open cover $\CV$ of $X$,  
	there is an open castle $\{(S_{i},V_{i}):i\in I\}$ such that 
	\begin{enumerate}[label=(\roman*)]
		\item\label{thm:AE-actions:refine} the collection $\left\{ sV_{i} \colon s\in S_{i}, i \in I \right\}$ of levels refines $\CV$,  
		\item\label{thm:AE-actions:extend} $\{(F S_{i},V_{i}):i\in I\}$ is again a castle, and
        
		\item\label{thm:AE-actions:remainder} there are sets $S_{i}'\subseteq S_{i}$ with $|S_{i}'|< \delta |S_{i}|$ for $i \in I$ such that 
		\[
		X\setminus\bigsqcup_{i\in I}S_{i}V_{i}\prec_{\alpha}\bigsqcup_{i\in I}S'_{i}V_{i}
		\]
		in the sense of Example~\ref{exa:comparison-actions}
	\end{enumerate}
\end{thm}

\begin{proof}
The proof is similar to the proof of Theorem \ref{thm:Kerr-AE} with the help of Lemma~\ref{lem:AE-relative-fortified}.

``$\Rightarrow$'': Without loss of generality, let $F\subseteq \Gamma$ be a finite symmetric set containing the identity $e_\Gamma$, a number $\delta>0$, and an open cover $\CV$ of $X$. Define a compact set $K=\bigsqcup_{\gamma\in F}B_\gamma$ and $L=X$. Note that $K$ is also symmetric in $\CG$. 
Lemma~\ref{lem:AE-relative-fortified} implies that for the compact sets $K$ and $L$ in $\CG$, the $\delta > 0$, and the finite open cover $\CV$ of $\CG^{(0)}$, there is a fortified castle $\CF$
satisfying 
\begin{enumerate}[label=(\roman*)]
		\item 
		for any $F_+ \in \CF_+$ with $F_+ \cap L \not= \varnothing$, there is $V \in \CV$ such that $F_+ \subseteq V$; 
		\item
		$\bigcup \CF_+$ is precompact in $\CG$;  
		\item\label{reminder3}
		there is a subcollection $\CW \subseteq \left( \CF^{-K} \right)_-^{(0)}$ such that $|\CW \cap \CT^{(0)}| \leq \delta |\CT^{(0)}|$ for every tower $\CT$ in $\left( \CF^{-K} \right)_-$ and
		$\CG^{(0)}\setminus\ \left( \bigcup \left( \CF^{-K} \right)_-^{(0)} \right) \prec_{\CG}  \bigcup \CW$.   
	\end{enumerate}
Then, the precompactness of $\CF_+$ implies that $\CF_-$ is precompact as well. Since $\CG=\bigcup_{\gamma\in \Gamma}B_\gamma$, there is a finite set $M\subseteq \Gamma$ such that $\bigcup \CF_-\subseteq \bigsqcup_{\gamma\in M}B_\gamma$.  Define $P_\gamma=B_\gamma\cap (\bigcup \CF_-)$ for $\gamma\in M$. Then the collection $\{P_\gamma: \gamma\in M\}$ forms an open partition of $\bigcup \CF_-$. Now Lemma~\ref{lem:castle-split} shows that there exists an open castle $\CE$ satisfying $\bigcup \CE=\bigcup \CF_-$ and refining $\CF_-$ and $\CP$. By Definition \ref{defn:castle-shrinking}, this further implies that $\bigcup \CF_-^{-K}=\bigcup \CE^{-K}$. Note that $\CE$ refines $\CP_\alpha$ and we decompose $\CE=\bigsqcup_{i\in I} \CT_i$ into towers. The construction in Remark~\ref{rem:dynamical-vs-groupoid-castles}\eqref{rem:dynamical-vs-groupoid-castles:surjective} shows that there exists a dynamical open castle $\{(T_i, V_i): i\in I\}$ such that $V_i\in \CT_i\unitSpace$ is a chosen level and
        \[T_i=\{\gamma\in \Gamma: \text{there exists }  E\in \CT_i \text{ such that } E\subseteq B_\gamma,\text{ and }s(E)=V_i\}.\] Moreover, one has $T_iV_i=\bigcup \CT_i\unitSpace$.
Then, as $\CE^{-K}$ is a subcastle of $\CE$, it corresponds to a dynamical subcastle $\{(S_i, V_i): i\in I\}$ of $\{(T_i, V_i): i\in I\}$ with $S_i\subseteq T_i$ and $\bigcup (\CE^{-K})^{(0)}=\bigcup_{i\in I}S_iV_i$. Then let $x\in sV_i$ for some $s\in S_i$. Then by the definition of $\CE^{-K}$, one has $Kx\subseteq \bigcup\CE$, which implies that $Fs\subseteq T_i$. Therefore, one has $FS_i\subseteq T_i$ for any $i\in I$, from which, the collection $\{(FS_i, V_i): i\in I\}$ is a subcastle of $\{(T_i, V_i): i\in I\}$. 

Finally, define $\CX=\{E\in \CE^{-K}: E\subseteq D\text{ for some }D\in \CW\}$, which corresponds to a subcastle $\{(S'_i, V_i): i\in I\}$ of $\{(S_i, V_i): i\in I\}$. Then using Proposition \ref{prop:castle-split-more}, the same argument as in proving Theorem \ref{thm:Kerr-AE} shows that $|S'_i|<\delta|S_i|$ holds for $i\in I$ and 
\[X\setminus\bigsqcup_{i\in I}S_{i}V_{i}\prec_{\alpha}\bigsqcup_{i\in I}S'_{i}V_{i}.\]

``$\Leftarrow$'': Let a compact set $K$ in $\CG$, a number $\varepsilon>0$, and an open cover $\CV$ of $X$ be given. Choose a finite symmetric set $F$ containing $e_\Gamma$ such that $K\subseteq \bigsqcup_{\gamma}B_\gamma$. Then for the $F$, $\varepsilon$ and the $\CV$, by assumption, there is an open castle $\{(S_{i},V_{i}):i\in I\}$ such that 
	\begin{enumerate}[label=(\roman*)]
		\item the collection $\left\{ sV_{i} \colon s\in S_{i}, i \in I \right\}$ of levels refines $\CV$,  
		\item $\{(F S_{i},V_{i}):i\in I\}$ is again a castle, and
        
		\item there are sets $S_{i}'\subseteq S_{i}$ with $|S_{i}'|< \delta |S_{i}|$ for $i \in I$ such that 
		\[
		X\setminus\bigsqcup_{i\in I}S_{i}V_{i}\prec_{\alpha}\bigsqcup_{i\in I}S'_{i}V_{i}
		\]
		in the sense of Example~\ref{exa:comparison-actions}
	\end{enumerate}
    Note that the castles $\{(FS_i, V_i): i\in I\}$ and $\{(S_i, V_i): i\in I\}$ provides open (groupoid) castles $\CD$ and $\CC$ in $\CG$, respectively, by Remark~\ref{rem:dynamical-vs-groupoid-castles}\eqref{rem:dynamical-vs-groupoid-castles:assignment}. First, $\CC\unitSpace=\{sV_i: s\in S_i, i\in I\}$ refines $\CV$.   Then the same argument in proving Theorem \ref{thm:Kerr-AE} shows that $\CC\subseteq \CD^{-K}$. Finally, choose $\CW\subseteq\CC\unitSpace$ to be $\CW=\{sV_i: s\in S'_i, i\in I\}$. Then the construction in  Remark~\ref{rem:dynamical-vs-groupoid-castles}\eqref{rem:dynamical-vs-groupoid-castles:assignment} shows that $|\CW\cap \CT_i\unitSpace|=|S'_i|\leq \delta|S_i|=\delta|\CT_i\unitSpace|$ for any tower $\CT_i$ in $\CC$. And the same argument in proving Theorem \ref{thm:Kerr-AE} also shows
    $\GU\setminus (\bigcup \CC\unitSpace)\prec_\CG \bigcup \CW$.
    \end{proof}

\section{Nesting of castles} \label{sec:nesting}
In this section, we introduce a notion called \textit{nesting} that involves a pair of castles.  Roughly speaking, this means a castle lives in another coarser castle with compatible dynamics, from which one can construct desired c.p.c.\ maps witnessing regularity properties such as the tracial $\CZ$-stability for the groupoid $C^*$-algebra. Therefore, this concept plays an important role in the proof of our main result in Section \ref{sec:Z-stable}.

\begin{defn}
	\label{defn:nesting} 
	Let $\CC, \CD$ be castles. We say $\CC$ is \emph{nested} in $\CD$ if 
	\[
		\left\{  DC \colon C \in \CC, \, D \in \CD \right\} \setminus \{ \varnothing \} = \CC \; .
	\]
	The nesting is \emph{faithful} if for any $D \in \CD$, there is $C \in \CC$ with $DC \not= \varnothing$. 
\end{defn}

\begin{rem}\label{rem: left-right-nesting}
	It is equivalent to define nesting of $\CC$ in $\CD$ by $\{CD  \colon C\in \CC, D\in \CD\}\setminus\{\varnothing\}=\CC$. This is mainly because one has $(CD)^{-1}=D^{-1}C^{-1}$ as well as $C^{-1}\in \CC$ and $D^{-1}\in \CD$ for any $C\in \CC$ and $D\in \CD$.

	An immediate consequence of the definition is that $\{CD \colon C\in \CC\unitSpace, D\in \CD\unitSpace\}\setminus\{\varnothing\}=\CC\unitSpace$. In particular, it follows that $\CC\unitSpace$ refines $\CD\unitSpace$. 
	
	It is also easy to see that the nesting is faithful if and only if for any $D \in \CD \unitSpace$, there is $C \in \CC \unitSpace$ with $C \subseteq D$. 
\end{rem}

\begin{prop}\label{prop: tower nest-abstract}
A castle $\CC$ is nested in a castle $\CD$ if and only if all towers in $\CC$ are nested in $\CD$.
\end{prop}
\begin{proof}
 Suppose $\CC$ is nested in $\CD$, which means $\{DC: C\in \CC, D\in \CD\}\setminus \{\varnothing\}=\CC$ by definition. First enumerate all towers in $\CC$ by $\CC_1, \dots, \CC_n$. Let $l\leq n$.  For any $C\in \CC_l$ and $D\in\CD$, we claim that if $DC\neq \varnothing$ then $DC\in \CC_l$.  Indeed, suppose $DC\neq \varnothing$, then $DC=C'\in \CC$.  Now, $s(C')\cap s(C)\neq \varnothing$ implies $s(C')=s(C)$ by Lemma \ref{lem:plurisection-basic}\eqref{lem:plurisection-basic:composability} and therefore $C'\in \CC_l$ since $C'\sim_\CC C$. This shows $\{DC: C\in \CC_l, D\in \CD\}\setminus \{\varnothing\}\subseteq \CC_l$. For the converse, for any $C\in \CC_l$, since $\CC$ is nested in $\CD$, there is a $C'\in \CC$ and a $D\in \CD$ such that $DC'=C$ and thus $D^{-1}C=s(D)C'=C'$ necessarily. This implies $C'\in \CC_l$ by the argument above and thus one has $\CC_l\subseteq \{DC': C'\in \CC_l, D\in \CD\}\setminus \{\varnothing\}$. Therefore, $\CC_l$ is nested in $D$. 
 
 For the converse direction, suppose all towers $\CC_l$ are nested in $\CD$. Then observe that 
 \[\{DC: C\in \CC, D\in \CD\}\setminus \{\varnothing\}=\bigcup_{l\leq n}(\{DC: C\in \CC_l, D\in \CD\}\setminus\{\varnothing\})=\bigcup_{l\leq n}\CC_l=\CC,\]
 which means that $\CC$ is nested in $\CD$.
\end{proof}

Let us reinterpret Definition~\ref{defn:nesting} through the lens of multisections in the spirit of Proposition~\ref{prop:tower-multisection} and Corollary~\ref{cor:castle-multisections}, which provide a ``visual'' description of nesting.  In light of Proposition \ref{prop: tower nest-abstract}, it suffices to describe how a tower/multisection is nested in a castle.

\begin{lem}
	\label{lem:nesting-towers-explicit} 
	Let $\CC=\{C_{i,j}:i,j\in F\}$ be a multisection and let $\CalD=\{D^l_{m,n}: m,n\in E_l, l\in I\}$ be a castle as in Corollary~\ref{cor:castle-multisections}. Then $\CC$ is nested in $\CD$ if and only if the following hold: 
	\begin{enumerate}[label=(\roman*)]
		\item \label{lem:nesting-towers-explicit::level} for any $\CC$-level $C_{i,i}\in\CUU$, there is a $\CalD$-level
		$D^l_{n,n}\in\DU$ with $C_{i,i}\subseteq D^l_{n,n}$;  
		\item \label{lem:nesting-towers-explicit::ladder} for any $l\in I$, $m,n\in E_l$ and $\CC$-level $C\subseteq D^l_{n, n}$, one has $D^l_{m, n}C\in \CC$.
	\end{enumerate}
\end{lem}
\begin{proof}
	Suppose $\{DC: C\in\CC, D\in\CD\}\setminus \{\varnothing\}=\CC$ holds. Now for any $\CC$-level $C_{i, i}$, there are $C\in \CC$ and $D\in \CD$ such that $DC=C_{i, i}$. Let $D^l_{n, n}$ be the $\CD$-level $r(D)$. One has \[C_{i, i}=r(C_{i, i})=r(DC)\subseteq r(D)= D^l_{n, n}.\]
    This has established (i). For (ii), let $m, n\in E_l$ be given. Note that for any $\CC$-level $C\subseteq D^l_{n, n}$, the product $D^l_{m,n}C\neq \varnothing$ and thus $D^l_{m, n}C\in \CC$ holds by our assumption.
	
	Now suppose (i) and (ii) hold. For any $C\in \CC$, (i) implies that there is a $\CD$-level $D^l_{n, n}$ such that $r(C)\subseteq D^l_{n, n}$. Now since $\CD$ is a castle, actually such the existence of the level is unique. This means if $D\in \CD$ such that $DC\neq \varnothing$, then $s(D)=D^l_{n, n}$ necessarily. Then (ii) above entails that $Dr(C)\in \CC$ and thus $DC=Dr(C)\cdot C\in \CC$. This shows that $\{DC: C\in \CC, D\in \CD\}\setminus \{\varnothing\}\subseteq \CC$. For the converse direction,  (i) implies that, for any $C\in \CC$,  there are an $l\in I$ and an $m\in E_l$ such that $r(C)\subseteq D^l_{m,m}$. Thus, one has $C=D^l_{m, m}r(C)C=D^l_{m, m}C\neq \varnothing$, which entails that $\{DC: C\in \CC, D\in \CD\}\setminus \{\varnothing\}= \CC$.
\end{proof}

The following schematic illustrates the concept of nesting by keeping track of the levels of $\CC$ and $\CD$. Here $\CC$ is a tower with height $7$, whose levels are indicated by the small circles, and $\CD$ is a castle with two towers of heights $2$ and $3$, whose levels are indicated by the larger ovals. The dotted lines and curves indicate how the levels group into towers. 

\begin{center}
\begin{tikzpicture}[scale=0.5]
    \draw[line width=0.04cm] (0,-6.0) rectangle (17.2,6.0);
    
    \draw[line width=0.04cm] (5.2,2.0) ellipse (4.0cm and 1.6cm);
    \draw[line width=0.04cm] (5.2,-2.4) ellipse (4.0cm and 1.6cm);
    \draw[line width=0.04cm] (12.2,4.4) ellipse (2.2cm and 1.2cm);
    \draw[line width=0.04cm] (12.2,0.8) ellipse (2.2cm and 1.2cm);
    \draw[line width=0.04cm] (12.2,-2.8) ellipse (2.2cm and 1.2cm);
    
    \draw[line width=0.04cm] (3.2,2.0) ellipse (0.8cm and 0.8cm);
    \draw[line width=0.04cm] (3.2,-2.4) circle (0.8cm);
    \draw[line width=0.04cm] (6.8,2.4) circle (0.8cm);
    \draw[line width=0.04cm] (6.8,-2.0) circle (0.8cm);
    \draw[line width=0.04cm] (11.6,4.4) circle (0.8cm);
    \draw[line width=0.04cm] (11.6,0.8) circle (0.8cm);
    \draw[line width=0.04cm] (11.6,-2.8) circle (0.8cm);
    
    \draw[line width=0.04cm, dotted] (1.2,2.0) -- (1.2,-2.4);
    \draw[line width=0.04cm, dotted] (9.2,2.0) -- (9.2,-2.4);
    \draw[line width=0.04cm, dotted] (10.0,4.4) -- (10.0,0.8);
    \draw[line width=0.04cm, dotted] (10.0,0.8) -- (10.0,-2.8);
    \draw[line width=0.04cm, dotted] (14.4,4.4) -- (14.4,0.8) -- (14.4,-2.8);
    \draw[line width=0.04cm, dotted] (2.4,2.0) -- (2.4,-2.4);
    \draw[line width=0.04cm, dotted] (4.0,2.0) -- (4.0,-2.4);
    \draw[line width=0.04cm, dotted] (6.0,2.4) -- (6.0,-2.0);
    \draw[line width=0.04cm, dotted] (7.6,2.4) -- (7.6,-2.0);
    \draw[line width=0.04cm, dotted] (10.8,4.4) -- (10.8,-2.8);
    \draw[line width=0.04cm, dotted] (12.4,4.4) -- (12.4,-2.8);
    
    \draw[line width=0.04cm, dotted] (2.4,2.0) .. controls (2.4,6.0) and (7.6,5.6) .. (7.6,2.4);
    \draw[line width=0.04cm, dotted] (4.0,2.0) .. controls (4.0,4.4) and (6.0,4.0) .. (6.0,2.4);
    \draw[line width=0.04cm, dotted] (7.6,-2.0) .. controls (7.6,-4.8) and (10.8,-4.8) .. (10.8,-2.8);
    \draw[line width=0.04cm, dotted] (6.0,-2.0) .. controls (6.4,-7.2) and (12.4,-6.0) .. (12.4,-2.8);
\end{tikzpicture}
\end{center}

Let us explore certain algebraic structures arising from a nesting of two castles. 

\begin{defn} \label{defn:nesting-comm-red} 
	Let $\CC, \CD$ be castles with $\CC$ {nested} in $\CD$. We define the following (collections of) collections of subsets in $\CG$: 
	\begin{enumerate}
		\item for any $D \in \CD$, $\RedCD{\CC}{D} := \left\{ C D  \colon C \in \CC\unitSpace \right\} \setminus \{\varnothing\} \subseteq \CC$; 
		\item $\RedCD{\CC}{\CD} := \left\{ \RedCD{\CC}{D} \colon D \in \CD \right\} \setminus \{\varnothing\}$; 
		\item for any $C \in \CC$, $\CommCD{C}{\CD} := \left\{ D C D^{-1} \colon D \in \CD \right\} \setminus \{\varnothing\} \subseteq \CC$; 
		\item $\CommCD{\CC}{\CD} := \left\{ \CommCD{C}{\CD} \colon C \in \CC \right\} \setminus \{\varnothing\}$;  
		\item for any $D \in \CD\unitSpace$, $\ResCD{\CC}{D} := \left\{ D C D \colon C \in \CC \right\} \setminus \{\varnothing\} \subseteq \CC$. 
	\end{enumerate}
	Here ``Red'' stands for ``reduction'', ``Comm'' stands for ``commutant'', and ``Res'' stands for ``restriction''
\end{defn}

\begin{rem} \label{rem:nesting-comm-red-basic} 
	Let $\CC, \CD$ be castles with $\CC$ {nested} in $\CD$. Recall from Lemma~\ref{lem:castle-basic}\eqref{lem:castle-basic:subquotient} that $\CC$ is a finite principal groupoid with the induced operations. We make the following observations: 
	\begin{enumerate}
		\item  \label{rem:nesting-comm-red-basic:Red} 
		For any $D \in \CD$, $\RedCD{\CC}{D} = \left\{ D C \colon C \in \CC\unitSpace
		 \right\} \setminus \{\varnothing\} = \left\{ C \in \CC \colon C \subseteq D \right\}$. 
		 
		 \item  \label{rem:nesting-comm-red-basic:Red-Res} 
		 For any $D \in \CD \unitSpace$, $\RedCD{\CC}{D} = \left\{ C \in \CC \unitSpace \colon C \subseteq D \right\} = \ResCD{\CC}{D} \cap \CC \unitSpace$. 
		 
		\item  \label{rem:nesting-comm-red-basic:Comm-empty} 
		For any $C \in \CC$, $\CommCD{C}{\CD} \not= \varnothing$ if and only if there is 
		$D \in \CD\unitSpace$ such that $r(C) \cup s(C) \subseteq D$ 
		if and only if there is 
		$D \in \CD\unitSpace$ such that $C \in \ResCD{\CC}{D} $. 
		
		\item  \label{rem:nesting-comm-red-basic:Comm-same} 
		For any $C, C' \in \CC$ with $\CommCD{C}{\CD} \not= \varnothing$ and $\CommCD{C'}{\CD} \not= \varnothing$, we have $\CommCD{C}{\CD} = \CommCD{C'}{\CD}$ 
		if and only if $\CommCD{C}{\CD} \cap \CommCD{C'}{\CD}$ is not empty 
		if and only if there is 
		$D \in \CD$ such that $D C D^{-1} = C'$. 
		
		\item  \label{rem:nesting-comm-red-basic:castles-C} 
		In the finite principal groupoid $\CC$, $\RedCD{\CC}{\CD}$ and $\CommCD{\CC}{\CD}$ are castles. 
		
		\item  \label{rem:nesting-comm-red-basic:castles-G} 
		For any $D \in \CD\unitSpace$, $\ResCD{\CC}{D}$ is a castle in $\CG$, whose towers are of the form $\ResCD{\CT}{D}$ where $\CT$'s are towers of $\CC$. It is nonempty if and only if $D \cap ( \bigcup \CC) \not= \varnothing$. 
		
		\item  \label{rem:nesting-comm-red-basic:isomorphism} 
		If the nesting is faithful, then the map $\rho \colon \CD \to \RedCD{\CC}{\CD}$ given by $D \mapsto \RedCD{\CC}{D}$ is an isomorphism of finite groupoids, where both groupoid structures are given by Lemma~\ref{lem:castle-basic}\eqref{lem:castle-basic:subquotient}.

		\item  \label{rem:nesting-comm-red-basic:Res-Comm} 
		For any $D \in \CD\unitSpace$, the map $\gamma_D \colon \ResCD{\CC}{D} \to \CommCD{\CC}{\CD}$ given by $C \mapsto \CommCD{C}{\CD}$ yields an injective groupoid homomorphism. Moreover, for any $D, D' \in \CD\unitSpace$, we have $\gamma_D \left( \ResCD{\CC}{D} \right) = \gamma_{D'} \left( \ResCD{\CC}{D'} \right)$ if and only if $D$ and $D'$ are in the same tower in $\CD$. Consequently, these maps induce a groupoid isomorphism 
		\[
			\CommCD{\CC}{\CD} \cong \bigsqcup_{[D] \in \CD\unitSpace / \sim_{\CD}} \ResCD{\CC}{D} 
			\cong \bigsqcup_{[D] \in \CD\unitSpace / \sim_{\CD}} \bigsqcup_{\CT \in \CC / \sim_{\CC}} \ResCD{\CT}{D} \; .
		\]
		
		\item \label{rem:nesting-comm-red-basic:commutation} 
		For any $C \in \CC$ and $D \in \CD$, the set product $\RedCD{\CC}{D} \CommCD{C}{\CD}$, if non-empty, in the groupoid $\CC$ is a singleton, and we have 
		\[
			\RedCD{\CC}{D} \CommCD{C}{\CD} = \CommCD{C}{\CD} \RedCD{\CC}{D} \; .
		\]
		Moreover, if $\CD$ is a tower, then every singleton in $\CC$ can be written as such a product. 
		
		\item \label{rem:nesting-comm-red-basic:algebras} 
		The linear map $C^*(\CD) \to C^*(\CC)$ defined by the prescription $\chi_{D} \mapsto \chi_{\RedCD{\CC}{D}}$ is an injective $*$-homomorphism, and the commutant of its image in $C^*(\CC)$ has $\left\{ \chi_{\CommCD{C}{\CD}} \colon C \in \CC \right\}$ as a linear basis. 
	\end{enumerate}
\end{rem}

In view of Remark~\ref{rem:nesting-comm-red-basic}\eqref{rem:nesting-comm-red-basic:commutation} and~\eqref{rem:nesting-comm-red-basic:algebras}, the following fact is not surprising. 

\begin{lem} \label{lem:nesting-comm-red-order-zero}
	Let $\CC, \CD$ be open castles with $\CC$ \emph{nested} in $\CD$. Let $f$ (respectively, $g$) be a positive contraction in $C_0 \left(\bigcup \CC \right) \cap C_{\operatorname{c}}(\CG)$ (respectively, $C_0 \left(\bigcup \CD \right) \cap C_{\operatorname{c}}(\CG)$) that is $\left(\bigcup \CC\right)$-invariant (respectively, $\left(\bigcup \CD\right)$-invariant). Denote by $\psi_{\chi_{\CommCD{\CC}{\CD}}}$ the canonical isomorphic embedding of $C^*$-algebra $C^*(\CommCD{\CC}{\CD})$ to $C^*(\CC)$. Let $\psi_f \colon C^*(\CC) \to C^*(\CG)$ and $\psi_g \colon C^*(\CD) \to C^*(\CG)$ be c.p.c.\ order zero maps as defined in Proposition~\ref{prop:castle-order-zero}. Suppose that $g(u) = 1$ for any $u \in \bigcup \CC\unitSpace$. Then the images of $\psi_f\circ \psi_{\chi_{\CommCD{\CC}{\CD}}}$ and $\psi_g$ commute. 
\end{lem}

\begin{proof}
Note that $\psi_{\chi_{\CommCD{\CC}{\CD}}}$ maps the Dirac function $\chi_{\CommCD{C}{\CD}} \mapsto \sum_{D\in \CD}\chi_{DCD^{-1}}$. Then $\psi_f\circ \psi_{\chi_{\CommCD{\CC}{\CD}}}$ maps  each  $\chi_{\CommCD{C}{\CD}} \mapsto \sum_{D\in \CD} f\chi_{DCD^{-1}}$. Recall that the function $\psi_g$ maps $\chi_D\mapsto g\chi_D$. Then it suffices to verify  
\[(\sum_{D\in \CD} f\chi_{DC_0D^{-1}})* g\chi_{D_0}=g\chi_{D_0}*(\sum_{D\in \CD} f\chi_{DC_0D^{-1}})\]
for any $C_0\in \CC, D_0\in \CD$ as $\{\CommCD{C}{\CD}: C\in \CC\}$ and $\CD$ form linear bases for $C^*(\CommCD{\CC}{\CD})$ and $C^*(\CD)$, respectively. To this end, let $C_0\in \CC$, and $D_0, D\in \CD$ be given. Then one has 
\begin{align*}
   &(f\chi_{DC_0D^{-1}}*g\chi_{D_0})(x)=\sum_{\{(y, z)\in \CG^2: x=yz\}}(f\chi_{DC_0D^{-1}})(y)\cdot (g\chi_{D_0})(z)\\
   & =\begin{cases}
			0 \, ,  & \text{if } x\notin DC_0C^{-1}\cdot D_0\\
			 f((r|_{DC_0D^{-1}})^{-1}(r(x)))g((s|_{D_0})^{-1}(s(x))\, , & \text{if } x\in DC_0D^{-1}\cdot D_0
		\end{cases}
        \end{align*}
and
\begin{align*}
 &(g\chi_{D_0}*f\chi_{DC_0D^{-1}})(x)=\sum_{\{(z, y)\in \CG^2: x=zy\}}(g\chi_{D_0})(z)\cdot(f\chi_{DC_0D^{-1}})(y)\\
   & =\begin{cases}
			0 \, ,  & \text{if } x\notin D_0\cdot DC_0C^{-1}\\
			 g((r|_{D_0})^{-1}(r(x)))f((s|_{DCD^{-1}})^{-1}(s(x)))\, , & \text{if } x\in D_0\cdot DC_0D^{-1}
		\end{cases}.
        \end{align*}
In addition, note that for $\CommCD{C_0}{\CD}=\{DC_0D^{-1}: D\in \CD\}\setminus \{\emptyset\}$,
there exists unique $D_1, D_2\in \CD$ such that $D_1C_0D^{-1}_1\cdot D_0\neq \emptyset$ and $D_0\cdot D_2C_0D_2^{-1}\neq \emptyset$ by Remark \ref{rem:nesting-comm-red-basic}\eqref{rem:nesting-comm-red-basic:Comm-empty}. 
This implies
\[  (\sum_{D\in \CD}f\chi_{DC_0D^{-1}})*g\chi_{D_0}=f\chi_{D_1C_0D^{-1}_1}*g\chi_{D_0}\]
and
\[(\sum_{D\in \CD}g\chi_{D_0}* f\chi_{DC_0D^{-1}})=g\chi_{D_0}*f\chi_{D_2C_0D_2^{-1}}.\]
Moreover, in this case, the equality  $D_1C_0D^{-1}_1\cdot D_0=D_0\cdot D_2C_0D_2^{-1}$ holds. Let $x\in D_1C_0D^{-1}_1\cdot D_0=D_0\cdot D_2C_0D_2^{-1}$ and denote by $y_1=(r|_{D_1C_0D_1^{-1}})^{-1}(r(x))$, $z_1=(s|_{D_0})^{-1}(s(x))$, $y_2=(s|_{D_2CD_2^{-1}})^{-1}(s(x))$ and $z_2=(r|_{D_0})^{-1}(r(x))$ for simplicity such that $x=y_1z_1=z_2y_2$. Since $\CC$ is nested in $\CD$, the bisection  $D_1C_0D^{-1}_1\cdot D_0=D_0\cdot D_2C_0D_2^{-1}\in \CC$. Consequently, one has that $x\in \bigcup \CC$ and thus $z_1, z_2\in \bigcup\CC$ as well. Furthermore, this implies $y_1=z_2y_2z^{-1}_1$, i.e. $y_1\sim_{\bigcup \CC}y_2$. Finally, By assumption, the $\bigcup \CD$-invariant function $g$ satisfies that  $g\equiv 1$ on $\bigcup \CC\unitSpace$. This shows that $g(z_1)=g(z_2)=1$. Combining these, one has
\[(f\chi_{D_1C_0D^{-1}_1}*g\chi_{D_0})(x)=f(y_1)=f(y_2)= (g\chi_{D_0}*f\chi_{D_2C_0D_2^{-1}})(x)\]
because $f$ is $\bigcup \CC$-invariant. Therefore, one indeed has
\[(\sum_{D\in \CD} f\chi_{DC_0D^{-1}})* g\chi_{D_0}=g\chi_{D_0}*(\sum_{D\in \CD} f\chi_{DC_0D^{-1}})\]
as for any $x\notin  D_1C_0D^{-1}_1\cdot D_0=D_0\cdot D_2C_0D_2^{-1}$, the functions on both sides have value $0$. 
\end{proof}

Finally, let us show that the definition of almost elementariness (with relative remainder control) for minimal groupoid can be upgraded to an equivalent characterization involving nestings. The following result should be compared with Definition~\ref{defn:AE} and Lemma ~\ref{lem:AE-fortified}. On the other hand, For the non-trivial ``only if'' direction of Proposition \ref{prop:AE-nesting},  we remark that the minimality is used only in the small reminder part \ref{prop:AE-nesting-minimal:remainder} and Proposition~\ref{lem:AE-relative-fortified}, a relative version of  Lemma ~\ref{lem:AE-fortified}, allows to establish the other three conditions. 

\begin{prop} \label{prop:AE-nesting}
	A minimal \'{e}tale groupoid $\CG$ with a compact unit space is {almost elementary}
	if and only if 
	for any nonempty compact set $K$ in $\CG$,  a non-empty open set $O_0$ in $\GU$,
	any fortified castle $\CE$ in $\CG$, 
	and any finite collection $\CV$ of open sets in $\CG^{(0)}$ covering $\CE_-$, 
	there are a fortified castle $\CF$ and an open castle $\CD$ 
	satisfying 
	\begin{enumerate}[label=(\roman*)]
		\item\label{prop:AE-nesting-minimal:nesting} $\CF_+$ is nested in $\CE_+$; 
		\item\label{prop:AE-nesting-minimal:extend} $\CF_+ \subseteq \CD^{-K}$ (in particular, $\CF_+$ is $K$-extendable);  
		\item\label{prop:AE-nesting-minimal:refine} $\CD\unitSpace$ refines $\CV$;    
		\item\label{prop:AE-nesting-minimal:remainder} 
		$\left( \overline{\bigcup \CE_-^{(0)}} \right)\setminus\ \left( \bigcup \CF_-^{(0)} \right) \prec_{\CG}  O_0$. 
	\end{enumerate}
\end{prop}

\begin{proof}
	For the ``if'' direction, we observe that $\left\{\CG\unitSpace\right\}$ is a clopen castle in which any open castle in $\CG$ is nested. Hence if we fix $\CE$ in the above to be the fortified castle $\left(\left\{\CG\unitSpace\right\}, \left\{\CG\unitSpace\right\}\right)$, we see the condition above verifies Definition~\ref{defn:AE}. 
	
	For the ``only if'' direction, we fix $K$, $\CV$ and $\CE$ as in the statement and seek a fortified castle $\CF$ with the desired properties. 
	Without loss of generality, we may assume that  $\CE_+$ is precompact, 
    since otherwise we may apply Remark~\ref{rem:fortified-castle-basics}\eqref{rem:fortified-castle-basics:sandwich-exist} and~\eqref{rem:fortified-castle-basics:sandwich} to replace $\CE$ by another fortified castle $\CE'$ with $\CE'_- = \CE_-$ and $\CE'_+$ sandwiched in $\CE$, to the effect that $\CE'_+$ is precompact by Remark~\ref{rem:fortified-castle-basics}\eqref{rem:fortified-castle-basics:precompact} and the compactness of $\CG\unitSpace$. In addition, we may also assume that $\GU\subseteq K$.
    
    Define compact sets $L_1=\overline{\bigcup \CE_-}\cup \GU$ and $L_2=\overline{\bigcup \CE_+}\cup \GU$.
 Note that $\CE_+$ is an open cover of $\overline{\bigcup \CE_-}$. Choose another open cover $\CU$ of $L_1$ refining $\CE_+\cup\CV$ and $(\CE_+\unitSpace\cup \{\GU\setminus \overline{\bigcup\CE_-\unitSpace}\})\vee \CV$. In addition, since $\CE_+$ is a castle, for any $U\in \CU$, if $U\cap \overline{E_-}\neq \emptyset$ then  $E_+\in \CE_+$ is the unique element in $\CE_+$ such that $U\subseteq E_+$.

Now apply Lemma~\ref{lem:AE-fortified} for the compact sets $KL_2, L_1$ and the open cover $\CU$ of $L_1$ to obtain a fortified castle $\CD$ satisfying: 
	\begin{enumerate}
		\item\label{bisection cover} for any $D_+ \in \CD_+$ with $D_+ \cap L_1 \not= \varnothing$, there is a $U\in \CU$ such that $D_+ \subseteq U$; 
		\item 
		$\bigcup \CD_+$ is precompact in $\CG$;  
		\item\label{small reminder nesting0}
		$\CG^{(0)}\setminus\ \left( \
        \bigcup \left( \CD^{-KL_2}_- \right)^{(0)} \right) \prec_{\CG} O_0$.   
	\end{enumerate}
Define
    \begin{align*}
        \CA:=\{&D_+\in \CD^{-K}_+: D_+\cap \overline{E_{-}}D_{+, 0}\overline{E'_-}\neq \emptyset 
    \text{ for some }
    D_{+, 0}\in (D^{-KL_2}_+)^{(0)}
    \text{ and }\\
    &(E_+, E_-), (E'_+, E'_-)\in \CE\}.  
    \end{align*}
We denote by 
 \[\CA\unitSpace:=\{s(D_+), r(D_+): D_+\in \CA\},\]
   from which one obtains a fortified open castle 
    \[\CF=\{(D_+, D_-)\in \CD^{-K}: s(D_+), r(D_+)\in \CA\unitSpace\}\]
    by Lemma \ref{lem:castle-basic}\eqref{lem:castle-basic:reduction}. 
We now claim that the collection $\CA$ coincides with the family
\begin{align*}
        \CB:=\{&E_{+}D_{+, 0}E'_+:
    D_{+, 0}\in (D^{-KL_2}_+)^{(0)}
    ,(E_+, E_-), (E'_+, E'_-)\in \CE \text{ such that }\\ 
    &\overline{E_-}D_{+, 0}\overline{E'_-}\neq \emptyset, D_{+, 0}\subseteq r(E'_+),\text{ and } D_{+, 0}\subseteq s(E_+)\}.  
    \end{align*}
 To this end, first let $D_+\in \CA$. Then there exists a $D_{+, 0}\in (\CD^{-KL_2}_+)^{(0)}$ and $(E_{+, 1}, E_{-,1})$, $(E_{+, 2}, E_{-, 2})$ in the fortified castle $\CE$ such that  $D_+\cap \overline{E_{-, 1}}D_{+, 0}\overline{E_{-, 2}}\neq \emptyset$, which implies that $D_+\cap \overline{E_{-}}\neq \emptyset$ where $E_-=E_{-, 1}E_{-, 2}$. Now the condition \eqref{bisection cover} for $\CD$ shows that  $D_+\subseteq E_+=E_{+, 1}E_{+, 2}$, which is the unique element in $\CE_+$ containing $D_+$ and $\overline{E_-}$, and is disjoint form $\overline{E'_-}$ for any other $E'_-$ in $\CE_-$. 
 Choose $\gamma\in D_+\cap \overline{E_{-,1}}D_{+, 0}\overline{E_{-,2}}$. Then $\gamma=xy$ for $x\in \overline{E_{-, 1}}$  and $y\in \overline{E_{-, 2}}$ with $s(x)=r(y)\in D_{+, 0}$.
 Then because $D_{+, 0}\in (\CD^{-KL_2}_+)\unitSpace$, by definition and $\GU\subseteq K$, one has 
 $x\in \overline{E_{-, 1}}\cdot D_{+, 0}
 \subseteq \bigcup \CD$ and therefore, there exists a $D_{+, 1}\in \CD$ such that $x\in D_{+, 1}$. This first implies that $s(x)\in s(D_{+, 1})\cap D_{+, 0}$ and $r(x)=r(\gamma)\in r(D_+)\cap r(D_{+, 1})$. Therefore, one has $s(D_{+, 1})=D_{+, 0}$ and $r(D_{+, 1})=r(D_+)$ as $\CD$ is a castle. Moreover, since $x\in D_{+, 1}\cap \overline{E_{-,1}}$, The condition \eqref{bisection cover} for $\CD$ shows that $D_{+, 1}\subseteq E_{+, 1}$, which further entails that $D_{+, 1}=E_{+, 1}D_{+, 0}$. Then for $y$, the similar argument shows that $D_{+, 2}=D_{+ ,0}E_{+, 2}$ and $s(D_{+, 2})=s(D_+)$. Consequently, one actually has 
 \[D_+=D_{+,1}D_{+, 2}=E_{+ ,1}D_{+, 0}E_{+, 2}.\]
 This shows that $\CA\subseteq \CB$.
 
Now, for the converse, let $E_{+,1}D_{+, 0}E_{+ ,2}\in \CB$ such that  $D_{+, 0}\in (\CD^{-KL_2}_+)\unitSpace$ and $(E_{+, 1}, E_{-,1}), (E_{+, 2}, E_{-, 2})\in \CE$ satisfy $\overline{E_{-, 1}}D_{+, 0}\overline{E_{-, 2}}\neq \emptyset$. Note that $D_{+, 0}\in\CD^{-KL_2}_+$ implies that $\overline{E_{-, 1}}D_{+, 0}\overline{E_{-, 2}}\subseteq \bigcup \CD$ and thus there exits a $D_+\in \CD$ such that the intersection $D_+\cap \overline{E_{-, 1}}D_{+, 0}\overline{E_{-, 2}}\neq \emptyset$. Then the same argument as above shows that $E_{+, 1}D_{+, 0}E_{+, 2}=D_+$. In addition, since $D_{+, 0}\in\CD^{-KL_2}_+$, one has $D_+=E_{+, 1}D_{+, 0}E_{+, 2}\in \CD_+^{-K}$, which shows that $\CB\subseteq \CA$ and thus $\CA=\CB$.

	We now show that $\CF$ is a fortified castle with the desired properties. We establish \ref{prop:AE-nesting-minimal:nesting} first. Let $D_+\in \CF_+$. Then by definition, $r(D_+)\in \CA\unitSpace$ and therefore, $r(D_+)\subseteq E_+$ for some $E_+\in \CE_+\unitSpace$ since $\CA=\CB$. This implies that $E_+D_+=D_+$ and thus $\CF_+\subseteq \CE_+\CF_+\setminus\{\emptyset\}$. For the converse, let $E_+\in \CE_+$ and $D_+\in \CF_+$ such that $E_+D_+\neq \emptyset$. Then because $r(D_+)\in \CA\unitSpace$, there exists a $D'_+=E_{+, 1}D_{+, 0}E_{2, +}\in \CA=\CB$ for some
    $D_{+, 0}\in (\CD^{-KL_2}_+)^{(0)}$ and $(E_{+, 1}, E_{-,1})$, $(E_{+, 2}, E_{-, 2})\in \CE$ with  $D_+\cap E_{-, 1}D_{+, 0}E_{-, 2}\neq \emptyset$ such that $r(D_+)=r(D'_+)\subseteq r(E_{+, 1})$. Since $\CE_+$ is a castle, one necessarily has $r(E_{+, 1})= s(E_+)$. Denote by $E'_{+, 1}=E_+E_{+, 1}$ and   this then implies 
    \[E_+D'_+=E_+E_{+, 1}D_{+, 0}E_{+, 2}=E'_{+, 1}D_{+ ,0}E_{+, 2}\in \CB=\CA\subseteq \CF_+\]
    by definition. Therefore, one has
    \[E_+D_+=E_+D'_+\cdot D'^{-1}_+\cdot D_+\in \CF_+\]
    and thus $\CE_+\CF_+\setminus \{\emptyset\}\subseteq\CF_+$. Consequently, $\CF_+$ is nested in $\CE_+$.

    Then, \ref{prop:AE-nesting-minimal:extend} and \ref{prop:AE-nesting-minimal:refine} for $\CF$ and $\CD$ are clear by the definition of $\CF$ and $\CU$. To establish \ref{prop:AE-nesting-minimal:remainder},  define 
    \[\CD_-|_{\CE}=\{D_-\in \CD_-: \overline{E_-}D_-\overline{E'_-}\neq \emptyset \text{ for some } (E_+, E_-), (E'_+, E'_-)\in \CE\}.\]
    We claim that the family $\CY=\{D_+\in \CD_+\unitSpace: \CD_-\in \CD^{-KL_2}_-\cap (\CD_-|_{\CE})\}\subseteq \CA\unitSpace$. Indeed, let $D_+\in \CD_+\unitSpace$ such that $D_-\in D^{-KL_2}_-\cap (\CD_-|_{\CE})$. There exists $(E_+, E_-), (E'_+, E'_-)\in \CE$ such that $\overline{E_-}D_-\overline{E'_-}\neq \emptyset$, which further implies that $\overline{E_-}D_+\overline{E'_-}\neq \emptyset$. Then note that by definition $D_+\in (\CD^{-KL_2}_+)\unitSpace$, one has $D_+\cap s(\overline{E_-})D_+r(\overline{E'_-})\neq \emptyset$. Using $\GU\subseteq L_2$, this further shows that $D_+\in \CA\unitSpace$ and thus has established the claim $\CY\subseteq \CA\unitSpace=\CF_+\unitSpace$.
    Then observe that by definition the family
    \[\CF_-\unitSpace=\{s(D_-), r(D_-): D_+\in \CA\}\]
     and thus
    \[(\CD^{-KL_2}_-)\unitSpace\cap (\CD_-|_{\CE})=\{D_-: D_+\in \CY\} \subseteq\CF_-\unitSpace.\]
    One the other hand, observe that
    \[(\CD^{-KL_2}_-)\unitSpace\cap (\CD_-|_{\CE})=\{D_-\in (\CD^{-KL_2}_-)\unitSpace: D_-\cap \overline{E_-}\neq \emptyset \text{ for some } E_-\in \CE_-\unitSpace\}.\]
 Write $\CY_-\unitSpace=(\CD^{-KL_2}_-)\unitSpace\cap (\CD_-|_{\CE})$ for simplicity. Then, by the condition \eqref{small reminder nesting0} for the castle $\CD$, one has 
    \[(\overline{\bigcup \CE_-\unitSpace})\setminus (\bigcup \CF_-\unitSpace)\subseteq (\overline{\bigcup \CE_-\unitSpace})  \setminus \bigcup \CY_-\unitSpace\subseteq \GU\setminus \bigcup (\CD^{-KL_2}_-)\unitSpace\prec_\CG O_0,\]
    which establishes the \ref{prop:AE-nesting-minimal:remainder}.
\end{proof}

Although we will not use this fact, we point out that Proposition~\ref{prop:AE-nesting} is formulated in a way that allows it to be applied recursively to obtain a sequence $\CF_0, \CF_1, \CF_2, \ldots$ of fortified castles with $\CF_0 = \left\{\left(\CG\unitSpace, \CG\unitSpace\right)\right\}$ that get increasingly ``fine'' and ``extendable'', with each step incurring a ``small'' remainder. Hence roughly speaking, a ``large part'' of $\CG$ may be approximated by a nested sequence of open castles. 

We also have a further strengthening of Proposition~\ref{prop:AE-nesting} (see Proposition~\ref{prop:AE-nesting-minimal} below), where the additional conclusion may be seen as a statement about recurrence, to be made precise in Definition~\ref{defn:castle-multiplicity}. 
This kind of recurrence will be important for us in Section~\ref{sec:Z-stable}, since it gives rise to nestings with ``large commutants'', which ultimately leads to large matrix algebras mapping into the central sequence algebra of a groupoid $C^*$-algebra.

\begin{defn} \label{defn:castle-multiplicity}
	Let $\CC, \CD$ be open castles with $\CC$ \emph{nested} in $\CD$. The \emph{minimal multiplicity} of the nesting is defined to be the minimal height of $\CommCD{\CC}{\CD}$ in the sense of Definition~\ref{defn:castle}. 
	
	By Remark~\ref{rem:nesting-comm-red-basic}\eqref{rem:nesting-comm-red-basic:Red-Res}, \eqref{rem:nesting-comm-red-basic:castles-G} and~\eqref{rem:nesting-comm-red-basic:Res-Comm}, this number is equal to 
	\[
	\min_{[D] \in \CD\unitSpace / \sim_{\CD}} \left( \min_{\CT \in \CC / \sim_{\CC} \colon D \cap ( \bigcup \CT) \not= \varnothing} \left| \left\{ C \in \CT \unitSpace \colon C \subseteq D \right\}  \right| \right) \; .
	\]
\end{defn}

\begin{lem}\label{lem: recurrence imply high multiplicity}
  Let $\CG$ be a minimal groupoid with a compact unit space.  Let $\CU$ be a finite disjoint collection of non-empty open sets in $\GU$ and $N>0$. Then there exists a family $\CW$ of precompact open bisections such that 
  \begin{enumerate}[label=(\roman*)]
  \item\label{lem: recurrence imply high multiplicity-range recurrence} the family $\{r(W): W\in \CW\}$ refines $\CU$;     
  \item\label{lem: recurrence imply high multiplicity-source fine cover} for each $W\in \CW$, there exists an open set $W'$ with $\overline{W'}\subseteq W$ and there exists an open cover $\CO$ of $\GU$ that refines $\{s(W'):  W\in \CW\}$;
      \item\label{lem: recurrence imply high multiplicity-multiplicity}for any $U\in \CU$, one has $|\{r(W)\subseteq U: W\in \CW\}|\geq N$.
  \end{enumerate}
\end{lem}
\begin{proof}
    For each $U\in \CU$, since $\GU$ is perfect by Remark \ref{rem: minimal-perfect}, choose $N$ disjoint open sets $V_{U, 1}, \cdots, V_{U, N}$ in $U$. Then for each $u\in \GU$, by minimality, there exists $\{x_{U, i}: 1\leq i\leq N, U\in \CU\}$ such that all $s(x_{U, i})=u$ and $r(x_{U, i})\in V_{U, i}$. This allows to choose precompact open bisections $ W'_{U, i}, W_{U, i}$ such that $\overline{W'_{U, i}}\subseteq W_{U, i}$ satisfying that $u\in \bigcap_{U\in \CU, 1\leq i\leq N}s(W'_{U, i}):= O_u$ and $r(W_{U, i})\subseteq V_{U, i}$. The collection of these $O_u$ forms an open cover of $\GU$ and by compactness of $\GU$, one obtains a finite open cover $\CO$ of $\GU$ by construction. There exists a finite family $\CW$ of precompact open bisections such that $\CO$ refines $\{s(W): W\in \CW\}$ and $\{r(W): W\in \CW\}$ refines $\CU$. Moreover, by our construction, for any $U\in \CU$ and $V_{U, i}$, there exists at least one $W\in \CW$ such that $r(W)\subseteq V_{U, i}$.  Therefore, one has $|\{r(W)\subseteq U: W\in \CW\}|\geq N$.
\end{proof}

Then we have the following strengthened version of Proposition~\ref{prop:AE-nesting} whose proof is similar to the Proposition~\ref{prop:AE-nesting} with additional adjustment using Lemma \ref{lem: recurrence imply high multiplicity}. The main idea is that using the recurrence of Lemma \ref{lem: recurrence imply high multiplicity} from the minimality and starting with a fortified castle $\CD$ with a stronger shrinkability, one may obtain a subset of $\CD_-^{-KL_2}\cap (\CD_-|_{\CE})$ in Proposition ~\ref{prop:AE-nesting} to play the exact role of the $\CD_-^{-KL_2}\cap (\CD_-|_{\CE})$ in Proposition ~\ref{prop:AE-nesting} and the same construction yields the desired $\CF$.

\begin{prop} \label{prop:AE-nesting-minimal}
	A minimal \'{e}tale groupoid $\CG$ with a compact unit space is {almost elementary}
	if and only if 
	for any $K$, $\CV$, $O_0$ and $\CE$ as in Proposition~\ref{prop:AE-nesting}, as well as any natural number $N>1$, there are a fortified castle $\CF$ and an open castle $\CC$ satisfying Proposition~\ref{prop:AE-nesting}\ref{prop:AE-nesting-minimal:extend}, \ref{prop:AE-nesting-minimal:refine}, and~\ref{prop:AE-nesting-minimal:remainder} together with: 
	\begin{enumerate}[label=(\roman*')]  
		\item\label{prop:AE-nesting-minimal:nesting-recurrent} $\CF_+$ is faithfully nested in $\CE_-$ with a minimal multiplicity at least $N$. 
	\end{enumerate}
\end{prop}
\begin{proof}
Indeed, let the compact set $K$, the open cover $\CV$ of $\GU$, the number $N>0$, the open set $O_0$, and the fortified castle $\CE$ be given. We may also assume $\CE_+$ is precompact and $\GU\subseteq K$. Then for $\CE_-$, a finite disjoint family of non-empty open sets in $\GU$ and $N>1$, Lemma \ref{lem: recurrence imply high multiplicity} implies that there exists  a collection $\CW$ of precompact open bisections satisfying conditions \begin{enumerate}
  \item\label{range recurrence} the family $\{r(W): W\in \CW\}$ refines $\CE_-$;     
  \item\label{source fine cover} for each $W\in \CW$, there exists an open $W'$ with $\overline{W'}\subseteq W$ and there exists an open cover $\CO$ of $\GU$ that refines $\{s(W'):  W\in \CW\}$;
      \item\label{multiplicity} one has $|\{r(W)\subseteq E_-: W\in \CW, O\subseteq s(W'), \}|\geq N$ for any $E_-\in \CE_-$ and $O\in \CO$.
     \end{enumerate}
 Denote by $\CW'$ the family consisting of all these $W'$. Then define a new compact set $K_0=\overline{\bigcup \CW'}\cup \GU$ and compact sets  $L_1=\overline{\bigcup \CE_-}\cup \GU$ and $L_2=\overline{\bigcup \CE_+}\cup \GU$ as in Proposition~\ref{prop:AE-nesting}. Note that $\CW$ is an open cover of $K_0$ and $\CE_+$ is an open cover of $\overline{\bigcup \CE_-}$. Note that $\CV\vee \CO$ is an open cover of $\GU$ and the condition \eqref{source fine cover} implies $s(\CW')$ is also an open cover of $\GU$, i.e., $s(K_0)=\GU$.  Then using $\GU$ is clopen in $\CG$, one chooses an open cover $\CU$
of $L_1\cup K_0$ that is fine enough such that the following hold.
\begin{enumerate}
    \item For any $W\in \CW$, if $U\cap \overline{W'}\neq \emptyset$ then $U\subseteq W$.
    \item For any $(E_+, E_-)\in \CE$, if $U\cap \overline{E_-}\neq \emptyset$ then $U\subseteq E_+$.
    \item $\{U\in \CW: U\subseteq \GU\}$ refines $\CV\vee \CO$.
\end{enumerate}

Then for compact sets $KL_2K_0$, $L_1\cup K_0$ and the open cover $\CU$, apply Lemma~\ref{lem:AE-fortified} to  obtain a fortified castle $\CD$ satisfying: 
	\begin{enumerate}[label=(\roman*)]
		\item\label{bisection cover1} for any $D_+ \in \CD_+$ with $D_+ \cap (L_1\cup K_0) \not= \varnothing$, there is a $U\in \CU$ such that $D_+ \subseteq U$; 
		\item 
		$\bigcup \CD_+$ is precompact in $\CG$;  
		\item\label{small reminder nesting}
		$\CG^{(0)}\setminus\ \left( \
        \bigcup \left( \CD^{-KL_2K_0}_- \right)^{(0)} \right) \prec_{\CG} O_0$.   
	\end{enumerate}
By the definition of the cover $\CU$ of $L_1\cup K_0$,  for any $D_+\in (\CD^{-KL_2K_0}_+)\unitSpace$ one has $D_+\subseteq O$ for some $O\in \CO$. Then by \eqref{multiplicity}, for any $E_-\in \CE_-$ there are at least $N$ open bisections $W\in \CW'$ such that $O\subseteq s(W')$ and $r(W')\subseteq r(W)\subseteq E_-$. Now, let $E_-\in \CE_-\unitSpace$ and $D_{+,0}\in (\CD^{-KL_2K_0}_+)\unitSpace$ be fixed. For $W'\in \CW'$ with $D_{+,0}\subseteq s(W')$ and $r(W)\subseteq E_-$, one has $W'D_{+,0}\subseteq \bigcup \CD_+$ because $W'\subseteq K$. Therefore there exists a $D_{+}\in \CD_+$ such that $D_{+}\cap W'D_{+,0}\neq \emptyset$, which implies that $D_{+}\subseteq W$ by the condition \ref{bisection cover1} for $\CD_+$ above. This further shows that $D_{+}=WD_{+, 0}$ and thus one has $D_{+}, r(D_{+})\in \CD^{-KL_2}$. Moreover, note that $r(D_{+})=r(WD_{+,0})\subseteq E_-$. This implies that the family 
\[\CX_{E_-, D_{+, 0}}:=\{D_+\in (\CD^{-KL_2}_+): r(D_+)\subseteq E_-, s(D_+)=D_{+, 0}\}\]
satisfies $|r(\CX_{E_-, D_{+, 0}})|\geq N$.
Note all $D_+\in \CX_{E_-, D_{+, 0}}$ come from the same tower in $\CD_+$. Define 
\[\CX=\bigcup_{D_{+,0}\in (\CD^{-KL_2K_0}_+)\unitSpace, E_-\in \CE_-\unitSpace}\CX_{E_-, D_{+, 0}}.\]
Then we follow exactly the approach in Proposition \ref{prop:AE-nesting} to define the family $\CA$ and $\CA\unitSpace$ using the current $\CX$ by 
 \begin{align*}
        \CA:=\{&D_+\in \CD^{-K}_+: D_+\cap \overline{E_{-}}D'_{+}\overline{E'_-}\neq \emptyset 
    \text{ for some }
    D'_{+}\in r(\CX)
    \text{ and }\\
    &(E_+, E_-), (E'_+, E'_-)\in \CE\}.  
    \end{align*}
from which we obtain a fortified castle $\CF$, as a fortified subcastle of $\CD$, such that 
\[\CF_+\unitSpace=\CA\unitSpace:=\{s(D_+), r(D_+): D_+\in \CA\}.\]
Then the same argument proving Proposition \ref{prop:AE-nesting} shows that actually
\begin{align*}
\CA=\{&D_+\in \CD^{-K}_+: D_+= E_{-}D'_{+}E'_- \text{ for some }
    D'_{+}\in r(\CX)
    ,(E_+, E_-),(E'_+, E'_-)\in \CE_+ \\
    &\text{such that } D'_+\subseteq s(E_-)\text{ and } D'_+\subseteq r(E'_-)\}.
\end{align*}

We now establish \ref{prop:AE-nesting-minimal:nesting-recurrent}. First, $\CF_+$ is nested in $\CE_-$ by the same argument for proving Proposition \ref{prop:AE-nesting}\ref{prop:AE-nesting-minimal:nesting}. Indeed, let $D_+\in \CF_+$ and then $r(D_+)\in \CA\unitSpace$ by definition and thus $ r(D_+)\subseteq E_-$ for some $E_-\in \CE_-\unitSpace$. This implies that $E_-D_+=D_+$ and thus $\CF_+\subseteq \CE_-\CF_+\setminus \{\emptyset\}$. For the converse direction, let $E_-D_+\neq \emptyset$ for some $(E_+,E_-)\in \CE$ and $D_+\in \CF_+$. This implies $r(D_+)\in \CA\unitSpace$ and $r(D_+)\subseteq s(E_-)$ by the construction. Then there exists a $D_{+, 1}\in \CA$ satisfying $D_{+, 1}=E_{-, 1}D'_+E_{-, 2}$ for some $D'_+\in r(\CX)$ and $E_{-, 1}, E_{-, 2}\in \CE_+$ such that $r(D_{+ ,1})=r(D_+)$. This implies that 
\[E_-D_{+, 1}=E_-E_{-, 1}D'_{+}E_{-, 2}\in \CA\subseteq \CF\]
and thus
\[E_-D_+=E_-D_{+, 1}\cdot D^{-1}_{+, 1}\cdot D_+\in \CF,\]
which shows that $\CE_-\CF_+\setminus \{\emptyset\}=\CF_+$. In addition, for any $E_-\in \CE_-$, recall that $|r(\CX_{E_-, D_{+, 0}})|\geq N$ and
\[r(\CX_{E_-, D_{+, 0}})\subseteq r(\CX)\subseteq\CA\unitSpace=\CF_+\unitSpace\]
holds for any $D_{+, 0}\in (\CD^{-KL_2K_0}_+)\unitSpace$.  This implies that $E_-D_+\neq \emptyset$ for any $D_+\in \CX_{E_-, D_{+, 0}}$ and thus the nesting is faithful.

To address the multiplicity part, we first claim that 
\begin{equation}\label{key formula}
   (\CD^{-KL_2K_0}_+)\unitSpace\subseteq r(\CX)\subseteq\CF_+\unitSpace\subseteq \CD_+\unitSpace. 
\end{equation}
It suffices to show $r(\CX)\subseteq\CF_+\unitSpace$ because the other two containment are obtained by definition and the fact $\GU\subseteq K_0$.  Let $D_+\in r(\CX)$, which implies that there exists a $E_-\in \CE_-\unitSpace$ such that $D_+\subseteq E_-$ by definition. Then one has $D_+\cap E_-D_+E_-=D_+\neq \emptyset$. Note additionally that $\GU\subseteq L_2$ implies $\CD_+^{-KL_2}\subseteq \CD_+^{-K}$. This implies that $r(\CX)\subseteq \CA\unitSpace=\CF_+\unitSpace$.

Now, let $E_-\in \CE_-\unitSpace$ and $\CT$ be a tower in $\CF_+$ such that $E_-\cap (\bigcup \CT\unitSpace)\neq \emptyset$. Using \eqref{key formula} and $\CD^{-KL_2K_0}_+$ intersects all possible tower in $\CD$, one has $\CT\unitSpace\cap (\CD^{-KL_2K_0}_+)\unitSpace\neq \emptyset$, which allows to choose one $D_{+, 0}\in (\CD^{-KL_2K_0}_+)\unitSpace\cap \CT\unitSpace$. Then by definition, one has 
\[r(\CX_{E_-, D_{+ ,0}})\subseteq E_-\cap \bigcup\CT\unitSpace.\]
and thus $|E_-\cap \CT\unitSpace|\geq N$. Consequently, by Definition \ref{defn:castle-multiplicity} One finally has that $\CF_+$ is nested in $\CE_-$ with a minimal multiplicity at least $N$. Therefore, we have established \ref{prop:AE-nesting-minimal:nesting-recurrent}.

Then, \ref{prop:AE-nesting-minimal:extend} and \ref{prop:AE-nesting-minimal:refine} for $\CF$ and $\CD$ follows from the definition of $\CF$ and $\CU$. We finally establish \ref{prop:AE-nesting-minimal:remainder}. First note that \ref{small reminder nesting} implies that $\mu(\GU\setminus \bigcup (\CD^{-KL_2K_0}_-)\unitSpace)\leq \mu(O_0)$ for any $\mu\in M(\CG)$ by Remark \ref{rem:subequivalence-basic}\eqref{rem:subequivalence-basic:measure-auxiliary}. Define
\[\CX_-=\{D_-\in \CD_-: D_+\in \CX\}.\]
Then from \eqref{key formula}, one also has
\begin{equation}\label{key formula2}
     (\CD^{-KL_2K_0}_-)\unitSpace\subseteq r(\CX_-)\subseteq\CF_-\unitSpace.
\end{equation}

Moreover, our construction implies that
\[\mu(\bigcup r(\CX_-))\geq N|\CE_-|\mu(\bigcup  (\CD^{-KL_2K_0}_-)^{(0)})>\mu(\bigcup  (\CD^{-KL_2K_0}_-)^{(0)})\]
for any $\mu\in M(\CG)$.  Using \eqref{key formula2}, this further entails
\[\mu(\overline{\bigcup \CE_-\unitSpace}\setminus \bigcup \CF_-\unitSpace)\leq \mu(\GU\setminus \bigcup r(\CX_-))<\mu(\GU\setminus \bigcup (\CD^{-KL_2K_0}_-)\unitSpace)\leq \mu(O_0)\]
for any $\mu\in M(\CG)$.
Therefore, Theorem \ref{thm:AE-DSC} shows that $(\overline{\bigcup \CE_-\unitSpace})\setminus (\bigcup \CF_-\unitSpace)\prec_\CG O_0$, which establishes  \ref{prop:AE-nesting-minimal:remainder}.
\end{proof}

\section{Tracial $\protect\CZ$-stability} \label{sec:Z-stable}

In this section, we investigate structural properties of the reduced groupoid
$C^{*}$-algebra $C_{\operatorname{r}}^{*}(\CG)$ associated to a minimal almost elementary \'{e}tale groupoid
$\CG$. In particular, we will show $C_{\mathrm{r}}^{*}(\CG)$ is tracially
$\CZ$-stable in the sense of Hirshberg and Orovitz \cite{HirshbergOrovitz2013Tracially}. 
We start by recalling the notion of Cuntz subequivalence, which is a source of motivation for the groupoid subequivalence relation $\precsim_{\CG}$ we gave in Definition~\ref{defn:subequivalence}. 

\begin{defn}[{\cite{Cuntz1978Dimension}}]
	Let $A$ be a $C^{*}$-algebra. 
	For $a,b\in A$, we say $a$ is
	\emph{Cuntz subequivalent} to $b$, and write $a\precsim b$,
	if there is a sequence $\{x_{n}\in A:n\in\N\}$ such that $\lim_{n\to\infty}\|a-x_{n}bx_{n}^{*}\|=0$.
	We say $a$ is
	\emph{Cuntz subequivalent} to $b$, and write $a\sim b$, if $a\precsim b$ and $b\precsim a$. 
\end{defn}

We will use the following elementary fact.

\begin{lem}\label{lem: Cuntz subequi}
   Let $A$ be a $C^*$-algebra and $x_1, \dots, x_n$ are elements in $A$. Then $\left(\sum_{i=1}^nx^*_i\right)\left(\sum_{i=1}^nx_i\right)\precsim \sum_{i=1}^n x^*_ix_i$ holds. 
\end{lem}

\begin{proof}
    We compute: 
    \begin{align*}
        & \left(\sum_{i=1}^nx_i \right)^* \, \left(\sum_{i=1}^nx_i \right) \\
        \leq&\ \sum_{F \subseteq \{1,\ldots,n\}} \left(\sum_{i=1}^n (-1)^{\chi_F(i)} \, x_i\right)^* \, \left(\sum_{i=1}^n  (-1)^{\chi_F(i)} \, x_i\right)  \\ 
        =&\ \sum_{i,j=1}^n \left( \sum_{F \subseteq \{1,\ldots,n\}} (-1)^{\chi_F(i) + \chi_F(j)} \right) x_i^* x_j \\
        =&\ \left( \sum_{i=1}^n 2^n x_i^* x_i \right) \\
        &+ \sum_{_{\quad i \not= j}^{i,j \in \{1,\ldots, n\}}} \left( \sum_{F \subseteq \{1,\ldots,n\} \setminus \{i\} } (-1)^{\chi_F(i) + \chi_F(j)} + (-1)^{\chi_{F \sqcup \{i\}}(i) + \chi_{F \sqcup \{i\}}(j)} \right) x_i^* x_j \\
        =&\ 2^n \sum_{i=1}^n x_i^* x_i \; ,
    \end{align*}
    which implies $\left(\sum_{i=1}^nx_i \right)^* \, \left(\sum_{i=1}^nx_i \right) \precsim 2^n \sum_{i=1}^n x_i^* x_i \sim \sum_{i=1}^n x_i^* x_i$. 
\end{proof}

In the following, we write $\osupp{f}$ for the \emph{open} support of a complex function $f$ on a topological space $X$, that is
\[
	\osupp{f} := \left\{ x \in X \colon f(x) \not= 0 \right\} \; .
\] 
We record an explicit relationship between Cuntz subequivalence and groupoid subequivalence.

\begin{lem}\label{lem:groupoid subequi to Cuntz subequi}
	Let $\CG$ be an \'{e}tale groupoid and let $f , g \in C_0(\CG\unitSpace)$. 
	If $\osupp{f} \precsim_{\CG} \osupp{g}$, then $f \precsim g$ in $C^{*}(\CG)$. 
\end{lem}
\begin{proof}
It suffices to show $(f-\varepsilon)_+\precsim g$ for any fixed $\varepsilon>0$. Define $K:=\supp(f-\varepsilon)_+$, which is a compact set contained in $O:=\supp^o(f)$. Then by Definition  \ref{defn:subequivalence}, there exists an open $s$-section $S$ such that $K\subseteq r(S)$ and $s(S)\subseteq O$. Using the compactness of $K$, we may choose a finite collection $\{U_i: i\in I\}$ of precompact open bisections $U_i$ satisfying 
\[K\subseteq \bigcup_{i\in I}r(U_i)  \quad \text{ and } \quad   \bigcup_{i\in I}\overline{U_i}\subseteq S  \; .\]
Now choose $h_i\in C_c(U_i)$ for each $i \in I$, such that the collection $\{h_ih^*_i: i\in I\}$ forms a partition of unity for the covering $K\subseteq \bigcup_{i\in I}r(U_i)$. Define $h=\sum_{i\in I}h_i$, which is supported on $\bigcup_{i\in I}U_i\subseteq S$. Then $hh^*$ is supported on $SS^{-1}$. Since $S$ is an $s$-section, it follows that $SS^{-1}=r(S)$ and $h h ^* = \sum_{i \in I} h_i h_i^*$, whence $hh^* \equiv 1$ on $K$. Therefore, one has 
\[(f-\varepsilon)_+\precsim hh^*\sim h^*h=(\sum_{i\in I}h^*_i)(\sum_{i\in I}h_i)\precsim \sum_{i\in I}h^*_ih_i\precsim g\]
by Lemma \ref{lem: Cuntz subequi} and the fact $\supp^o(\sum_{i\in I}h^*_ih_i)\subseteq \bigcup_{i\in I} s(U_i)\subseteq s(S)$. 
\end{proof}

\begin{defn}[{\cite[Definition 2.1]{HirshbergOrovitz2013Tracially}}]
	\label{8.1} \label{defn:tracial-Z-stable}
	A unital $C^{*}$-algebra $A$
	is said to be \textit{tracially} $\CZ$\textit{-stable} if $A\neq\C$ and
	for any finite set $F\subseteq A$, $\varepsilon>0$, any nonzero positive
	element $a\in A_{+}$ and $n\in\N$, there is an order zero c.p.c.\ map
	$\psi:M_{n}(\C)\to A$ such that the following hold: 
	\begin{enumerate}
		\item $1_{A}-\psi(1_{M_{n}(\C)})\precsim a$; 
		\item $\|[\psi(x),y]\|<\varepsilon$ for any $x\in M_{n}(\C)$ with $\|x\|=1$ and any $y\in F$. 
		 
	\end{enumerate}
\end{defn}
In addition, it was proved in \cite{HirshbergOrovitz2013Tracially}
that tracial $\CZ$-stability is equivalent to $\CZ$-stability
when the $C^{*}$-algebra $A$ under consideration is
unital simple separable and nuclear. 

We will verify this property for groupoid $C^*$-algebras under the assumption of almost elementariness. 
We begin with the following lemma, which is a groupoid version of \cite[Lemma 7.9]{Phillips-large}. 
It will help us get more control over the nonzero positive element $a\in A_{+}$ in Definition~\ref{defn:tracial-Z-stable}. 

\begin{lem}\label{lem: small function}
	 Let $\CG$ be an effective \'{e}tale
	groupoid on a compact sapce. For any non-zero element $a\in C_{\operatorname{r}}^{*}(\CG)_{+}$
	there is a non-zero function $g\in C(\GU)_{+}$ such that $g\precsim a$. If $\CG$ is assumed ample, then $g$ can be chosen to be supported on a compact open set in $\GU$.\end{lem}
\begin{proof}
	Let $a\in C_{\operatorname{r}}^{*}(\CG)_{+}\setminus\{0\}$. Without loss of generality,
	one may assume $\|a\|=1$. Let $\varepsilon = (1/6)\|E(a)\|$ where
	$E:C_{\operatorname{r}}^{*}(\CG)\to C(\GU)$ is the canonical faithful expectation.
	Then there is an $h\in C_{c}(\CG)_{+}$ such that $\|a^{1/2}-h\|<\varepsilon$
	and $\|h\|\leq1$. Then one has $\|h*h^{*}-a\|<2\varepsilon$ and $\|h^{*}*h-a\|<2\varepsilon$.
	Then for $h^{*}*h\in C_{c}(\CG)_{+}$, Lemma 4.2.5 in \cite{Sims-groupoids}
	shows that there is a function $f\in C(\GU)_{+}$ such that $\|f\|=1$,
	$f*h^{*}*h*f=fE(h^{*}*h)f$ and $\|f*h^{*}*h*f\|\geq\|E(h^{*}*h)\|-\varepsilon$.
	Then one has 
	\[
	\|f*h^{*}*h*f\| \geq \|E(h^{*}*h)\|-\varepsilon\geq\|E(a)\|-3\varepsilon = 3\varepsilon.
	\]
	Then define $g=(f*h^{*}*h*f-2\varepsilon)_{+}=(fE(h^{*}*h)f-2\varepsilon)_{+}\in C(\GU)_{+}\setminus\{0\}$.
	Then using Lemma 1.6 and 1.7 in \cite{Phillips-large}, one has 
	\[
	g\sim(h*f^{2}*h^{*}-2\varepsilon)_{+}\precsim(h*h^{*}-2\varepsilon)_{+}\precsim a
	\]
	as desired. If $\CG$ is ample then choose a compact open set $V\subseteq \supp(g)$ and then $1_V\precsim g\precsim a$. Then $1_V$ is what we want.
\end{proof}

As hinted in Proposition~\ref{prop:castle-order-zero}, we will construct order zero maps from open castles. Let us first see how the combination of fineness and extendability (see Definition~\ref{defn:castle-extendable}) of an open castle allows one to obtain a sort of local approximation of functions on $\CG$ supported not too far away from $\CG\unitSpace$. 

\begin{lem} \label{lem:castle-local-approximation}
	Let $\CG$ be an \'{e}tale groupoid. 
	Let $\CC$ be an open castle in $\CG$, let $f$ be a positive contraction in $C_0 \left(\bigcup \CC \right) \cap C_{\operatorname{c}}(\CG)$ that is $\left(\bigcup \CC\right)$-invariant, and let $\psi_f \colon C^*(\CC) \to  C^*(\CG)$ be the c.p.c.\ order zero map given by Proposition~\ref{prop:castle-order-zero}. 
	Choose $x_C \in C$ for each $C \in \CC$, and define a map $\varphi \colon C_{c}(\CG) \to C^*(\CC)$ by $g \mapsto \left( g(x_C) \right)_{C \in \CC}$.  
    Let $K$ be a compact bisection in $\CG$ contained in an open bisection $B$. Suppose $\CC$ refines the cover $\CV=\{B, \CG\setminus K\}$. Then for any $g \in C_{c}(\CG)$ supported in $K$ and any contraction $e \in C\left(\left(\CC^{-K}\right)\unitSpace\right) \subseteq C^*(\CC)$,  one has
    \[\max \left\{ \left\| \psi_f(e) g - \psi_f \left( e \varphi (g) \right) \right\|_r  ,  \left\| g \psi_f(e) - \psi_f \left( \varphi (g) e \right) \right\|_r \right\}\leq \|e\|_\infty\cdot \max_{C\in \CC}\operatorname{osc}_C(g).\]
    Here, the notion $\operatorname{osc}_{C} (\cdot)$ denotes the  oscillation of the function on $C$, i.e., $\operatorname{osc}_C(f):= \sup_{x, y \in C} |f(x) - f(y)|$ for a function $f\in C_c(\CG)$.
    
    More generally, if $K$ is a compact set that can be covered by $n$ open bisections $\{B_1, \dots, B_n\}$ in $\CG$ and let $\{h_1,\dots, h_n\}\subseteq C_c(\CG)$ be a partition of unity subordinate to $\{B_1, \dots, B_n\}$. Suppose $\CC$ refines the cover $\CV=\{B_1, \dots, B_n, \CG\setminus K\}$ such that if $C\in \CC$ satisfies $C\cap K\cap B_i\neq \emptyset$ then $C\subseteq B_i$ holds for any $1\leq i\leq n$. Then for any $g \in C_{c}(\CG)$ supported in $K$ 
	and any contraction $e \in C\left(\left(\CC^{-K}\right)\unitSpace\right) \subseteq C^*(\CC)$, 
	we have 
	\begin{align*}
	&\max \left\{ \left\| \psi_f(e) g - \psi_f \left( e \varphi (g) \right) \right\|_r  ,  \left\| g \psi_f(e) - \psi_f \left( \varphi (g) e \right) \right\|_r \right\}\\
 &\leq\|e\|_\infty\cdot\sum_{i=1}^n(\max_{C\in \CC}(\operatorname{osc}_{C}(g)+\|g\|_\infty\operatorname{osc}_C(h_i)).
	\end{align*}

\end{lem}
\begin{proof}
We proceed to prove the following more general fact to establish the two claims together. Let $K$ be a compact set. First, suppose $P$ is a precompact subset of $K$ and $P$ is contained in an open bisection $B$ (for the first claim, take $P=K$). Then suppose if $C\in \CC$ satisfies $C\cap P\neq \emptyset$ then $C\subseteq B$ holds. In addition, suppose $g$ is supported on $P$. 

Then by the  assumption, note that $P$ is also a bisection in $\CG$. Now, let $e\in C((\CC^{-K})\unitSpace)$ be a contraction, which means there exists a subfamily $\CE\subseteq (\CC^{-K})\unitSpace$ such that  $e$ can be written as
\[e=\sum_{C\in \CE}a_C\chi_C\]
for some $a_C\in \C$ for each $C\in \CE$. Then a direct computation shows 
\[\psi_f(e*\varphi(g))=\sum_{C\in \CC, r(C)\in \CE}f|_C\cdot a_{r(C)}g(x_C)=\sum_{C\in \CC, r(C)\in \CE, C\cap P\neq \emptyset}f|_C\cdot a_{r(C)}g(x_C),\]
where the last equality is due to the fact that $g$ is supported on  $P$. In addition, one has 
\[\psi_f(e)*g=\sum_{C\in \CE}(a_Cf|_C)* g.\]
Note that each $f|_C*g$ is supported on the bisection $CP\subseteq CK\subseteq \bigcup \CC$ as $C\in \CE\subseteq \CC^{-K}$. Then observe
\begin{align*}
    \psi_f(e*\varphi(g))-\psi_f(e)*g=&\sum_{C_0\in \CE}\sum_{\substack{C\in \CC, r(C)=C_0\\ C\cap P\neq \emptyset}}f|_C\cdot a_{r(C)}g(x_C)-\sum_{C_0\in \CE}a_{C_0}f|_{C_0}*g\\
    =&\sum_{C_0\in \CE}(\sum_{\substack{C\in \CC, r(C)=C_0\\ C\cap P\neq \emptyset}}f|_C\cdot a_{C_0}g(x_C)-a_{C_0}f|_{C_0}*g).
    \end{align*}
Now for $z\in \CG$, note that 
\begin{align*}
  &( \sum_{\substack{C\in \CC, r(C)=C_0\\ C\cap P\neq \emptyset}}f|_C\cdot a_{C_0}g(x_C))(z)\\
  &=\begin{cases}
			 f|_C(z)\cdot a_{C_0}g(x_C)\, ,  & \text{if } z\in C \text{ for some } C\in \CC \text{ with } r(C)=C_0 \text{ and }C\cap P\neq \emptyset  \\
			0\, , & \text{otherwise} 
		\end{cases}\\
&=\begin{cases}
			 f|_{C_0}(r(z))\cdot a_{C_0}g(x_C)\, ,  & \text{if } z\in C \text{ for some } C\in \CC \text{ with } r(C)=C_0 \text{ and }C\cap P\neq \emptyset  \\
			0\, , & \text{otherwise} 
		\end{cases}
        \end{align*}
   because $f$ is $\bigcup \CC$-invariant.  In addition, this implies that the function \[\sum_{\substack{C\in \CC, r(C)=C_0\\ C\cap P\neq \emptyset}}f|_C\cdot a_{C_0}g(x_C)\]
   is supported on $\bigcup \{C\in \CC: r(C)=C_0, C\cap P\neq \emptyset\}\subseteq B$. On the other hand, note that 
      \begin{align*}
          (a_{C_0}f|_{C_0}*g)(z)=
        \begin{cases}
            f|_{C_0}(r(z))a_{C_0}g(z)\, , & z\in C_0\cdot P\\
            0\, , &\text{otherwise}      \end{cases}
               \end{align*}  
Then recall that $C_0\cdot P\subseteq C_0\cdot K\subseteq \bigcup \CC$, which implies that 
\begin{align*}
          (a_{C_0}f|_{C_0}*g)(z)=
        \begin{cases}
            f|_{C_0}(r(z))a_{C_0}g(z)\, , & z\in C\cap P\text{ for some } C\in \CC \text{ with } r(C)=C_0\\
            0\, , &\text{otherwise}      \end{cases}
               \end{align*} 
and thus the function $a_{C_0}f|_{C_0}*g$ is also supported on $P\subseteq B$. This further implies that 
\[T_{C_0}:=\sum_{\substack{C\in \CC, r(C)=C_0\\ C\cap P\neq \emptyset}}f|_C\cdot a_{C_0}g(x_C)-a_{C_0}f|_{C_0}*g\]
is also supported on the bisection $B$ and therefore one has
\begin{align*}  \|T_{C_0}\|_r&=\|T_{C_0}\|_\infty=\sup_{\substack{C\in \CC, r(C)=C_0\\ C\cap P\neq \emptyset}}\sup_{z\in C}|f|_{C_0}(r(z))a_{C_0}g(x_C)-f|_{C_0}(r(z))a_{C_0}g(z)|\\
&\leq \|f|_{C_0}\|_\infty\cdot \|e\|_\infty\cdot  \max_{C \in \CC} \operatorname{osc}_{C} (g)\leq\|e\|_\infty\cdot  \max_{C \in \CC} \operatorname{osc}_{C} (g)
\end{align*}
because $f$ is a contraction. Finally, to be more precise, note that each $T_{C_0}$ is supported on the bisection  $C_0\cdot B$. Then note that for $C_0\neq C'_0\in \CE$, actually $C_0$ is disjoint from $C'_0$.  This implies that $T_{C_0}*T^*_{C'_0}=0$ because $T_{C_0}*T^*_{C'_0}$ is supported on $C_0BB^{-1}C'_0=C_0r(B)C'_0=\emptyset$. Moreover, $T^*_{C_0}*T_{C'_0}=0$ also holds because $T^*_{C_0}*T_{C'_0}$ is supported on $B^{-1}C_0C'_0B=\emptyset$. Therefore, one has
\begin{align*}
     \|\psi_f(e*\varphi(g))-\psi_f(e)*g\|_r=\|\sum_{C_0\in \CE}T_{C_0}\|_r=\sup_{C_0\in \CE}\|T_{C_0}\|_r\leq \|e\|_\infty\cdot  \max_{C \in \CC} \operatorname{osc}_{C} (g).
     \end{align*}
This has establish the first claim by setting $P=K$.

Now, in general, let $K$ be a compact set covered  by $n$ open bisections $B_1, \dots, B_n$ and denote by $\CV=\{B_1,\dots ,B_n, \CG\setminus K\}$ the collection of open sets that $\CC$ refines. Let $g$ be the function supported on $K$ and let $\{h_1, \dots, h_n\}$ be the unity of partition subordinate to the open cover $\{B_1, \dots, B_n\}$. Then the function $gh_i$ is supported on $P_i:=B_i\cap K\subseteq B_i$. Then by our assumption of $\CC$, if $C\in \CC$ satisfies $C\cap P_i\neq \emptyset$ then $C\subseteq B_i$. Then apply the above argument to $f$, $e$, $gh_i$, one obtains
\begin{align*}
     \|\psi_f(e*\varphi(gh_i))-\psi_f(e)*(gh_i)\|_r\leq \|e\|_\infty\cdot\max_{C\in \CC}\operatorname{osc}_{C}(gh_i). 
     \end{align*}
 Moreover, note that 
 \[|g(x)h_i(x)-g(y)h_i(y)|\leq |g(x)||h_i(x)-h_i(y)|+|h_i(y)||g(x)-g(y)|,\]
 which further implies that 
 \begin{align*}
     \|\psi_f(e*\varphi(gh_i))-\psi_f(e)*(gh_i)\|_r\leq \|e\|_\infty\cdot(\max_{C\in \CC}(\operatorname{osc}_{C}(g)+\|g\|_\infty\operatorname{osc}_C(h_i)). 
     \end{align*} 
Finally, since $g=\sum_{i=1}^ngh_i$,  observes $\varphi(g)=\sum_{i=1}^n\varphi(gh_i)$. This thus implies that
\[ \psi_f(e*\varphi(g))-\psi_f(e)*g=\sum_{i=1}^n(\psi_f(e*\varphi(gh_i))-\psi_f(e)*(gh_i)).\]
Therefore, one has
\[\|\psi_f(e*\varphi(g))-\psi_f(e)*g\|_r\leq \|e\|_\infty\cdot\sum_{i=1}^n(\max_{C\in \CC}(\operatorname{osc}_{C}(g)+\|g\|_\infty\operatorname{osc}_C(h_i)).\]
Using the same method, one also obtains
\[\left\| g *\psi_f(e) - \psi_f \left( \varphi (g)* e \right) \right\|_r\leq \|e\|_\infty\cdot\sum_{i=1}^n(\max_{C\in \CC}(\operatorname{osc}_{C}(g)+\|g\|_\infty\operatorname{osc}_C(h_i)).\]
This finishes the proof of the second claim.
\end{proof}

We shall need to use ``approximately flat'' functions on the unit space of a groupoid. 

\begin{defn} \label{defn:approximately-flat-functions}
	Let $\CG$ be a groupoid. Let $K \subseteq \CG$ and let $\varepsilon \geq 0$. Then a function $f \colon \CG\unitSpace \to \mathbb{C}$ is \emph{$(K,\varepsilon)$-flat} if for any $x \in K$, we have $|f(s(x)) - f(r(x))| \leq \varepsilon$. 
\end{defn}

These functions are ``approximately central'', to be made precise below. 

\begin{lem} \label{lem:approximately-flat-central}
	Let $\CG$ be a groupoid, let $K \subseteq \CG$ and let $\varepsilon \geq 0$.  
	Suppose $K$ can be covered by $n$ bisections in $\CG$. Then for any $(K,\varepsilon)$-flat function $f \colon \CG\unitSpace \to \mathbb{C}$ and any function $g \colon \CG \to \mathbb{C}$ supported on $K$, we have $\|[f,g]\| \leq n  \,\varepsilon \,\|g\|_\infty$ in $C^*(\CG)$. 
\end{lem}

\begin{proof}
First suppose $K$ is a set covered by one open bisection $O$. Then necessarily $K\subseteq O$ and thus $K$ itself is a bisection.  Then since $f$ is supported on $\GU$, which is also a bisection. Then Proposition \ref{prop: product of functions supported on bisections} implies that $f*g$ and $g*f$ are supported on $K$ since $K=\GU\cdot K=K\cdot \GU$. Moreover,  $(f*g)(x)=f(r(x))g(x)$ and $(g*f)(x)=g(x)f(s(x))$ hold for all $x\in K$. Note that $[f, g]=f*g-g*f$ is also supported on $K$ and therefore Proposition \ref{prop: norms on bisections} entails that 
\[\|[f, g]\|_r=\|[f, g]\|_\infty=\sup_{x\in K}|f(r(x))g(x)-g(x)f(s(x))|\leq \varepsilon \,\|g\|_\infty\]
because $f$ is $(K, \varepsilon)$-flat.

Now, suppose $K$ is a compact set covered by $n$ open bisections, say, $O_1, \dots, O_n$. Then one may choose a partition of unity $\{h_1, \dots, h_n\}$ subordinate to the open cover $\{O_1, \dots, O_n\}$. Then for each $1\leq i\leq n$, the function $gh_i$ is supported on $O_i\cap K$, which is a bisection contained in $O_i$. Then the argument above shows that 
\[\left\|[f,gh_i]\right\|_r\leq \sup_{x\in K\cap O_i}\left|f(r(x))gh_i(x)-gh_i(x)f(s(x))\right|\leq \varepsilon \, \|g\|_\infty \]
since $K\cap O_i\subseteq O_i$. Then because $g=\sum_{i=1}^ngh_i$, one has
\[ \left\|[f, g]\right\|_r=\left\|\sum_{i=1}^n[f, gh_i]\right\|_r\leq \sum_{i=1}^n\left\|[f, gh_i]\right\|_r\leq n \,\varepsilon \, \|g\|_\infty \]
as desired.
\end{proof}

Next we give a simple construction of ``approximate flat'' functions on a castle. The idea is to build a ``staircase'' by repeating the shrinking procedure in \Cref{defn:castle-shrinking}. 
To simplify our exposition, we shall restrict ourselves to the case of a finite principal groupoid $\CC$, which suffices for our purposes. 
Note that each subgroupoid $\CD$ of $\CC$ corresponds to a castle of $\CC$ via the mapping $x \mapsto \{x\}$. 
This correspondence allows us to define a shrinking $\CD^{-K}$ (with $K \subseteq \CC$) as another subgroupoid of $\CC$. 

\begin{lem}\label{lem: decay function}
    Let $\CC$ be a finite principal groupoid, let $n \in \mathbb{N}$, and let $K$ be a subset in $\CC$ containing $\CC^{(0)}$. 
    Let $\CC_0, \CC_1, \ldots, \CC_n$ be subgroupoids in $\CC$ with $\CC_0 = \CC$ and $\CC_k \subseteq (\CC_{k-1})^{-K}$ for $k = 1, \ldots, n$. 
    Then there exists a $(K, 1/n)$-flat function $e: \CC\unitSpace\to [0, 1]$ such that $e\equiv 1$ on $(\CC_n)\unitSpace$ and $e\equiv 0$ on $\CC\unitSpace\setminus (\CC_1)\unitSpace$. 
\end{lem}
\begin{proof}
    Define $e: \CC\unitSpace\to [0, 1]$ by $e\equiv 1$ on $\CP_n:=(\CC_n)\unitSpace$, and $e\equiv m/n$ on $\CP_{m}:=(\CC_m)\unitSpace\setminus (\CC_{m+1})\unitSpace$ for $m=0, 1, \dots, n-1$. 
    Now for any $x\in K$, observe that 
    \begin{enumerate}[label=(\roman*)]
        \item  if $s(x)\in \CP_0$, then $r(x)\in \CP_1\cup \CP_0$;
        \item if $s(x)\in \CP_m$ for some $m=1,\dots, n-1$, then $r(x)\in \CP_{m-1}\cup \CP_m \cup \CP_{m+1}$;
         \item if $s(x)\in \CP_n$, then $r(x)\in \CP_n\cup \CP_{n-1}$.
    \end{enumerate}
    This implies that $|e(s(x))-e(r(x))|\leq 1/n$ and thus $e$ is $(K, 1/n)$-flat.
\end{proof}

\begin{lem}\label{lem: matrix-embedding}
    Let $\CG$ be a groupoid and $n, N\in \N_+$. Let $\CC, \CD$ be open castles in $\CG$ such that $\CC$ is faithfully nested in $\CD$  with a minimal multiplicity of at least $nN$. Then there exists a canonical $*$-embedding $\phi: M_n(\C)\to C^*(\CommCD{\CC}{\CD})$ such that 
    \[\sup_{\mu\in M(\CG)}\mu \left(\bigcup\{C\in \CC\unitSpace: \phi(1_n)(\CommCD{C}{\CD})=0\} \right)\le1/N.\]
    \end{lem}
    \begin{proof}
        Let $D\in \CD\unitSpace$ and a tower $\CT\subseteq \CC$. Because $\CC$ is faithfully nested in $\CD$  with a minimal multiplicity at least $nN$, there are at least $nN$ levels $C\in \CT\unitSpace$ contained in $D$, i.e., 
        \[|\RedCD{\CT}{D}|=|\{C\in \CT\unitSpace: C\subseteq D\}|\geq nN\]
        by Remark \ref{rem:nesting-comm-red-basic}(\ref{rem:nesting-comm-red-basic:Red-Res}). Then divide the set $\RedCD{\CT}{D}$ into $n$ subsets with the same cardinality $\lfloor\frac{|\RedCD{\CT}{D}|}{n}\rfloor$. Denote by $Q_{\CT, D, 1}, \dots, Q_{\CT, D, n}$ for these families and note that 
        \[|\RedCD{\CT}{D}\setminus \bigsqcup_{i=1}^nQ_i|\leq n-1.\]
        Now, choose bijective maps $f_{\CT, D, i}: Q_{\CT, D, 1}\to Q_{\CT, D, i}$ for each $i=1,\dots, n$ such that $f_{\CT, D, 1}=\id_{Q_{\CT, D, 1}}$. Define 
        \[Q_{\CT, D, i, 1}=\{C\in \CT: s(C)\in Q_{\CT, D, 1}, r(C)=f_{\CT, D, i}(s(C))\}\]
        and
        \[Q_{\CT, D, i, j}=Q_{\CT, D, i, 1}\cdot Q^{-1}_{\CT, D, j, 1}.\]
        Moreover, denote by $\CS$ a representative transversal set for $\CD\unitSpace/\sim_\CD$.   Define 
        \[Q_{i, j}=\bigcup\{Q_{\CT, D, i, j}: D\in \CS, \CT\in \CC/\sim_\CC\}\]
        and an injective $*$-homomorphism $\phi: M_n(\C)\to C^*(\CommCD{\CC}{\CD})$ spanned by
        \[\phi: e_{ij}\mapsto \sum_{C\in Q_{i,j}}\chi_{\CommCD{C}{\CD}}\]
        where $e_{ij}$ denotes the $(i, j)$-generator of $M_n(\C)$.

        Then let $\mu\in M(\CG)$. Then by definition of $\phi$, for each $D\in \CD\unitSpace$ and $\CC$-tower $\CT$, one has
        \[\{C\in \CC\unitSpace: C\subseteq D, C\in \CT\unitSpace, \phi(1_n)(\CommCD{C}{\CD})=0\}=\RedCD{\CT}{D}\setminus \bigcup_{i=1}^nQ_{\CT, D, i}\]
        which has cardinality bounded by $n-1$. Choosing a $C_\CT\in \CT\unitSpace$, this thus implies that
        \begin{equation}\label{key for1}
            \mu(\bigcup\{C\in \CC\unitSpace: \phi(1_n)(\CommCD{C}{\CD})=0\} \leq \sum_{\CT\in \CC/\sim_\CC}\sum_{D\in \CD\unitSpace}(n-1)\mu(C_\CT).
        \end{equation}
    On the other hand, for each $\CC$-tower $\CT$,  note that 
\begin{equation}\label{kery for2}
    \sum_{\CT\in \CC/\sim_\CC}\sum_{D\in \CD\unitSpace}nN\mu(C_\CT)\leq \mu(\bigcup \CC\unitSpace)\leq 1
\end{equation}
       because $\CC$ is faithfully nested in $\CD$ with multiplicity at least $nN$. Combining inequalities \eqref{key for1}, \eqref{kery for2}, and the fact that $\mu\in M(\CG)$ is arbitrary, one finally has
 \[\sup_{\mu\in M(\CG)}\mu\left(\bigcup\{C\in \CC\unitSpace: \phi(1_n)(\CommCD{C}{\CD})=0\}\right)\le1/N. \qedhere \]
 \end{proof}

\begin{thm}
	\label{thm: tracial Z stable} Let $\CG$ be a minimal
	\'{e}tale groupoid with a compact infinite unit space. If $\CG$ is almost
	elementary, then $C_{\operatorname{r}}^{*}(\CG)$ is simple and tracially $\CZ$-stable. 
\end{thm}
\begin{proof}
	By Theorem~\ref{thm:effective-AE}, the groupoid $\CG$ is effective. Combining this with the minimality of $\CG$, it is a standard fact that $C_r^{*}(\CG)$ is simple (see, e.g., \cite[Proposition 10.3.7]{Sims-groupoids}). 
	
	It remains to show $\CG$ is tracially $\CZ$-stable. To this end, let $F\subseteq C_{\operatorname{r}}^{*}(\CG)$, $\varepsilon>0$, $a\in C_{\operatorname{r}}^{*}(\CG)_{+} \setminus \{0\}$ and $n\in\N$ be given as in Definition~\ref{defn:tracial-Z-stable}. 
	Observe that since $C_{\operatorname{c}}(\CG)$ is dense in $C_{\operatorname{r}}^{*}(\CG)$, we may assume without loss of generality that $F = \left\{ f_1, \ldots, f_m \right\} \subseteq C_{\operatorname{c}}(\CG)$.
    One may further assume that each $K_i:=\supp(f_i)$ is a bisection and contained in an open bisection $B_i$ for $i=1,\dots, m$. Moreover, one may also assume $a\in C(\CG\unitSpace)_+ \setminus \{0\}$ by Lemma \ref{lem: small function}. Define $O:=\osupp{a}$ and choose three disjoint non-empty open sets $O_0, O_1, O_2$ in $O$ such that $\min_{\mu\in M(\CG)}\mu(O_2)>\varepsilon_0$ for some $\varepsilon_0>0$ by minimality. Define $K=\bigcup_{i=1}^mK_i$ and another compact set 
    \[L := (\bigcup_{i=1}^{m} \overline{B_i})\cup (\bigcup_{i=1}^{m} \overline{B_i})^{-1}\cup \GU\]
    satisfying $K\subseteq L$. Choose a finite open cover $\CV$ of $L$ refining the open cover $\bigvee_{i=1}^m\{B_i, \CG\setminus K_i\}$ of $\CG$ and such that for any $V \in \CV$ and any $i \in \{1, \ldots, m\}$, we have $\sup_{x, y \in V} | f_i(x) - f_i(y) |< \varepsilon/3$.
	
	Choose an integer $k\geq3m/\varepsilon$. By Lemma~\ref{lem:AE-fortified}, there is a fortified castle $\CE$
	satisfying 
	\begin{enumerate}
		\item \label{small E-level}
		for any $E_+ \in \CE_+$ with $E_+ \cap L \not= \varnothing$, there is $V \in \CV$ such that $E_+ \subseteq V$ (this further implies that if $E_+\cap K_i\neq \emptyset$ then $E_+\subseteq B_i$ for each $i=1,\dots, m$ and in particular $E_+\subseteq L$); 
		\item
		$\bigcup \CE_+$ is precompact in $\CG$;  
		\item\label{small E-reminder}
        $\GU\setminus \bigcup (\CE^{-L^k}_-)\unitSpace\prec_\CG O_0$.
			\end{enumerate}
Remark \ref{rem:subequivalence-basic}\eqref{rem:subequivalence-basic: compact to open} shows that there exists an open set $Q_0$ containing the compact set $\GU\setminus \bigcup (\CE^{-L^k}_-)\unitSpace$  such that $Q_0\precsim_\CG O_0$.
	
	Then, choose $N$ large enough such that $N>1/\varepsilon_0$. For this $N$, the fortified castle $\CE=(\CE_+, \CE_-)$, the cover $\CV$,  and the open set $O_1$, Proposition~\ref{prop:AE-nesting-minimal} shows that there exists a fortified open castle $\CF$ and a castle $\CB$ such that 
	\begin{enumerate}[label=(\roman*)]
           \item\label{AE-nesting-minimal:nesting-recurrent} $\CF_+$ is faithfully nested in $\CE_-$ with a minimal multiplicity at least $nN$.     
		\item\label{AE-nesting-minimal:extend} $\CF_+ \subseteq \CB^{-L}$ (in particular, $\CF_+$ is $L$-extendable);  
		\item\label{AE-nesting-minimal:refine} $\CD\unitSpace$ refines $\CV$;    
		\item\label{AE-nesting-minimal:remainder} 
		$\left( \overline{\bigcup \CE_-^{(0)}} \right)\setminus\ \left( \bigcup \CF_-^{(0)} \right) \prec_{\CG}  O_1$. 
	\end{enumerate}	
Remark \ref{rem:subequivalence-basic}\eqref{rem:subequivalence-basic: compact to open} shows that there exists an open set $Q_1$ containing the compact set $\left( \overline{\bigcup \CE_-^{(0)}} \right)\setminus\ \left( \bigcup \CF_-^{(0)} \right)$  such that $Q_1\precsim_\CG O_1$. Then Lemma \ref{lem: matrix-embedding} implies there exists an injective $*$-homomorphism $\phi: M_n(\C)\to C^*(\CommCD{\CF_+}{\CE_+})$ such that  
\[\sup_{\mu\in M(\CG)}\mu(\bigcup\{F_+\in \CF_+\unitSpace: \phi(1_n)(\CommCD{\CF_+}{\CE_+})=0\})\le1/N.\]

By the condition~\eqref{small E-level} for $\CE$, every member of the collection $K':=\{E_+\in \CE_+: E_+\cap K\neq \emptyset\}$ is a subset of $L$, whence the sequence $(\CE_+^{-L^i})_{i = 0}^n$ of castles, when considered as subgroupoids of $\CE_+$, satisfy $\CE_+^{-L^i} \subseteq (\CE_+^{-L^{i-1}})^{K'}$ for $i = 1,2,\ldots, k$. Applying \Cref{lem: decay function}, with $K$ replaced by $K'$, to the sequence $(\CE_+^{-L^i})_{i = 0}^n$  yields a $(K', 1/k)$-flat function $e\in C(\CE_+\unitSpace)$ satisfying $e\equiv 1$ on $(\CE_+^{-L^k})\unitSpace$ and $e\equiv 0$ on $\CE_+\unitSpace\setminus (\CE_+^{-L})\unitSpace$.

Choose a positive contraction $f\in C_0(\bigcup \CF_+)\cap C_c(\CG)$ that is $(\bigcup \CF_+)$-invariant and a positive contraction $g\in C_0(\bigcup \CE_+)\cap C_c(\CG)$ that is $(\bigcup \CE_+)$-invariant. Moreover, we can make $f\equiv1$ on $\overline{\bigcup \CF_-}$ and $g\equiv1$ on $\overline{\bigcup \CE_-}$. Note that this also implies $g\equiv 1$ on $\bigcup \CF_+\unitSpace$ because $\CF_+$ is nested in $\CE_-$.

Let $\psi_f \colon C^*(\CF_+) \to C^*_r(\CG)$, $\psi_g \colon C^*(\CE_+) \to C^*_r(\CG)$ be c.p.c.\ order zero maps and canonical embedding $\psi_{\chi_{\CommCD{\CF_+}{\CE_+}}}: C^*(\CommCD{\CF_+}{\CE_+})\to C^*(\CF_+)$ obtained in Proposition~\ref{prop:castle-order-zero} and Lemma \ref{lem:nesting-comm-red-order-zero}, respectively. Now we define a c.p.c order zero map $\psi: M_n(\C)\to C^*_r(\CG)$ by 
\[\psi(\cdot)=(\psi_f\circ\psi_{\chi_{\CommCD{\CF_+}{\CE_+}}}\circ \phi (\cdot))* \psi_g(e).\]
Now, let $f_i$ and $b\in M_n(\C)$ with $\|b\|=1$ be given as in the beginning.	Note that actually the function $e$ belongs to $C((\CE^{-L}_+)\unitSpace)\subseteq C((\CE^{-K}_+)\unitSpace)\subseteq C((\CE^{-K_i}_+)\unitSpace)$.
	Then applying Lemma~\ref{lem:castle-local-approximation} to $K_i\subseteq B_i$, the function $e$, $\psi_g$ above, and a choice of $x_{E_+}\in E_+$ for each $E_+\in \CE_+$, together with condition \eqref{small E-level} for $\CE_+$, one has
\begin{align}\label{ineq1}
  & \max \left\{ \left\| \psi_g(e) f_i - \psi_g \left( e \varphi (f_i) \right) \right\|  ,  \left\| f_i \psi_g(e) - \psi_g \left( \varphi (f_i) e \right) \right\| \right\}\\
  &\leq \|e\|_\infty \, \max_{E_+ \in \CE_+} \operatorname{osc}_{E_+} (f_i) 
    \leq\varepsilon/3\nonumber, 
\end{align}
in which $\varphi: C_c(\CG)\to C^*(\CE_+)$ is defined to yield $\varphi(f_i)=(f_i(x_{E_+}))_{E_+\in \CE_+}$. Therefore, $\varphi(f_i)\in C^*(\CE_+)$ is supported on $K'=\{E_+\in \CE_+: E_+\cap K\neq \emptyset\}$. Moreover, by the condition \eqref{small E-level} for $\CE$, the set $K'$ is covered by all $\{E_+\in \CE_+: E_+\subseteq B_i\}$ for $i=1,\dots, m$, which are bisections in the groupoid $\CE_+$,. Then, Lemma \ref{lem:approximately-flat-central} shows that 
\begin{equation}\label{ineq2}
    \|[\varphi(f_i), e]\|\leq m/k\leq \varepsilon/3.
    \end{equation}
Then, recall that the images of $\psi_f\circ \psi_{\chi_{\CommCD{\CF_+}{\CE_+}}}\circ \phi$ and  of $\psi_g$ commute by Lemma~\ref{lem:nesting-comm-red-order-zero}. Combining this with inequalities \eqref{ineq1} and \eqref{ineq2}, one has
\begin{align*}
    \psi(b)*f_i&=\psi_f(\psi_{\chi_{\CommCD{\CF_+}{\CE_+}}}(\phi(b)))*\psi_g(e)*f_i\\
    &\approx_{\frac{\varepsilon}{3}} \psi_f(\psi_{\chi_{\CommCD{\CF_+}{\CE_+}}}(\phi(b)))*\psi_g(e\varphi(f_i))\\
    &\approx_{\frac{\varepsilon}{3}}\psi_f(\psi_{\chi_{\CommCD{\CF_+}{\CE_+}}}(\phi(b)))*\psi_g(\varphi(f_i)e)\\ 
&=\psi_g(\varphi(f_i)e)*\psi_f(\psi_{\chi_{\CommCD{\CF_+}{\CE_+}}}(\phi(b)))\\
&\approx_{\frac{\varepsilon}{3}} f_i*\psi_g(e)*\psi_f(\psi_{\chi_{\CommCD{\CF_+}{\CE_+}}}(\phi(b)))\\
&=f_i*\psi(b),
\end{align*}
which implies that $\|[\psi(b), f_i]\|<\varepsilon$.

Now, it is left to show $1_{C^*_r(\CG)}-\psi(1_{M_n(\C)})\precsim a$. To do so, first write
\[Q_2 :=\bigcup\{F_+\in \CF_+\unitSpace: \phi(1_n)(\CommCD{F_+}{\CE_+})=0\}\]
for simplicity. Since 
\[\sup_{\mu\in M(\CG)}\mu(Q)<1/N<\varepsilon_0<\min_{\mu\in M(\CG)}\mu(O_2),\]
Theorem \ref{thm:AE-DSC} implies that $Q_2\precsim_\CG O_2$.

Then because $f\equiv 1$ on $\bigcup \CF_-\unitSpace$, observe $h:=1_{C^*_r(\CG)}-\psi(1_{M_n(\C)})\in C(\GU)_+$ is supported on 
\[S:=(\GU\setminus \bigcup (\CE^{-L^k}_-)\unitSpace)\cup (\left( \overline{\bigcup \CE_-^{(0)}} \right)\setminus\ \left( \bigcup \CF_-^{(0)} \right))\cup Q_2\]
which is a subset of $Q_0\cup Q_1\cup Q_2$. This further implies that 
\[\osupp{h}\subseteq Q_0\cup Q_1\cup Q_2\precsim_\CG O_0\sqcup O_1\sqcup O_2\subseteq O=\osupp{a}\]
by Remark \ref{rem:subequivalence-basic}\eqref{rem:subequivalence-basic:additive-precsim}. Thus, Lemma \ref{lem:groupoid subequi to Cuntz subequi} entails 
\[1_{C^*_r(\CG)}-\psi(1_{M_n(\C)})=h\precsim a.\]
Therefore, the $C^*$-algebra $C^*_r(\CG)$ is tracial $\CZ$-stable.
\end{proof}

Applying Theorem \ref{thm: tracial Z stable} directly, we have the following natural
corollaries.
\begin{cor}
	\label{cor: strict comparison} 
	Let $\CG$ be a minimal second countable \'{e}tale groupoid with a compact unit space. 
	If $\CG$ is almost elementary, then $C_{\operatorname{r}}^{*}(\CG)$ has strict comparison
	for positive elements. 
\end{cor}
\begin{proof}
	If $\CG\unitSpace$ is an infinite set, then 
	Theorem~\ref{thm: tracial Z stable} implies that $C_{\operatorname{r}}^{*}(\CG)$ is a simple unital
	tracially $\CZ$-stable $C^{*}$-algebra, whence it follows from Theorem~3.3 in \cite{HirshbergOrovitz2013Tracially}
	that $C_{\operatorname{r}}^{*}(\CG)$ has strict comparison of positive
	elements. 
	
	If $\CG\unitSpace$ is a finite set, then by Corollary~\ref{cor:AE-finite-pair}, $C_{\operatorname{r}}^{*}(\CG)$ is a matrix algebra, which also has strict comparison of positive elements. 
\end{proof}

We remark that Corollary~\ref{cor: strict comparison} does not assume nuclearity
of $C_{\operatorname{r}}^{*}(\CG)$. However, given the current technologies and interests in the classification and structure theory of $C^*$-algebras, the most interesting case appears to be where $\CG$ is (topologically) amenable.

\begin{cor}
	\label{8.12} Let $\CG$ be a minimal (topologically) amenable second
	countable \'{e}tale groupoid on a compact space. Suppose $\CG$
	is almost elementary. Then $C_{\operatorname{r}}^{*}(\CG)$ is unital simple separable
	nuclear and $\CZ$-stable and thus has nuclear dimension one. In addition,
	in this case $C_{\operatorname{r}}^{*}(\CG)$ is classified by its Elliott invariant.
	Finally, If $M(\CG)\neq\emptyset$, then $C_{\operatorname{r}}^{*}(\CG)$ is quasidiagonal
	and if $M(\CG)=\emptyset$ then $C_{\operatorname{r}}^{*}(\CG)$ is a unital Kirchberg
	algebra. \end{cor}
\begin{proof}
	Since $\CG$ is assumed to be amenable, Theorem \ref{thm: tracial Z stable} implies
	$C_{\operatorname{r}}^{*}(\CG)$ is unital simple separable nuclear and tracially
	$\CZ$-stable and thus $\CZ$-stable by Theorem 4.1 in \cite{HirshbergOrovitz2013Tracially}.
	In this case, the nuclear dimension $\dimnuc(C_{\operatorname{r}}^{*}(\CG))=1$ by
	Theorem A and Corollary C in \cite{CastillejosEvingtonTikuisisWhiteWinter2019Nuclear}.
	Therefore, $C_{\operatorname{r}}^{*}(\CG)$ is classified by Elliott invariant by
	the recent progress of classification theorem for unital simple nuclear
	separable $C^{*}$-algebras having finite nuclear dimendion and satisfying
	the UCT via combining  \cite{GLN2020a, GLN2020b,ElliottGongLinNiu2015classification,TikuisisWhiteWinter2017Quasidiagonality,CGSTW,
Phillips2000classification}.
	Finally, if $M(\CG)\neq\emptyset$ then there is a non-zero tracial
	state on $C_{\operatorname{r}}^{*}(\CG)$. This implies that $C_{\operatorname{r}}^{*}(\CG)$ is
	stably finite and thus quasidiagonal by Corollary 6.1 in \cite{TikuisisWhiteWinter2017Quasidiagonality}.
	On the other hand, if $M(\CG)=\emptyset$, then $C_{\operatorname{r}}^{*}(\CG)$
	is traceless and thus $C_{\operatorname{r}}^{*}(\CG)$ is purely infinite by Corollary
	5.1 in \cite{Rordam2004stable}. Therefore, in this case $C_{\operatorname{r}}^{*}(\CG)$
	is a unital Kirchberg algebra. 
\end{proof}

Then, combining Theorem \ref{thm:Matui-AE} and \ref{thm: tracial Z stable}, we have the following result.

\begin{cor}
	\label{cor: Matui AF tracial Z} Let $\CG$ be an ample minimal $\sigma$-compact groupoid on a compact space. Suppose $\CG$ is
	almost finite in Matui's sense. Then $C_{\operatorname{r}}^{*}(\CG)$ is tracially
	$\CZ$-stable and thus has the strict comparison for positive elements.
	If we assume $\CG$ is also amenable and second countable then $C_{\operatorname{r}}^{*}(\CG)$ is $\CZ$-stable
	and quasidiagonal. \end{cor}

In addition, by combining Theorem \ref{thm:Kerr-AE} and \ref{thm: tracial Z stable}, we may recover the following result due to Kerr in \cite{Kerr2020Dimension}.
\begin{cor}
	\label{8.14} Let $\alpha:\Gamma\curvearrowright X$ be a minimal
	action of a countable discrete amenable group $\Gamma$ on a
	compact metrizable space $X$. Suppose $\alpha$ is almost finite
	in Kerr's sense. Then the crossed product $C(X)\rtimes_{r}\Gamma$
	is $\CZ$-stable and quasidiagonal. \end{cor}

We finally provide several applications of our result on (tracial) $\CZ$-stability
of almost finite groupoids.
\begin{example}
	\label{8.15} Ito, Whittaker and Zacharias \cite{ItoWhittakerZacharias-tiling} established $\CZ$-stability for Kellendonk's
	$C^{*}$-algebra associated to an aperiodic and repetitive tiling on $\mathbb{R}^d$ with finite
	local complexity. 
	This $C^{*}$-algebra can be presented as the reduced groupoid $C^{*}$-algebra
	of a minimal, second countable, principle, \'{e}tale groupoid called the tiling groupoid. 
	Directly applying a result of Matui \cite[Remark~6.4]{Matui2012Homology}, the authors of \cite{ItoWhittakerZacharias-tiling} observed that tiling groupoids are almost finite in the sense of Matui. 
	They then obtained tracial $\CZ$-stability by carefully adapting the techniques in \cite{Kerr2020Dimension} for group actions, which are based on Ornstein-Weiss quasitilings on groups, 
	to the setting of tiling groupoids acting on their unit spaces. 
	Now, in view of Theorem~\ref{thm:Matui-AE}, the main result in \cite{ItoWhittakerZacharias-tiling} can be subsumed in Corollary~\ref{cor: Matui AF tracial Z} above. 
\end{example}

\begin{example}
	In \cite[Theorem 6]{Elek-qualitative}, Elek constructed
	a class of groupoids called \textit{geometric groupoids}, that are minimal, 
	principal, almost finite, ample, second countable, \'{e}tale, but not amenable.
	By Corollary~\ref{cor: Matui AF tracial Z}, it follows that for any such geometric groupoid $\CG$, its $C^*$-algebra $C_{\operatorname{r}}^{*}(\CG)$ is a
	non-nuclear tracially $\CZ$-stable $C^*$-algebra. 
\end{example}

\begin{example}
    It was recently shown in \cite{ELNZ2026} that for a countable amenable discrete group $\Gamma$, the natural action of $\Gamma$ on the universal minimal set $M\subset \beta\Gamma$ (which is typically non-metrizable) gives rise to a transformation groupoid that is almost finite in Matui's sense (see \cite[Corollary 3.6]{ELNZ2026}). Therefore, its $C^*$-algebra $C(M)\rtimes \Gamma$  is non-separable, tracial $\CZ$-stable and thus has the strict comparison by Corollary \ref{cor: Matui AF tracial Z}. This result should be compared to \cite[Theorem 1.2]{ELNZ2026}.
\end{example}

\section*{Acknowledgements}
This research is supported by National Key R\&D Program of China 2022YFA100700. 
The first author is supported by NSFC No.~12571133. 
The second author is partially supported by NSFC Key Program No.~12231005. 

The authors would like to thank David Kerr for his helpful comments during a number of discussions, as well as Guoliang Yu and Zhizhang Xie for their hospitality during the authors' visits to Texas A\&M University, where a part of the work was completed.

\bibliographystyle{alpha}
\bibliography{Almost_elementary}

\end{document}